\documentclass[12pt,reqno]{amsart}

\usepackage{amssymb,amscd,upref,xspace,varioref}

\makeatletter
\ifcase \@ptsize \or
\or
\fi
\makeatother

\newtheorem{theorem}{Theorem}[section]
\newtheorem{lemma}[theorem]{Lemma}
\newtheorem{corollary}[theorem]{Corollary}
\newtheorem{proposition}[theorem]{Proposition}

\theoremstyle{definition}
\newtheorem{definition}{Definition}[section]
\newtheorem{example}{Example}[section]

\theoremstyle{remark}
\newtheorem{remark}{Remark}[section]
\newtheorem{xca}{Exercise}[section]

\numberwithin{equation}{section}
\numberwithin{figure}{section}

\newcommand{\cprime}{\/{\mathsurround=0pt$'$}}

\DeclareFontFamily{OML}{cyr}{}
\DeclareFontShape{OML}{cyr}{m}{n}{
   <5> <6> <7> <8> <9> gen * wncyr
   <10> <10.95> <12> <14.4> <17.28> <20.74> <24.88> wncyr10
  }{}
\DeclareSymbolFont{rusletters}{OML}{cyr}{m}{n}
\DeclareSymbolFontAlphabet{\rusmath}{rusletters}
\DeclareMathSymbol\re{\rusmath}{rusletters}{"03}

\makeatletter
\let\reF=\ref
\newcommand{\eqreF}[1]{\textup{\tagform@{\reF{#1}}}}
\DeclareRobustCommand\veqref[1]{%
  \unskip~\eqreF{#1}%
  \@vpageref[\unskip]{#1}%
}
\makeatother

\newcommand{\al}{{\alpha}}      \newcommand{\om}{{\omega}}
\newcommand{\be}{{\beta}}       \newcommand{\Om}{{\Omega}}
\newcommand{\ga}{{\gamma}}      \newcommand{\si}{{\sigma}}
\newcommand{\Ga}{{\Gamma}}      
\newcommand{\de}{{\delta}}      \newcommand{\ta}{{\theta}}
\newcommand{\De}{{\Delta}}      
\newcommand{\ep}{{\epsilon}}    
      \newcommand{\vf}{{\varphi}}
\newcommand{\la}{{\lambda}}     \newcommand{\ve}{{\varepsilon}}
\newcommand{\La}{{\Lambda}}     
\newcommand{\na}{{\nabla}}

\newcommand{\bD}{{\mathbf{D}}}
\newcommand{\bE}{\mathbf{E}}

\newcommand{\CA}{\mathcal{A}}   \newcommand{\CN}{\mathcal{N}}
\newcommand{\CB}{\mathcal{B}}   
\newcommand{\CC}{\mathcal{C}}   \newcommand{\CP}{\mathcal{P}}
\newcommand{\cD}{\mathcal{D}}   
\newcommand{\CE}{\mathcal{E}}   \newcommand{\CR}{\mathcal{R}}
\newcommand{\CF}{\mathcal{F}}   
\newcommand{\CG}{\mathcal{G}}   \newcommand{\CT}{\mathcal{T}}
\newcommand{\CH}{\mathcal{H}}   \newcommand{\CU}{\mathcal{U}}
\newcommand{\CI}{\mathcal{I}}   \newcommand{\CV}{\mathcal{V}}
\newcommand{\CJ}{\mathcal{J}}   
\newcommand{\CK}{\mathcal{K}}   
\newcommand{\CL}{\mathcal{L}}   
\newcommand{\CM}{\mathcal{M}}
\let\E=\CE
\let\F=\CF
\let\R=\CR

\newcommand{\BBR}{\mathbb{R}}    
    \newcommand{\BBZ}{\mathbb{Z}}
\newcommand{\BB}[1]{{\mathbb{{#1}}}}

\newcommand{\bex}{\begin{example}}
\newcommand{\eex}{\end{example}}
\newcommand{\bpr}{\begin{proposition}}
\newcommand{\epr}{\end{proposition}}
\newcommand{\bre}{\begin{remark}}
\newcommand{\ere}{\end{remark}}
\newcommand{\bde}{\begin{definition}}
\newcommand{\ede}{\end{definition}}
\newcommand{\bth}{\begin{theorem}}
\newcommand{\ethm}{\end{theorem}}
\newcommand{\ble}{\begin{lemma}}
\newcommand{\ele}{\end{lemma}}
\newcommand{\bco}{\begin{corollary}}
\newcommand{\eco}{\end{corollary}}
\newcommand{\eq}[1]{\begin{equation}\label{#1}}

\newcommand{\abs}[1]{\lvert#1\rvert}
\newcommand{\norm}[1]{\lVert#1\rVert}
\newcommand{\dd}[2]{\dfrac{\partial#1}{\partial#2}}
\newcommand{\id}{\mathrm{id}}
\newcommand{\gl}{\mathfrak{gl}}
\newcommand{\Hom}{\mathrm{Hom}}
\newcommand{\hatotimes}{\mathbin{\widehat{\otimes}}}

\newcommand{\pdr}[2]{\dfrac{\partial#1}{\partial#2}}
\newcommand{\Gal}{\Gamma_{\mathrm{loc}}}
\newcommand{\ldb}{{\mathrm{[\![}}}
\newcommand{\rdb}{{\mathrm{]\!]}}}
\newcommand{\rnij}[2]{\ldb{#1},{#2}\rdb^{\rm rn}}
\newcommand{\fnij}[2]{\ldb{#1},{#2}\rdb^{\rm fn}}
\newcommand{\rest}[1]{\left|_{{#1}}\right.}

\let\ot=\otimes
\let\pat=\partial
\let\Ra=\Rightarrow

\let\sbs=\subset
\let\sps=\supset

\let\wg=\wedge
\let\vk=\varkappa

\newcommand{\cnst}{\mathrm{const}}

\newcommand{\xra}{\xrightarrow}
\newcommand{\xla}{\xleftarrow}
\newcommand{\hra}{\hookrightarrow}

\newcommand{\Smbl}{\mathrm{Smbl}}
\newcommand{\Diff}{\mathrm{Dif{}f}}
\newcommand{\hH}{\bar{H}}
\newcommand{\hS}{\bar{\mathrm{S}}}
\newcommand{\CDiff}{\mathcal{C}\mathrm{Dif{}f}}
\newcommand{\CDer}{\mathcal{C}\mathrm{D}}
\newcommand{\Dr}{\mathrm{D}}
\newcommand{\Dv}{{\mathrm{D}^v}}
\newcommand{\Ber}{{\mathrm{Ber}}}
\newcommand{\ip}{\mathrm{i}}
\newcommand{\Ld}{\mathrm{L}}
\newcommand{\tp}{\mathrm{tp}}
\newcommand{\J}{\bar{\mathcal{J}}}
\newcommand{\hL}{\bar{\Lambda}}
\newcommand{\hj}{\bar{\jmath}}
\newcommand{\ji}{{\bar{\jmath}_{\infty}}}
\newcommand{\hd}{\bar{d}}
\newcommand{\Ei}{{\mathcal{E}^{\infty}}}
\newcommand{\Ji}{{J^{\infty}}}
\newcommand{\Ci}{{C^{\infty}}}
\newcommand{\sltwo}{\mathfrak{sl}_2}
\newcommand{\CLa}[1]{\mathcal{C}^{#1}\Lambda}
\newcommand{\hJi}{\bar{\mathcal{J}}^{\infty}}

\DeclareMathOperator{\ad}{ad}

\DeclareMathOperator{\alt}{alt}
\DeclareMathOperator{\Ann}{Ann}

\DeclareMathOperator{\cl}{cl}
\DeclareMathOperator{\codim}{codim}
\DeclareMathOperator{\coker}{coker}

\DeclareMathOperator{\im}{im}

\DeclareMathOperator{\ord}{ord}
\DeclareMathOperator{\rank}{rank}

\DeclareMathOperator{\smbl}{smbl}

\DeclareMathOperator{\sym}{sym}

\newcommand{\cm}{$\mathcal{C}$-mod\-ule\xspace}
\newcommand{\cms}{$\mathcal{C}$-mod\-ules\xspace}
\newcommand{\fm}{$\mathcal{F}$-mod\-ule\xspace}

\newcommand{\cd}{$\mathcal{C}$-dif\-fer\-en\-tial\xspace}
\newcommand{\css}{the Vinogradov $\mathcal{C}$-spec\-tral
                  sequence\xspace}
\newcommand{\kcc}{the $\mathcal{C}$-co\-ho\-mol\-ogy\xspace}
\newcommand{\ldot}{\,\dot{}}
\newcommand{\FN}{Fr\"o\-li\-cher--Nij\-en\-huis\xspace}

\begin{document}

\begin{titlepage}

\hfill Preprint DIPS 7/98

\hfill math.DG/9808130

\vspace*{2cm}
\begin{center}
\large{\textbf{\uppercase{
Homological methods \\[9pt]
in equations of mathematical physics}}}\footnote{Lectures given in
August 1998 at the International Summer School in Levo\v{c}a,
\hbox{Slovakia}.\\ This work was supported in part by RFBR grant
97-01-00462 and INTAS grant 96-0793}
\end{center}

\vspace{3.5cm}
\begin{flushleft}
\textsc{Joseph KRASIL{\cprime}SHCHIK}\footnote{Correspondence
to: J. Krasil{\cprime}shchik, 1st Tverskoy-Yamskoy per., 14, apt.\
45, \\
125047 Moscow, Russia \\
\textit{E-mail}: \texttt{josephk@glasnet.ru}}\\[5pt]
\textit{Independent University of Moscow and \\
The Diffiety Institute, \\
Moscow, Russia}\\[7pt]
and \\[7pt]
\textsc{Alexander VERBOVETSKY}\footnote{Correspondence to:
A. Verbovetsky, Profsoyuznaya 98-9-132, 117485 Moscow, Russia \\
\textit{E-mail}: \texttt{verbovet@mail.ecfor.rssi.ru}}\\[5pt]
\textit{Moscow State Technical University and \\
The Diffiety Institute, \\
Moscow, Russia}
\end{flushleft}
\end{titlepage}

\setcounter{page}{2}

\tableofcontents
\newpage

\section*{Introduction}
Mentioning (co)homology theory in the context of differential equations would
sound a bit ridiculous some 30--40 years ago: what could be in common
between the essentially analytical, dealing with functional spaces theory of
partial differential equations (PDE) and rather abstract and algebraic
cohomologies?

Nevertheless, the first meeting of the theories took place in the papers by
D.~Spencer and his school (\cite{Spencer1,Goldsch1}), where cohomologies were
applied to analysis of overdetermined systems of linear PDE generalizing
classical works by Cartan \cite{Cartan1}. Homology operators and groups
introduced by Spencer (and called the \emph{Spencer operators} and
\emph{Spencer homology} nowadays) play a basic role in all computations
related to modern homological applications to PDE (see below).

Further achievements became possible in the framework of the
geometrical approach to PDE. Originating in classical works by Lie,
B\"acklund, Darboux, this approach was developed by A.~Vinogradov and
his co-workers (see \cite{KLV,Vin5}). Treating a differential equation
as a submanifold in a suitable jet bundle and using a nontrivial
geometrical structure of the latter allows one to apply powerful tools
of modern differential geometry to analysis of nonlinear PDE of a general
nature. And not only this: speaking the geometrical language makes it
possible to clarify underlying algebraic structures, the latter giving
better and deeper understanding of the whole picture, \cite[Ch.~1]{KLV}
and \cite{Vin6,Kras4}.

It was also A.~Vinogradov to whom the next homological application to
PDE belongs. In fact, it was even more than an application: in a series
of papers \cite{Vin3,Vin2,Vin1}, he has demonstrated that the adequate
language for Lagrangian formalism is a special spectral sequence (the
so-called \emph{Vinogradov $\CC$-spectral sequence}) and obtained first
spectacular results using this language. As it happened, the area of
the $\CC$-spectral sequence applications is much wider and extends to
scalar differential invariants of geometric structures
\cite{VerbVinGess}, modern field theory
\cite{BarnBrandtHenn1,BarnBrandtHenn2,Barnich1,Brandt1,Henneaux1}, etc.
A lot of work was also done to specify and generalize Vinogradov's
initial results, and here one could mention those by I.~M.~Anderson
\cite{Ander1,Ander2}, R.~L.~Bryant and P.~A.~Griffiths \cite{BryantGriff1},
D.~M.~Gessler \cite{Gessler1,Gessler2}, M.~Marvan \cite{Marvan4,Marvan1},
T.~Tsujishita \cite{Tsuji2,Tsuji3,Tsuji1}, W.~M.~Tulczyjew
\cite{Tulcz3,Tulcz1,Tulcz2}.

Later, one of the authors found out that another cohomology theory
($\CC$-cohomologies) is naturally related to any PDE \cite{Kras1}. The
construction uses the fact that the infinite prolongation of any
equation is naturally endowed with a flat connection (the \emph{Cartan
connection}). To such a connection, one puts into correspondence a
differential complex based on the \emph{\FN bracket}
\cite{Nijenh,FrolNijenh}. The group $H^0$ for this complex coincides
with the \emph{symmetry algebra} of the equation at hand, the group
$H^1$ consists of equivalence classes of \emph{deformations} of the
equation structure.  Deformations of a special type are identified with
\emph{recursion operators} \cite{Olver} for symmetries. On the other
hand, this theory seems to be dual to the term $E_1$ of the Vinogradov
$\CC$-spectral sequence, while special cochain maps relating the former
to the latter are Poisson structures on the equation \cite{Kras5}.

Not long ago, the second author noticed (\cite{Verb6}) that both
theories may be understood as \emph{horizontal cohomologies} with
suitable coefficients.  Using this observation combined with the fact
that the horizontal de Rham cohomology is equal to the cohomology of
the \emph{compatibility complex} for the universal linearization
operator, he found a simple proof of the vanishing theorem for the term
$E_1$ (the ``$k$-line theorem'') and gave a complete description of
$\CC$-cohomology in the ``2-line situation''.

Our short review will not be complete, if we do not mention applications of
cohomologies to the singularity theory of solutions of nonlinear PDE
(\cite{Lychagin1}), though this topics is far beyond the scope of these lecture
notes.

$$\star\ \star\ \star$$

The idea to expose the above mentioned material in a lecture course at the
Summer School in Levo\v{c}a belongs to Prof. D.~Krupka to whom we are
extremely grateful.

We tried to give here a complete and self-contained picture which was not
easy under natural time and volume limitations.
To make reading easier, we included the Appendix containing basic facts
and definitions from homological algebra.
In fact, the material needs
not 5 days, but 3--4 semester course at the university level, and we
really do hope that these lecture notes will help to those who became
interested during the lectures. For further details (in the geometry of PDE
especially) we refer the reader to the books \cite{KLV} and \cite{Symm} (an
English translation of the latter is to be published by the American
Mathematical Society in 1999). For advanced reading we also strongly
recommend the collection \cite{zzConfMoscow97}, where one will find a lot of
cohomological applications to modern physics.

\vspace{1cm}

\hfill J.~Krasil{\cprime}shchik

\hfill A.~Verbovetsky

\hfill Moscow, 1998

\newpage

\section{Differential calculus over commutative algebras}\label{sec:calc}
Throughout this section we shall deal with a commutative algebra $A$ over
a field $\Bbbk$ of zero characteristic. For further details we refer the
reader to \cite[Ch.~I]{KLV} and \cite{Kras4}.
\subsection{Linear differential operators}
Consider two $A$-modules $P$ and
$Q$ and the group $\Hom_\Bbbk(P,Q)$. Two $A$-module structures can be
introduced into this group:
\begin{equation}\label{sec1:eq:str}
(a\De)(p)=a\De(p),\quad (a^+\De)(p)=\De(ap),
\end{equation}
where $a\in A$, $p\in P$, $\De\in\Hom_\Bbbk(P,Q)$. We also set
$$\de_a(\De)=a^+\De-a\De,\quad \de_{a_0,\dots,a_k}=\de_{a_0}\circ\dots\circ
\de_{a_k},$$
$a_0,\dots,a_k\in A$. Obviously, $\de_{a,b}=\de_{b,a}$ and
$\de_{ab}=a^+\de_b+b\de_a$ for any $a,b\in A$.
\bde\label{sec1:df:do}
A $\Bbbk$-homomorphism $\De\colon P\to Q$ is called a \emph{linear
differential operator} of order $\le k$ over the algebra $A$, if
$\de_{a_0,\dots,a_k}(\De)=0$ for all $a_0,\dots,a_k\in A$.
\ede
\bpr\label{sec1:pr:equiv}
If $M$ is a smooth manifold, $\xi,\zeta$ are smooth locally trivial
vector bundles over $M$, $A=\Ci(M)$ and $P=\Ga(\xi),Q=\Ga(\zeta)$ are the
modules of smooth sections, then any linear differential operator acting from
$\xi$ to $\zeta$ is an operator in the sense of Definition \vref{sec1:df:do}
and vice versa.
\epr
\begin{xca}
Prove this fact.
\end{xca}
Obviously, the set of all differential operators of order $\le k$
acting from $P$ to $Q$ is a subgroup in $\Hom_\Bbbk(P,Q)$ closed with
respect to both multiplications \veqref{sec1:eq:str}. Thus we obtain two
modules denoted by $\Diff_k(P,Q)$ and $\Diff_k^+(P,Q)$ respectively. Since
$a(b^+\De)=b^+(a\De)$ for any $a,b\in A$ and $\De\in\Hom_\Bbbk(P,Q)$, this group
also carries the structure of an $A$-bimodule denoted by
$\Diff_k^{(+)}(P,Q)$. Evidently, $\Diff_0(P,Q)=\Diff_0^+(P,Q)=\Hom_A(P,Q)$.

It follows from Definition \vref{sec1:df:do} that any differential operator
of order $\le k$ is an operator of order $\le l$ for all $l
\ge k$ and consequently we obtain the embeddings $\Diff_k^{(+)}(P,Q)
\sbs\Diff_l^{(+)}(P,Q)$, which allow us to define the filtered bimodule
$\Diff^{(+)}(P,Q)=\bigcup_{k\ge0}\Diff_k^{(+)}(P,Q)$.

We can also consider the $\mathbb{Z}$-graded module associated to the filtered
module $\Diff^{(+)}(P,Q)$:
$\Smbl(P,Q)=\bigoplus_{k\ge0}\Smbl_k(P,Q)$, where $\Smbl_k(P,Q)=
\Diff_k^{(+)}(P,Q)/\Diff_{k-1}^{(+)}(P,Q)$, which is called the \emph{module
of symbols}. The elements of $\Smbl(P,Q)$ are called \emph{symbols} of
operators acting from $P$ to $Q$. It easily seen that two module structures
defined by \veqref{sec1:eq:str} become identical in $\Smbl(P,Q)$.

The following properties of linear differential operator are directly
implied by the definition:
\bpr
Let $P,Q$ and $R$ be $A$-modules. Then:
\begin{enumerate}
\item If $\De_1\in\Diff_k(P,Q)$ and $\De_2\in\Diff_l(Q,R)$ are two
differential operators, then their composition
$\De_2\circ\De_1$ lies in $\Diff_{k+l}(P,R)$.
\item The maps
$$\mathrm{i}^{\cdot,+}\colon\Diff_k(P,Q)\to\Diff_k^+(P,Q),\
\mathrm{i}^{+,\cdot}\colon\Diff_k^+(P,Q)\to\Diff_k(P,Q)$$
generated by the identical map of $\Hom_\Bbbk(P,Q)$ are
differential operators of order $\le k$.
\end{enumerate}
\epr
\begin{corollary}\label{sec1:co:iso}
There exists an isomorphism
$$\Diff^+(P,\Diff^+(Q,R))=\Diff^+(P,\Diff(Q,R))$$
generated by the operators $\mathrm{i}^{\cdot,+}$ and $\mathrm{i}^{+,\cdot}$.
\end{corollary}
Introduce the notation $\Diff_k^{(+)}(Q)=\Diff_k^{(+)}(A,Q)$ and define the
map $\cD_k\colon\Diff_k^+(Q)\to Q$ by setting $\cD_k(\De)=\De(1)$.
Obviously, $\cD_k$ is an operator of order $\le k$. Let also
\begin{equation}\label{sec1:eq:rep1}
\psi\colon\Diff_k^+(P,Q)\to\Hom_A(P,\Diff_k^+(Q)),\quad
\De\mapsto\psi_\De,
\end{equation}
be the map defined by $(\psi_\De(p))(a)=\De(ap)$, $p\in P$, $a\in A$.
\bpr\label{sec1:pr:dorep1}
The map \veqref{sec1:eq:rep1} is an isomorphism of $A$-modules.
\epr
\begin{proof}
Compatibility of $\psi$ with $A$-module structures is obvious. To complete
the proof it suffices to note that the correspondence
$$
\Hom_A(P,\Diff_k^+(Q))\ni\vf\mapsto\cD_k\circ\vf\in\Diff_k^+(P,Q)
$$
is inverse to $\psi$.
\end{proof}
The homomorphism $\psi_\De$ is called $\Diff$-\emph{associated} to $\De$.
\begin{remark}
Consider the correspondence $P\Ra\Diff_k^+(P,Q)$ and for any $A$-homomorphism
$f\colon P\to R$ define the homomorphism
\[
\Diff_k^+(f,Q)\colon \Diff_k^+(R,Q)\to\Diff_k^+(P,Q)
\]
by setting
$\Diff_k^+(f,Q)(\De)=\De\circ f$. Thus, $\Diff_k^+(\cdot,Q)$ is a
contravariant functor from the category of all $A$-modules to itself.
Proposition \vref{sec1:pr:dorep1} means that this functor is
representable and the module $\Diff_k^+(Q)$ is its representative
object. Obviously, the same is valid for the functor $\Diff^+(\cdot,Q)$
and the module $\Diff^+(Q)$.
\end{remark}
From Proposition \vref{sec1:pr:dorep1} we also obtain the following
\begin{corollary}
There exists a unique homomorphism
$$c_{k,l}=c_{k,l}(P)\colon\Diff_k^+(\Diff_l^+(P))\to\Diff_{k+l}(P)$$
such that the diagram
$$\label{sec1:eq:cD}\begin{CD}
\Diff_k^+(\Diff_l^+(P))@>\cD_k>>\Diff_l^+(P)\\
@Vc_{k,l}VV                     @VV\cD_lV\\
\Diff_{k+l}^+(P)@>\cD_{k+l}>>P
\end{CD}$$
is commutative.
\end{corollary}
\begin{proof}
It suffices to use the fact that the composition
$$\cD_l\circ\cD_k:\Diff_k(\Diff_l(P))\xra{}P$$
is an operator of order $\le k+l$ and to set $c_{k,l}=
\psi_{\cD_l\circ\cD_k}$.
\end{proof}

The map $c_{k,l}$ is called the \emph{gluing homomorphism} and from the
definition it follows that $(c_{k,l}(\De))(a)=(\De(a))(1)$, $\De\in
\Diff_k^+(\Diff_l^+(P))$, $a\in A$.
\begin{remark}\label{sec1:re:nat}
The correspondence $P\Ra\Diff_k^+(P)$ also becomes a (covariant) functor,
if for a homomorphism $f\colon P\to Q$ we define the homomorphism
$\Diff_k^+(f)\colon\Diff_k^+(P)\to \Diff_k^+(Q)$ by
$\Diff_k^+(f)(\De)=f\circ\De$. Then the correspondence $P\Ra
c_{k,l}(P)$ is a natural transformation of functors
$\Diff_k^+(\Diff_l^+(\cdot))$ and $\Diff_{k+l}^+(\cdot)$ which means that
for any $A$-homomorphism $f\colon P\to Q$ the diagram
$$\begin{CD}
\Diff_k^+(\Diff_l^+(P))@>\Diff_k^+(\Diff_l^+(f))>>\Diff_k^+(\Diff_l^+(Q))\\
@Vc_{k,l}(P)VV                     @VVc_{k,l}(Q)V\\
\Diff_{k+l}^+(P)@>\Diff_{k+l}^+(f)>>\Diff_{k+l}^+(Q)
\end{CD}$$
is commutative.

Note also that the maps $c_{k,l}$ are compatible with the natural embeddings
$\Diff_k^+(P)\to\Diff_s^+(P),\,k\le s$, and thus we can define the
gluing $c_{*,*}\colon \Diff^+(\Diff^+(\cdot))\to\Diff^+(\cdot)$.
\end{remark}

\subsection{Multiderivations and the Dif{}f-Spencer complex}
Let $A^{\ot k}=A\ot_\Bbbk\dotsm\ot_\Bbbk A$, $k$ times.
\begin{definition}\label{sec1:df:poly}
A $\Bbbk$-linear map $\na\colon A^{\ot k}\to P$ is called a
\emph{skew-symmetric multiderivation} of $A$ with values in an
$A$-module $P$, if the following conditions hold:
\begin{align*}
\mathrm{(1)}\ &\na(a_1,\dots,a_i,a_{i+1},\dots,a_k)+
\na(a_1,\dots,a_{i+1},a_i,\dots,a_k)=0,
\\
\mathrm{(2)}\ &\na(a_1,\dots,a_{i-1},ab,a_{i+1},\dots,a_k)=\\
&a\na(a_1,\dots,a_{i-1},b,a_{i+1},\dots,a_k)+
b\na(a_1,\dots,a_{i-1},a,a_{i+1},\dots,a_k)
\end{align*}
for all $a,b,a_1,\dots,a_k\in A$ and any $i$, $1\le i\le k$.
\end{definition}
The set of all skew-symmetric $k$-derivations forms an $A$-module denoted
by $\Dr_k(P)$. By definition, $\Dr_0(P)=P$. In particular, elements of
$\Dr_1(P)$ are called
\emph{$P$-valued derivations} and form a submodule in $\Diff_1(P)$ (but not
in the module $\Diff_1^+(P)$!).

There is another, functorial definition of the modules $\Dr_k(P)$: for any
$\na\in\Dr_k(P)$ and $a\in A$ we set $(a\na)(a_1,\dots,a_k)=
a\na(a_1,\dots,a_k)$. Note
first that the composition $\ga_1\colon \Dr_1(P)\hra\Diff_1(P)
\xra{\mathrm{i}^{\cdot,+}}\Diff_1^+(P)$ is a monomorphic differential
operator of order $\le 1$. Assume now that the first-order monomorphic
operators $\ga_i=\ga_i(P)\colon\Dr_i(P)\to\Dr_{i-1}(\Diff_1^+(P))$ were
defined for all $i\le k$. Assume also that all the maps $\ga_i$ are
natural\footnote{This means that for any $A$-homomorphism $f\colon P\to
Q$ one has $\ga_i(Q)\circ
\Dr_i(f)=\Dr_{i-1}(\Diff_1^+(f))\circ\ga_i(P)$.} operators. Consider
the composition
\begin{equation}
\label{sec1:eq:ker}
\Dr_k(\Diff_1^+(P))\xra{\ga_k}\Dr_{k-1}(\Diff_1^+(\Diff_1^+(P)))
\xra{\Dr_{k-1}(c_{1,1})}\Dr_{k-1}(\Diff_2^+(P)).
\end{equation}
\bpr\label{sec1:pr:polyf}
The following facts are valid:
\begin{enumerate}
\item $\Dr_{k+1}(P)$ coincides with the kernel of the composition
\eqref{sec1:eq:ker}.
\item The embedding $\ga_{k+1}\colon \Dr_{k+1}(P)\hra\Dr_k(\Diff_1^+(P))$ is
a first-order differential operator.
\item The operator $\ga_{k+1}$ is natural.
\end{enumerate}
\epr
The proof reduces to checking the definitions.

\bre\label{sec1:re:dot}
We saw above that the $A$-module $\Dr_{k+1}(P)$ is the kernel of the
map $\Dr_{k-1}(c_{1,1})\circ\ga_k$, the latter being not an $A$-module
homomorphism but a differential operator. Such an effect arises in the
following general situation.

Let $\mathrm{F}$ be a functor acting on a subcategory of the category
of $A$-modules. We say that $\mathrm{F}$ is \emph{$\Bbbk$-linear}, if
the corresponding map
$\mathrm{F}_{P,Q}\colon\Hom_\Bbbk(P,Q)\to\Hom_\Bbbk(P,Q)$ is linear
over $\Bbbk$ for all $P$ and $Q$ from our subcategory. Then we can
introduce a new $A$-module structure in the the $\Bbbk$-module
$\mathrm{F}(P)$ by setting $a\ldot q= (\mathrm{F}(a))(q)$, where
$q\in\mathrm{F}(P)$ and
$\mathrm{F}(a)\colon\mathrm{F}(P)\to\mathrm{F}(P)$ is the homomorphism
corresponding to the multiplication by $a$: $p\mapsto ap$, $p\in P$.
Denote the module arising in such a way by $\mathrm{F}\ldot(P)$.

Consider two $\Bbbk$-linear functors $\mathrm{F}$ and $\mathrm{G}$ and
a natural transformation $\De$:
$P\Ra\De(P)\in\Hom_\Bbbk(\mathrm{F}(P),\mathrm{G}(P))$.
\begin{xca}
Prove that the natural transformation
$\De$ induces a natural homomorphism of $A$-modules $\De\ldot\colon\mathrm{F}
\ldot(P)\to
\mathrm{G}\ldot(P)$ and thus its kernel is always an $A$-module.
\end{xca}
\ere

From Definition \vref{sec1:df:poly} it also follows that elements of the modules
$\Dr_k(P),\,k\ge 2$, may be understood as derivations $\De\colon A\to
\Dr_{k-1}(P)$ satisfying $(\De(a))(b)=-(\De(b))(a)$.
We call $\De(a)$ the \emph{evaluation}\label{sec1:p:eval} of the
multiderivation $\De$ at the element $a\in A$.
Using this interpretation, define by induction on $k+l$ the operation
$\wg\colon \Dr_k(A)\ot_A\Dr_l(P)\to\Dr_{k+l}(P)$ by setting
$$
a\wg p=ap,\ a\in\Dr_0(A)=A,p\in\Dr_0(P)=P,
$$
and
\begin{equation}\label{sec1:eq:wedge}
(\De\wg\na)(a)=\De\wg\na(a)+(-1)^l\De(a)\wg\na.
\end{equation}
Using elementary induction on $k+l$, one can easily prove the following
\bpr
The operation $\wg$ is well defined and satisfies the following properties:
\begin{align*}
\mathrm{(1)}\ &\De\wg(\De'\wg\na)=(\De\wg\De')\wg\na,\\
\mathrm{(2)}\ &(a\De+a'\De')\wg\na=a\De\wg\na+a'\De'\wg\na,\\
\mathrm{(3)}\ &\De\wg(a\na+a'\na')=a\De\wg\na+a'\De\wg\na',\\
\mathrm{(4)}\ &\De\wg\De'=(-1)^{kk'}\De'\wg\De
\end{align*}
for any elements $a,a'\in A$ and multiderivations
$\De\in\Dr_k(A)$, $\De'\in\Dr_{k'}(A)$, $\na\in\Dr_l(P)$, $\na'\in
\Dr_{l'}(P)$.
\epr
Thus, $\Dr_*(A)=\bigoplus_{k\ge0}\Dr_k(A)$ becomes a $\mathbb{Z}$-graded
commutative algebra and $\Dr_*(P)=\bigoplus_{k\ge0}\Dr_k(P)$ is a
graded $\Dr_*(A)$-module. The correspondence $P\Ra\Dr_*(P)$ is a functor
from the category of $A$-modules to the category of graded
$\Dr_*(A)$-modules.

Let now $\na\in\Dr_k(\Diff_l^+(P))$ be a multiderivation. Define
\begin{equation}\label{sec1:eq:spen}
(S(\na)(a_1,\dots,a_{k-1}))(a)=(\na(a_1,\dots,a_{k-1},a)(1)),
\end{equation}
$a,a_1,\dots,a_{k-1}\in A$. Thus we obtain the map
$$
S\colon \Dr_k(\Diff_l^+(P))\to\Dr_{k-1}(\Diff_{l+1}^+(P))
$$
which can be represented as the composition
\begin{equation}\label{sec1:eq:spen1}
\Dr_k(\Diff_l^+(P))\xra{\ga_k}\Dr_{k-1}(\Diff_1^+(\Diff_l^+(P)))
\xra{\Dr_{k-1}(c_{1,l})}\Dr_{k-1}(\Diff_{l+1}^+(P)).
\end{equation}
\bpr
The maps $S\colon \Dr_k(\Diff_l^+(P))\to\Dr_{k-1}(\Diff_{l+1}^+(P))$ possess the
following properties:
\begin{enumerate}
\item $S$ is a differential operator of order $\le 1$.
\item $S\circ S=0$.
\end{enumerate}
\epr
\begin{proof}
The first statement follows from \veqref{sec1:eq:spen1}, the second one is
implied by  \veqref{sec1:eq:spen}.
\end{proof}
\bde
The operator $S$ is called the \emph{$\Diff$-Spencer operator}. The sequence
of operators
$$
0\xla{}P\xla{\cD}\Diff^+(P)\xla{S}\Diff^+(P)\xla{S}
\Dr_2(\Diff^+(P))\xla{}\dotsb
$$
is called the \emph{$\Diff$-Spencer complex}.
\ede
\subsection{Jets}
Now we shall deal with the functors $Q\Ra\Diff_k(P,Q)$ and their
representability.

Consider an $A$-module $P$ and the tensor product $A\ot_\Bbbk P$. Introduce
an $A$-module structure in this tensor product by setting
$$a(b\ot p)=(ab)\ot p,\ a,b\in A,\ p\in P,$$
and consider the $\Bbbk$-linear map $\epsilon\colon P\to
A\ot_\Bbbk P$ defined by $\epsilon(p)=1\ot p$. Denote by $\mu^k$ the
submodule in $A\ot_\Bbbk P$ generated by the elements of the form
$(\de_{a_0,\dots,a_k}(\epsilon))(p)$ for all $a_0,\dots,a_k\in A$
and $p\in P$.
\bde
The quotient module $(A\ot_\Bbbk P)/\mu^k$ is called the \emph{module of
$k$-jets} for $P$ and is denoted by $\CJ^k(P)$.
\ede
We also define the map $j_k\colon P\to\CJ^k(P)$ by setting
$j_k(p)=\epsilon(p)\bmod\mu^k$.
Directly from the definition of $\mu^k$ it follows that $j_k$ is a
differential operator of order $\le k$.
\bpr\label{sec1:pr:jetrep}
There exists a canonical isomorphism
\begin{equation}
\psi\colon \Diff_k(P,Q)\to\Hom_A(\CJ^k(P),Q),\quad\De\mapsto\psi^\De,
\end{equation}
defined by the equality $\De=\psi^\De\circ j_k$ and called
\emph{Jet}-associated to $\De$.
\epr
\begin{proof}
Note first that since the module $\CJ^k(P)$ is generated by the elements of
the form $j_k(p),\ p\in P$, the homomorphism $\psi^\De$, if defined, is
unique. To establish existence of $\psi^\De$, consider the homomorphism
$$
\eta\colon \Hom_A(A\ot_\Bbbk P,Q)\to\Hom_\Bbbk(P,Q),\quad
\eta(\vf)=\vf\circ\epsilon.
$$
Since $\vf$ is an $A$-homomorphism, one has
$$\de_a(\eta(\vf))=\de_a(\vf\circ\epsilon)=\vf\circ\de_a(\epsilon)
=\eta(\de_a(\vf)),\quad a\in A.$$
Consequently, the element $\eta(\vf)$ is an operator of order $\le k$
if and only if $\vf(\mu^k)=0$, i.e., restricting $\eta$ to $\Diff_k(P,Q)
\sbs\Hom_\Bbbk(P,Q)$ we obtain the desired isomorphism $\psi$.
\end{proof}
The proposition proved means that the functor $Q\Ra\Diff_k(P,Q)$ is
representable and the module $\CJ^k(P)$ is its representative object.

Note that the correspondence $P\Ra\CJ^k(P)$ is a functor itself: if $\vf\colon P\to
Q$ is an $A$-module homomorphism, we are able to define the homomorphism
$\CJ^k(\vf)\colon
\CJ^k(P)\to\CJ^k(Q)$ by the commutativity condition
$$\begin{CD}
P@>{j_k}>>\CJ^k(P)\\
@V{\vf}VV @VV\CJ^k(\vf)V\\
Q@>{j_k}>>\CJ^k(Q)
\end{CD}
$$

The universal property of the operator $j_k$ allows us to introduce the
natural transformation $c^{k,l}$ of the functors $\CJ^{k+l}(\cdot)$ and
$\CJ^k(\CJ^l(\cdot))$ defined by the commutative diagram
$$\begin{CD}
P@>j_l>>\CJ^l(P)\\
@Vj_{k+l}VV @VVj_kV\\
\CJ^{k+l}(P)@>c^{k,l}>>\CJ^k(\CJ^l(P))
\end{CD}$$
It is called the \emph{co-gluing homomorphism} and is dual to the gluing one
discussed in Remark \vref{sec1:re:nat}.

Another natural transformation related to functors $\CJ^k(\cdot)$ arises
from the embeddings $\mu^l\hra\mu^k,\, l\ge k$, which generate the
projections $\nu_{l,k}\colon \CJ^l(P)\to\CJ^k(P)$ dual to the embeddings
$\Diff_k(P,Q)\hra\Diff_l(P,Q)$. One can easily see that if $f\colon P\to P'$ is
an $A$-module homomorphism, then $\CJ^k(f)\circ\nu_{l,k}=\nu_{l,k}\circ
\CJ^l(f)$. Thus we obtain the sequence of projections
$$
\dotsb\xra{}\CJ^k(P)\xra{\nu_{k,k-1}}\CJ^{k-1}(P)\xra{}\dotsb\xra{}\CJ^1(P)
\xra{\nu_{1,0}}\CJ^0(P)=P
$$
and set $\CJ^\infty(P)=\projlim\CJ^k(P)$. Since $\nu_{l,k}\circ
j_l=j_k$, we can also set $j_\infty=\projlim j_k\colon P\to
\CJ^\infty(P)$.

Let $\De\colon P\to Q$ be an operator of order $\le k$. Then for any $l
\ge 0$ we have the commutative diagram
$$\begin{CD}
P@>\De>>Q\\
@Vj_{k+l}VV@VVj_lV\\
\CJ^{k+l}(P)@>\psi_l^\De>>\CJ^l(Q)
\end{CD}$$
where $\psi_l^\De=\psi^{j_l\circ\De}$. Moreover, if $l'\ge l$, then
$\nu_{l',l}\circ\psi_{l'}^\De=\psi_{l}^\De\circ\nu_{k+l',k+l}$ and we
obtain the homomorphism $\psi_\infty^\De\colon \CJ^\infty(P)\to\CJ^\infty(Q)$.

Note that the co-gluing homomorphism is a particular case of the above
construction: $c^{k,l}=\psi_k^{j_l}$. Thus, passing to the inverse limits,
we obtain the co-gluing $c^{\infty,\infty}$:
$$\begin{CD}
P@>{j_\infty}>>\CJ^\infty(P)\\
@V{j_\infty}VV @VV{j_\infty}V\\
\CJ^\infty(P)@>{c^{\infty,\infty}}>>\CJ^\infty(\CJ^\infty(P))
\end{CD}$$

\subsection{Compatibility complex}\label{sub:comp}
The following construction will play an important role below.

Consider a differential operator $\De\colon Q\to Q_1$ of order $\le k$.
Without loss of generality we may assume that its Jet-associated
homomorphism $\psi^\De\colon \CJ^k(Q)\to Q_1$ is epimorphic. Choose an integer
$k_1\ge 0$ and define $Q_2$ as the cokernel of the homomorphism
$\psi_{k_1}^\De\colon \CJ^{k+k_1}(Q)\to\CJ^k(Q_1)$,
$$
0\to\CJ^{k+k_1}(Q)\xra{\psi_{k_1}^\De}\CJ^{k_1}(Q_1)\to Q_2\to 0.
$$
Denote the composition of the operator $j_{k_1}\colon Q_1\to\CJ^{k_1}(Q_1)$ with
the natural projection $\CJ^{k_1}(Q_1)\to Q_2$ by $\De_1\colon Q_1\to Q_2$.
By construction, we have
$$
\De_1\circ\De=\psi^{\De_1}\circ j_{k_1}\circ\De=\psi^{\De_1}\circ
\psi_{k_1}^{\De}\circ j_{k+k_1}.
$$
\begin{xca}
\label{sec1:ex:comop}
Prove that $\De_1$ is a compatibility operator for the operator $\De$,
i.e., for any operator $\nabla$ such that $\nabla\circ\De=0$ and
$\ord\nabla\ge k_1$, there exists an operator $\square$ such that
$\nabla=\square\circ\Delta_1$.
\end{xca}

We can now apply the procedure to the operator $\De_1$ and some integer
$k_2$ obtaining $\De_2\colon Q_2\to Q_3$, etc. Eventually, we obtain the
complex
$$
0\xra{} Q\xra{\De}Q_1\xra{\De_1}Q_2\xra{\De_2}\dotsb\xra{}
Q_i\xra{\De_i}Q_{i+1}\xra{}\dotsb
$$
which is called the \emph{compatibility complex} of the operator $\De$.

\subsection{Differential forms and the de Rham complex}
Consider the embedding $\be\colon A\to\CJ^1(A)$ defined by
$\be(a)=aj_1(1)$ and define the module $\La^1=\CJ^1(A)/\im\be$. Let $d$
be the composition of $j_1$ and the natural projection
$\CJ^1(A)\to\La^1$. Then $d\colon A\to\La^1$ is a differential operator
of order $\le 1$ (and, moreover, lies in $\Dr_1(\La^1)$).

Let us now apply the construction of the previous subsection to the
operator $d$ setting all $k_i$ equal to $1$ and preserving the notation $d$
for the operators $d_i$. Then we get the compatibility complex
$$
0\xra{}A\xra{d}\La^1\xra{d}\La^2\xra{}\dotsb\xra{}\La^k\xra{d}\La^{k+1}
\xra{}\dotsb
$$
which is called the \emph{de Rham complex} of the algebra $A$. The elements
of $\La^k$ are called \emph{$k$-forms} over $A$.
\bpr\label{sec1:pr:forms}
For any $k\ge 0$, the module $\La^k$ is the representative object
for the functor $\Dr_k(\cdot)$.
\epr
\begin{proof}
It suffices to compare the definition of $\La^k$ with the description
of $\Dr_k(P)$ given by Proposition \vref{sec1:pr:polyf}.
\end{proof}
\bre
In the case $k=1$, the isomorphism between $\Hom_A(\La^1,\cdot)$ and
$\Dr_1(\cdot)$ can be described more exactly. Namely, from the definition
of the operator $d\colon A\to\La^1$ and from Proposition
\vref{sec1:pr:jetrep} it follows that
any derivation $\na\colon A\to P$ is uniquely represented as the composition
$\na=\vf^\na\circ d$ for some homomorphism $\vf^\na\colon \La^1\to P$.
\ere

As a consequence Proposition \vref{sec1:pr:forms}, we obtain the following
\bco
The module $\La^k$ is the $k$-th exterior power of $\La^1$.
\eco
\begin{xca}\label{sec1:es:dual}
Since $\Dr_k(P)=\Hom_A(\La^k,P)$, one can introduce the pairing $\langle
\cdot,\cdot\rangle\colon\Dr_k(P)\ot\La^k\xra{}P$. Prove that the evaluation
operation (see p.~\pageref{sec1:p:eval}) and the wedge product are mutually
dual with respect to this pairing, i.e.,
$$\langle X,da\wg\om\rangle=\langle X(a),\om\rangle$$
for all $X\in\Dr_{k+1}(P)$, $\om\in\La^k$, and $a\in A$.
\end{xca}
The following proposition establishes the relation of the de Rham
differential to the wedge product.
\bpr[the Leibniz rule]
For any $\om\in\La^k$ and $\ta\in\La^l$ one has
$$
d(\om\wg\ta)=d\om\wg\ta+(-1)^k\om\wg d\ta.
$$
\epr
\begin{proof}
We first consider the case $l=0$, i.e., $\ta=a\in A$. To do it, note that
the wedge product $\wg\colon \La^k\ot_A\La^l\to\La^{k+l}$,
due to Proposition \vref{sec1:pr:forms}, induces the natural embeddings
of modules
$\Dr_{k+l}(P)\to\Dr_k(\Dr_l(P))$. In particular, the embedding $\Dr_{k+1}(P)
\to\Dr_k(\Dr_1(P))$ can be represented as the composition
$$
\Dr_{k+l}(P)\xra{\ga_{k+1}}\Dr_k(\Diff_1^+(P))\xra{\la}\Dr_k(\Dr_1(P)),
$$
where $(\la(\na))(a_1,\dots,a_k)=\na(a_1,\dots,a_k)-(\na(a_1,\dots,a_k))(1)$.
In a dual way, the wedge product is represented as
$$
\La^k\ot_A\La^1\xra{\la'}\CJ^1(\La^k)\xra{\psi^d}\La^{k+1},
$$
where $\la'(\om\ot da)=(-1)^k(j_1(\om a)-j_1(\om)a)$. Then
\begin{multline*}
(-1)^k\wg\om da=(-1)^k\psi^d(\la'(\om\ot da))\\
=\psi^d(j_1(\om a)-j_1(\om)a)=d(\om a)-d(\om)a.
\end{multline*}
The general case is implied by the identity
$$d(\om\wg da)=(-1)^kd(d(\om a)-d\om\cdot a)=(-1)^{k+1}d(d\om\cdot a).$$
\end{proof}

Let us return back to Proposition \vref{sec1:pr:forms} and consider
the $A$-bilinear pairing
$$
\langle\cdot,\cdot\rangle\colon \Dr_k(P)\ot_A\La^k\to P
$$
again. Take a form $\om\in\La^k$ and a derivation $X\in\Dr_1(A)$. Using the
definition of the wedge product in $\Dr_*(P)$ (see equality
\veqref{sec1:eq:wedge}), we can set
\begin{equation}\label{sec1:eq:inter}
\langle\De,i_X\om\rangle=(-1)^{k-1}\langle X\wg\De,\om\rangle
\end{equation}
for an arbitrary $\De\in\Dr_{k-1}(P)$.
\bde
The operation $\ip_X\colon \La^k\to\La^{k-1}$ defined by \veqref{sec1:eq:inter}
is called the \emph{internal product}, or \emph{contraction}.
\ede
\bpr
For any $X,Y\in\Dr_1(A)$ and $\om\in\La^k,\ta\in\La^l$ one has
\begin{align*}
\mathrm{(1)}\ &\ip_X(\om\wg\ta)=\ip_X(\om)\wg\ta+(-1)^k\om\wg \ip_X(\ta),\\
\mathrm{(2)}\ &\ip_X\circ \ip_Y=-\ip_Y\circ \ip_X
\end{align*}
\epr
In other words, internal product is a derivation of the $\mathbb{Z}$-graded
algebra $\La^*=\bigoplus_{k\ge0}\La^k$ of degree $-1$ and $\ip_X,\ip_Y$
commute as graded maps.

Consider a derivation $X\in\Dr_1(A)$ and set
\begin{equation}\label{sec1:eq:lie}
\Ld_X(\om)=[\ip_X,d](\om)=\ip_X(d(\om))+d(\ip_X(\om)),\ \om\in\La^*.
\end{equation}
\bde\label{sec1:df:lie}
The operation $\Ld_X\colon \La^*\to\La^*$ defined by \ref{sec1:eq:lie}
is called the \emph{Lie derivative}.
\ede
Directly from the definition one obtains the following properties of Lie
derivatives:
\bpr\label{sec1:pr:lie}
Let $X,Y\in\Dr_1(A)$, $\om,\ta\in\La^*$, $a\in A$, $\al,\be\in\Bbbk$. Then
the following identities are valid:
\begin{align*}
\mathrm{(1)}\ &\Ld_{\al X+\be Y}=\al \Ld_X+\be \Ld_Y,\\
\mathrm{(2)}\ &\Ld_{aX}=a\Ld_X+da\wg \ip_X,\\
\mathrm{(3)}\ &\Ld_X(\om\wg\ta)=\Ld_X(\om)\wg\ta+\om\wg \Ld_X(\ta),\\
\mathrm{(4)}\ &[d,\Ld_X]=d\circ \Ld_X-\Ld_X\circ d=0,\\
\mathrm{(5)}\ &\Ld_{[X,Y]}=[\Ld_X,\Ld_Y],\mbox{ where } [X,Y]=X\circ Y-Y\circ X,\\
\mathrm{(6)}\ &\ip_{[X,Y]}=[\Ld_X,\ip_Y]=[\ip_X,\Ld_Y].
\end{align*}
\epr

To conclude this subsection, we present another description of the
$\Diff$-Spencer complex. Recall Remark \vref{sec1:re:dot} and introduce
the ``dotted'' structure into the modules $\Dr_k(\Diff_l^+(P))$ and note
that $\Diff_l^+(P)\ldot=\Diff_l(P)$. Define
the isomorphism
$$\zeta\colon (\Dr_k(\Diff^+))\ldot(P)
=\Hom_A(\La^k,\Diff^+)\ldot
=\Diff^+(\La^k,P)\ldot=\Diff(\La^k,P).$$
Then we have
\bpr\label{sec1:pr:spenew}
The above defined map $\zeta$ generates the isomorphism of complexes
$$\begin{CD}
\dotsb@<<<(\Dr_{k-1}(\Diff^+))\ldot(P)@<S\ldot<<(\Dr_k(\Diff^+))\ldot(P)
@<<<\dotsb\\
@.@V{\zeta}VV@V{\zeta}VV@.\\
\dotsb@<<<\Diff(\La^{k-1},P)@<{v}<<\Diff(\La^k,P)@<<<\dotsb
\end{CD}$$
where $S\ldot$ is the operator induced on ``dotted'' modules by the
$\Diff$-Spencer operator, while $v(\na)=\na\circ d$.
\epr

\subsection{Left and right differential modules}\label{sub:dmod}
From now on till the end of this section we shall assume the modules under
consideration to be projective.
\bde\label{sec1:df:lmod}
An $A$-module $P$ is called a \emph{left differential module}, if there
exists an $A$-module homomorphism $\la\colon P\to\CJ^\infty(P)$ satisfying
$\nu_{\infty,0}\circ\la=\id_P$ and such that the diagram
$$\begin{CD}
P@>\la>>\CJ^\infty(P)\\
@V{\la}VV @VV{\CJ^\infty(\la)}V
\\
\CJ^\infty(P)@>c^{\infty,\infty}>>\CJ^\infty(\CJ^\infty(P))
\end{CD}$$
is commutative.
\ede
\ble\label{sec1:le:lmod}
Let $P$ be a left differential module. Then for any differential operator
$\De\colon Q_1\to Q_2$ there exists an operator $\De_P\colon
Q_1\ot_AP\to Q_2\ot_AP$ satisfying $(\id_Q)_P=\id_{Q\ot_AP}$ for
$Q=Q_1=Q_2$ and
$$(\De_2\circ\De_1)_P =(\De_2)_P\circ(\De_1)_P$$
for any operators $\De_1\colon Q_1\to Q_2$, $\De_2\colon Q_2\to Q_3$.
\ele
\begin{proof}
Consider the map
$$
\overline{\De}\colon Q_1\ot_A(A\ot_\Bbbk P)\to Q_2\ot_AP,\quad
q\ot a\ot p\mapsto\De(aq)\ot p.
$$
Since
$$
\overline{\De}(q\ot\de_a(\epsilon)(p))=\overline{\de_a\De}(q\ot 1\ot p),
\quad
p\in P,\quad q\in Q_1,\quad a\in A,
$$
the map
$$
\xi_P(\De)\colon Q_1\ot_A\CJ^\infty(P)\to Q_2\ot_A P
$$
is well defined. Set now the operator $\De_P$ to be the composition
$$
Q_1\ot_A P\xra{\id\ot\la}Q_1\ot_A\CJ^\infty(P)\xra{\xi_P(\De)}Q_2\ot_AP,
$$
which is a $k$-th order differential operator in an obvious way. Evidently,
$(\id_Q)_P=\id_{Q\ot_AP}$.

Now,
\begin{align*}
(\De_2\circ\De_1)_P&=\,\xi_P(\De_2\circ\De_1)\circ(\id\ot\la)\\
&=\,\xi_P(\De_2)\circ\xi_{\CJ^\infty(P)}(\De_1)\circ(\id\ot c^{\infty,\infty})
\circ(\id\ot\la)\\
&=\,\xi_P(\De_2)\circ\xi_{\CJ^\infty(P)}(\De_1)\circ(\id\ot\CJ^\infty(\la))\circ
(\id\circ\la)\\
&=\,\xi_P(\De_2)\circ(\id\ot\la)\circ\xi_P(\De_1)\circ(\id\ot\la)=
(\De_2)_P\circ(\De_1)_P,
\end{align*}
which proves the second statement.
\end{proof}

Note that the lemma proved shows in particular that any left differential
module is a left module over the algebra $\Diff(A)$ which justifies our
terminology.

Due to the above result, any complex of differential operators
$\dotsb\xra{} Q_i\xra{} Q_{i+1}\xra{}\dotsb$ and a left differential
module $P$ generate the complex $\dotsb\xra{} Q_i\ot_AP\xra{}
Q_{i+1}\ot_AP\xra{}\dotsb$ ``with coefficients'' in $P$. In particular,
since the co-gluing $c^{\infty,\infty}$ is in an obvious way
co-associative, i.e., the diagram
$$\begin{CD}
\CJ^\infty(P)@>{c^{\infty,\infty}(P)}>>\CJ^\infty(\CJ^\infty(P))\\
@V{c^{\infty,\infty}(P)}VV@VV\CJ^\infty({c^{\infty,\infty}(P))}V\\
\CJ^\infty(\CJ^\infty(P))@>{c^{\infty,\infty}(\CJ^\infty(P)})>>\CJ^\infty(\CJ^\infty(\CJ^\infty(P)))
\end{CD}$$
is commutative, $\CJ^\infty(P)$ is a left differential module with $\la=
c^{\infty,\infty}$. Consequently, we can consider the de Rham complex with
coefficients in $\CJ^\infty(P)$:
\begin{multline*}
0\xra{}
P\xra{j_\infty}\CJ^\infty(P)\xra{}\La^1\ot_A\CJ^\infty(P)\xra{}\dotsb\\
\dotsb\xra{}\La^i\ot_A\CJ^\infty(P)\xra{}\La^{i+1}\ot_A\CJ^\infty(P)
\xra{}\dotsb
\end{multline*}
which is the inverse limit for the Jet-\emph{Spencer complexes} of $P$
\begin{multline*}
0\xra{}
P\xra{j_k}\CJ^k(P)\xra{S}\La^1\ot_A\CJ^{k-1}(P)\xra{S}\dotsb\\
\dotsb\xra{S}\La^i\ot_A\CJ^{k-i}(P)\xra{S}\La^{i+1}\ot_A\CJ^{k-i-1}(P)
\xra{}\dotsb,
\end{multline*}
where $S(\omega\otimes j_{k-i}(p))=d\omega\otimes j_{k-i-1}(p)$.

Let $\De\colon P\to Q$ be a differential operator and $\psi_\infty^\De\colon
\CJ^\infty(P)\to\CJ^\infty(Q)$ be the corresponding homomorphism. The
kernel $E_\De=\ker\psi_\infty^\De$ inherits the left differential module
structure of $\CJ^\infty(P)$ and we can consider the de Rham complex with
coefficients in $E_\De$:
\begin{equation}\label{sec1:eq:jspcd}
0\xra{}E_\De\xra{}\La^1\ot_AE_\De\xra{}\dotsb\xra{}
\La^i\ot_AE_\De\xra{}\La^{i+1}\ot_AE_\De\xra{}\dotsb
\end{equation}
which is called the Jet-\emph{Spencer complex} of the operator $\De$.

Now we shall introduce the concept dual to that of left differential
modules.
\bde
An $A$-module $P$ is called a \emph{right differential module}, if there
exists an $A$-module homomorphism $\rho\colon \Diff^+(P)\to P$ that
satisfies the equality
$\rho\left|_{\Diff_0^+(P)}\right.=\id_P$ and makes the diagram
$$\begin{CD}
\Diff^+(\Diff^+(P))@>{c_{\infty,\infty}}>>\Diff^+(P)\\
@V{\Diff^+(\rho)}VV@VV{\rho}V\\
\Diff^+(P)@>{\rho}>>P
\end{CD}$$
commutative.
\ede
\ble\label{sec1:le:rmod}
Let $P$ be a right differential module. Then for any differential operator
$\De\colon Q_1\to Q_2$ of order $\le k$ there exists an operator
$$\De^P\colon\Hom_A(Q_2,P)\to\Hom_A(Q_1,P)$$
of order $\le k$ satisfying $\id_Q^P
=\id_{\Hom_A(Q,P)}$ for $Q=Q_1=Q_2$ and
$$(\De_2\circ\De_1)^P=\De_1^P\circ\De_2^P$$
for any operators $\De_1\colon Q_1\to Q_2$, $\De_2\colon Q_2\to Q_3$.
\ele
\begin{proof}
Let us define the action of $\De^P$ by setting $\De^P(f)=
\rho\circ\psi_{f\circ\De}$, where $f\in\Hom_A(Q_2,P)$. Obviously,
this is a $k$-th order differential operator
and $\id_Q^P=\id_{\Hom_A(Q,P)}$. Now,
\begin{align*}
(\De_2\circ\De_1)^P&=\rho\circ\psi_{f\circ\De_2\circ\De_1}=
\rho\circ c_{\infty,\infty}\circ\Diff^+(\psi_{f\circ\De_2})\circ\psi_{\De_1}
\\
&=\rho\circ\Diff^+(\rho\circ\psi_{f\circ\De_2})\circ\psi_{\De_1}=
\rho\circ\Diff^+(\De_2^P(f))\circ\psi_{\De_1}\\
&=\De_1^P(\De_2^P(f)).
\end{align*}
Hence, $(\cdot)^P$ preserves composition.
\end{proof}

From the lemma proved it follows that any right differential module is
a right module over the algebra $\Diff(A)$.

Let $\dotsb\to Q_i\xra{\De_i}Q_{i+1}\to\dotsb$ be a complex of
differential operators and $P$ be a right differential module. Then, by
Lemma \vref{sec1:le:rmod}, we can construct the dual complex
$\dotsb\xla{}
\Hom_A(Q_i,P)\xla{\De_i^P}\Hom_A(Q_{i+1},P)\xla{}\dotsb$ with
coefficients in $P$. Note that the $\Diff$-Spencer complex is a
particular case of this construction. In fact, due to properties of the
homomorphism $c_{\infty,\infty}$ the module $\Diff^+(P)$ is a right
differential module with $\rho=c_{\infty,\infty}$. Applying Lemma
\vref{sec1:le:rmod} to the de Rham complex, we obtain the
$\Diff$-Spencer complex.

Note also that if $\De\colon P\to Q$ is a differential operator, then the cokernel
$C_\De$ of the homomorphism $\psi_\De^\infty\colon \Diff^+(P)\to\Diff^+(Q)$
inherits the right differential module structure of $\Diff^+(Q)$. Thus we
can consider the complex
$$
0\xla{}\coker\De\xla{\cD}C_\De\xla{}\Dr_1(C_\De)\xla{}\dotsb
\xla{}\Dr_i(C_\De)\xla{}\Dr_{i+1}(C_\De)\xla{}\dotsb
$$
dual to the de Rham complex with coefficients in $C_\De$. It is called
the \emph{$\Diff$-Spencer complex} of the operator $\De$.

\subsection{The Spencer cohomology}\label{sub:exact}

Consider an important class of commutative algebras.
\bde\label{sec1:df:smootha}
An algebra $A$ is called \emph{smooth}, if the module $\La^1$ is projective
and of finite type.
\ede

In this section we shall work over a smooth algebra $A$.

Take two $\Diff$-Spencer complexes, of orders $k$ and $k-1$, and consider
their embedding
$$\begin{CD}
0@<<<P@<<<\Diff_k^+(P))@<<<\Dr_1(\Diff_{k-1}^+(P))@<<<\dotsb\\
@.   @|     @AAA             @AAA\\
0@<<<P@<<<\Diff_{k-1}^+(P))@<<<\Dr_1(\Diff_{k-2}^+(P))@<<<\dotsb
\end{CD}$$
Then, if the algebra $A$ is smooth, the direct sum of the corresponding
quotient complexes is of the form
$$
0\xla{}\Smbl(A,P)\xla{\de}\Dr_1(\Smbl(A,P))\xla{\de}\Dr_2(\Smbl(A,P))\xla{}
\dotsb
$$
By standard reasoning, exactness of this complex implies that of
$\Diff$-complexes.

\begin{xca}
Prove that the operators $\de$ are $A$-homomorphisms.
\end{xca}

Let us describe the structure of the modules $\Smbl(A,P)$. For the time
being, use the notation $D=\Dr_1(A)$. Consider the homomorphism $\al_k\colon
P\ot_AS^k(D)\to\Smbl_k(A,P)$ defined by
$$
\al_k(p\ot\na_1\cdot\dots\cdot\na_k)=\smbl_k(\De),\quad
\De(a)=(\na_1\circ\dots\circ\na_k)(a)p,
$$
where $a\in A$, $p\in P$, and $\smbl_k\colon\Diff_k(A,P)\xra{}\Smbl_k(A,P)$
is the natural projection.
\ble\label{sec1:le:symbols}
If $A$ is a smooth algebra, the homomorphism $\al_k$ is an isomorphism.
\ele
\begin{proof}
Consider a differential operator $\De\colon A\to P$ of order $\le k$. Then
the map $s_\De\colon A^{\ot k}\to P$ defined by $s_\De(a_1,\dots,a_k)=
\de_{a_1,\dots,a_k}(\De)$ is a symmetric multiderivation and thus the
correspondence $\De\mapsto s_\De$ generates a homomorphism
\begin{equation}\label{sec1:eq:iso}
\Smbl_k(A,P)\to\Hom_A(S^k(\La^1),P)=S^k(D)\ot_AP,
\end{equation}
which, as it can be checked by direct computation, is inverse to $\al_k$.
Note that the second equality in \veqref{sec1:eq:iso} is valid because $A$ is
a smooth algebra.
\end{proof}

\begin{xca}
Prove that the module $\Smbl_k(P,Q)$ is isomorphic to the module $S^k(D)
\ot_A\Hom_A(P,Q)$.
\end{xca}

\begin{xca}\label{sec1:es:Kosz}
Dualize Lemma \vref{sec1:le:symbols}. Namely, prove that the kernel of
the natural projection $\nu_{k,k-1}\colon \CJ^k(P)\to\CJ^{k-1}(P)$ is
isomorphic to $S^k(\La^1)\ot_AP$, with the isomorphism $\al^k\colon
S^k(\La^1)\ot_AP\to\ker\nu_{k,k-1}$ given by
\[
\al^k(da_1\cdot\ldots\cdot da_k\ot p)=\delta_{a_1,\dots,a_k}(j_k)(p),
\qquad p\in P.
\]
\end{xca}

Thus we obtain:
$$\Dr_i(\Smbl_k(P))=\Hom_A(\La^i,P\ot_AS^k(D))=P\ot_AS^k(D)\ot_A\La^i(D).$$
But from the definition of the Spencer operator it easily follows
that the action of the operator
$$\de\colon P\ot_AS^k(D)\ot_A\La^i(D)\to
P\ot_AS^{k+1}(D)\ot_A\La^{i-1}(D)$$
is expressed by
\begin{multline*}
\de(p\ot\si\ot\na_1\wg\dots\wg\na_i) \\
=\sum_{l=1}^i(-1)^{l+1}p\ot\si\cdot\na_l\ot\na_1\wg\dots\wg\hat{\na}_l\wg
\dots\wg\na_i
\end{multline*}
where $p\in P$, $\si\in S^k(D)$, $\na_l\in D$ and the ``hat'' means that the
corresponding term is omitted. Thus we see that the operator $\de$ coincides
with the Koszul differential (see the Appendix) which implies exactness of
$\Diff$-Spencer complexes.

The Jet-Spencer complexes are dual to them and consequently, in the
situation under consideration, are exact as well. This can also be proved
independently by considering
two Jet-Spencer complexes of orders $k$ and $k-1$ and their projection
\[
\begin{CD}
0@>>>P@>>>\CJ^k(P))@>>>\La^1\ot_A\CJ^{k-1}(P)@>>>\dotsb\\
@.   @|     @VVV             @VVV\\
0@>>>P@>>>\CJ^{k-1}(P))@>>>\La^1\ot_A\CJ^{k-2}(P)@>>>\dotsb\\
\end{CD}
\]
\label{sec1:pg:jetsp}
Then the corresponding kernel complexes are of the form
\begin{multline*}
0\xra{}S^k(\La^1)\ot_A P \xra{\de}\La^1\ot_AS^{k-1}(\La^1)\ot_A P\\
\xra{\de}\La^2\ot_AS^{k-2}(\La^1)\ot_A P\xra{} \dotsb
\end{multline*}
and are called the \emph{$\de$-Spencer complexes} of $P$.
These are complexes of $A$-homomorphisms. The operator
$$\de\colon\La^s\ot_AS^{k-s}(\La^1)\ot_A
P\to\La^{s+1}\ot_AS^{k-s-1}(\La^1)\ot_A P$$
is defined by $\de(\omega\ot u\ot p)=(-1)^s\omega\wedge i(u)\ot p$,
where $i\colon
S^{k-s}(\La^1)\to\La^1\ot S^{k-s-1}(\La^1)$ is the natural inclusion.
Dropping the multiplier $P$ we get the de~Rham complexes with polynomial
coefficients. This proves that the $\de$-Spencer complexes and,
therefore, the Jet-Spencer complexes are exact.

Thus we have the following
\begin{theorem}\label{sec1:th:exactspen}
If $A$ is a smooth algebra, then all $\Diff$-Spencer complexes and
\textup{Jet}-Spencer complexes are exact.
\end{theorem}

Now, let us consider an operator $\De\colon P\to P_1$ of order $\le k$.
Our aim is to compute the Jet-Spencer cohomology of
$\De$, i.e., the cohomology of the complex \veqref{sec1:eq:jspcd}.
\begin{definition}
A complex of \cd operators $\dotsb\xra{}P_{i-1}\xra{\Delta_i}P_i
\xra{\Delta_{i+1}}P_{i+1}\xra{}\dotsb$ is called \emph{formally exact}, if
the complex
\[
\dotsb \xra{} \J^{k_i+k_{i+1}+l}(P_{i-1})
\xra{\varphi_{\Delta_i}^{k_i+k_{i+1}+l}} \J^{k_{i+1}+l}(P_i)
\xra{\varphi_{\Delta_{i+1}}^{k_{i+1}+l}}\J^l(P_{i+1}) \xra{} \dotsb,
\]
with $\ord\Delta_j\le k_j$, is exact for any $l$.
\end{definition}
\begin{theorem}
\label{sec1:th:jetcomp}
Jet-Spencer cohomology of $\De$ coincides with the cohomology of any
formally exact complex of the form
\[
0\xra{}P\xra{\De}P_1\xra{}P_2\xra{}P_3\xra{}\dotsb
\]
\end{theorem}
\begin{proof}
Consider the following commutative diagram
\[\minCDarrowwidth=20pt
\begin{CD}
@. \vdots @. \vdots @. \vdots @. \\
@.  @AAA  @AAA  @AAA  @. \\
0 @>>> \La^2\otimes\CJ^{\infty}(P) @>>> \La^2\otimes\CJ^{\infty}(P_1)
@>>> \La^2\otimes\CJ^{\infty}(P_2) @>>> \dotsb \\
@.  @AA\hd A  @AA\hd A  @AA\hd A  @. \\
0 @>>> \La^1\otimes\CJ^{\infty}(P) @>>> \La^1\otimes\CJ^{\infty}(P_1)
@>>> \La^1\otimes\CJ^{\infty}(P_2) @>>> \dotsb \\
@.  @AA\hd A  @AA\hd A  @AA\hd A  @. \\
0 @>>> \CJ^{\infty}(P) @>>> \CJ^{\infty}(P_1) @>>> \CJ^{\infty}(P_2) @>>> \dotsb \\
@.  @AAA  @AAA  @AAA  @. \\
@.  0  @. 0 @. 0  @.
\end{CD}
\]
where the $i$-th column is the de Rham complex with coefficients in the
left differential module $\CJ^\infty(P_i)$.
The horizontal maps are induced by the operators $\Delta_i$. All the
sequences are exact except for the terms in the left column and the
bottom row.  Now the standard spectral sequence arguments (see the Appendix)
completes the proof.
\end{proof}

Our aim now is to prove that in a sense all compatibility complexes are
formally exact. To this end, let us discuss the notion of
involutiveness of a differential operator.

The map $\psi_l^{\De}\colon\CJ^{k+l}(P)\to\CJ^l(P_1)$ gives rise to the
map
$$\smbl_{k,l}(\De)\colon S^{k+l}(\La^1)\ot P\to S^l(\La^1)\ot P_1$$
called the \emph{$l$-th prolongation of the symbol of $\De$}.
\begin{xca}
Check that 0-th prolongation map
$\smbl_{k,0}\colon\Diff_k(P,P_1)\to\Hom(S^k(\La^1)\ot P,P_1)$
coincides with the natural projection of differential operators to their
symbols,
$\smbl_k\colon\Diff_k(P,P_1)\to\Smbl_k(P,P_1)$.
\end{xca}

Consider the \emph{symbolic module} $g^{k+l}=\ker\smbl_{k,l}(\De)\subset
S^{k+l}(\La^1)\ot P$ of the operator $\De$. It is easily
shown that the subcomplex of the $\de$-Spencer complex
\begin{equation}
0 \xra{} g^{k+l} \xra{\delta} \La^1\otimes g^{k+l-1}
\xra{\delta} \La^2\otimes g^{k+l-2} \xra{\delta} \dotsb
\end{equation}
is well defined. The cohomology of this complex in the term $\La^i\otimes
g^{k+l -i}$ is denoted by $H^{k+l,i}(\Delta)$ and is said to be
\emph{$\delta$-Spencer cohomology of the operator $\Delta$}.
\begin{xca}
Prove that $H^{k+l,0}(\Delta)=H^{k+l,1}(\Delta)=0$.
\end{xca}
The operator $\Delta$ is called \emph{involutive} (in the sense of
Cartan), if $H^{k+l,i}(\Delta)=0$ for all $i\ge 0$.

\bde
An operator $\De$ is called \emph{formally integrable}, if for all $l$
modules
$E_{\De}^l=\ker\psi_{\De}^l\subset\CJ^{k+l}(P)$ and $g^{k+l}$ are
projective and the natural mappings $E_{\De}^l\to E_{\De}^{l-1}$
are surjections.
\ede
Till the end of this section we shall assume all the operators under
consideration to be formally integrable.
\begin{theorem}\label{sec1:th:fecc}
If the operator $\Delta$ is involutive, then the compatibility complex
of $\De$ is formally exact for all positive integers $k_1$, $k_2$,
$k_3$, \dots.
\end{theorem}
\begin{proof}
Suppose that the compatibility complex of $\De$
\[
P\xra{\De}P_1\xra{\De_1}P_2\xra{\De_2}\dotsb
\]
is formally exact in terms $P_1$, $P_2$, \dots, $P_{i-1}$. The
commutative diagram
\begin{multline*}
\setlength{\multlinegap}{0pt}
\begin{CD}
@. 0@.0@.0@. @. \\
@. @VVV @VVV @VVV  \\
0@>>>g^K@>>>S^K\ot P@>>>S^{K-k}\ot P_1@>>>\dotsb \\
@. @VVV @VVV @VVV \\
0@>>>E_{\De}^{K-k}@>>>\CJ^K(P)@>>>\CJ^{K-k}(P_1)@>>>\dotsb \\
@. @VVV @VVV @VVV \\
0@>>>E_{\De}^{K-k-1}@>>>\CJ^{K-1}(P)@>>>\CJ^{K-k-1}(P_1)@>>>\dotsb \\
@. @VVV @VVV @VVV \\
@. 0@.0@.0 @.
\end{CD} \\[7pt]
\begin{CD}
@. 0 @. 0 \\
@. @VVV @VVV \\
\dotsb@>>>S^{k_i}\ot P_i@>>>P_{i+1}@>>>0 \\
@. @VVV @VVV \\
\dotsb@>>>\CJ^{k_i}(P_i)@>>>P_{i+1}@>>>0 \\
@. @VVV @VVV \\
\dotsb@>>>\CJ^{k_i-1}(P_i)@>>>0 \\
@. @VVV \\
@.0
\end{CD}
\end{multline*}
where $S^j=S^j(\La^1)$, $K=k+k_1+k_2+\dots+k_i$, shows that the complex
\[
0\xra{}g^K\xra{}S^K\ot P\xra{}S^{K-k}\ot P_1\xra{}\dotsb
\xra{}S^{k_i}\ot P_i
\]
is exact.

What we must to prove is that the sequences
\[
S^{k_{i-1}+k_i+l}\otimes P_{i-1} \xra{} S^{k_i+l}\otimes P_i \xra{}
S^l\otimes P_{i+1}
\]
are exact for all $l\ge 1$. The proof is by induction on $l$, with
the inductive step involving the standard spectral sequence
arguments applied to the commutative diagram
\[
\minCDarrowwidth=20pt
\begin{CD}
0 @>>> S^l\otimes P_{i+1} @>\delta >>
\La^1\otimes S^{l-1}\otimes P_{i+1} @>\delta >>
\La^2\otimes S^{l-2}\otimes P_{i+1} @>\delta >> \dotsb \\
@. @AAA  @AAA  @AAA @. \\
0 @>>> S^{k_i+l}\otimes P_i @>\delta >>
\La^1\otimes S^{k_i+l-1}\otimes P_i @>\delta >>
\La^2\otimes S^{k_i+l-2}\otimes P_i @>\delta >> \dotsb \\
@. @AAA  @AAA  @AAA @. \\
@. \vdots @. \vdots @. \vdots  @.  \\
@. @AAA  @AAA  @AAA @. \\
0 @>>> S^{K+l}\otimes P_0 @>\delta >>
\La^1\otimes S^{K+l-1}\otimes P_0 @>\delta >>
\La^2\otimes S^{K+l-2}\otimes P_0 @>\delta >> \dotsb \\
@. @AAA  @AAA  @AAA @. \\
0 @>>> g^{K+l} @>\delta >>
\La^1\otimes g^{K+l-1} @>\delta >>
\La^2\otimes g^{K+l-2} @>\delta >> \dotsb \\
@. @AAA  @AAA  @AAA @. \\
@.  0  @.  0  @.  0  @.
\end{CD}\]
\end{proof}

\begin{example}
For the de~Rham differential $d\colon A\to\Lambda^1$ the symbolic
modules $g^l$ are trivial. Hence, the de~Rham differential is
involutive and, therefore, the de~Rham complex is formally exact.
\end{example}

\begin{example}
\label{sec1:ex:pf}
Consider the geometric situation and suppose that the manifold $M$ is
a (pseudo-)Riemannian manifold. For an integer $p$ consider the operator
$\Delta=d{*}d\colon\Lambda^p\to\Lambda^{n-p}$, where $*$ is the Hodge
star operator on the modules of differential forms.
Let us show that the complex
\[
\hL^p \xra{\Delta} \hL^{n-p} \xra{d} \hL^{n-p+1} \xra{d}
\Lambda^{n-p+2} \xra{d} \dotsb \xra{d} \Lambda^n \xra{} 0
\]
is formally exact and, thus, is the compatibility complex for the
operator $\Delta$. In view of the previous example we must prove that
the image of the map $\smbl(\Delta)\colon S^{l+2}\otimes\Lambda^p \to
S^l\otimes\Lambda^{n-p}$ coincides with the image of the map
$\smbl(d)\colon S^{l+1}\otimes\Lambda^{n-p-1}\to
S^l\otimes\Lambda^{n-p}$ for all $l\ge 0$. Since
$\Delta*=d{*}d*=d(*d{*}+d)$, it is sufficient to show that the map
$\smbl(*d{*}+d)\colon
S^{l+1}\otimes(\Lambda^{n-p+1}\oplus\Lambda^{n-p-1}) \to
S^l\otimes\Lambda^{n-p}$ is an epimorphism. Consider $\smbl(L)\colon
S^l\otimes\Lambda^{n-p}\to S^l\otimes\Lambda^{n-p}$, where
$L=(*d{*}+d)(*d{*}\pm d)$ is the Laplace operator. From coordinate
considerations it easily follows that the symbol of the Laplace operator
is epimorphic, and so the symbol of the operator $*d{*}+d$ is also
epimorphic.
\end{example}

The condition of involutiveness is not necessary for the
formal exactness of the compatibility complex due to the following

\begin{theorem}[$\delta$-Poincar\'e lemma]
If the algebra $A$ is Noetherian, then for any operator $\Delta$
there exists an integer $l_0=l_0(m,n,k)$, where $m=\rank P$, such that
$H^{k+l,i}(\Delta)=0$ for $l\ge l_0$ and $i\ge 0$.
\end{theorem}

Proof can be found, e.g., in \cite{KLV,BCGGG}.
\label{sec1:pg:ecc}
Thus, from the proof of Theorem \vref{sec1:th:fecc} we see that for
sufficiently large integer $k_1$ the compatibility complex is formally
exact for any operator $\De$.

We shall always assume that compatibility complexes are formally exact.

\subsection{Geometrical modules}\label{sub:geommod}
There are several directions to generalize or specialize the above described
theory. Probably, the most important one, giving rise to various interesting
specializations, is associated with the following concept.
\bde
An abelian subcategory $\CM(A)$ of the category of all $A$-modules is
said to be \emph{differentially closed}, if
\begin{enumerate}
\item it is closed under tensor product over $A$,
\item it is closed under the action of the functors
$\Diff^{(+)}_k(\cdot,\cdot)$ and $\Dr_i(\cdot)$,
\item the functors
$\Diff^{(+)}_k(P,\cdot),\Diff^{(+)}_k(\cdot,Q)$ and $\Dr_i(\cdot)$
are representable in $\CM(A)$, whenever $P$, $Q$ are objects of $\CM(A)$.
\end{enumerate}
\ede

As an example consider the following situation. Let $M$ be a smooth
(i.e., $\Ci$-class) finite-dimensional manifold and set $A=\Ci(M)$. Let
$\pi:E\to M$, $\xi:F\to M$ be two smooth locally trivial
finite-dimensional vector bundles over $M$ and $P=\Ga(\pi),Q=\Ga(\xi)$
be the corresponding $A$-modules of smooth sections.

One can prove that the module $\Diff^{(+)}_k(P,Q)$ coincides with the
module of $k$-th order differential operators acting from the bundle
$\pi$ to $\xi$ (see Proposition \vref{sec1:pr:equiv}). Further, the
module $\Dr(A)$ coincides with the module of vector fields on the
manifold $M$.

However if one constructs representative objects for the functors such as
$\Diff_k(P,\cdot)$ and $\Dr_i(\cdot)$ in the category of all
$A$-modules, the modules $\CJ^k(P)$ and $\La^i$ will not coincide with
``geometrical'' jets and differential forms.

\begin{xca}
Show that in the case $M=\mathbb{R}$ the form $d(\sin x)-\cos x\,dx$ is
nonzero.
\end{xca}

\bde
A module $P$ over $\Ci(M)$ is called \emph{geometrical}, if
$$\bigcap_{x\in M}\mu_xP=0,$$
where $\mu_x$ is the ideal in $\Ci(M)$ consisting of functions vanishing
at point $x\in M$.
\ede
Denote by $\CG(M)$ the full subcategory of the category of all modules
whose objects are geometrical $\Ci(M)$-modules. Let $P$ be an
$A$-module and set
$$
\CG(P)=P\!\Big/\!\!\bigcap_{x\in M}\mu_xP.
$$
Evidently,
$\CG(P)$ is a geometrical module while the correspondence $P\Ra \CG(P)$
is a functor from the category of all $\Ci(M)$-modules to the category
$\CG(M)$ of geometrical modules.

\bpr
Let $M$ be a smooth finite-dimensional manifold and $A=C^\infty(M)$. Then
\begin{enumerate}
\item The category $\CG(A)$ of geometrical $A$-modules is
differentially closed.
\item The representative objects for the functors
$\Diff_k(P,\cdot)$ and $\Dr_i(\cdot)$ in $\CG(A)$ coincide with
$\CG(\CJ^k(P))$ and $\CG(\La^i)$ respectively.
\item The module $\CG(\La^i)$ coincides with the module of
differential $i$-forms on $M$.
\item If $P=\Ga(\pi)$ for a smooth locally trivial
finite-dimensional vector bundle $\pi:E\to M$, then the module
$\CG(\CJ^k(P))$ coincides with the module $\Ga(\pi_k)$, where $\pi_k:
J^k(\pi)\to M$
is the bundle of $k$-jets for the bundle $\pi$ \textup{(}see Section
\ref{sub:fjets}\textup{)}.
\end{enumerate}
\epr
\begin{xca}
Prove (1), (2), and (3) above.
\end{xca}
The situation described in this Proposition will be referred to as the
\emph{geometrical} one.

Another example of a differentially closed category is the category of
filtered geometrical modules over a filtered algebra. This category is
essential to construct differential calculus over manifolds of infinite jets
and infinitely prolonged differential equations (see Sections
\ref{sub:ijets} and \ref{sub:basstr} respectively).

\begin{remark}
The logical structure of the above described theory is obviously
generalized to the supercommutative case. For a noncommutative
generalization see \cite{Verb3,Verb4}.
\end{remark}

\newpage

\section{Algebraic model for Lagrangian formalism}\label{sec:lagr}
Using the above introduced algebraic concepts, we shall construct now
an algebraic model for Lagrangian formalism; see also \cite{Verb1}. For
geometric motivations, we refer the reader to Section \ref{css:sec} and to
Subsection \ref{css.genfun:subsec} especially.

\subsection{Adjoint operators}
Consider an $A$-module $P$ and the complex of $A$-homomorphisms
\begin{equation}\label{sec2:eq:adjcom}
0\xra{}\Diff^+(P,A)\xra{w}\Diff^+(P,\La^1)\xra{w}\Diff^+(P,\La^2)\xra{w}\dotsb,
\end{equation}
where, by definition, $w(\na)=d\circ\na\in\Diff^+(P,\La^{i+1})$ for the
operator $\na\in\Diff^+(P,\La^i)$.
Let $\hat{P}_n,\,n\ge 0$, be the cohomology module of this complex at
the term $\Diff^+(P,\La^n)$.

Any operator $\De\colon P\to Q$ determines the natural cochain map
$$\begin{CD}
\dotsb@>>>\Diff^+(Q,\La^{i-1})@>{w}>>\Diff^+(Q,\La^i)@>>>\dotsb\\
@. @V\tilde{\De}VV@V\tilde{\De}VV\\
\dotsb@>>>\Diff^+(P,\La^{i-1})@>{w}>>\Diff^+(P,\La^i)@>>>\dotsb
\end{CD}$$
where $\tilde{\De}(\na)=\na\circ\De\in\Diff^+(P,\La^i)$ for $\na\in
\Diff^+(Q,\La^i)$.
\bde\label{sec2:df:adj}
The cohomology map $\De_n^*\colon \hat{Q}_n\to\hat{P}_n$ induced by $\tilde\De$
is called the ($n$-th) \emph{adjoint operator} for $\De$.
\ede
Below we assume $n$ to be fixed and omit the corresponding subscript.
The main properties of the adjoint operator are described by
\bpr\label{sec2:pr:adjprop}
Let $P,Q$ and $R$ be $A$-modules. Then
\begin{enumerate}
\item\label{sec2:it:ord}
If $\De\in\Diff_k(P,Q)$, then $\De^*\in\Diff_k(\hat{Q},\hat{P})$.
\item\label{sec2:it:comp}
If $\De_1\in\Diff(P,Q)$ and $\De_2\in\Diff(Q,R)$, then $(\De_2\circ\De_1)^*=
\De_1^*\circ\De_2^*$.
\end{enumerate}
\epr
\begin{proof}
Let $[\na]$ denote the cohomology class of $\na\in\Diff^+(P,\La^n)$, where
$w(\na)=0$.

(1) Let $a\in A$. Then
\begin{align*}
\de_a(\De^*)([\na])&=\De^*([\na])-\De^*(a[\na])=
[\na\circ a\circ\De]-[\na\circ\De\circ a]\\
&=(a\circ\De)^*([\na])-(\De\circ a)^*([\na])=-\de_a(\De^*)([\na]).
\end{align*}
Consequently, $\de_{a_0,\dots,a_k}(\De^*)=(-1)^{k+1}
(\de_{a_0,\dots,a_k}(\De))^*$ for any $a_0,\dots,a_k\in A$.

(2) The second statement is implied by the following identities:
$$(\De_2\circ\De_1)^*([\na])=[\na\circ\De_2\circ\De_1]=
\De_1^*([\na\circ\De_2])=\De_1^*(\De_2^*([\na])),$$
which concludes the proof.
\end{proof}
\begin{example}
Let $a\in A$ and $a=a_P\colon P\to P$ be the operator of multiplication by
$a$: $p\mapsto ap$. Then obviously $a_P^*=a_{\hat{P}}$.
\end{example}
\begin{example}
Let $p\in P$ and $p\colon A\to P$ be the operator acting by $a\mapsto
ap$. Then, by Proposition \ref{sec2:pr:adjprop}~\veqref{sec2:it:ord},
$p^*\in\Hom_A(\hat{P},\hat{A})$.  Thus there exists a natural paring
$\langle\cdot,\cdot\rangle\colon P\ot_A\hat{P} \to\hat{A}$ defined by
$\langle p,\hat{p}\rangle=p^*(\hat{p})$, $\hat{p}\in \hat{P}$.
\end{example}

\subsection{Berezinian and integration}
Consider a complex of differential operators $\dotsb\xra{}
P_k\xra{\De_k}P_{k+1} \xra{}\dotsb$. Then, by Proposition
\vref{sec2:pr:adjprop}, $\dotsb\xla{}\hat{P}_k\xla{\De_k^*}\hat{P}_{k+1}
\xla{}\dotsb$ is a complex of differential operators as well. This complex
called \emph{adjoint} to the initial one.
\bde
The complex adjoint to the de Rham complex of the algebra $A$ is called
the \emph{complex of integral forms} and is denoted by
$$
0\xla{}\Sigma_0\xla{\de}\Sigma_1\xla{\de}\dotsb,
$$
where $\Sigma_i=\hat{\La}^i, \de=d^*$. The module $\Sigma_0=\hat{A}$ is
called the \emph{Berezinian} (or the \emph{module of the volume forms}) and
is denoted by $\CB$.
\ede
Assume that the modules under consideration are projective and of finite type.
Then we have
$\hat{P}=\Hom_A(P,\CB)$. In particular, $\Sigma_i=\hat{\La}^i=\Dr_i(\CB)$.

Let us calculate the Berezinian in the geometrical situation (see Subsection
\ref{sub:geommod}), when $A=\Ci(M)$.
\begin{theorem}
\label{sec2:th:ber}
If $A=\Ci(M)$, $M$ being a smooth finite-dimensional manifold, then
\begin{enumerate}
\item $\hat{A}_s=0$ for $s\ne n=\dim M$.
\item $\hat{A}_n=\CB=\La^n$, i.e., the Berezinian coincides with the
module of forms of maximal degree. This isomorphism takes each form
$\omega\in\La^n$ to the cohomology class of the zero-order operator
$\omega\colon A\to\La^n$, $f\mapsto f\omega$.
\end{enumerate}
\end{theorem}
The proof is similar to that of Theorem \vref{sec1:th:exactspen} and is
left to the reader.

In the geometrical situation there exists a natural isomorphism
$\La^i\to\Dr_{n-i}(\La^n)=\Sigma_i$ which takes $\omega\in\La^i$ to the
homomorphism $\omega\colon\La^{n-i}\to\La^n$ defined by
$\omega(\eta)=\eta\wedge\omega$, $\eta\in\La^{n-i}$.

\begin{xca}
Show that $\langle\omega_1,\omega_2\rangle=\omega_1\wedge\omega_2$,
$\omega_1\in\La^i$, $\omega_2\in\La^{n-i}$.
\end{xca}

\begin{xca}
Prove that $d_i^*=(-1)^{i+1}d_{n-i-1}$, where
$d_i\colon\La^i\to\La^{i+1}$ is the de~Rham differential.
\end{xca}
Thus, in the geometrical situation the complex of integral forms
coincides (up to a sign) with the de~Rham complex.

\begin{xca}
Prove the coordinate formula for the adjoint operator:
\begin{enumerate}
\item if $\Delta=\sum_{\sigma}a_{\sigma}
\dd{^{\abs{\sigma}}}{x_{\sigma}}$ is a scalar operator, then
$\Delta^*=\sum_{\sigma}(-1)^{\abs{\sigma}}
\dd{^{\abs{\sigma}}}{x_{\sigma}}\circ a_{\sigma}$;
\item if $\Delta=\norm{\Delta_{ij}}$ is a matrix operator, then
$\Delta^*=\norm{\Delta_{ji}^*}$.
\end{enumerate}
\end{xca}

The operator $\cD\colon \Diff^+(\La^k)\to\La^k$ defined
\vpageref{sec1:eq:cD} generates the map $\int\colon \CB\to H^*(\La^\bullet)$
from the
Berezinian to the de Rham cohomology group of $A$. Namely, for any operator
$\na\in\Diff(A,\La^n)$ satisfying $d\circ\na=0$ we set $\int[\na]=[\na(1)]$,
where $[\cdot]$ denotes the cohomology class.
\bpr
The map $\int\colon \CB\to H^*(\La^\bullet)$ possesses the following
properties:
\begin{enumerate}
\item
If $\om\in\Sigma_1$, then $\int\de\om=0$.
\item For any differential operator $\De\colon P\to Q$ and elements $p\in P$,
$\hat{q}\in\hat{Q}$ the identity
$$
\int\langle\De(p),\hat{q}\rangle=\int\langle p,\De^*(\hat{q})\rangle
$$
holds.
\end{enumerate}
\epr
\begin{proof}
(1) Let $\om=[\na]\in\Sigma_1$. Then $\de\om=[\na\circ d]$ and consequently
$\int\om=[\na d(1)]=0$.

(2) Let $\hat{q}=[\na]$ for some operator $\na\colon Q\to\La^n$. Then
\begin{multline*}
\int\langle\De(p),\hat{q}\rangle=\int[\na\De(p)]=\int\na\circ\De\circ p\\
=\int\langle p,[\na\circ\De]\rangle=\int\langle p,\De^*(\hat{q})\rangle,
\end{multline*}
which completes the proof.
\end{proof}

\bre
Note that the Berezinian $\CB$ is a differential right module (see Subsection
\ref{sub:dmod}) and the complex of integral forms may be understood as
the complex dual to the de Rham complex with coefficients in $\CB$.
\ere

\begin{xca}
Show that in the geometrical situation the right action of vector
fields can also be defined via $X(\omega)=-\Ld_X(\omega)$, where $\Ld_X$ is
the Lie derivative.
\end{xca}

Now we establish a relationship between the de~Rham cohomology and the
homology of the complex of integral forms.

\begin{proposition}[algebraic Poincar\'e duality]
There exists a spectral sequence $(E^r_{p,q},d^r_{p,q})$ with
$$E^2_{p,q}=H_p((\Sigma_\bullet)_{-q}),$$
the homology of complexes of integral forms, and converging to the
de~Rham cohomology $H(\Lambda^\bullet)$.
\end{proposition}
\begin{proof}
Consider the commutative diagram
\[
\minCDarrowwidth=23pt
\begin{CD}
@.  0 @.  0    @.   0   @.\\
@.    @AAA    @AAA    @AAA\\
0@>>> \Diff^+(A,A)@>w>>\Diff^+(A,\La^1)@>w>>\Diff^+(A,\La^2)@>>>\dotsb\\
@.   @AA{\tilde{d}}A @AA{\tilde{d}}A      @AA{\tilde{d}}A\\
0@>>>
\Diff^+(\La^1,A)@>w>>\Diff^+(\La^1,\La^1)@>w>>\Diff^+(\La^1,\La^2)@>>>\dotsb\\
@.    @AA{\tilde{d}}A @AA{\tilde{d}}A      @AA{\tilde{d}}A\\
0@>>>
\Diff^+(\La^2,A)@>w>>\Diff^+(\La^2,\La^1)@>w>>\Diff^+(\La^2,\La^2)@>>>\dotsb\\
@.    @AA{\tilde{d}}A @AA{\tilde{d}}A      @AA{\tilde{d}}A\\
@.  \vdots@.    \vdots@.      \vdots
\end{CD}
\]
where the differential $\tilde{d}\colon\Diff^+(\La^{k+1},P)\to
\Diff^+(\La^k,P)$ is defined by
$\tilde{d}(\Delta)=\Delta\circ d$. The statement follows easily from
the standard spectral sequence arguments.
\end{proof}

\subsection{Green's formula}
Let $Q$ be an $A$-module. Then a natural homomorphism $\xi_Q\colon
Q\to\Hat{\Hat Q}$ defined by $\xi_Q(q)(\hat q)=\langle
q,\hat{q}\rangle$ exists. Consequently, to any operator $\De\colon P\to
\hat{Q}$ there
corresponds the operator $\De^\circ\colon Q\to\hat{P}$, where
$\De^\circ=\De^*\circ\xi_Q$. This operator will also be called
\emph{adjoint} to $\De$.

\begin{remark}
In the geometrical situation the two notions of adjointness coincide.
\end{remark}

\begin{example}
Let $\hat{q}\in\hat{Q}$ and $\hat{q}\colon A\to\hat{Q}$ be the
zero-order operator defined by $a\mapsto a\hat{q}$. The adjoint operator
is $\hat{q}$ itself understood as an element of $\Hom_A(Q,\CB)$.
\end{example}

\bpr\label{sec2:pr:adjpr}
The correspondence $\De\mapsto\De^\circ$ possesses the following properties:
\begin{enumerate}
\item Let $\De\in\Diff(P,\hat{Q})$ and $\De(p)=[\na_p]$, where
$\na_p\in\Diff(Q,\La^i)$. Then $\De^\circ(q)=[\square_q]$, where $\square_q
\in\Diff(P,\La^i)$ and $\square_q(p)=\na_p(q)$.
\item For any $\De\in\Diff(P,\hat{Q})$, one has $(\De^\circ)^\circ=\De$.
\item For any $a\in A$, one has $(a\De)^\circ=\De^\circ\circ a$.
\item If $\De\in\Diff_k(P,\CB)$, then $\De^\circ=j_k^*\circ(a\De)$.
\item If $X\in\Dr_1(\CB)$, then $X+X^\circ=\de X\in\Diff_0(A,\CB)=\CB$.
\end{enumerate}
\epr
\begin{proof}
Statements (1), (3), and (4) are the direct consequences of the definition.
Statement (2) is implied by (1). Let us prove (5).

Evidently, $\de_a(j_1)=j_1(a)-aj_1(1)\in\CJ^1(A)$. Hence for an operator
$\De\in
\Diff_1(A,P)$ one has $(\de_a(j_1))^*(\De)=\De(a)-a\De(1)=(\de_a\De)(1)$.
Consequently,
$$\de_a(X+X^\circ)(1)=(\de_aX)(1)+(\de_a(j_1^*))(X)=
(\de_aX)(1)-\de_a(j_1)^*(X)=0$$
and finally $\de X=j_1^*(X)=X^\circ(1)=X+X^\circ$.
\end{proof}

Note that Statements (1) and (4) of Proposition \vref{sec2:pr:adjpr} can be
taken for the definition of $\De^\circ$.

Note now that from Proposition \vref{sec1:pr:spenew} it follows that the
modules $\Dr_i(P),\,i\ge 2$, can be described as
$$
\Dr_i(P)=\{\,\na\in\Diff_1(\La^{i-1},P)\mid\na\circ d=0\,\}.
$$
Taking $\CB$ for $P$, one can easily show that $\de\na=\na^\circ(1)$ and the
last equality holds for $i=1$ as well. Proposition \vref{sec2:pr:adjpr}
shows that the correspondence $\De\mapsto\De^\circ$ establishes an
isomorphism between the modules $\Diff(P,\hat{Q})$ and $\Diff^+(Q,\hat{P})$
which, taking into account Proposition \vref{sec1:pr:spenew}, means that
the $\Diff$-Spencer complex of the module $\hat{P}$ is isomorphic to the
complex
\begin{equation}\label{sec2:eq:spenber}
0\xla{}\hat{P}\xla{\mu}\Diff(P,\CB)\xla{\om}\Diff(P,\Sigma_1)\xla{\om}
\Diff(P,\Sigma_2)\xla{}\dotsb,
\end{equation}
where $\om(\na)=\de\circ\na,\,\mu(\na)=\na^\circ(1)$. From Theorem
\vref{sec1:th:exactspen} one immediately obtains
\begin{theorem}\label{sec2:th:exspenber}
Complex \veqref{sec2:eq:spenber} is exact.
\end{theorem}
\bre
Let $\De\colon P\to Q$ be a differential operator. Then obviously the following
commutative diagram takes place:
$$\begin{CD}
0@<<<\hat{Q}@<\mu<<\Diff(Q,\CB)@<\om<<\Diff(Q,\Sigma_1)@<\om<<\dotsb\\
@.@V{\De^*}VV@VVV@VVV\\
0@<<<\hat{P}@<\mu<<\Diff(P,\CB)@<\om<<\Diff(P,\Sigma_1)@<\om<<\dotsb\\
\end{CD}$$
\ere

As a corollary of Theorem \vref{sec2:th:exspenber} we obtain
\begin{theorem}[Green's formula]
If $\De\in\Diff(P,\hat{Q}),p\in P,q\in Q$, then
$$
\langle q,\De(p)\rangle-\langle\De^\circ(q),p\rangle=\de G
$$
for some integral 1-form $G\in\Sigma_1$.
\end{theorem}
\begin{proof}
Consider an operator $\na\in\Diff(A,\CB)$. Then $\na-\na^\circ(1)$ lies in
$\ker\mu$ and consequently there exists an operator $\square\in
\Diff(A,\Sigma_1)$ satisfying $\na-\na^\circ(1)=\om(\square)=\de\circ
\square$. Hence, $\na(1)-\na^\circ(1)=\de G$, where $G=\square(1)$. Setting
$\na(a)=\langle q,\De(ap)\rangle$ we obtain the result.
\end{proof}
\begin{remark}
The integral $1$-form $G$ is dependent on $p$ and $q$. Let us show
that we can choose $G$ in such a way that the map $p\times q\mapsto
G(p,q)$ is a bidifferential operator. Note first that the map
$\omega\colon\Diff^+(A,\Sigma_1)\to\Diff^+(A,\CB)$ is an
$A$-homomorphism. Since the module $\Diff^+(A,\CB)$ is projective,
there exists an $A$-homomorphism
$\vk\colon\im\omega\to\Diff^+(A,\Sigma_1)$ such that
$\omega\circ\vk=\id$. We can put $\square=\vk(\na-\na(1))$. Thus
$G=\vk(\na-\na(1))(1)$. This proves the required statement.
\end{remark}

\begin{remark}
From algebraic point of view, we see that in the geometrical situations
there is the multitude of misleading isomorphisms, e.q.,
$\CB=\La^n$, $\De^\circ=\De^*$, etc. In generalized
settings, for example, in supercommutative situation (see
Subsection~\vref{sub:supeq}), these isomorphisms disappear.
\end{remark}

\subsection{The Euler operator}\label{sub:Euler}
Let $P$ and $Q$ be $A$-modules. Introduce the notation
$$
\Diff_{(k)}(P,Q)=
\underbrace{\Diff(P,\dots,\Diff(P}_{k\ \mathrm{times}},Q)\dots)
$$
and set $\Diff_{(*)}(P,Q)=\bigoplus_{k=0}^\infty\Diff_{(k)}(P,Q)$. A
differential operator $\na\in\Diff_{(k)}(P,Q)$ satisfying the condition
$$\na(p_1,\dots,p_i,p_{i+1},\dots,p_k)=\si\na(p_1,\dots,p_{i+1},p_i,\dots,p_k)
$$
is called \emph{symmetric}, if $\si=1$, and \emph{skew-symmetric}, if $\si=
-1$ for all $i$. The modules of symmetric and skew-symmetric operators will
be denoted by $\Diff_{(k)}^{\sym}(P,Q)$ and $\Diff_{(k)}^{\alt}(P,Q)$,
respectively. From Theorem \vref{sec2:th:exspenber} and Corollary
\vref{sec1:co:iso} it follows that for any $k$ the complex
\begin{equation}\label{sec2:eq:polysigma}
0\xla{}\Diff_{(k)}(P,\CB)\xla{\om}\Diff_{(k)}(P,\Sigma_1)\xla{\om}
\Diff_{(k)}(P,\Sigma_2)\xla{\om}\dotsb,
\end{equation}
where $\om(\na)=\de\circ\na$, is exact in all positive degrees, while
its $0$-homology is of the form
$H_0(\Diff_{(k)}(P,\Sigma_\bullet))=\Diff_{(k-1)}(P,\hat{P})$. This result
can be refined in the following way.
\begin{theorem}\label{sec2:th:symskew}
The symmetric
\begin{equation}\label{sec2:eq:sym}
0\xla{}\Diff_{(k)}^{\sym}(P,\CB)\xla{\om}\Diff_{(k)}^{\sym}(P,\Sigma_1)\xla{\om}
\Diff_{(k)}^{\sym}(P,\Sigma_2)\xla{\om}\dotsb
\end{equation}
and skew-symmetric
\begin{equation}\label{sec2:eq:skew}
0\xla{}\Diff_{(k)}^{\alt}(P,\CB)\xla{\om}\Diff_{(k)}^{\alt}(P,\Sigma_1)\xla{\om}
\Diff_{(k)}^{\alt}(P,\Sigma_2)\xla{\om}\dotsb
\end{equation}
are acyclic complexes in all positive degrees, while the $0$-homologies
denoted by $L_k^{\sym}(P)$ and $L_k^{\alt}(P)$ respectively are of the form
\begin{align*}
L_k^{\sym}&=\{\,\na\in\Diff_{(k-1)}^{\sym}(P,\hat{P})\mid
(\na(p_1,\dots,p_{k-2}))^\circ=\na(p_1,\dots,p_{k-2})\,\},\\
L_k^{\alt}&=\{\,\na\in\Diff_{(k-1)}^{\alt}(P,\hat{P})\mid
(\na(p_1,\dots,p_{k-2}))^\circ=-\na(p_1,\dots,p_{k-2})\,\}
\end{align*}
for $k>1$ and
$$L_1^{\sym}(P)=L_1^{\alt}(P)=\hat{P}.$$
\end{theorem}
\begin{proof}
We shall consider the case of symmetric operators only, since the case of
skew-symmetric ones is proved in the same way exactly.

Obviously, the complex \eqref{sec2:eq:sym} is a direct summand in
\veqref{sec2:eq:polysigma} and due to this fact the only thing we need to
prove is that the diagram
$$\begin{CD}
\Diff_{(k-1)}(P,\hat{P})@<\mu_{(k-1)}<<\Diff_{(k)}(P,\CB)\\
@V{\rho'}VV @VV{\rho}V\\
\Diff_{(k-1)}(P,\hat{P})@<\mu_{(k-1)}<<\Diff_{(k)}(P,\CB)
\end{CD}$$
is commutative. Here
\begin{align*}
\mu_{(k-1)}(\na)(p_1,\dots,p_{k-1})=&\,(\na(p_1,\dots,p_{k-1}))^\circ(1),\\
\rho(\na)(p_1,\dots,p_{k-1},p_k)=&\,\na(p_1,\dots,p_k,p_{k-1}),\\
\rho'(\na)(p_1,\dots,p_{k-2})=&\,(\na(p_1,\dots,p_{k-2}))^\circ.
\end{align*}\label{sec2:pg:mu}
Note that $\mu_{(k-1)}=\Diff_{(k-1)}(\mu)$, where $\mu$ is defined in
\veqref{sec2:eq:spenber}. To prove commutativity, it suffices to consider
the case $k=2$. Let $\na\in\Diff_{(2)}(P,\CB)$ and $\na(p_1,p_2)=
[\De_{p_1,p_2}]$. Then $\mu_{(1)}(\na)(p_1)=[\De'_{p_1}]$, where
$\De'_{p_1}(p_2)=\De_{p_1,p_2}(1)$. Further, $\rho'(\mu_{(1)}(\na))=
[\De''_{p_1}]$, where
$$\De''_{p_1}(p_2)=\De'_{p_2}(p_1)=\De_{p_2,p_1}(1).$$
On the other hand, one has $\rho(\na)(p_1,p_2)=\na(p_2,p_1)$ and
$\mu_{(1)}(\rho(\na))(p_1)=[\square_{p_1}]$, where $\square_{p_1}(p_2)=
\De_{p_2,p_1}(1)$.
\end{proof}
\bde
The elements of the space $\CL ag(P)=\bigoplus_{k=1}^\infty L_k^{\sym}(P)$ are
called \emph{Lagrangians} of the module $P$. An operator $L\in
\Diff_{(*)}^{\sym}(P,\CB)$ is called a \emph{density} of a Lagrangian $\CL$,
if $\CL=L\bmod\im\om$. The natural correspondence $\mathbf{E}\colon
\Diff_{(*)}^{\sym}(P,\CB)\to\Diff_{(*)}^{\sym}(P,\hat{P}),\,L\mapsto\CL$ is
called the \emph{Euler operator}, while operators of the form $\De=
\mathbf{E}(L)$ are said to be \emph{Euler--Lagrange operators}.
\ede
Theorem \vref{sec2:th:symskew} implies the following
\begin{corollary}\label{sec2:co:eulag}
For any projective $A$-module $P$ one has:
\begin{enumerate}
\item An operator $\De\in\Diff_{(*)}^{\sym}(P,\hat{P})$ is an
Euler--Lagrange operator if and only if $\De$ is self-adjoint, i.e.,
if $\De\in L_*^{\sym}(P)$.
\item A density $L\in\Diff_{(*)}^{\sym}(P,\CB)$ corresponds to a
trivial Lagrangian, i.e., $\mathbf{E}(L)=0$, if and only if $L$ is a total
divergence, i.e., $L\in\im\om$.
\end{enumerate}
\end{corollary}
\subsection{Conservation laws}
Denote by $\CF$ the commutative algebra of nonlinear operators\footnote{In
geometrical situation, this algebra is identified with the algebra of
polynomial functions on infinite jets (see the next section).}
$\Diff_{(*)}^{\sym}(P,A)$. Then for
any $A$-module $Q$ one has
$$\Diff_{(*)}^{\sym}(P,Q)=\CF\ot_AQ.$$
Let $\De\in
\CF\ot_AQ$ be a differential operator and let us set $\CF_\De=\CF/
\mathfrak{a}$, where $\mathfrak{a}$ denotes the ideal in $\CF$ generated by
the operators of the form $\square\circ\De$, $\square\in\Diff(Q,A)$.

Thus, fixing $P$, we obtain the functor $Q\Ra\CF\ot_AQ$ and fixing an
operator
$\De\in\Diff_{(*)}(P,Q)$ we get the functor $Q\Ra\CF_\De\ot_AQ$ acting from
the category $\CM_A$ to $\CM_\CF$ and to $\CM_{\CF_\De}$ respectively, where
$\CM$ denotes the category of all modules over the corresponding algebra.
These functors in an obvious way generate natural transformations of the
functors $\Diff_k^{(+)}(\cdot)$, $\Dr_k(\cdot)$, etc., and of their
representative objects
$\CJ^k(P)$, $\La^k$, etc. For example, to any operator $\na\colon Q_1\to Q_2$
there
correspond operators $\CF\ot\na\colon \CF\ot_AQ_1\to\CF\ot_AQ_2$ and
$\CF_\De\ot\na\colon \CF_\De\ot_AQ_1\to\CF_\De\ot_AQ_2$.

These natural transformations allow us to lift the theory of linear
differential operators from $A$ to $\CF$ and to restrict the lifted theory
to $\CF_\De$. They are in parallel to the theory of $\CC$-differential
operators (see the next section).

\label{sec2:p:ulin}
The natural embeddings
$$\Diff_{(k)}^{\sym}(P,R)\hra\Diff_{(k-1)}^{\sym}(P,\Diff(P,R))$$
generate the map $\ell\colon \CF\ot_AR\to
\CF\ot_A\Diff(P,R)$, $\vf\mapsto\ell_\vf$, which is called the
\emph{universal linearization}. Using this map, we can rewrite Corollary
\ref{sec2:co:eulag}~(1) on page \pageref{sec2:co:eulag} in the form
$\ell_\De=\ell_\De^\circ$ while the
Euler operator is written as $\mathbf{E}(L)=\ell_L^\circ(1)$. Note also
that $\ell_{\vf\psi}=\vf\ell_\psi+\psi\ell_\vf$ for any $\vf,\psi\in
\CF\ot_AR$.
\begin{definition}
The group of \emph{conservation laws} for the algebra $\CF_\De$ (or for the
operator $\De$) is the first homology group of the complex of integral
forms
\begin{equation}
0\xla{}\CF_\De\ot_A\CB\xla{}\CF_\De\ot_A\Sigma_1\xla{}
\CF_\De\ot_A\Sigma_2\xla{}\dotsb
\end{equation}
with coefficients in $\CF_\De$.
\end{definition}

\newpage

\section{Jets and nonlinear differential equations. Symmetries}
\label{sec:geom}
We expose here main facts concerning geometrical approach to jets (finite
and infinite) and to nonlinear differential operators. We shall confine
ourselves with the case of vector bundles, though all constructions below
can be carried out---with natural modifications---for an arbitrary
locally trivial bundle $\pi$ (and even in more general settings). For further
reading, the books \cite{KLV,Symm} together with the paper \cite{Vin7} are
recommended.

\subsection{Finite jets}\label{sub:fjets}
Let $M$ be an $n$-dimensional smooth, i.e., of the class $\Ci$, manifold and
$\pi\colon E\to M$ be
a smooth $m$-dimensional vector bundle over $M$. Denote by $\Ga(\pi)$
the $\Ci(M)$-module of sections of the bundle $\pi$. For any point
$x\in M$ we shall also consider the module $\Gal(\pi;x)$ of all {\em
local} sections at $x$.

For a section $\vf\in\Gal(\pi;x)$ satisfying $\vf(x)=\ta\in E$, consider its
graph $\Ga_\vf\sbs E$ and all sections $\vf'\in\Gal(\pi;x)$ such that
\begin{enumerate}
\item[(a)]$\vf(x)=\vf'(x)$;
\item[(b)]the graph $\Ga_{\vf'}$ is tangent to $\Ga_{\vf}$ with order
$k$ at $\ta$.
\end{enumerate}
Conditions (a) and (b) determine equivalence relation $\sim_x^k$ on
$\Gal(\pi;x)$ and we denote the equivalence class of $\vf$ by
$[\vf]_x^k$. The quotient set $\Gal(\pi;x)/\!\sim_x^k$ becomes an
$\BBR$-vector space, if we put
$$[\vf]_x^k+[\psi]_x^k=[\vf+\psi]_x^k,\ a[\vf]_x^k=[a\vf]_x^k,\quad
\vf,\psi\in\Gal(\pi;x),\quad a\in\BBR,$$
while the natural projection $\Gal(\pi;x)\to\Gal(\pi;x)/\!\sim_x^k$
becomes a linear map. We denote this quotient space by $J_x^k(\pi)$.
Obviously, $J_x^0(\pi)$ coincides with $E_x=\pi^{-1}(x)$.

The tangency class $[\vf]_x^k$ is completely
determined by the point $x$ and partial derivatives at $x$ of the section
$\vf$ up to order $k$. From here it follows that $J_x^k(\pi)$ is
finite-dimensional with
\begin{equation}\label{sec3:eq:dim}
\dim J_x^k(\pi)=m\sum_{i=0}^k\binom{n+i-1}{n-1}=m\binom{n+k}{k}.
\end{equation}

\bde
The element $[\vf]_x^k\in J_x^k(\pi)$ is called the
{\em $k$-jet} of the section $\vf\in\Gal(\pi;x)$ at the point $x$.
\index{jet of section at a point}
\ede

The $k$-jet of $\vf$ at $x$ can be identified with the $k$-th order Taylor
expansion of the section $\vf$. From the definition it follows that it
is independent of coordinate choice.

Consider now the set
\begin{equation}\label{sec3:eq:jetm}
J^k(\pi)=\bigcup_{x\in M} J_x^k(\pi)
\end{equation}
and introduce a smooth manifold structure on $J^k(\pi)$ in the following
way. Let $\left\{\CU_\al\right\}_\al$ be an atlas in $M$ such that the
bundle $\pi$ becomes trivial over each $\CU_\al$, i.e., $\pi^{-1}(\CU_\al)
\simeq\CU_\al\times F$, where $F$ is the ``typical fiber''. Choose a
basis $e_1^\al,\dots,e_m^\al$ of local sections of $\pi$ over $\CU_\al$.
Then any section of $\pi\left|_{\CU_\al}\right.$ is representable in
the form $\vf=u^1e_1^\al+\dots+u^me_m^\al$ and the functions $x_1,\dots,
x_n,u^1,\dots,u^m$, where $x_1,\dots,x_n$ are local coordinates in
$\CU_\al$, constitute a local coordinate system in $\pi^{-1}(\CU_\al)$.
Let us define the functions $u_\si^j\colon
\bigcup_{x\in\CU_\al}J_x^k(\pi)\to\BBR$, where
$\si=i_1\dotsc i_r$, $|\si|=r\le k$, by
\begin{equation}\label{sec3:eq:coor}
u_\si^j([\vf]_x^k)=
\left.\dd{^{\abs{\si}} u^j}{x_{i_1}\dotsm\partial x_{i_r}}\right|_x.
\end{equation}
Then these functions, together with local coordinates $x_1,\dots,x_n$,
determine the map $f_\al\colon \bigcup_{x\in\CU_\al}J_x^k(\pi)\to\CU_\al
\times\BBR^N$, where $N$ is the number defined by \veqref{sec3:eq:dim}.
Due to computation rules for partial derivatives under coordinate
transformations, the map
$$(f_\al\circ f_\be^{-1})\left|_{\CU_\al\cap\CU_\be}\right.\colon
\left(\CU_\al\cap\CU_\be\right)\times\BBR^N\to
\left(\CU_\al\cap\CU_\be\right)\times\BBR^N$$
is a diffeomorphism preserving the natural projection
$\left(\CU_\al\cap\CU_\be\right)\times\BBR^N\to
\left(\CU_\al\cap\CU_\be\right)$. Thus we have proved the following
result:
\bpr\label{sec3:pr:smooth}
The set $J^k(\pi)$ defined by \eqref{sec3:eq:jetm} is a
smooth manifold while the projection $\pi_k\colon J^k(\pi)\to M$,
$\pi_k\colon [\vf]_x^k\mapsto x$, is a smooth vector bundle.
\epr
\bde\label{sec3:df:jets}
Let $\pi\colon E\to M$ be a smooth vector bundle, $\dim M=n$, $\dim E=n+m$.
\begin{enumerate}
\item The manifold $J^k(\pi)$ is called the {\em manifold of $k$-jets}
for $\pi$;\index{manifold of $k$-jets}
\item The bundle $\pi_k\colon J^k(\pi)\to M$ is called the {\em bundle of
$k$-jets} for $\pi$;\index{bundle of $k$-jets}
\item The above constructed local coordinates
$\left\{x_i,u_\si^j\right\}$, $i=1,\dots,n$, $j=1,\dots,m$,
$|\si|\le k$, are called the {\em special coordinate system} on
$J^k(\pi)$ associated to the trivialization
$\left\{\CU_\al\right\}_\al$ of the bundle $\pi$.\index{special
coordinate system}
\end{enumerate}
\ede
Obviously, the bundle $\pi_0$
coincides with $\pi$.

Since tangency of two manifolds with order $k$ implies
tangency with less order, there exists a map
$$\pi_{k,l}\colon J^k(\pi)\to J^l(\pi),\quad
[\vf]_x^k\mapsto[\vf]_x^l,\quad k\ge l,$$
which is a smooth fiber
bundle. If $k\ge l\ge s$, then obviously
\begin{equation}\label{sec3:eq:comm}
\pi_{l,s}\circ\pi_{k,l}=\pi_{k,s},\quad\pi_l\circ\pi_{k,l}=\pi_k.
\end{equation}
On the other hand, for any section $\vf\in\Ga(\pi)$ (or $\in\Gal(\pi;x)$)
we can define the map $j_k(\vf)\colon M\to J^k(\pi)$ by setting
$j_k(\vf)(x)=[\vf]_x^k$. Obviously, $j_k(\vf)\in\Ga(\pi_k)$ (respectively,
$j_k(\vf)\in\Gal(\pi_k;x)$).
\bde\label{sec3:df:jk}
The section $j_k(\vf)$ is called the {\em $k$-jet of the
section} $\vf$.\index{$k$-jet of the section} The correspondence
$j_k\colon \Ga(\pi)\to\Ga(\pi_k)$ is called the {\em $k$-jet operator}.
\index{operator of $k$-jet}
\ede
From the definition it follows that
\begin{equation}\label{sec3:eq:comm1}
\pi_{k,l}\circ j_k(\vf)=j_l(\vf),\quad j_0(\vf)=\vf,\quad k\ge l,
\end{equation}
for any $\vf\in\Ga(\pi)$.

Let $\vf,\psi\in\Ga(\pi)$ be two sections, $x\in M$ and $\vf(x)=\psi(x)=\ta
\in E$. It is a tautology to say that the manifolds $\Ga_\vf$ and $\Ga_\psi$
are tangent to each other with order $k+l$ at $\ta$ or that the manifolds
$\Ga_{j_k(\vf)},\Ga_{j_k(\psi)}\sbs J^k(\pi)$ are tangent with order $l$
at the point $\ta_k=j_k(\vf)(x)=j_k(\psi)(x)$.
\bde\label{sec3:df:Rpl}
Let $\ta_k\in J^k(\pi)$. An {\em $R$-plane}\index{$R$-plane}
at $\ta_k$ is an $n$-dimen\-si\-on\-al plane tangent to a manifold of
the form $\Ga_{j_k(\vf)}$ such that $[\vf]_x^k=\ta_k$.
\ede
Immediately from definitions we obtain the following result.
\bpr\label{sec3:pr:fiber}
Consider a point $\ta_k\in J^k(\pi)$. Then the fiber of the bundle $\pi_{k+1,k}\colon
J^{k+1}(\pi)\to J^k(\pi)$ over $\ta_k$ coincides with the set of all
$R$-planes at $\ta_k$.
\epr
For $\ta_{k+1}\in J^{k+1}(\pi)$, we shall denote the corresponding
$R$-plane at $\ta_k=\pi_{k+1,k}(\ta_{k+1})$ by $L_{\ta_{k+1}}\sbs
T_{\ta_k}(J^k(\pi))$.

\subsection{Nonlinear differential operators}\label{sub:ndo}
Let us consider now the algebra of smooth functions on $J^k(\pi)$ and denote it by
$\CF_k=\CF_k(\pi)$. Take another vector bundle $\pi'\colon E'\to M$ and consider
the pull-back $\pi_k^*(\pi')$. Then the set of sections of $\pi_k^*(\pi')$ is
a module over $\CF_k(\pi)$ and we denote this module by $\CF_k(\pi,\pi')$. In
particular, $\CF_k(\pi)=\CF_k(\pi,\mathbf{1}_M)$, where $\mathbf{1}_M$ is the
trivial one-dimensional bundle over $M$.

The surjections $\pi_{k,l}$ and $\pi_k$ for all $k\ge l\ge0$  generate
the natural embeddings
$\nu_{k,l}=\pi_{k,l}^*\colon \CF_l(\pi,\pi')\to\CF_k(\pi,\pi')$ and
$\nu_k=\pi_k^*\colon \Ga(\pi')\to\CF_k(\pi,\pi')$. Due to \veqref{sec3:eq:comm},
we have the equalities
\begin{equation}\label{sec3:eq:comm3}
\nu_{k,l}\circ\nu_{l,s}=\nu_{k,s},\quad\nu_{k,l}\circ\nu_l=\nu_k,\quad
k\ge l\ge s.
\end{equation}
Identifying $\CF_l(\pi,\pi')$ with its image in $\CF_k(\pi,\pi')$
under $\nu_{k,l}$, we can consider $\CF_k(\pi,\pi')$ as a filtered
module,
\begin{equation}\label{sec3:eq:filt}
\Ga(\pi')\hra\CF_0(\pi,\pi')\hra\dotsb\hra\CF_{k-1}(\pi,\pi')\hra
\CF_k(\pi,\pi'),
\end{equation}
over the filtered algebra $\Ci(M)\hra\CF_0\hra\dotsb\hra\CF_{k-1}\hra\CF_k$
with the embeddings $\CF_k\cdot\CF_l(\pi,\pi')\sbs\CF_{\max(k\,,\,l)}
(\pi,\pi')$. Let $F\in\CF_k(\pi,\pi')$. Then we have the correspondence
\begin{equation}\label{sec3:eq:ndo}
\De=\De_F\colon \Ga(\pi)\to\Ga(\pi'),\quad
\De(\vf)= j_k(\vf)^*(F),\quad\vf\in\Ga(\pi).
\end{equation}
\bde\label{sec3:df:ndo}
A correspondence $\De$ of the form \veqref{sec3:eq:ndo} is called a ({\em
nonlinear}) {\em differential operator}\index{nonlinear differential
operator} of order $\le k$ acting from the bundle $\pi$ to the bundle
$\pi'$. In particular, when $\De(a\vf+b\psi)=a\De(\vf)+b\De(\psi)$, $a,b\in
\mathbb{R}$, the operator $\De$ is said to be {\em linear}.\index{linear
differential operator}
\ede
\begin{example}\label{sec3:ex:jk}
Let us show that the $k$-jet operator $j_k\colon \Ga(\pi)\to
\Ga(\pi_k)$ (Definition \vref{sec3:df:jk}) is differential. To do this,
recall that the total space of the pull-back $\pi_k^*(\pi_k)$ consists
of points $(\ta_k,\ta'_k)\in J^k(\pi)\times J^k(\pi)$ such that
$\pi_k(\ta_k)=\pi_k(\ta'_k)$. Consequently, we may define the {\em diagonal
section} $\rho_k\in\CF_k(\pi,\pi_k)$ of the bundle $\pi_k^*(\pi_k)$ by
setting $\rho_k(\ta_k)=\ta_k$. Obviously, $j_k=\De_{\rho_k}$, i.e.,
$$j_k(\vf)^*(\rho_k)=j_k(\vf),\quad\vf\in\Ga(\pi).$$
The operator $j_k$ is linear.
\end{example}
\begin{example}\label{sec3:ex:deRham}
Let $\tau^*\colon T^*M\to M$ be the cotangent bundle of $M$
and $\tau_p^*\colon \bigwedge^pT^*M\to M$ be its $p$-th exterior power. Then
the {\em de Rham differential}\index{de Rham differential} $d$ is a
first order linear differential operator acting from $\tau_p^*$ to
$\tau_{p+1}^*,\ p\ge0$.
\end{example}

Let us prove now that composition of nonlinear differential operators is a
differential operator again. Let $\De\colon \Ga(\pi)\to\Ga(\pi')$ be a differential
operator of order $\le k$. For any $\ta_k=[\vf]_x^k\in J^k(\pi)$, set
\begin{equation}\label{sec3:eq:Psi}
\Phi_\De(\ta_k)=
[\De(\vf)]_x^0=(\De(\vf))(x).
\end{equation}
Evidently, the map $\Phi_\De$ is a morphism of fiber bundles (but not of
vector bundles!), i.e., $\pi'\circ\Phi_\De=\pi_k$.
\bde\label{sec3:df:Phi}
The map $\Phi_\De$ is called the {\em representative
morphism}\index{representative morphism} of the operator $\De$.
\ede
For example, for $\De=j_k$ we have $\Phi_{j_k}=\id_{J^k(\pi)}$.
Note that there exists a one-to-one correspondence between nonlinear
differential operators and their representative morphisms: one
can easily see it just by inverting equality \veqref{sec3:eq:Psi}. In fact,
if $\Phi\colon J^k(\pi)\to E'$ is a morphism of $\pi$ to $\pi'$, a
section $\vf\in\CF(\pi,\pi')$ can be defined by setting $\vf(\ta_k)=
(\ta_k,\Phi(\ta_k))\in J^k(\pi)\times E'$. Then, obviously, $\Phi$
is the representative morphism for $\De=\De_\vf$.

\bde\label{sec3:df:prol}
Let $\De\colon \Ga(\pi)\to\Ga(\pi')$ be a $k$-th order differential operator. Its
{\em $l$-th prolongation}\index{prolongation of differential operator} is the
composition $\De^{(l)}=j_l\circ\De\colon \Ga(\pi)\to\Ga(\pi_l)$.
\ede
\ble\label{sec3:le:prol}
For any $k$-th order differential operator $\De$,
its $l$-th prolongation is a $(k+l)$-th order operator.
\ele
\begin{proof}
In fact, for any $\ta_{k+l}=[\vf]_x^{k+l}\in J^{k+l}(\pi)$ set
$\Phi_\De^{(l)}(\ta_{k+l})=[\De(\vf)]_x^l\in J^l(\pi)$. Then the
operator, for which the morphism $\Phi_\De^{(l)}$ is representative,
coincides with $\De^{(l)}$.
\end{proof}
\bco\label{sec3:co:comp}
The composition $\De'\circ\De$ of two nonlinear differential operators
$\De\colon \Ga(\pi)\to
\Ga(\pi')$ and $\De'\colon \Ga(\pi')\to\Ga(\pi'')$ of orders $\le k$ and
$\le k'$ respectively is a $(k+k')$-th order differential operator.
\eco
\begin{proof}
Let $\Phi_\De^{(k')}\colon J^{k+k'}(\pi)\to J^{k'}(\pi')$ be the representative
morphism for $\De^{(k')}$.  Then the operator $\square$, for which the
composition $\Phi_{\De'}\circ\Phi_\De^{(k')}$ is the representative
morphism, coincides with $\De'\circ\De$.
\end{proof}

The following obvious proposition describes main properties of prolongations
and representative morphisms.
\bpr\label{sec3:pr:propprol}
Let $\De\colon \Ga(\pi)\to\Ga(\pi')$ and $\De'\colon \Ga(\pi')\to\Ga(\pi'')$
be two differential operators of orders $k$ and $k'$ respectively.
Then:
\begin{enumerate}
\item $\Phi_{\De'\circ\De}=\Phi_{\De'}\circ\Phi_\De^{(k')}$,
\item $\Phi_\De^{(l)}\circ j_{k+l}(\vf)=\De^{(l)}(\vf)$ for
any $\vf\in\Ga(\pi)$ and $l\ge0$,
\item $\pi_{l,l'}\circ\Phi_\De^{(l)}=\Phi_\De^{(l')}\circ
\pi_{k+l,k+l'}$, i.e., the diagram
\begin{equation}\label{sec3:eq:propprol}
\begin{CD}
J^{k+l}(\pi)@>{\Phi_\De^{(l)}}>>J^l(\pi')\\
@V{\pi_{k+l,k+l'}}VV@VV{\pi'_{l,l'}}V\\
J^{k+l'}(\pi)@>{\Phi_\De^{(l')}}>>J^{l'}(\pi')
\end{CD}
\end{equation}
is commutative for all $l\ge l'\ge0$.
\end{enumerate}
\epr

\subsection{Infinite jets}\label{sub:ijets}
We now pass to infinite limit in all previous constructions.
\bde\label{sec3:df:ijets}
The {\em space of infinite jets}\index{space of infinite
jets} $\Ji(\pi)$ of the fiber bundle $\pi\colon E\to M$ is the inverse limit
of the sequence
$$\dotsb\xra{} J^{k+1}(\pi)\xra{\pi_{k+1,k}}J^k(\pi)\xra{}\dotsb\xra{}
J^1(\pi) \xra{\pi_{1,0}}E\xra{\pi}M, $$ i.e.,
$\Ji(\pi)=\projlim_{\{\pi_{k,l}\}}J^k(\pi)$.
\ede

Thus a \emph{point} $\ta$ of $\Ji(\pi)$ is a sequence of points
$\{\ta_k\}_{k\ge0}$, $\ta_k\in J^k(\pi)$, such that
$\pi_{k,l}(\ta_k)=\ta_l,\ k\ge l$. Points of
$\Ji(\pi)$ can be understood as $m$-dimensional formal series and can be
represented in the form $\ta=[\vf]_x^\infty,\ \vf\in\Gal(\pi)$.

A \emph{special coordinate system} associated to a trivialization
$\{\CU_\al\}_\al$ is given by the functions $x_1,\dots,x_n,
\dots,u_\si^j,\dotsc$.

A \emph{tangent vector} to $\Ji(\pi)$ at a point $\ta$ is defined as
a system of vectors $\{w,v_k\}_{k\ge0}$ tangent to $M$ and to
$J^k(\pi)$ respectively such that $(\pi_k)_*v_k=w$, $(\pi_{k,l})_*v_k=v_l$
for all $k\ge l\ge0$.

A \emph{smooth bundle} $\xi$ over $\Ji(\pi)$ is a system of bundles $\eta
\colon Q\to M$, $\xi_k\colon P_k\to J^k(\pi)$ together with smooth maps
$\Psi_k\colon P_k\to Q$, $\Psi_{k,l}\colon P_k\to P_l$, $k\ge l\ge0$, such
that
$$
\Psi_l\circ\Psi_{k,l}=\Psi_k,\quad\Psi_{k,l}\circ\Psi_{l,s}=\Psi_{k,s},\quad
k\ge l\ge s\ge0
$$
For example, if $\eta\colon Q\to M$ is a bundle, then the
pull-backs $\pi_k^*(\eta)\colon \pi_k^*(Q)\to J^k(\pi)$ together with
natural projections $\pi_k^*(Q)\to\pi_l^*(Q)$ and $\pi_k^*(Q)\to Q$
form a bundle over $\Ji(\pi)$. We say that $\xi$ is a {\em vector bundle}
over $\Ji(\pi)$, if $\eta$ and all $\xi_k$ are vector bundles while the
maps $\Psi_k$ and $\Psi_{k,l}$ are fiberwise linear.

A \emph{smooth map} of $\Ji(\pi)$ to $\Ji(\pi')$, where $\pi\colon E\to M$,
$\pi'\colon E'\to M'$, is defined as a system $F$ of maps $F_{-\infty}\colon
M\to M'$, $F_k\colon J^k(\pi)\to J^{k-s}(\pi')$, $k\ge s$, where $s\in\BBZ$
is a fixed integer called the {\em degree of $F$}, such that
$$
\pi_{k-r,k-s-1}\circ F_k=F_{k-1}\circ\pi_{k,k-1},\quad k\ge s+1,
$$
and
$$
\pi_{k-s}\circ F_k=F_{-\infty}\circ\pi_k,\quad k\ge s.
$$
For example, if $\De\colon\Ga(\pi)\to\Ga(\pi')$ is a differential operator
of order $s$, then the system of maps $F_{-\infty}=\id_M$, $F_k=
\Phi_\De^{(k-s)}$, $k\ge s$ (see the previous subsection), is a smooth map
of $\Ji(\pi)$ to $\Ji(\pi')$.

A \emph{smooth function}\index{smooth function on $\Ji(\pi)$}
on $\Ji(\pi)$ is an element of the direct limit $\CF=\CF(\pi)=
\injlim_{\{\pi_{k,l}^*\}}\CF_k(\pi)$, where $\CF_k(\pi)$ is the
algebra of smooth functions on $J^k(\pi)$. Thus, a smooth function on
$\Ji(\pi)$ is a function on
$J^k(\pi)$ for some finite but an arbitrary $k$. The set $\CF=\CF(\pi)$ of
such functions is identified with $\bigcup_{k=0}^\infty\CF_k(\pi)$ and forms
a commutative filtered algebra. Using duality between
smooth manifolds and algebras of smooth functions on these manifolds, we
deal in what follows with the algebra $\CF(\pi)$ rather than with the
manifold $\Ji(\pi)$ itself.

From this point of view,  a \emph{vector field}\index{vector field on
$\Ji(\pi)$} on $\Ji(\pi)$ is a filtered derivation of $\CF(\pi)$, i.e.,
an $\BBR$-linear map $X\colon \CF(\pi)\to\CF(\pi)$ such that
$$
X(fg)=fX(g)+gX(f),\quad f,g\in\CF(\pi),\quad X(\CF_k(\pi))\sbs\CF_{k+l}(\pi)
$$
for all $k$ and some $l=l(X)$. The latter is called the {\em filtration
degree}
of the field $X$. The set of all vector fields is a filtered Lie algebra
over $\BBR$ with respect to commutator $[X,Y]$ and is denoted by
$\Dr(\pi)=\bigcup_{l\ge0}\Dr^{(l)}(\pi)$.

\emph{Differential forms} of degree $i$ on $\Ji(\pi)$ are defined as
elements of the filtered $\CF(\pi)$-module $\La^i=\La^i(\pi)=
\bigcup_{k\ge0}{\La^i(\pi_k)}$, where $\La^i(\pi_k)=\La^i(J^k(\pi))$
and the module $\La^i(\pi_k)$ is considered to be embedded into
$\La^i(\pi_{k+1})$ by the map $\pi_{k+1,k}^*$. Defined in such a way, these forms
possess all basic properties
of differential forms on finite-dimensional manifolds. Let us mention the most
important ones:
\begin{enumerate}
\item The module $\La^i(\pi)$ is the $i$-th exterior power of
$\La^1(\pi)$, $\La^i(\pi)=\bigwedge^i\La^1(\pi)$. Respectively, the
operation of {\em wedge product} $\wg\colon \La^p(\pi)\ot\La^q(\pi)\to
\La^{p+q}(\pi)$ is defined and $\La^*(\pi)=\bigoplus_{i\ge0}\La^i(\pi)$
becomes a supercommutative $\mathbb{Z}$-graded algebra.
\item The module $\Dr(\pi)$ is dual to $\La^1(\pi)$, i.e.,
\begin{equation}\label{sec3:eq:gradrep}
\Dr(\pi)=\Hom_{\CF(\pi)}^{\phi}(\La^1(\pi),\CF(\pi)),
\end{equation}
where $\Hom_{\CF(\pi)}^{\phi}(\cdot,\cdot)$ denotes the module of filtered
homomorphisms over
$\CF(\pi)$. Moreover, equality \veqref{sec3:eq:gradrep} is established in the
following way: there is a derivation $d\colon \CF(\pi)\to\La^1(\pi)$
(the \emph{de Rham differential} on $\Ji(\pi)$) such
that for any vector field $X$ there exists a uniquely defined
filtered homomorphism $f_X$ satisfying $f_X\circ d=X$.
\item The operator $d $ is extended up to maps $d \colon \La^i(\pi)\to
\La^{i+1}(\pi)$ in such a way that the sequence
$$
0\xra{}\CF(\pi)\xra{d}\La^1(\pi)\xra{}\dotsb\xra{}\La^i(\pi)\xra{d}
\La^{i+1}(\pi) \xra{}\dotsb
$$
becomes a complex, i.e., $d\circ d=0$. This
complex is called the {\em de Rham complex} on $\Ji(\pi)$. The latter is a derivation of
the superalgebra $\La^*(\pi)$.
\end{enumerate}

Using algebraic techniques (see Section \ref{sec:calc}), we can introduce the
notions of \emph{inner product} and \emph{Lie derivative} and to prove
their basic properties (cf.\ Proposition \vref{sec1:pr:lie}). We can also
define \emph{linear differential operators} over $\Ji(\pi)$ as follows. Let
$P$ and $Q$ be two filtered
$\CF(\pi)$-modules and $\De\in\Hom_{\CF(\pi)}^\phi(P,Q)$. Then $\De$ is
called a linear differential operator of order $\le k$ acting from $P$
to $Q$, if
$$(\de_{f_0}\circ\de_{f_1}\circ\dots\circ\de_{f_k})\De=0$$
for all $f_0,\dots,f_k\in\CF(\pi)$, where, as in Section \ref{sec:calc},
$(\de_f\De)p=\De(fp)-f\De(p)$. We write $k=\ord(\De)$.

Due to existence of filtrations in the algebra $\CF(\pi)$, as well as in
modules $P$ and $Q$, one can define
differential operators of \emph{infinite order}\index{differential
operators of infinite order} acting from $P$ to $Q$.
Namely, let $P=\{P_l\}$, $Q=\{Q_l\}$, $P_l\sbs P_{l+1}$, $Q_l\sbs Q_{l+1}$,
$P_l$ and $Q_l$ being $\CF_l(\pi)$-modules. Let $\De\in\Hom_{\CF(\pi)}^\phi(P,Q)$
and $s$ be filtration of $\De$, i.e., $\De(P_l)\sbs Q_{l+s}$. We can
always assume that $s\ge0$. Suppose now that $\De_l=\De
\left|_{P_l}\right.\colon P_l\to Q_l$ is a linear differential operator of
order $o_l$ over $\CF_l(\pi)$ for any $l$. Then we say that $\De$ is a linear
differential operator of order growth $o_l$. In particular,
if $o_{l}=\al l+\be$, $\al,\be\in\BBR$, we say that $\De$ is
of constant growth $\al$.

\emph{Distributions}. Let $\ta\in\Ji(\pi)$. The tangent plane to $\Ji(\pi)$
at the point $\ta$ is the set of all tangent vectors to $\Ji(\pi)$ at
this point (see above). Denote such a plane by $T_\ta=T_\ta(\Ji(\pi))$.
Let $\ta=\{x,\ta_k\}$, $x\in M$, $\ta_k\in J^k(\pi)$ and $v=\{w,v_k\}$, $v'=
\{w',v'_k\}\in T_\ta$. Then the linear combination $\la v+\mu v'=
\{\la w+\mu w',\la v_k+\mu v'_k\}$ is again an element of $T_\ta$ and thus
$T_\ta$ is a vector space. A correspondence $\CT\colon \ta\mapsto \CT_\ta\sbs
T_\ta$, where $\CT_\ta$ is a linear subspace, is called a \emph{distribution}
on $\Ji(\pi)$. Denote by $\CT \Dr(\pi)\sbs \Dr(\pi)$ the submodule of vector
fields lying in $\CT$, i.e., a vector field $X$ belongs to $\CT \Dr(\pi)$ if
and
only if $X_\ta\in\CT_\ta$ for all $\ta\in\Ji(\pi)$. We say that the
distribution $\CT$ is {\em integrable}, if it satisfies the formal Frobenius
condition: for any vector fields
$X,Y\in\CT \Dr(\pi)$ their commutator lies in $\CT \Dr(\pi)$
as well, or $[\CT \Dr(\pi),\CT \Dr(\pi)]\sbs\CT \Dr(\pi)$.

This condition can expressed in a dual way as follows. Let us set
$$
\CT^1\La(\pi)=\{\,\om\in\La^1(\pi)\mid \ip_X\om=0,\ X\in\CT \Dr(\pi)\,\}
$$
and consider the ideal $\CT\La^*(\pi)$ generated in $\La^*(\pi)$ by
$\CT^1\La(\pi)$. Then the distribution $\CT$ is integrable if and only if
the ideal $\CT\La^*(\pi)$ is differentially closed: $d (\CT\La^*(\pi))\sbs
\CT\La^*(\pi)$.

Finally, we say that a submanifold $N\sbs\Ji(\pi)$ is an {\em integral
manifold} of $\CT$, if $T_\ta N\sbs\CT_\ta$ for any point $\ta\in N$. An
integral manifold $N$ is called {\em locally maximal} at a point $\ta\in
N$, if there no neighborhood $\CU\sbs N$ of $\ta$ is embedded to other
integral manifold $N'$ such that $\dim N\le \dim N'$.

\subsection{Nonlinear equations and their solutions}\label{sub:nde}
Let $\pi\colon E\to M$ be a vector bundle.
\bde\label{sec3:df:nde}
A submanifold $\CE\sbs J^k(\pi)$ is called a ({\em nonlinear})
{\em differential equation}\index{nonlinear differential equation} of order
$k$ in the bundle $\pi$. We say that $\CE$ is a {\em linear equation},
if $\CE\cap\pi_x^{-1}(x)$ is a linear subspace in $\pi_x^{-1}(x)$ for all
$x\in M$. In other words, $\CE$ is a linear subbundle in the bundle $\pi_k$.
\ede
We shall always assume that $\CE$ is projected surjectively to $E$ under
$\pi_{k,0}$.

\bde\label{sec3:df:sol}
A (local) section $f$ of the bundle $\pi$ is called a (local)
{\em solution}\index{solution} of the equation $\CE$, if its graph lies in
$\CE$: $j_k(f)(M)\sbs\CE$.
\ede
We say that the equation $\CE$ is {\em determined}, if $\codim\,\CE=\dim\pi$,
that it is {\em overdetermined}, if $\codim\,\CE>\dim\pi$, and that it is
{\em underdetermined}, if $\codim\,\CE<\dim\pi$.

Obviously, in a special coordinate system these definitions coincide with
``usual'' ones.

One of the ways to represent differential equations is as follows. Let
$\pi'\colon \BBR^r
\times\bar{\CU}\to\bar{\CU}$ be the trivial $r$-dimensional bundle. Then
the set of functions $(F^1,\dots,F^r)$ can be understood as a section $\vf$
of the pull-back $\left(\pi_k\left|_\CU\right.\right)^*(\pi')$, or as
a nonlinear operator $\De=\De_\vf$ defined in $\CU$, while the equation
$\CE$ is characterized by the condition
\begin{equation}\label{sec3:eq:eqop}
\CE\cap\CU=\{\,\ta_k\in\CU\mid \vf(\ta_k)=0\,\}.
\end{equation}
More generally, any equation $\CE\sbs J^k(\pi)$ can be represented in the form
similar to \veqref{sec3:eq:eqop}. Namely, for any equation $\CE$ there exists
a fiber bundle $\pi'\colon E'\to M$ and a section $\vf\in\CF_k(\pi,\pi)$ such
that $\CE$ coincides with the set of zeroes for $\vf$: $\CE=\{\vf=0\}$.
In this case we say that $\CE$ is {\em associated to the operator}
$\De=\De_\vf\colon \Ga(\pi)\to\Ga(\pi')$\index{equation associated to operator}
and use the notation $\CE=\CE_\De$.

\begin{example}\label{sec3:ex1:deRhameq}
Let $\pi=\tau_p^*\colon\bigwedge^pT^*M\to M$, $\pi'=\tau_{p+1}^*\colon
\bigwedge^{p+1}
T^*M\to M$ and $d\colon \Ga(\pi)=\La^p(M)\to
\Ga(\pi')=\La^{p+1}(M)$ be the de Rham differential (see Example
\vref{sec3:ex:deRham}). Thus we obtain a first-order equation $\CE_{d }$ in the
bundle $\tau_p^*$. Consider the case $p=1,n\ge2$ and choose local
coordinates $x_1,\dots,x_n$ in $M$. Then any form $\om\in\La^1(M)$ is
represented as $\om=u^1d  x_1+\dots+u^nd  x_n$ and we have
$\CE_{d }=\{\,u^j_i=u^i_j\mid i<j\,\}$. This equation is underdetermined
when $n=2$, determined for $n=3$ and overdetermined for $n>3$.
\end{example}
\begin{example}\label{sec3:ex:monge}
Consider an arbitrary vector bundle $\pi\colon E\to M$ and a differential form
$\om\in\La^p(J^k(\pi))$, $p\le\dim M$. Then the condition $j_k(\vf)^*
(\om)=0,\ \vf\in\Ga(\pi)$, determines a $(k+1)$-st order equation $\CE_\om$
in the bundle $\pi$. Consider the case $p=\dim M=2$, $k=1$ and choose a special
coordinate system $x,y,u,u_x,u_y$ in $J^1(\pi)$. Let $\vf=\vf(x,y)$ be a
local section and
\begin{multline*}
\om=Ad  u_x\wg d  u_y+(B_1d  u_x+B_2 d  u_y)\wg d  u+
d  u_x\wg(B_{11}d  x+B_{12}d  y)\\
+d  u_y\wg(B_{21}d  x+B_{22}d  y)+d  u\wg(C_1d  x+C_2d  y)+Dd  x\wg d  y,
\end{multline*}
where $A,B_i,B_{ij},C_i,D$ are functions of $x,y,u,u_x,u_y$. Then we have
\begin{multline*}
j_1(\vf)^*\om=\Bigl(A^\vf(\vf_{xx}\vf_{yy}-\vf_{xy}^2)+
(\vf_yB_1^\vf+B_{12}^\vf)\vf_{xx}\\
-(\vf_xB_2^\vf+B_{12}^\vf)\vf_{yy}+
(\vf_yB_2^\vf-\vf_xB_1^\vf+B_{22}^\vf-B_{11}^\vf)\vf_{xy}\\
+\vf_xC_2^\vf-\vf_yC_1^\vf+D^\vf)\Bigr)d  x\wg dy,
\end{multline*}
where $F^\vf= j_1(\vf)^*F$ for any $F\in\CF_1(\pi)$.
Hence, the equation $\CE_\om$ is of the form
\begin{equation}\label{sec3:eq:monge}
a(u_{xx}u_{yy}-u_{xy}^2)+b_{11}u_{xx}+b_{12}u_{xy}+b_{22}u_{yy}+c=0,
\end{equation}
where $a=A,b_{11}=u_yB_1+B_{12}$, $b_{12}=u_yB_2-u_xB_1+B_{22}-B_{11}$,
$b_{22}=u_xB_2+B_{12}$, $c=u_xC_2-u_yC_1+D$ are functions on $J^1(\pi)$.
Equation \veqref{sec3:eq:monge} is the so-called two-dimensional {\em
Monge--Ampere equation}
\index{Monge--Ampere equation} and obviously any such an equation
can be represented as $\CE_\om$ for some $\om\in\La^1(J^1(\pi))$
(see \cite{LychaginRub} for more details).
\end{example}
\begin{example}\label{sec3:ex:conn}
Consider again a bundle $\pi\colon E\to M$ and a section $\na\colon E\to J^1(\pi)$ of the
bundle $\pi_{1,0}\colon J^1(\pi)\to E$. Then the graph
$\CE_\na=\na(E)\sbs J^1(\pi)$ is a first-order equation in the bundle
$\pi$. Let $\ta_1\in\CE_\na$. Then, due to Proposition \vref{sec3:pr:fiber},
the point $\ta_1$ is identified with the pair $(\ta_0,L_{\ta_1})$, where $\ta_0=
\pi_{1,0}(\ta_1)\in E$, while $L_{\ta_1}$ is the $R$-plane at $\ta_0$
corresponding to $\ta_1$. Hence, the section $\na$ (or the equation
$\CE_\na$) may be understood as a distribution of horizontal
(i.e., nondegenerately projected to
$T_xM$ under $(\pi_k)_*$, where $x=\pi_k(\ta_k)$)
$n$-dimensional planes on $E$: $\CT_\na\colon E\ni\ta\mapsto\ta_1=L_{\na(\ta)}$.
In other words, $\na$ is a \emph{connection} in the bundle $\pi$.
A solution of the equation $\CE_\na$, by definition, is a section $\vf\in
\Ga(\pi)$ such that $j_1(\vf)(M)\sbs\na(E)$. It means that at any point
$\ta=\vf(x)\in\vf(M)$ the plane $\CT_\na(\ta)$ is tangent to the graph
of the section $\vf$. Thus, solutions of $\CE_\na$ coincide with integral
manifolds of $\CT_\na$.

In a local coordinate system $(x_1,\dots,x_n,u^1,\dots,u^m,\dots,u_i^j,\dots)$,
$i=1,\dots,n$, $j=1,\dots,m$, the equation
$\CE_\na$ is represented as
\begin{equation}\label{sec3:eq:conn}
u_i^j=\na_i^j(x_1,\dots,x_n,u^1,\dots,u^m),\ i=1,\dots,n,j=1,\dots,m,
\end{equation}
$\na_i^j$ being smooth functions.
\end{example}
\begin{example}\label{sec3:ex:flat}
As we saw in the previous example, to solve the equation
$\CE_\na$ is the same as to find integral $n$-dimensional manifolds
of the distribution $\CT_\na$. Hence, the former to be solvable, the
latter is to satisfy the Frobenius theorem. Thus, for solvable $\CE_\na$
we obtain conditions on the section $\na\in\Ga(\pi_{1,0})$. Let write
down these conditions in local coordinates.

Using representation \veqref{sec3:eq:conn}, note that $\CT_\na$ is given by
1-forms
$$\om^j=d  u^j-\sum_{i=1}^n\na_i^jd  x_i,\ j=1,\dots,m.$$
Hence, the integrability conditions may be expressed as
$$
d \om^j=\sum_{i=1}^m\rho_i^j\wg\om_i,\ j=1,\dots,m,
$$
for some 1-forms $\rho_i^i$. After elementary computations, we obtain
that the functions $\na_i^j$ must satisfy the following relations:
\begin{equation}\label{sec3:eq:flat}
\pdr{\na_\al^j}{x_\be}+\sum_{\ga=1}^m\na_\al^\ga\pdr{\na_\be^j}{u^\ga}=
\pdr{\na_\be^j}{x_\al}+\sum_{\ga=1}^m\na_\be^\ga\pdr{\na_\al^j}{u^\ga}
\end{equation}
for all $j=1,\dots,m$, $1\le\al<\be\le m$. Thus we got a naturally
constructed first-order equation $\CI(\pi)\sbs J^1(\pi_{1,0})$, whose
solutions are horizontal $n$-dimensional distributions in $E=J^0(\pi)$.
\end{example}

\subsection{Cartan distribution on $J^k(\pi)$}\label{sub:cardis}
We shall now introduce a very important structure on $J^k(\pi)$ responsible
for ``individuality'' of these manifolds.
\bde\label{sec3:df:cardis}
Let $\pi\colon E\to M$ be a vector bundle. Consider a point $\ta_k\in J^k(\pi)$ and
the span $\CC_{\ta_k}^k\sbs T_{\ta_k}(J^k(\pi))$ of all $R$-planes at the
point $\ta_k$.
\begin{enumerate}
\item The correspondence $\CC^k=\CC^k(\pi)\colon \ta_k\mapsto\CC_{\ta_k}^k$
is called the {\em Cartan distribution} on $J^k(\pi)$.\index{Cartan
distribution on $J^k(\pi)$}
\item Let $\CE\sbs J^k(\pi)$ be a differential equation. Then the
correspondence $\CC^k(\CE)\colon \CE\ni\ta_k\mapsto\CC_{\ta_k}^k\cap
T_{\ta_k}\CE\sbs T_{\ta_k}\CE$ is called the {\em Cartan distribution} on
$\CE$.\index{Cartan distribution on $\CE$}
We call elements of the Cartan distributions {\em Cartan
planes}.\index{Cartan plane}
\item A point $\ta_k\in\CE$ is called {\em regular},\index{regular
point} if the Cartan plane $\CC_{\ta_k}^k(\CE)$ is of maximal dimension.
We say that $\CE$ is a {\em regular equation}\index{regular equation},
if all its points are regular.
\end{enumerate}
\ede
In what follows, we deal with regular equations or with neighborhoods of
regular points. As it can be easily seen, for any regular point there
exists a neighborhood of this point all points of which are regular.

Let $\ta_k\in J^k(\pi)$ be represented in the form
\begin{equation}\label{sec3:eq:point}
\ta_k=[\vf]_x^k,\quad\vf\in\Ga(\pi),\quad x=\pi_k(\ta_k).
\end{equation}
Then, by definition, the Cartan plane $\CC_{\ta_k}^k$ is spanned by the
vectors
\begin{equation}\label{sec3:eq:vect}
j_k(\vf)_{*,x}(v),\quad v\in T_xM,
\end{equation}
for all $\vf\in\Gal(\pi)$ satisfying \veqref{sec3:eq:point}.

Let $x_1,\dots,x_n,\dots,u_\si^j,\dots$, $j=1,\dots,m$, $|\si|\le k$, be
a special coordinate system in a neighborhood of $\ta_k$. The vectors of
the form \veqref{sec3:eq:vect} can be expressed as linear combinations
of the vectors
\begin{equation}\label{sec3:eq:vectco}
\pdr{}{x_i}+\sum_{|\si|\le k}\sum_{j=1}^m
\frac{\pat^{|\si|+1}\vf^j}{\pat x_\si\pat x_i}\pdr{}{u_\si^j},
\end{equation}
where $i=1,\dots,n$. Using this representation, we prove the following
result:
\bpr\label{sec3:pr:cplane}
For any point $\ta_k\in J^k(\pi)$, $k\ge1$, the Cartan plane
$\CC_{\ta_k}^k$ is of the form $\CC_{\ta_k}^k=
(\pi_{k,k-1})_*^{-1}(L_{\ta_k})$, where $L_{\ta_k}$ is the $R$-plane at
the point $\pi_{k,k-1}(\ta_k)\in J^{k-1}(\pi)$ determined by the point
$\ta_k$.
\epr
\begin{proof}
Denote the vector \veqref{sec3:eq:vectco} by $v_i^{k,\vf}$. It is obvious that for
any two sections $\vf$ and $\vf'$ satisfying \veqref{sec3:eq:point} the difference
$v_i^{k,\vf}-v_i^{k,\vf'}$ is a $\pi_{k,k-1}$-vertical vector and any such
a vector can be obtained in this way. On the other hand, the vectors
$v_i^{k-1,\vf}$ do not depend on section $\vf$ satisfying
\veqref{sec3:eq:point} and form a
basis in the space $L_{\ta_k}$.
\end{proof}
\bre\label{sec3:re:notint}
From the result proved it follows that the Cartan
distribution on $J^k(\pi)$ can be locally considered as generated by the
vector fields
\begin{equation}\label{sec3:eq:cbasis}
D_i^{(k-1)}=\pdr{}{x_i}+\sum_{|\si|\le k-1}\sum_{j=1}^m
u_{\si i}^j\pdr{}{u_\si^j},\quad V_\tau^s=\pdr{}{u_\tau^s},\quad
|\tau|=k,s=1,\dots,m.
\end{equation}
From here, by direct computations, it follows that $[V_\tau^s,D_i^{(k-1)}]
=V_{\tau-i}^s$, where
$$
V_{\tau-i}^s=\left\{\begin{array}{ll}
V_{\tau'},&\mbox{ if }\tau=\tau'i,\\
0,&\mbox{ if }\tau\mbox{ does not contain }i.
\end{array}\right.
$$
But, as it follows from Proposition \vref{sec3:pr:cplane}, vector fields
$V_\si^j$ for $|\si|\le k$ do not lie in $\CC^k$. Thus, the Cartan
distribution on $J^k(\pi)$ is not integrable.
\ere

Introduce 1-forms in special coordinates on $J^{k+1}(\pi)$:
\begin{equation}\label{sec3:eq:cforms}
\om_\si^j=d  u_\si^j-\sum_{i=1}^nu_{\si i}^jd  x_i,
\end{equation}
where $j=1,\dots,m$, $|\si|<k$. From the representation \veqref{sec3:eq:cbasis}
we immediately obtain the following important property of the forms
introduced:
\bpr\label{sec3:pr:cforms}
The system of forms \veqref{sec3:eq:cforms} annihilates the Cartan
distribution on $J^k(\pi)$, i.e., a vector field $X$ lies in $\CC^k$ if and
only if $\ip_X\om_\si^j=0$ for all $j=1,\dots,m,|\si|<k$.
\epr
\bde\label{sec3:df:cforms}
The forms \veqref{sec3:eq:cforms} are called the {\em Cartan forms}
\index{Cartan forms} on $J^k(\pi)$ associated to the special coordinate
system $x_i,u_\si^j$.
\ede
Note that the $\CF_k(\pi)$-submodule generated in $\La^1(J^k(\pi)$ by
the forms \veqref{sec3:eq:cforms} is independent of the choice of coordinates.
\bde\label{sec3:df:cmodule}
The $\CF_k(\pi)$-submodule generated in $\La^1(J^k(\pi))$ by the Cartan forms
is called the {\em Cartan submodule}\index{Cartan submodule}. We denote this
submodule by $\CC\La^1(J^k(\pi))$.
\ede

We shall now describe maximal integral manifolds of the Cartan distribution
on $J^k(\pi)$. Let $N\sbs J^k(\pi)$ be an integral manifold of the Cartan
distribution. Then from Proposition \vref{sec3:pr:cforms} it follows that the
restriction of any Cartan form $\om$ to $N$ vanishes. Similarly, the
differential $d \om$ vanishes on $N$. Therefore, if vector fields $X,Y$ are
tangent to $N$, then $d\om\left|_N\right.(X,Y)=0$.
\bde\label{sec3:df:invo}
Let $\CC_{\ta_k}^k$ be the Cartan plane at $\ta\in J^(\pi)$.
\begin{enumerate}
\item Two vectors $v,w\in\CC_{\ta_k}^k$ are said to be {\em in
involution}\index{vectors in involution}, if $d \om_{\ta_k}(v,w)=0$
for any $\om\in\CC\La^1(J^k(\pi))$.
\item A subspace $W\sbs\CC_{\ta_k}^k$ is said to be {\em involutive},
\index{involutive subspace} if any two vectors $v,w\in W$ are in involution.
\item An involutive subspace is called {\em maximal}\index{maximal
involutive subspace}, if it cannot be embedded into any other involutive
subspace.
\end{enumerate}
\ede

Consider a point $\ta_k=[\vf]_x^k\in J^k(\pi)$. Then from Proposition
\vref{sec3:pr:cforms} it follows that the direct sum decomposition
$\CC_{\ta_k}^k=T_{\ta_k}^v\oplus T_{\ta_k}^\vf$ takes place, where
$T_{\ta_k}^v$
denotes the tangent plane to the fiber of
the projection $\pi_{k,k-1}$ passing through the point $\ta_k$, while
$T_{\ta_k}^\vf$ is the tangent plane to the graph of $j_k(\vf)$. Hence,
the involutiveness is sufficient to be checked for the following pairs
of vectors $v,w\in\CC_{\ta_k}^k$:
\begin{enumerate}
\item $v,w\in T_{\ta_k}^v$;
\item $v,w\in T_{\ta_k}^\vf$;
\item $v\in T_{\ta_k}^v$, $w\in T_{\ta_k}^\vf$.
\end{enumerate}

Note now that the tangent space $T_{\ta_k}^v$ is identified with the
tensor product $S^k(T_x^*)\ot E_x,\ x=\pi_k(\ta_k)\in M$, where $T_x^*$ is
the fiber of the cotangent bundle to $M$ at $x$, $E_x$ is
the fiber of the bundle $\pi$ at the same point while $S^k$ denotes the
$k$-th symmetric power. Then any vector $w\in T_xM$ determines
the map $\de_w\colon S^k(T_x^*)\ot E_x\to S^{k-1}(T_x^*)\ot E_x$ by
$$
\de_w(\rho_1\cdot\ldots\cdot\rho_k)\ot e=\sum_{i=1}^k
\rho_1\cdot\ldots\cdot\langle\rho_i,w\rangle\cdot\ldots\cdot\rho_k\ot e,
$$
where the dot ``$\,\cdot\,$'' denotes multiplication in $S^k(T_x^*)$,
$\rho_i\in T_x^*$, $e\in E_x$ while $\langle\cdot,\cdot\rangle$ is the natural
pairing between $T_x^*$ and $T_x$.
\bpr\label{sec3:pr:invo} Let
$v,w\in\CC_{\ta_k}^k$. Then:
\begin{enumerate}
\item All pairs $v,w\in T_{\ta_k}^v$ are in involution.
\item All pairs $v,w\in T_{\ta_k}^\vf$
are in involution too.  If $v\in T_{\ta_k}^v$ and $w\in T_{\ta_k}^\vf$,
then they are in involution if and only if $\de_{(\pi_k)_*(w)}v=0$.
\end{enumerate}
\epr
\begin{proof}
Note first that the involutiveness conditions are sufficient to check for the
Cartan forms \veqref{sec3:eq:cforms} only. The all three results follow from the
representation \veqref{sec3:eq:cbasis} by straightforward computations.
\end{proof}

Let $\ta_k\in J^k(\pi)$ and $F_{\ta_k}$ be the fiber of the bundle $\pi_{k,k-1}$
passing through the point $\ta_k$ while $H\sbs T_xM$ be a linear subspace.
Using the
linear structure, we identify the fiber $F_{\ta_k}$ of the bundle
$\pi_{k,k-1}$ with its tangent space and define the space
$$\Ann(H)=\{\,v\in F_{\ta_k}\mid \de_wv=0,\ \forall w\in H\,\}.$$
Then, as it follows from Proposition \vref{sec3:pr:invo}, the following
description of maximal involutive subspaces takes place:
\bco\label{sec3:co:invsum}
Let $\ta_k=[\vf]_x^k$, $\vf\in\Gal(\pi)$. Then any maximal
involutive subspace $V\sbs\CC_{\ta_k}^k(\pi)$ is of the form
$V=j_k(\vf)_*(H)\oplus\Ann(H)$ for some $H\sbs T_xM$.
\eco
If $V$ is a maximal involutive subspace, then the corresponding space
$H$ is obviously $\pi_{k,*}(V)$. We call the dimension of $H$ the {\em type}
of the maximal involutive subspace $V$\index{type of maximal involutive
subspace} and denote it by $\tp(V)$.
\bpr\label{sec3:pr:invdim}
Let $V$ be a maximal involutive subspace. Then
$$\dim V=m\binom{n-r+k-1}{k}+r,$$
where $n=\dim M$, $m=\dim\pi$, $r=\tp(V)$.
\epr
\begin{proof}
Let us choose local coordinates in $M$ in such a way that the vectors $\pat/
\pat x_1,\dots,\pat/\pat x_r$ form a basis in $H$. Then, in the corresponding
special system in $J^k(\pi)$, coordinates along $\Ann(H)$ will consist
of those functions $u_\si^j,\ |\si|=k$, for which multi-index $\si$ does not
contain indices $1,\dots,r$.
\end{proof}

Let $N\sbs J^k(\pi)$ be a maximal integral manifold of the Cartan
distribution and $\ta_k\in N$. Then the tangent plane to $N$ at $\ta_k$
is a maximal involutive plane. Let its type be equal to $r(\ta_k)$.
\bde\label{sec3:df:invtp}
The number $\tp(N)=\max\limits_{\ta_k\in N}r(\ta_k)$
is called the {\em type} of the maximal integral manifold $N$ of the
Cartan distribution.\index{type of maximal integral manifold of the
Cartan distribution}
\ede
Obviously, the set $g(N)=\{\,\ta_k\in N\mid r(\ta_k)=\tp(N)\,\}$
is everywhere dense in $N$. We call the points $\ta_k\in g(N)$ {\em
generic}.\index{generic point of maximal integral manifold of the
Cartan distribution} Let $\ta_k$ be such a point and $\CU$ be its
neighborhood in $N$ consisting of generic points.  Then:
\begin{enumerate}
\item $N'=\pi_{k,k-1}(N)$ is an integral manifold of the Cartan distribution on
$J^{k-1}(\pi)$;
\item $\dim(N')=\tp(N)$ and
\item $\pi_{k-1}\left|_{N'}
\right.\colon N'\to M$ is an immersion.
\end{enumerate}

\begin{theorem}\label{sec3:th:invstr}
Let $N\sbs J^{k-1}(\pi)$ be an integral manifold of the
Cartan distribution on  $J^k(\pi)$ and $\CU\sbs N$ be an open domain
consisting of generic points. Then
$$\CU=\{\,\ta_k\in J^k(\pi)\mid L_{\ta_k}\sps T_{\ta_{k-1}}\CU'\,\},$$
where $\ta_{k-1}=\pi_{k,k-1}(\ta_k),\ \CU'=\pi_{k,k-1}(\CU)$.
\end{theorem}
\begin{proof}
Let $M'=\pi_{k-1}(\CU')\sbs M$. Denote its dimension (which equals to
$\tp(N)$) by $r$ and choose local coordinates in $M$ in such a way that
the submanifold $\CV'$ is determined by the equations $x_{r+1}=\dots=x_n=0$
in these coordinates. Then, since $\CU'\sbs J^{k-1}(\pi)$ is an integral
manifold and $\pi_{k-1}\left|_{\CU'}\right.\colon \CU'\to\CV'$ is a diffeomorphism,
in corresponding special coordinates the manifold $\CU'$ is given by the
equations
$$
u_\si^j=\left\{\begin{array}{ll}
\dfrac{{\partial}^{|\si|} \vf^j}{\partial x_\si},&\mbox{ if }\si
\mbox{ does not contain }r+1,\dots,n,\\
0&\mbox{ otherwise,}
\end{array}\right.
$$
for all $j=1,\dots,m$, $|\si|\le k-1$ and some smooth function $\vf=
\vf(x_1,\dots,x_r)$. Hence, the tangent plane $H$ to $\CU'$ at $\ta_{k-1}$
is spanned by the vectors of the form \veqref{sec3:eq:vectco} with $i=1,\dots,r$.
Consequently, a point $\ta_k$, such that $L_{\ta_k}\sps H$, is determined
by the coordinates
$$
u_\si^j=\left\{\begin{array}{ll}
\dfrac{{\partial}^{|\si|} \vf^j}{\partial x_\si},&\mbox{ if }\si
\mbox{ does not contain }r+1,\dots,n,\\
\mbox{arbitrary real numbers}&\mbox{ otherwise,}
\end{array}\right.
$$
where $j=1,\dots,m$, $|\si|\le k$. Hence, if $\ta_k$, $\ta'_k$ are two such
points, then the vector $\ta_k-\ta'_k$ lies in $\Ann(H)$, as it follows
from the proof of Proposition \vref{sec3:pr:invdim}. As it can be easily
seen, any integral manifold of the Cartan distribution projecting to
$\CU'$ is contained in $\CU$, which concludes the proof.
\end{proof}
\bre\label{sec3:re:invmax}
Note that maximal integral manifolds $N$ of type $\dim M$
are exactly graphs of jets $j_k(\vf),\vf\in\Gal(\pi)$. On the other hand,
if $\tp(N)=0$, then $N$ coincides with a fiber of the projection
$\pi_{k,k-1}\colon J^k(\pi)\to J^{k-1}(\pi)$.
\ere

\subsection{Classical symmetries}\label{sub:csym}
Having the basic structure on $J^k(\pi)$, we can now introduce
transformations preserving this structure.
\bde\label{sec3:df:sym}
Let $\CU$, $\CU'\sbs J^k(\pi)$ be open domains.
\begin{enumerate}
\item A diffeomorphism $F\colon \CU\to\CU'$ is called a {\em Lie
transformation},\index{Lie transformation} if it preserves the Cartan
distribution, i.e., $F_*(\CC_{\ta_k}^k)=\CC_{F(\ta_k)}^k$
for any point $\ta_k\in\CU$.
\end{enumerate}
Let $\CE$, $\CE'\sbs J^k(\pi)$ be differential equations.
\begin{enumerate}
\item[(2)] A Lie transformation $F\colon \CU\to\CU$ is called a ({\em local})
{\em equivalence},\index{local equivalence} if $F(\CU\cap\CE)=
\CU'\cap\CE'$.
\item[(3)] A (local) equivalence is called a ({\em local}) {\em symmetry},
\index{symmetry} if $\CE=\CE'$ and $\CU=\CU'$.
\end{enumerate}
\ede
Below we shall not distinguish between local and global versions of the
concepts introduced.
\begin{example}\label{sec3:ex:j0}
Consider the case $J^0(\pi)=E$. Then, since any
$n$-di\-men\-sional horizontal plane in $T_\ta E$ is tangent to some section
of the bundle $\pi$, the Cartan plane $\CC_\ta^0$ coincides with the
whole space $T_\ta E$. Thus the Cartan distribution is trivial in this case
and any diffeomorphism of $E$ is a Lie transformation.
\end{example}
\begin{example}\label{sec3:ex:cont}
Since the Cartan distribution on $J^k(\pi)$ is locally determined by the
Cartan forms (see \veqref{sec3:eq:cforms}), the condition of $F$ to be
a Lie transformation can be reformulated as
\begin{equation}\label{sec3:eq:liecond}
F^*\om_\si^j=
\sum_{\al=1}^m\sum_{|\tau|<k}\la_{\si,\tau}^{j,\al}\om_\tau^\al,\quad
j=1,\dots,m,\quad |\si|<k,
\end{equation}
where $\la_{\si,\tau}^{j,\al}$ are smooth functions on $J^k(\pi)$. Equations
\veqref{sec3:eq:liecond} are the base for computations in local coordinates.

In particular, if $\dim\pi=1$ and $k=1$, equations \veqref{sec3:eq:liecond} reduce to
the only condition $F^*\om=\la\om$, where $\om=d  u-\sum_{i=1}^nu_i\,d
x_i$. Hence, Lie transformations in this case are just {\em contact
transformations}\index{contact transformations} of the natural contact
structure in $J^1(\pi)$.
\end{example}

Let $F\colon J^k(\pi)\to J^k(\pi)$ be a Lie transformation. Then graphs of $k$-jets
are taken by $F$ to $n$-dimensional
maximal manifolds. Let $\ta_{k+1}$ be a point of $J^{k+1}(\pi)$ and
represent $\ta_{k+1}$ as a pair $(\ta_k,L_{\ta_{k+1}})$, or, which is the
same, as a class of graphs of $k$-jets tangent to each other at $\ta_k$.
Then the image $F_*(L_{\ta_{k+1}})$
will almost always be an $R$-plane at $F(\ta_k)$.
Denote the corresponding point in $J^{k+1}(\pi)$ by $F^{(1)}(\ta_{k+1})$.
\bde\label{sec3:df:lift}
Let $F\colon J^k(\pi)\to J^k(\pi)$ be a Lie transformation. The
above defined map $F^{(1)}\colon J^{k+1}(\pi)\to J^{k+1}(\pi)$ is called the
1-{\em lifting} of $F$.\index{1-lifting of Lie transformation}
\ede
The map $F^{(1)}$ is a Lie transformation at the domain of its
definition, since almost everywhere it takes graphs of $(k+1)$-jets to
graphs of the same kind. Hence, for any $l\ge1$ we can define $F^{(l)}=
(F^{(l-1)})^{(1)}$ and call this map the $l$-{\em lifting} of
$F$.\index{$l$-lifting of Lie transformation}
\begin{theorem}\label{sec3:th:lietr}
\index{structure of Lie transformations}%
Let $\pi\colon E\to M$ be an $m$-dimensional vector bundle over
an $n$-dimensional manifold $M$ and $F\colon J^k(\pi)\to J^k(\pi)$ be a Lie
transformation. Then:
\begin{enumerate}
\item If $m>1$ and $k>0$, then the map $F$ is of the form $F=
G^{(k)}$ for some diffeomorphism $G\colon J^0(\pi)\to J^0(\pi)$\textup{;}
\item If $m=1$ and $k>1$, then the map $F$ is of the form $F=
G^{(k-1)}$ for some contact transformation $G\colon J^1(\pi)\to J^1(\pi)$.
\end{enumerate}
\end{theorem}
\begin{proof}
Recall that fibers of the projection $\pi_{k,k-1}\colon J^k(\pi)\to
J^{k-1}(\pi)$ for $k\ge1$ are maximal integral manifolds of the Cartan
distribution of type 0 (see Remark \vref{sec3:re:invmax}). Further, from
Proposition \vref{sec3:pr:invdim} it follows in the cases $m>1$, $k>0$
and $m=1$, $k>1$ that
they are integral manifolds of maximal dimension, provided $n>1$. Therefore,
the map $F$ is $\pi_{k,\ep}$-fiberwise, where $\ep=0$ for $m>1$ and
$\ep=1$ for $m=1$.

Thus there exists a map $G\colon J^\ep(\pi)\to J^\ep(\pi)$ such that
$\pi_{k,\ep}\circ F=G\circ\pi_{k,\ep}$ and $G$ is a Lie transformation
in an obvious way. Let us show that $F=G^{(k-\ep)}$. To do this, note first
that in fact, by the same reasons, the transformation $F$ generates a
series of Lie transformations $G_l\colon J^l(\pi)\to J^l(\pi),\ l=\ep,\dots,k$,
satisfying $\pi_{l,l-1}\circ G_l=G_{l-1}\circ\pi_{l,l-1}$ and $G_k=F, G_\ep=
G$. Let us compare the maps $F$ and $G_{k-1}^{(1)}$.

From Proposition \vref{sec3:pr:cplane} and the definition of Lie transformations
we obtain
$$F_*((\pi_{k,k-1})_*^{-1}(L_{\ta_k}))=F_*(\CC_{\ta_k}^k)=
\CC_{F(\ta_k)}=(\pi_{k,k-1})_*^{-1}(L_{F(\ta_k)})$$
for any $\ta_k\in J^k(\pi)$. But
$$F_*((\pi_{k,k-1})_*^{-1}(L_{\ta_k}))=
(\pi_{k,k-1})_*^{-1}(G_{k-1,*}(L_{\ta_k}))$$
and consequently
$G_{k-1,*}(L_{\ta_k})=L_{F(\ta_k)}$. Hence, by the definition of 1-lifting
we have $F=G_{k-1}^{(1)}$. Using this fact as a base of elementary
induction, we obtain the result of the theorem for $\dim M>1$.

Consider the case $n=1$, $m=1$ now. Since all maximal integral manifolds are
one-dimensional in this case, it should be treated in a special way.
Denote by $\CV$ the distribution consisting of vector fields tangent to
the fibers of the projection $\pi_{k,k-1}$. We must show that
\begin{equation}\label{sec3:eq:pres}
F_*\CV=\CV
\end{equation}
for any Lie transformation $F$, which is equivalent to $F$ being
$\pi_{k,k-1}$-fiberwise.

Let us prove \veqref{sec3:eq:pres}. To do it, consider an arbitrary distribution $\CP$
on a manifold $N$ and introduce the notation
\begin{equation}\label{sec3:eq:CD}
\CP D=\{\,X\in \Dr(N)\mid X\text{ lies in }\CP\,\}
\end{equation}
and
\begin{equation}\label{sec3:eq:DC}
D_\CP=\{\,X\in \Dr(N)\mid [X,Y]\in\CP,\ \forall Y\in \CP D\,\}.
\end{equation}
Then one can show (using coordinate representation, for example) that
$D\CV=D\CC^k\cap D_{[D\CC^k,D\CC^k]}$
for $k\ge2$. But Lie transformations preserve the distributions at the
right-hand side of the last equality and consequently preserve $D\CV$.
\end{proof}
\bde\label{sec3:df:infsym}
Let $\pi\colon E\to M$ be a vector bundle and $\CE\sbs J^k(\pi)$
be a $k$-th order differential equation.
\begin{enumerate}
\item A vector field $X$ on $J^k(\pi)$ is called a {\em Lie
field}\index{Lie field}, if the corresponding one-parameter group
consists of Lie transformations.
\item A Lie field is called an {\em infinitesimal
symmetry}\index{infinitesimal symmetry} of the equation $\CE$, if it is
tangent to $\CE$.
\end{enumerate}
\ede

Since in the sequel we shall deal with infinitesimal symmetries only, we
shall call them just {\em symmetries}. By definition, one-parameter groups of
transformations corresponding to symmetries preserve generalized
solutions\footnote{A \emph{generalized solution} of an equation $\CE$
is a maximal integral manifold $N\sbs\CE$ of the Cartan distribution on
$\CE$; see \cite{Lychagin1}.}.

Let $X$ be a Lie field on $J^k(\pi)$ and $F_t\colon J^k(\pi)\to J^k(\pi)$ be its
one-parameter group. Then we can construct the $l$-lifting $F_t^{(l)}\colon
J^{k+l}(\pi)\to J^{k+l}(\pi)$ and the corresponding Lie field $X^{(l)}$ on
$J^{k+l}(\pi)$. This field is called the {\em $l$-lifting}\index{$l$-lifting
of Lie field} of the field $X$. As we shall see a bit later, liftings of Lie
fields are defined globally and can be described explicitly. An immediate
consequence of the definition and of Theorem \vref{sec3:th:lietr} is
\begin{theorem}\label{sec3:th:liefld}
Let $\pi\colon E\to M$ be an $m$-dimensional vector bundle over
an $n$-dimensional manifold $M$ and $X$ be a Lie field on $J^k(\pi)$. Then:
\begin{enumerate}
\item If $m>1$ and $k>0$, the field $X$ is of the form $X=
Y^{(k)}$ for some vector field $Y$ on $J^0(\pi)$\textup{;}
\item If $m=1$ and $k>1$, the field $X$ is of the form $X=
Y^{(k-1)}$ for some contact vector field $Y$ on $J^1(\pi)$.
\end{enumerate}
\end{theorem}

To finish this subsection, we describe coordinate expressions for Lie fields.
Let $(x_1,\dots,x_n,\dots,u_\si^j,\dots)$ be a special coordinate system in
$J^k(\pi)$.
Then from \veqref{sec3:eq:liecond} it follows that
$$X=\sum_{i=1}^nX_i\pdr{}{x_i}+
\sum_{j=1}^m\sum_{|\si|\le k}X_\si^j\pdr{}{u_\si^j}$$
is a Lie field if and only if
\begin{equation}\label{sec3:eq:lift}
X_{\si i}^j=D_i(X_\si^j)-\sum_{\al=1}^nu_{\si \al}^jD_i(X_\al),
\end{equation}
where
\begin{equation}\label{sec3:eq:totder}
D_i=
\pdr{}{x_i}+\sum_{j=1}^m\sum_{|\si|\ge0}u_{\si i}^j\pdr{}{u_\si^j}
\end{equation}
are the so-called {\em total derivatives}.\index{total derivatives}
\begin{xca}
It is easily seen that the operators \eqref{sec3:eq:totder} do not preserve
the algebras $\CF_k$: they are derivations acting from $\CF_k$ to
$\CF_{k+1}$. Prove that nevertheless for any contact field on $J^1(\pi)$,
$\dim\pi=1$, or for an arbitrary vector field on $J^0(\pi)$ (regardless
of the dimension of $\pi$) the formulas above determine a vector field on
$J^k(\pi)$.
\end{xca}

Recall now that a contact field $X$ on $J^1(\pi)$ is completely determined
by its {\em generating function}\index{generating function of contact
field} $f=\ip_X\om$, where $\om=du-\sum_iu_i\,dx_i$ is the Cartan
(contact) form on $J^1(\pi)$. The contact field corresponding to
$f\in\CF_1(\pi)$ is denoted by $X_f$ and is given by the expression
\begin{gather}\label{sec3:eq:cont}
\begin{split}
X_f=-\sum_{i=1}^n\pdr{f}{u_{1_i}}\pdr{}{x_i}+
\Biggl(f-&\sum_{i=1}^nu_i\pdr{f}{u_i}\Biggr)\pdr{}{u}\\
&+\sum_{i=1}^n\left(\pdr{f}{x_i}+u_i\pdr{f}{u}\right)\pdr{}{u_i}.
\end{split}
\end{gather}

Thus, starting with a field \veqref{sec3:eq:cont} in the case $\dim\pi=1$ or with
an arbitrary field on $J^0(\pi)$ for $\dim\pi>1$ and using
\veqref{sec3:eq:lift}, we can obtain efficient expressions for Lie fields.
\bre\label{sec3:re:gensec}
Note that in the multi-dimensional case $\dim\pi>1$ we can introduce the
functions $f^j=\ip_X\om^j$, where
$\om^j=d  u^j-\sum_iu_i^j\,dx_i$ are the Cartan forms on $J^1(\pi)$.
Such a function may be understood as an element of the module
$\CF_1(\pi,\pi)$. The local conditions of a section $f\in\CF_1(\pi,\pi)$ to
generate a Lie field is as follows:
$$\pdr{f^\al}{u_i^\al}=\pdr{f^\be}{u_i^\be},\quad
\pdr{f^\al}{u_i^\be}=0,\quad\al\neq\be.$$
We call $f$ the \emph{generating function} (though, strictly speaking,
the term \emph{generating section} should be used)
\index{generating section of Lie field}\index{generating function of Lie
field} of the Lie field $X$, if $X$ is a lifting of the field $X_f$.
\ere

Let us write down the conditions of a Lie field to be a symmetry.
Assume that an equation $\CE$ is given by the relations $F^1=0,\dots,F^r=0$,
where $F^j\in\CF_k(\pi)$. Then $X$ is a symmetry of $\CE$ if and only if
$$X(F^j)=\sum_{\al=1}^r\la_\al^jF^\al,\quad j=1,\dots,r,$$
where $\la_\al^j$ are smooth functions, or
$X(F^j)\left|_{\CE}\right.=0,\ j=1,\dots,r$.
These conditions can be rewritten in terms of generating sections and we
shall do it later in a more general situation.

\subsection{Prolongations of differential equations}\label{sub:proleq}
Prolongations are differential consequences of a given differential
equation. Let us give a formal definition.
\bde\label{sec3:df:proleq}
Let $\CE\sbs J^k(\pi)$ be a differential equation of order $k$. Define the
set
$$\CE^1=\{\,\ta_{k+1}\in J^{k+1}(\pi)\mid \pi_{k+1,k}(\ta_{k+1})\in\CE,\
L_{\ta_{k+1}}\sbs T_{\pi_{k+1,k}(\ta_{k+1})}\CE\,\}$$
and call it the {\em first prolongation} of the equation $\CE$.\index{first
prolongation of equation}
\ede
If the first prolongation $\CE^1$ is a submanifold in $J^{k+1}(\pi)$, we
define the second prolongation of $\CE$ as $(\CE^1)^1\sbs J^{k+2}(\pi)$, etc.
Thus the {\em $l$-th prolongation}\index{$l$-th prolongation of equation} is
a subset $\CE^l\sbs J^{k+l}(\pi)$.

Let us redefine the notion of $l$-th prolongation directly. Namely, take a
point $\ta_k\in\CE$ and consider a section $\vf\in\Gal(\pi)$ such that the
graph of $j_k(\vf)$ is tangent to $\CE$ with order $l$. Let $\pi_k(\ta_k)=
x\in M$. Then $[\vf]_x^{k+l}$ is a point of $J^{k+l}(\pi)$ and the set
of all points obtained in such a way obviously coincides with $\CE^l$,
provided all intermediate prolongations $\CE^1,\dots,\CE^{l-1}$ be well
defined in the sense of Definition \vref{sec3:df:proleq}.

Assume now that locally $\CE$ is given by the equations
$F^1=0,\dots,F^r=0$, $F^j\in\CF_k(\pi)$
and $\ta_k\in\CE$ is the origin of the chosen special coordinate system.
Let $u^1=\vf^1(x_1,\dots,x_n),\dots,u^m=\vf^m(x_1,\dots,x_n)$ be a local
section of the bundle $\pi$. Then the equations of the first prolongation are
$$\pdr{F^j}{x_i}+\sum_{\al,\si}
\pdr{F^j}{u_\si^\al}u_{\si i}^\al=0,\quad i=1,\dots,n,\quad j=1,\dots,r,$$
combined with the initial equations $F^r=0$.
From here and by comparison with the coordinate representation of
prolongations for nonlinear differential operators (see Subsection
\ref{sub:ndo}), we obtain the following result:
\bpr\label{sec3:pr:proleq}
Let $\CE\sbs J^k(\pi)$ be a differential equation. Then
\begin{enumerate}
\item If the equation $\CE$ is determined by a differential
operator $\De\colon \Ga(\pi)\to\Ga(\pi')$, then its $l$-th prolongation is
given by the $l$-th prolongation $\De^{(l)}\colon \Ga(\pi)\to\Ga(\pi'_l)$ of
the operator $\De$.
\item If $\CE$ is locally described by the system
$F^1=0,\dots,F^r=0,\ F^j\in\CF_k(\pi)$, then the system
\begin{equation}\label{sec3:eq:prolco}
D_\si F^j=0,\quad |\si|\le l,j=1,\dots,r,
\end{equation}
where $D_\si= D_{i_1}\circ\dots\circ D_{i_{\abs{\si}}}$,
$\si=i_1\dots i_{|\si|}$, corresponds to
$\CE^l$. Here $D_i$ stands for the $i$-th total derivative \textup{(}see
\veqref{sec3:eq:totder}\textup{)}.
\end{enumerate}
\epr

From the definition it follows that for any $l\ge l'\ge0$ one has
$\pi_{k+l,k+l'}(\CE^l)\sbs\CE^{l'}$ and consequently one has the maps
$\pi_{k+l,k+l'}\colon \CE^l\to\CE^{l'}$.
\bde\label{sec3:df:formint}
An equation $\CE\sbs J^k(\pi)$ is called {\em formally
integrable}\index{formally integrable equation}, if
\begin{enumerate}
\item all prolongations
$\CE^l$ are smooth manifolds and
\item all the maps $\pi_{k+l+1,k+l}\colon \CE^{l+1}
\to\CE^l$ are smooth fiber bundles.
\end{enumerate}
\ede
\bde\label{sec3:df:infprol}
The inverse limit $\projlim_{l\to\infty}\CE^l$ with respect
to projections $\pi_{l+1,l}$ is called the {\em infinite prolongation}
\index{infinite prolongation of equation} of the equation $\CE$ and is
denoted by $\Ei\sbs\Ji(\pi)$.
\ede
\subsection{Basic structures on infinite prolongations}\label{sub:basstr}
Let $\pi\colon E\to M$ be a vector bundle and $\CE\sbs J^k(\pi)$ be a $k$-th
order differential equation. Then we have embeddings $\ve_l\colon \CE^l\sbs
J^{k+l}(\pi)$ for all $l\ge0$. Since, in general, the sets $\CE^l$ may
not be smooth manifolds, we define a {\em function} on $\CE^l$ as the
restriction $f\left|_{\CE^l}\right.$ of a smooth function $f\in
\CF_{k+l}(\pi)$. The set $\CF_l(\CE)$ of all functions on $\CE^l$ forms
an $\BBR$-algebra in a natural way and $\ve_l^*\colon \CF_{k+l}(\pi)\to\CF_l(\CE)$
is a homomorphism of algebras. In the case of formally integrable equations, the
algebra $\CF_l(\CE)$ coincides with $\Ci(\CE^l)$. Let $I_l=\ker\ve_l^*$.
Evidently, $I_l(\CE)\sbs I_{l+1}(\CE)$. Then $I(\CE)=
\bigcup_{l\ge0}I_l(\CE)$
is an ideal in $\CF(\pi)$ which is called the {\em ideal of the equation}
\label{sec3:ideal}
$\CE$.\index{ideal of equation} The {\em function algebra}\index{algebra of
functions on $\Ei$} on $\Ei$ is the quotient algebra
$\CF(\CE)=\CF(\pi)/I(\CE)$ and
coincides with $\injlim_{l\to\infty}\CF_l(\CE)$ with respect to the
system of homomorphisms $\pi_{k+l+1,k+l}^*$. For all $l\ge0$, we have
the homomorphisms $\ve_l^*\colon \CF_l(\CE)\to\CF(\CE)$. When $\CE$ is
formally integrable, they are monomorphic, but in any case the algebra
$\CF(\CE)$ is filtered by the images of $\ve_l^*$.

To construct differential calculus on $\Ei$, one needs the general
algebraic scheme exposed in Section \ref{sec:calc} and applied to the
filtered algebra $\CF(\CE)$. However, in the case of formally integrable
equations, due to the fact that all $\CE^l$ are smooth manifolds, this
scheme may be simplified and combined with a purely geometrical approach
(cf.\ with similar constructions of Subsection \ref{sub:ijets}).

In special coordinates the infinite prolongation of the equation $\CE$ is
determined by the system similar to \veqref{sec3:eq:prolco} with the only
difference
that $|\si|$ is unlimited now. Thus, the ideal $I(\CE)$ is generated by
the functions $D_\si F^j$, $|\si|\ge0$, $j=1,\dots,m$. From these remarks
we obtain the following important fact.
\begin{remark}\label{sec3:re:tottan}
Let $\CE$ be a formally integrable equations. Then from the
above said it follows that the ideal $I(\CE)$ is stable with respect to the
action of the total derivatives $D_i,\ i=1,\dots,n$. Consequently, the
restrictions $D_i^\CE=D_i\rest{\CE}
\colon \CF(\CE)\to\CF(\CE)$ are well defined and $D_i^\CE$ are filtered
derivations. In other words, we obtain that the vector
fields $D_i$ are tangent to any infinite prolongation and thus determine
vector fields on $\Ei$. We shall often skip the superscript $\CE$ in the
notation of the above defined restrictions.
\end{remark}
\begin{example}\label{sec3:ex:evol}
Consider a system of evolution equations of the form
$$u_t^j=f^j(x,t,\dots,u^\al,\dots,u_x^\al,\dots),\quad
j,\al=1,\dots,m.$$
Then the set of functions $x_1,\dots,x_n,t,\dots,u_{i_1,\dots,i_r,0}^j$
with $1\le i_k\le n$, $j=1,\dots,m$, where $t=x_{n+1}$, may
be taken for \emph{internal coordinates} on $\Ei$. The total derivatives
restricted to $\Ei$ are expressed as
\begin{align}\label{sec3:eq:evtot}
\begin{split}
D_i=&\pdr{}{x_i}+\sum_{j=1}^n\sum_{|\si|\ge0}u_{\si i}^j
\pdr{}{u_\si^j},\ i=1,\dots,n,\\
D_t=&\pdr{}{t}+\sum_{j=1}^n\sum_{|\si|\ge0}D_\si(f^j)\pdr{}{u_\si^j}
\end{split}
\end{align}
in these coordinates, while the Cartan forms restricted to $\Ei$ are written
down as
\begin{equation}\label{sec3:eq:evcar}
\om_\si^j=du_\si^j-\sum_{i=1}^nu_{\si i}^j\,dx_i-D_\si(f^j)\,dt.
\end{equation}
\end{example}

Let $\pi\colon E\to M$ be a vector bundle and $\CE\sbs J^k(\pi)$ be a formally
integrable equation.
\bde\label{sec3:df:cartji}
Let $\ta\in\Ji(\pi)$. Then
\begin{enumerate}
\item The {\em Cartan plane}\index{Cartan plane on $\Ji(\pi)$} $\CC_\ta
=\CC_\ta(\pi)\sbs T_\ta\Ji(\pi)$ at
$\ta$ is the linear envelope of tangent planes to all manifolds
$j_\infty(\vf)(M),\ \vf\in\Ga(\pi)$, passing through $\ta$.
\item If $\ta\in\Ei$, then the intersection $\CC_\ta(\CE)=\CC_\ta(\pi)
\cap T_\ta\Ei$ is called {\em Cartan plane}\index{Cartan plane on $\Ei$}
of $\Ei$ at $\ta$.
\end{enumerate}
The correspondence $\ta\mapsto\CC_\ta(\pi),\ \ta\in\Ji(\pi)$ (respectively,
$\ta\mapsto\CC_\ta(\Ei)$, $\ta\in\Ei$) is called the {\em Cartan distribution}
on $\Ji(\pi)$ (respectively, on $\Ei$).\index{Cartan distribution on
$\Ji(\pi)$}\index{Cartan distribution on $\Ei$}
\ede
\bpr\label{sec3:pr:cartji}
For any vector bundle $\pi\colon E\to M$ and a formally integrable
equation $\CE\sbs J^k(\pi)$ one has:
\begin{enumerate}
\item The Cartan plane $\CC_\ta(\pi)$ is $n$-dimensional at any
point $\ta\in\Ji(\pi)$.
\item Any point $\ta\in\Ei$ is generic, i.e., $\CC_\ta(\pi)\sbs
T_\ta\Ei$ and thus one has $\CC_\ta(\Ei)=\CC_\ta(\pi)$.
\item Both distributions, $\CC(\pi)$ and $\CC(\Ei)$, are
integrable.
\end{enumerate}
\epr
\begin{proof}
Let $\ta\in\Ji(\pi)$ and $\pi_\infty(\ta)=x\in M$. Then the point $\ta$
completely determines all partial derivatives of any section $\vf\in
\Gal(\pi)$ such that its graph passes through $\ta$. Consequently, all
such graphs have a common tangent plane at this point, which coincides with
$\CC_\ta(\pi)$. This proves the first statement.

To prove the second one, recall Remark \vref{sec3:re:tottan}: locally,
any vector field $D_i$ is tangent to $\Ei$. But as it follows from
\veqref{sec3:eq:cforms},
$\ip_{D_i}\om_\si^j=0$ for any $D_i$ and for any Cartan form $\om_\si^j$.
Hence, linear independent vector fields $D_1,\dots,D_n$ locally lie both in
$\CC(\pi)$ and in $\CC(\Ei)$ which gives the result.

Finally, as it follows from the above said, the module
\begin{equation}\label{sec3:eq:CDpi}
\CC \Dr(\pi)=\{\,X\in \Dr(\pi)\mid X\mbox{ lies in } \CC(\pi)\,\}
\end{equation}
is locally generated by the fields $D_1,\dots,D_n$. But it is easily seen
that $[D_\al,D_\be]=0$ for all $\al,\be=1,\dots,n$ and consequently
$[\CC \Dr(\pi),\CC \Dr(\pi)]\sbs\CC \Dr(\pi)$. The same reasoning is valid for
\begin{equation}\label{sec3:eq:CDei}
\CC \Dr(\CE)=\{\,X\in \Dr(\Ei)\mid X\mbox{ lies in } \CC(\Ei)\,\}.
\end{equation}
This completes the proof of the proposition.
\end{proof}
\bpr\label{sec3:pr:cartmax}
Maximal integral manifolds of the Cartan distribution
$\CC(\pi)$ are graph of infinite jets of sections $j_\infty(\vf)$, $\vf\in
\Gal(\pi)$.
\index{maximal integral manifolds of the Cartan distribution on $\Ji(\pi)$}
\epr
\begin{proof}
Note first that graphs of infinite jets are integral manifolds of the
Cartan distribution of maximal dimension (equaling to $n$) and that any
integral manifold projects to $J^k(\pi)$ and $M$ without singularities.

Let now $N\sbs\Ji(\pi)$ be an integral manifold and $N^k=\pi_{\infty,k}N
\sbs J^k(\pi)$, $N'=\pi_\infty N\sbs M$. Hence, there exists a
diffeomorphism $\vf'\colon N'\to N^0$ such that $\pi\circ\vf'=\id_{N'}$. Then by
the Whitney theorem on extension for smooth functions (see \cite{Malgr1}),
there exists a local
section $\vf\colon M\to E$ satisfying $\vf\left|_{N'}\right.=\vf'$ and $j_k(\vf)(M)
\sps N^k$ for all $k>0$. Consequently, $j_\infty(\vf)(M)\sps N$.
\end{proof}
\bco\label{sec3:co:cartmaxeq}
Maximal integral manifolds of the Cartan distribution on $\Ei$
coincide locally with graphs of infinite jets of solutions.
\eco

Consider a point $\ta\in\Ji(\pi)$ and let $x=\pi_\infty(\ta)\in M$
be its projection to $M$. Let $v$
be a tangent vector to $M$ at the point $x$. Then, since the Cartan plane
$\CC_\ta$ isomorphically projects to $T_xM$, there exists a unique
tangent vector $\CC v\in T_\ta\Ji(\pi)$ such that $(\pi_\infty)_*(\CC v)=v$.
Hence, for any vector field $X\in \Dr(M)$ we can define a vector field $\CC X
\in \Dr(\pi)$ by setting $(\CC X)_\ta=\CC(X_{\pi_\infty(\ta)})$. Then, by
construction, the field $\CC X$ is projected by $(\pi_\infty)_*$ to $X$ while
the correspondence $\CC\colon \Dr(M)\to \Dr(\pi)$ is a $\Ci(M)$-linear one.
In other words, this correspondence is a linear connection in the bundle
$\pi_\infty\colon \Ji(\pi)\to M$.
\bde\label{sec3:df:cartcon}
The connection $\CC\colon \Dr(M)\to \Dr(\pi)$ defined above
is called the {\em Cartan connection} in $\Ji(\pi)$.\index{Cartan connection
in $\Ji(\pi)$}
\ede

For any formally integrable equation, the Cartan connection is obviously
restricted to the bundle $\pi_\infty\colon \Ei\to M$\index{Cartan connection in
$\Ei$} and we preserve the same notation $\CC$ for this restriction.

Let $(x_1,\dots,x_n,\dots,u_\si^j,\dots)$ be a special coordinate system in
$\Ji(\pi)$
and $X=X_1\pat/\pat x_1+\dots+X_n\pat/\pat x_n$ be a vector field on $M$
represented in this coordinate system. Then the field $\CC X$ is to be of
the form $\CC X=X+X^v$, where $X^v=\sum_{j,\si}X_\si^j\pat/\pat u_\si^j$ is
a $\pi_\infty$-vertical field. The defining conditions $\ip_{\CC X}
\om_\si^j=0$, where $\om_\si^j$ are the Cartan forms on $\Ji(\pi)$, imply
\begin{equation}\label{sec3:eq:cconcoo}
\CC X=\sum_{i=1}^nX_i\left(\pdr{}{x_i}+\sum_{j,\si}u_{\si i}^j
\pdr{}{u_\si^j}\right)=\sum_{i=1}^nX_iD_i.
\end{equation}
In particular, $\CC(\pat/\pat x_i)=D_i$, i.e., total derivatives
are the liftings to $\Ji(\pi)$ of the corresponding partial derivatives by
the Cartan connection.

Let now $V$ be a vector field on $\Ei$
and $\ta\in\Ei$ be a point. Then the vector $V_\ta$ can be projected
parallel to the Cartan plane $\CC_\ta$ to the fiber of the projection
$\pi_\infty\colon \Ei\to M$ passing through $\ta$. Thus we get a vertical vector
field $V^v$. Hence, for any $f\in\CF(\CE)$ a differential one-form $U_\CC(f)
\in \La^1(\CE)$ is defined by
\begin{equation}\label{sec3:eq:ue}
\ip_V(U_\CC(f))= V^v(f),\quad V\in \Dr(\CE).
\end{equation}
The correspondence $f\mapsto U_\CC(f)$ is a derivation of the algebra
$\CF(\CE)$ with the values in the $\CF(\CE)$-module $\La^1(\CE)$, i.e.,
$U_\CC(fg)=fU_\CC(g)+gU_\CC(f)$ for all $f,g\in\CF(\CE)$.
\bde\label{sec3:df:ue}
The derivation $U_\CC\colon \CF(\CE)\to\La^1(\CE)$ is called the
{\em structural element} of the equation\index{structural element of
equation} $\CE$.
\ede

In the case $\Ei=\Ji(\pi)$ the structural element $U_\CC$ is locally
represented in the form
\begin{equation}\label{sec3:eq:upi}
U_\CC=\sum_{j,\si}\om_\si^j\ot\pdr{}{u_\si^j},
\end{equation}
where $\om_\si^j$ are the Cartan forms on $\Ji(\pi)$. To obtain the
expression in the general case, one needs to rewrite \veqref{sec3:eq:upi}
in local coordinates.

The following result is a consequence of definitions:
\bpr\label{sec3:pr:ccom}
For any vector field $X\in \Dr(M)$ the equality
\begin{equation}\label{sec3:eq:ccom}
j_\infty(\vf)^*(\CC X(f))=X(j_\infty(\vf)^*(f)),\quad
f\in\CF(\pi),\quad\vf\in\Gal(\pi),
\end{equation}
takes place. Equality {\em \veqref{sec3:eq:ccom}} uniquely determines
the Cartan connection in $\Ji(\pi)$.
\epr
\bco\label{sec3:co:cflat}
The Cartan connection in $\Ei$ is flat, i.e.,
$\CC [X,Y]=[\CC X,\CC Y]$ for any $X,Y\in \Dr(M)$.
\eco
\begin{proof}
Consider the case $\Ei=\Ji(\pi)$. Then from Proposition \vref{sec3:pr:ccom} we
have
\begin{multline*}
j_\infty(\vf)^*(\CC[X,Y](f))=[X,Y](j_\infty(\vf)^*(f))\\
=X(Y(j_\infty(\vf)^*(f)))-Y(X(j_\infty(\vf)^*(f)))
\end{multline*}
for any $\vf\in\Gal(\pi),f\in\CF(\pi)$. On the other hand,
\begin{multline*}
j_\infty(\vf)^*([\CC X,\CC Y](f))=j_\infty(\vf)^*(\CC X(\CC Y(f))-
\CC Y(\CC X(f)))\\
=X(j_\infty(\vf)^*(Y(f)))-Y(j_\infty(\vf)^*(\CC X(f)))\\
=X(Y(j_\infty(\vf)^*(f)))-Y(X(j_\infty(\vf)^*(f))).
\end{multline*}
To prove the statement for an arbitrary formally integrable equation
$\CE$, it suffices to note that the Cartan connection in $\Ei$ is obtained
by restricting the fields $\CC X$ to infinite prolongation of $\CE$.
\end{proof}

The construction of Proposition \vref{sec3:pr:ccom} can be generalized.
Let $\pi\colon E\to M$ be a vector bundle and $\xi_1\colon E_1\to M$,
$\xi_2\colon E_2\to M$ be
another two vector bundles.
\bde\label{sec3:df:dolift}
Let $\De\colon \Ga(\xi_1)\to\Ga(\xi_2)$ be a linear differential
operator. The {\em lifting} $\CC\De\colon \CF(\pi,\xi_1)\to\CF(\pi,\xi_2)$ of
the operator $\De$\index{lifting of linear operator} is defined by
\begin{equation}\label{sec3:eq:dolift}
j_\infty(\vf)^*(\CC\De(f))=\De(j_\infty(\vf)^*(f)),
\end{equation}
where $\vf\in\Gal(\pi)$, $f\in\CF(\pi,\xi_1)$ are arbitrary sections.
\ede
\bpr\label{sec3:pr:dolift}
Let $\pi\colon E\to M$, $\xi_i\colon E_i\to M$, $i=1,2,3$, be vector bundles. Then
\begin{enumerate}
\item For any $\Ci(M)$-linear differential operator $\De\colon
\Ga(\xi_1)\to\Ga(\xi_2)$, the operator $\CC\De$ is an $\CF(\pi)$-linear
differential operator of the same order.
\item For any $\De,\square\colon \Ga(\xi_1)\to\Ga(\xi_2)$ and
$f,g\in\CF(\pi)$, one has
$$\CC(f\De+g\square)=f\CC\De+g\CC\square.$$
\item\label{sec3:it:cc} For $\De_1\colon
\Ga(\xi_1)\to\Ga(\xi_2)$ and $\De_2\colon\Ga(\xi_2)\to\Ga(\xi_3)$, one has
$$\CC(\De_2\circ\De_1)=\CC\De_2\circ\CC\De_1.$$
\end{enumerate}
\epr

From this proposition and from Proposition \vref{sec3:pr:ccom} it follows that
if $\De$ is a scalar differential operator $\Ci(M)\to\Ci(M)$ locally
represented as $\De=\sum_\si a_\si\pat^{|\si|}/\pat x_\si,\ a_\si\in\Ci(M)$,
then $\CC\De=\sum_\si a_\si D_\si$ in the corresponding special coordinates.
If $\De=\|\De_{ij}\|$ is a matrix operator, then $\CC\De=\|\CC\De_{ij}\|$.
Obviously, $\CC\De$ may be understood as a constant
differential operator acting from sections of the bundle $\pi$ to linear
differential operators from $\Ga(\xi_1)$ to $\Ga(\xi_2)$. This observation
is generalized as follows.
\bde\label{sec3:df:cdif}
An $\CF(\pi)$-linear differential operator $\De$ acting from
the module $\CF(\pi,\xi_1)$
to $\CF(\pi,\xi_2)$ is called a \emph{$\CC$-differential operator},
\index{$\CC$-differential operator} if it admits restriction to graphs of
infinite jets, i.e., if for any section $\vf\in\Ga(\pi)$ there exists
an operator $\De_\vf\colon \Ga(\xi_1)\to\Ga(\xi_2)$ such that
\begin{equation}\label{sec3:eq:cdif}
j_\infty(\vf)^*(\De(f))=\De_\vf(j_\infty(\vf)^*(f))
\end{equation}
for all $f\in\CF(\pi,\xi_1)$.
\ede
Thus, $\CC$-differential operators are nonlinear differential operators
taking their values in $\Ci(M)$-modules of linear differential operators.
\begin{xca}\label{sec3:es:Cuniq}
Consider a \cd operator $\Delta\colon \F(\pi,\xi_1)\to\F(\pi,\xi_2)$.
Prove that if $\Delta(\pi^*(f))=0$ for all $f\in\Gamma(\xi_1)$, then
$\Delta=0$.
\end{xca}
\bpr\label{sec3:pr:cdif}
Let $\pi$, $\xi_1$, and $\xi_2$ be vector bundles over $M$. Then any
$\CC$-differential operator $\De\colon \CF(\pi,\xi_1)\to\CF(\pi,\xi_2)$ can be
presented in the form $\De=\sum_\al a_\al\CC\De_\al$, $a_\al\in\CF(\pi)$,
where $\De_\al$ are linear differential operators acting from $\Ga(\xi_1)$
to $\Ga(\xi_2)$.
\epr
\begin{proof}
Recall that we consider the filtered theory; in particular, there
exists an integer $l$ such that $\De(\CF_k(\pi,\xi_1))\sbs
\CF_{k+l}(\pi,\xi_2)$ for all $k$. Consequently, since $\Ga(\xi_1)$ is
embedded into $\CF_0(\pi,\xi_1)$, we have $\De(\Ga(\xi_1))\sbs
\CF_l(\pi,\xi_2)$ and the restriction $\bar{\De}=\De\left|_{\Ga(\xi_1)}
\right.$ is a $\Ci(M)$-differential operator taking its values in
$\CF_l(\pi,\xi_2)$.

On the other hand, the operator $\bar{\De}$ is represented in the form
$\bar{\De}=\sum_\al a_\al\De_\al$, $a_\al\in\CF_l(\pi)$, with $\De_\al\colon
\Ga(\xi_1)
\to\Ga(\xi_2)$ being $\Ci(M)$-linear differential operators. Define $\CC
\bar{\De}=\sum_\al a_\al\CC\De_\al$. Then the difference $\De-\CC
\bar{\De}$ is a $\CC$-differential operator such that its restriction to
$\Ga(\xi_1)$ vanishes. Therefore, by Exercise \vref{sec3:es:Cuniq}
$\De=\CC\bar{\De}$.
\end{proof}
\bco\label{sec3:co:cdif}
$\CC$-differential operators admit restrictions to infinite
prolongations: if $\De\colon \CF(\pi,\xi_1)\to\CF(\pi,\xi_2)$ is a
$\CC$-differential operator and $\CE\sbs J^k(\pi)$ is a $k$-th order
equation, then there exists a linear differential operator $\De_\CE\colon
\CF(\CE,\xi_1)\to\CF(\CE,\xi_2)$ such that $\ve^*\circ\De=\De_\CE\circ\ve^*$,
where $\ve\colon \Ei\hra\Ji(\pi)$ is the natural embedding.
\eco
\begin{proof}The result immediately follows from Remark
\vref{sec3:re:tottan} and from Proposition \vref{sec3:pr:cdif}.
\end{proof}
\begin{example}\label{sec3:ex:hordR}\label{sec3:p:start}
Let $\xi_1=\tau_i^*$, $\xi_2=\tau_{i+1}^*$, where $\tau_p^*\colon \bigwedge^p
T^*M\to M$ (see Example \vref{sec3:ex:deRham}), and $\De=d\colon\La^i(M)\to
\La^{i+1}(M)$ be the de Rham differential. Then we obtain the first-order
operator $\hd=\CC d\colon
\hL^i(\pi)\to\hL^{i+1}(\pi)$, where $\hL^p(\pi)$ denotes the module $\CF(\pi,
\tau_p^*)$. Due Corollary \vref{sec3:co:cdif} the operators $\hd\colon \hL^i(\CE)\to
\hL^{i+1}(\CE)$ are also defined, where $\hL^p(\CE)=\CF(\CE,\tau_p^*)$.
\end{example}
\bde\label{sec3:df:hordR}
Let $\CE\sbs J^k(\pi)$ be an equation.
\begin{enumerate}
\item Elements of the module $\hL^i(\CE)$ are called {\em horizontal
$i$-forms}\index{horizontal forms} on the infinite prolongation $\Ei$.
\item The operator $\hd\colon \hL^i(\CE)\to\hL^{i+1}(\CE)$ is called
the {\em horizontal de Rham differential}\index{horizontal de Rham
differential} on $\Ei$.
\end{enumerate}
\ede
From Proposition \ref{sec3:pr:dolift}~\veqref{sec3:it:cc} it follows
that $\hd\circ\hd=0$. The sequence
$$
0\xra{}\CF(\CE)\xra{\hd}\hL^1(\CE)\xra{}\dotsb\xra{}\hL^i(\CE)\xra{\hd}
\hL^{i+1}(\CE)\to\dotsb
$$
is called the {\em horizontal de Rham complex}\index{horizontal de Rham
complex} of the equation $\CE$. Its
cohomology is called the {\em horizontal de Rham
cohomology}\index{horizontal de Rham cohomology} of $\CE$ and is
denoted by $\hH^*(\CE)=\bigoplus_{i\ge0}\hH^i(\CE)$.

In local coordinates, horizontal forms of degree $p$ on $\Ei$ are represented
as $\om=\sum_{i_1<\dots<i_p}a_{i_1\dots i_p}d  x_{i_1}\wg\dots\wg dx_{i_p}$,
where $a_{i_1\dots i_p}\in\CF(\CE)$, while the horizontal de Rham
differential acts as
\begin{equation}\label{sec3:eq:hdcoor}
\hd(\om)=\sum_{i=1}^n\sum_{i_1<\dots<i_p}
D_i(a_{i_1\dots i_p})\, d  x_i\wg d x_{i_1}\wg\dots\wg d x_{i_p}.
\end{equation}
In particular, we see that both $\hL^i(\CE)$ and $\hH^i(\CE)$ vanish for
$i>\dim M$.

Consider the algebra $\La^*(\CE)$ of all differential forms on $\Ei$ and
note that one has the embedding $\hL^*(\CE)\hra\La^*(\CE)$. Let us extend the
horizontal de Rham differential to this algebra as follows:
$$\hd(d\om)=-d(\hd(\om)),\quad
\hd(\om\wg\ta)=\hd(\om)\wg\ta+(-1)^p\om\wg\hd(\ta)\quad
\om\in\La^p(\CE).$$
Obviously, these conditions define the differential $\hd\colon \La^i(\CE)\to
\La^{i+1}(\CE)$ and its restriction to $\hL^*(\CE)$ coincides with
the horizontal de Rham differential.

Let us also set $d_\CC=d-\hd\colon \La^*(\CE)\to\La^*(\CE)$ and call
$d_\CC$ the \emph{Cartan} (or \emph{vertical}) \emph{differential} on $\Ei$.
Then from definitions we obtain
$$
d=\hd+d_\CC,\quad \hd\circ\hd=d_\CC\circ d_\CC=0,\quad
d_\CC\circ\hd+\hd\circ d_\CC=0,
$$
i.e.,  the pair $(\hd,d _\CC)$ forms a bicomplex in $\La^*(\CE)$
with the total differential $d$.
It is called the \emph{variational bicomplex}\label{sec3:pg:varbic} and
will be discussed in more details in Section \ref{css:sec}.

Denote by $\CC\La^1(\CE)$ the \label{sec3:pg:catrsubm}\emph{Cartan
submodule} in $\La^1(\CE)$, i.e.,  the module of 1-forms vanishing on
the Cartan distribution on $\Ei$ (cf.\ with Definition
\vref{sec3:df:cforms}). Then the splitting $d=\hd+d_\CC$ implies the
direct sum decomposition
$$\La^1(\CE)=\hL^1(\CE)\oplus\CC\La^1(\CE),$$
which gives
\begin{equation}\label{sec3:eq:split}
\La^i(\CE)=\bigoplus_{p+q=i}\hL^q(\CE)\ot_{\CF(\CE)}\CC^p\La(\CE),
\end{equation}
where
$\CC^p\La(\CE)=\underbrace{\CC\La^1(\CE)\wg\dots\wg
\CC\La^1(\CE)}_{p\ \mathrm{times}}$.

To conclude this section, we shall write down the coordinate representation
for the Cartan differential $d_\CC$ and the extended differential $\hd$.
First note that
by definition and due to representation \veqref{sec3:eq:hdcoor}, one has
\begin{equation}\label{sec3:eq:dc}
d_\CC(f)=\sum_{j,\si}\pdr{f}{u_\si^j}\om_\si^j,\quad f\in\CF(\pi).
\end{equation}
In particular, $d _\CC$ takes coordinate functions $u_\si^j$ to the
corresponding Cartan forms. This is reason why we called $d_\CC$ the
Cartan differential on $\Ei$. It is easily
seen that $d_\CC\left|_{\CF(\CE)}\right.=U_\CC(\CE)$ (see Definition
\vref{sec3:df:ue}).
To finish computations, it suffices to compute $\hd(\om_\si^j)$:
$$
\hd(\om_\si^j)=\hd d_\CC(u_\si^j)=-d_\CC\hd(u_\si^j)
$$
and thus
\begin{equation}\label{seq3:eq:hd}
\hd(\om_\si^j)=-\sum_{i=1}^n\om_{\si i}^j\wg d x_i.
\end{equation}
Note that from the results obtained it follows that
\begin{gather*}
\hd(\hL^q(\CE)\ot\CC^p\La(\CE))\sbs\hL^{q+1}(\CE)\ot\CC^p\La(\CE),\\
d_\CC(\hL^q(\CE)\ot\CC^p\La(\CE))\sbs\hL^q(\CE)\ot\CC^{p+1}\La(\CE).
\label{sec3:p:end}
\end{gather*}

Now let us define the \emph{module of horizontal jets}. Let $\xi$ be a
vector bundle over $M$. Say that two elements of $P=\CF(\CE,\xi)$ are
horizontally equivalent up to order $k\le\infty$ at point $\theta\in
\Ei$, if their total derivatives up to order $k$ coincide at $\theta$.
The horizontal jet space $\bar{J}_\theta^k (P)$ is $P$ modulo this
relation, and the collection $\bar{J}^k(P)=\bigcup
_{\theta\,\in\,\Ei}\bar{J}^k_\theta(P)$ constitutes the
\emph{horizontal jet bundle} $\bar{J}^k(P)\to\Ei$. We denote
the module of sections of horizontal jet bundle by $\J^k(P)$.

As with the usual jet bundles, there exist the natural \cd operators
\[
\hj_k\colon P\to\J^k(P),
\]
and the natural projections
$\nu_{k,l}\colon\J^k(P)\to\J^l(P)$ such that
$\nu_{k,l}\circ\hj_k=\hj_l$. The operators $\hj_k$ and $\nu_{k,l}$ are
restrictions of the operators $\CC j_k$ and $\CC\pi^*_{k,l}$ to $\Ei$.

\cd operators, horizontal forms and jets constitute a ``subtheory'' in
the differential calculus on an infinitely prolonged equation. It is,
roughly speaking, ``the total derivatives calculus'' and is called
\emph{\cd calculus}. It is easily shown that all components of usual
calculus and the Lagrangian formalism discussed above have their
counterparts in the framework of \cd calculus. All constructions of
Sections \ref{sec:calc} and \ref{sec:lagr} are carried over into \cd
calculus word for word as long as the operators, jets, and forms in
them are assumed respectively \cd and horizontal.

\subsection{Higher symmetries}\label{sub:hsym}
Let $\pi\colon E\to M$ be a vector bundle and $\CE\sbs J^k(\pi)$ be a differential
equation. We shall still assume $\CE$ to be formally integrable, though it
not restrictive in this context.

Consider a symmetry $F:J^k(\pi)\xra{}J^k(\pi)$ of the equation $\CE$ and let
$\ta_{k+1}$ be a point of the first prolongation
$\CE^1$ such that $\pi_{k+1,k}(\ta_{k+1})=\ta_k\in\CE$. Then the $R$-plane
$L_{\ta_{k+1}}$ is taken to the $R$-plane $F_*(L_{\ta_{k+1}})$, since $F$ is a
Lie transformation, and $F_*(L_{\ta_{k+1}})\sbs T_{F(\ta_k)}$, since $F$ is a
symmetry. Consequently, the lifting $F^{(1)}\colon J^{k+1}(\pi)\to J^{k+1}(\pi)$ is
a symmetry of $\CE^1$. By the same reasons, $F^{(l)}$ is a symmetry of the
$l$-th prolongation of $\CE$. From here it also follows that for any
infinitesimal symmetry $X$ of the equation $\CE$, its $l$-th lifting is
is a symmetry of $\CE^l$ as well.
\bpr\label{sec3:pr:symprol}
Symmetries of a formally integrable equation $\CE\sbs J^k(\pi)$ coincide with
symmetries of any prolongation of this equation. The same is valid for
infinitesimal symmetries.
\epr
\begin{proof}
We showed already that to any (infinitesimal) symmetry of $\CE$ there
corresponds an (infinitesimal) symmetry of $\CE^l$. Consider an
(infinitesimal) symmetry of $\CE^l$. By Theorems \vref{sec3:th:lietr} and
\vref{sec3:th:liefld}, it is $\pi_{k+l,k}$-fiberwise and therefore generates
an (infinitesimal) symmetry of the equation $\CE$.
\end{proof}

The result proved means that a symmetry of $\CE$ generates a symmetry of
$\Ei$ which preserves every prolongation of finite order. A natural step
to generalize the concept of symmetry is to consider ``all symmetries''
of $\Ei$. Recall the notation
$$\CC \Dr(\pi)=\{\,X\in \Dr(\pi)\mid X\mbox{ lies in } \CC(\pi)\,\}$$
(cf.\ with \veqref{sec3:eq:CD}).
\bde\label{sec3:df:symcar}
Let $\pi$ be a vector bundle. A vector field $X\in \Dr(\pi)$ is called a {\em
symmetry} of the Cartan distribution $\CC(\pi)$ on $\Ji(\pi)$,
\index{symmetry of the Cartan distribution on $\Ji(pi)$} if $[X,\CC \Dr(\pi)]
\sbs\CC \Dr(\pi)$.
\ede
Thus, the set of symmetries coincides with $\Dr_\CC(\pi)$ (see
\veqref{sec3:eq:DC})
and forms a Lie algebra over $\BBR$ and a module over $\CF(\pi)$. Note that
since the Cartan distribution on $\Ji(\pi)$ is integrable, one has
$\CC \Dr(\pi)\sbs\Dr_\CC(\pi)$ and, moreover, $\CC \Dr(\pi)$ is an ideal in the
Lie algebra $\Dr_\CC(\pi)$.

Note also that symmetries belonging to $\CC \Dr(\pi)$ are tangent to any
integral manifold of the Cartan distribution. By
this reason, we call such symmetries {\em trivial}.\index{trivial symmetry
of the Cartan distribution on $\Ji(pi)$} Respectively, the elements of
the quotient Lie algebra
$$\sym(\pi)=\Dr_\CC(\pi)/\CC \Dr(\pi)$$
are called {\em nontrivial symmetries} of the Cartan distribution on
$\Ji(\pi)$.\index{nontrivial symmetry of the Cartan distribution on
$\Ji(pi)$}

Let now $\Ei$ be the infinite prolongation of an equation $\CE\sbs J^k(\pi)$.
Then, since $\CC \Dr(\pi)$ is spanned by the fields of the form $\CC Y$,
where $Y\in
\Dr(M)$ (see Remark \vref{sec3:re:tottan}), any vector field from $\CC \Dr(\pi)$ is
tangent to $\Ei$. Consequently, either all elements of the coset $[X]=X\bmod
\CC\Dr(\pi)$, $X\in\Dr(\pi)$, are tangent to $\Ei$ or neither of them do.
In the first case we say that the coset $[X]$ is {\em tangent} to $\Ei$.
\bde\label{sec3:df:symEi}
An element $[X]=X\bmod\CC \Dr(\pi)$, $X\in\Dr(\pi)$, is called a
{\em higher symmetry}\index{higher symmetry} of $\CE$, if it is tangent to
$\Ei$.
\ede
The set of all higher symmetries forms a Lie algebra over $\BBR$ and is
denoted by $\sym(\CE)$. We shall usually omit the adjective {\em higher} in
the sequel.

Consider a vector field $X\in \Dr(\pi)$. Then, substituting $X$ into the
structural element $U_\CC$ (see \veqref{sec3:eq:upi}), we obtain a field
$X^v\in \Dr(\pi)$. The correspondence $U_\CC\colon X\mapsto X^v=\ip_X U_\CC$
possesses the following properties:
\begin{enumerate}
\item The field $X^v$ is vertical, i.e., $X^v(\Ci(M))=0$.
\item $X^v=X$ for any vertical field.
\item $X^v=0$ if and only if the field $X$ lies in $\CC \Dr(\pi)$.
\end{enumerate}
Therefore, we obtain the direct sum decomposition of $\CF(\pi)$-modules
$$\Dr(\pi)=\Dr^v(\pi)\oplus\CC \Dr(\pi),$$
where $\Dr^v(\pi)$ denotes the Lie algebra of vertical fields. A direct
corollary of these properties is the following result.
\bpr\label{sec3:pr:vert}
For any coset $[X]\in\sym(\CE)$ there exists a unique vertical
representative and thus
\begin{equation}\label{sec3:eq:vert}
\sym(\CE)=\{\,X\in\Dr^v(\CE)\mid [X,\CC \Dr(\CE)]\sbs\CC \Dr(\CE)\,\},
\end{equation}
where $\CC \Dr(\CE)$ is spanned by the fields $\CC Y$, $Y\in\Dr(M)$.
\epr
\ble\label{sec3:le:vert}
Let $X\in\sym(\pi)$ be a vertical vector field. Then it is
completely determined by its restriction to $\CF_0(\pi)\sbs\CF(\pi)$.
\ele
\begin{proof}
Let $X$ satisfy the conditions of the lemma and $Y\in \Dr(M)$. Then for any
$f\in\Ci(M)$ one has
$$[X,\CC Y](f)=X(\CC Y(f))-\CC Y(X(f))=X(Y(f))=0$$
and hence the commutator $[X,\CC Y]$ is a vertical vector field. On the
other hand, $[X,\CC Y]\in\CC \Dr(\pi)$ because $\CC \Dr(\pi)$ is a Lie algebra
ideal. Consequently, $[X,\CC Y]=0$.
Note now that in special coordinates we have $D_i(u_\si^j)=u_{\si i}^j$
for all $\si$ and $j$. From the above said it follows that
\begin{equation}\label{sec3:eq:Xu}
X(u_{\si i}^j)=D_i(X(u_\si^j)).
\end{equation}
But a vertical derivation is determined by its values at the coordinate
functions $u_\si^j$.
\end{proof}

Let now $X_0\colon \CF_0(\pi)\to\CF(\pi)$ be a derivation. Then equalities
\veqref{sec3:eq:Xu} allow one to reconstruct locally a vertical derivation
$X\in \Dr(\pi)$ satisfying $X\left|_{\CF_0(\pi)}\right.=X_0$. Obviously, the
derivation $X$ lies in $\sym(\pi)$ over the neighborhood
under consideration. Consider two neighborhoods $\CU_1$, $\CU_2\sbs\Ji(\pi)$
with the corresponding special coordinates in each of them and two symmetries
$X^i\in\sym(\pi\left|_{\CU_i}\right.)$, $i=1,2$, arising by the described
procedure. But the restrictions of $X^1$ and $X^2$ to $\CF_0(\pi
\left|_{\CU_1\cap\CU_2}\right.)$ coincide. Hence, by Lemma \vref{sec3:le:vert},
the field $X^1$ coincides with $X^2$ over the intersection $\CU_1\cap
\CU_2$. Hence, the reconstruction procedure $X_0 \mapsto X$ is a
global one. So we have established a one-to-one correspondence between
elements of $\sym(\pi)$ and derivations $\CF_0(\pi)\to\CF(\pi)$.

Note now that due to vector bundle
structure in $\pi\colon E\to M$, derivations $\CF_0(\pi)\to\CF(\pi)$ are
identified with sections of $\pi_\infty^*(\pi)$, or with
elements of $\CF(\pi,\pi)$.
\begin{theorem}\label{sec3:th:sympi}
Let $\pi\colon E\to M$ be a vector bundle. Then:
\begin{enumerate}
\item The $\CF(\pi)$-module $\sym(\pi)$ is in one-to-one
correspondence with elements of the module $\CF(\pi,\pi)$.
\item In special coordinates the correspondence $\CF(\pi,\pi)\to
\sym(\pi)$ is expressed by the formula\footnote{To denote evolutionary
vector fields (see Definition \vref{sec3:df:evfld}), we use the Cyrillic
letter $\re$, which is pronounced like ``e'' in ``ten''.}
\begin{equation}\label{sec3:eq:evol}
\vf\mapsto\re_\vf=\sum_{j,\si}D_\si(\vf^j)\pdr{}{u_\si^j},
\end{equation}
where $\vf=(\vf^1,\dots,\vf^m)$ is the component-wise representation of
the section $\vf\in\CF(\pi,\pi)$.
\end{enumerate}
\end{theorem}
\begin{proof}
The first part of the theorem has already been proved. To prove the second
one, it suffices to use equality \veqref{sec3:eq:Xu}.
\end{proof}
\bde\label{sec3:df:evfld}{Let $\pi\colon E\to M$ be a vector bundle.
\begin{enumerate}
\item The field $\re_\vf$ of the form \veqref{sec3:eq:evol} is called an
\emph{evolutionary vector field}\index{evolutionary vector field} on
$\Ji(\pi)$.
\item The section $\vf\in\CF(\pi,\pi)$ is called the {\em generating
function}\index{generating function} of the field $\re_\vf$.
\end{enumerate}}
\ede
\bre\label{sec3:re:evol}
Let $\zeta\colon N\to M$ be an arbitrary smooth fiber bundle and
$\xi\colon P\to M$ be a vector bundle. Then it easy to show that any
$\zeta$-vertical vector field $X$ on $N$ can be uniquely lifted up to an
$\BBR$-linear map $X^\xi\colon \Ga(\zeta^*(\xi))\to\Ga(\zeta^*(\xi))$ such that
\begin{equation}\label{sec3:eq:Xxi}
X^\xi(f\psi)=X(f)\psi+fX^\xi(\psi),\quad f\in\Ci(N),
\quad\psi\in\Ga(\zeta^*(\xi)).
\end{equation}
In particular, taking $\pi_\infty$ for $\zeta$, for any evolutionary
vector field
$\re_\vf$ we obtain the family of maps $\re_\vf^\xi\colon \CF(\pi,\xi)\to
\CF(\pi,\xi)$ satisfying \veqref{sec3:eq:Xxi}.
\ere

Consider the map $\re_\vf^\pi\colon \CF(\pi,\pi)\to\CF(\pi,\pi)$ and recall
the element $\rho_0\in\CF_0(\pi,\pi)\sbs\CF(\pi,\pi)$
(see Example \vref{sec3:ex:jk}). It can be easily seen that
\begin{equation}\label{sec3:eq:genf}
\re_\vf^\pi(\rho_0)=\vf
\end{equation}
which can be taken for the definition of the generating section.

Let $\re_\vf,\re_\psi$ be two evolutionary vector fields. Then, since
$\sym(\pi)$ is a Lie algebra and by Theorem \vref{sec3:th:sympi}, there exists
a unique section $\{\vf,\psi\}$ satisfying $[\re_\vf,\re_\psi]=
\re_{\{\vf,\psi\}}$.
\bde\label{sec3:df:highj}
The section $\{\vf,\psi\}\in\CF(\pi,\pi)$ is called the ({\em higher})
{\em Jacobi bracket}\index{higher Jacobi bracket} of the sections $\vf,\psi
\in\CF(\pi)$.
\ede
\bpr\label{sec3:pr:highj}
Let $\vf,\psi\in\CF(\pi,\pi)$ be two sections. Then:
\begin{enumerate}
\item $\{\vf,\psi\}=\re_\vf^\pi(\psi)-\re_\psi^\pi(\vf)$.
\item In special coordinates, the Jacobi bracket of $\vf$ and
$\psi$ is expressed by the formula
\begin{equation}\label{sec3:eq:highj}
\{\vf,\psi\}^j=\sum_{\al,\si}\left(D_\si(\vf^\al)\pdr{\psi^j}{u_\si^\al}
-D_\si(\psi^\al)\pdr{\vf^j}{u_\si^\al}\right),
\end{equation}
where superscript $j$ denotes the $j$-th component of the corresponding
section.
\end{enumerate}
\epr
\begin{proof}
To prove (1), let us use \veqref{sec3:eq:genf}:
$$
\{\vf,\psi\}=\re_{\{\vf,\psi\}}^\pi(\rho_0)=\re_\vf^\pi(\re_\psi^\pi(\rho_0))
-\re_\psi^\pi(\re_\vf^\pi(\rho_0))=\re_\vf^\pi(\psi)-\re_\psi^\pi(\vf).
$$
The second statement follows from the first one and from equality
\veqref{sec3:eq:evol}.
\end{proof}

Consider now a nonlinear operator $\De\colon \Ga(\pi)\to\Ga(\xi)$ and let
$\psi_\De\in\CF(\pi,\xi)$
be the corresponding section. Then for any $\vf\in\CF(\pi,\pi)$ the section
$\re_\vf^\xi(\psi_\De)\in\CF(\pi,\xi)$ is defined and we can set
\begin{equation}\label{sec3:eq:ulin}
\ell_\De(\vf)=\re_\vf^\xi(\psi_\De).
\end{equation}
\bde\label{sec3:df:ulin}
The operator $\ell_\De\colon \CF(\pi,\pi)\to\CF(\pi,\xi)$ defined
by \veqref{sec3:eq:ulin} is called the {\em universal linearization operator}
\index{universal linearization operator}\footnote{Cf.\ with the algebraic
definition on page \pageref{sec2:p:ulin}.} of the operator $\De\colon
\Ga(\pi)\to\Ga(\xi)$.
\ede
From the definition and equality \veqref{sec3:eq:evol} we obtain that for a scalar
differential operator
$$\De\colon\vf\mapsto F(x_1,\dots,x_n,\dots,
\frac{\partial^{|\si|}\vf^j}{\partial x_\si},
\dots)$$
one has $\ell_\De=(\ell_\De^1,\dots,\ell_\De^m)$, $m=\dim\pi$, where
\begin{equation}\label{sec3:eq:ulinco}
\ell_\De^\al=\sum_{\si}\pdr{F}{u_\si^\al}D_\si.
\end{equation}
If $\dim\xi=r>1$ and $\De=(\De_1,\dots,\De_r)$, then
\begin{equation}\label{sec3:eq:ulinmat}
\ell_\De=\left\Vert\ell_{\De^\al}^\be\right\Vert,\quad
\al=1,\dots,m,\quad\be=1,\dots,r.
\end{equation}
In particular, we see that the following statement is valid.
\bpr\label{sec3:pr:ulincd}
For any differential operator $\De$, its universal
linearization is a $\CC$-differential operator.
\epr

Now we can describe the algebra $\sym(\CE)$, $\CE\sbs J^k(\pi)$ being a
formally integrable equation. Let $I(\CE)\sbs\CF(\pi)$ be the ideal of the
equation $\CE$ (see page \pageref{sec3:ideal}). Then, by definition,
$\re_\vf$ is a symmetry of $\CE$ if and only if
\begin{equation}\label{sec3:eq:idpres}
\re_\vf(I(\CE))\sbs I(\CE).
\end{equation}
Assume now that $\CE$ is given by a differential operator $\De\colon \Ga(\pi)\to
\Ga(\xi)$ and locally is described by the system of equations
$F^1=0,\dots,F^r=0,\ F^j\in\CF(\pi)$.
Then the functions $F^1,\dots,F^r$ are differential generators of the
ideal $I(\CE)$ and condition \veqref{sec3:eq:idpres} may be rewritten as
\begin{equation}\label{sec3:eq:evsym}
\re_\vf(F^j)=\sum_{\al,\si}a_{\si,j}^\al D_\si(F^\al),\quad
j=1,\dots,m,\quad a_\si^\al\in\CF(\pi).
\end{equation}
Using of \veqref{sec3:eq:ulin}, the last equation acquires the
form\footnote{We use the notation $\ell_F$,
$F\in\CF(\pi,\xi)$, as a synonym for $\ell_\De$, where $\De\colon
\Ga(\pi)\to\Ga(\xi)$ is the operator corresponding to the
section $F$.}
\begin{equation}\label{sec3:eq:evsym1}
\ell_{F^j}(\vf)=\sum_{\al,\si}a_{\si,j}^\al D_\si(F^\al),\quad j=1,\dots,m,
\quad a_\si^\al\in\CF(\pi).
\end{equation}
But by Proposition \vref{sec3:pr:ulincd}, the universal linearization is a
$\CC$-differential operator and consequently can be restricted to $\Ei$
(see Corollary \vref{sec3:co:cdif}). It means that we can rewrite
\veqref{sec3:eq:evsym1} as
\begin{equation}\label{sec3:eq:evsym2}
\ell_{F^j}\left|_{\Ei}\right.(\vf\left|_{\Ei}\right.)=0,\quad j=1,\dots,m.
\end{equation}
Combining these equations with \eqref{sec3:eq:ulinco} and
\veqref{sec3:eq:ulinmat}, we obtain
the following fundamental result:
\begin{theorem}\label{sec3:th:highsym}
\index{structure of $\sym(\Ei)$}Let $\CE\sbs J^k(\pi)$ be a
formally integrable equation
and $\De=\De_\CE\colon \Ga(\pi)\to\Ga(\xi)$ be the operator corresponding to
$\CE$. Then an evolutionary vector field $\re_\vf$, $\vf\in\CF(\pi,\pi)$ is
a symmetry of $\CE$ if and only if
\begin{equation}\label{sec3:eq:evsym3}
\ell_\CE(\bar{\vf})=0,
\end{equation}
where $\ell_\CE$ and $\bar{\vf}$ denote restrictions of $\ell_\De$ and $\vf$
on $\Ei$ respectively. In other words, $\sym(\CE)=\ker\ell_\CE$.
\end{theorem}
\begin{xca}
Show that classical symmetries (see Subsection \ref{sub:csym}) are
embedded in $\sym\CE$ as a Lie subalgebra. Describe their generating
functions.
\end{xca}
\bre\label{sec3:re:symsec}
From the result obtained it follows that higher symmetries
of $\CE$ can be identified with elements of $\CF(\CE,\pi)$ satisfying
equation \veqref{sec3:eq:evsym3}. Below we shall say that a {\em section
$\vf\in
\CF(\CE,\pi)$ is a symmetry} of $\CE$ keeping in mind this identification.
Note that due to the fact that $\sym(\CE)$ is a Lie algebra, for any two
symmetries $\vf,\psi\in\CF(\CE,\pi)$ their Jacobi bracket $\{\vf,\psi\}_\CE
\in\CF(\CE,\pi)$ is well defined and is a symmetry as well. If no confusion
arises, we shall omit the subscript $\CE$ in the notation of the Jacobi
bracket.
\ere

Finally, we give a useful description of the modules $\Dv(\E)$ and
$\CLa{k}(\E)$. Denote $\vk=\F(\E,\pi)$.

First consider the case $\Ei=\Ji(\pi)$. From the coordinate expression
\veqref{sec3:eq:evol} for an evolutionary vector field it immediately
follows that any vertical tangent vector at point $\theta\in\Ji(\pi)$
can be realized in the form $\left.\re_{\vf}\right|_{\theta}$ for some
$\vf$. This shows that the map $\vf\mapsto\re_{\vf}$ yields the
canonical isomorphism
\[
\Dv(\pi)=\hJi(\vk).
\]
The dual isomorphism reads
\[
\CLa{1}(\pi)=\CDiff(\varkappa,\F).
\]
In coordinates, this isomorphism takes the form $\omega_{\sigma}^j$ to
the operator
\[
(0,\dots,0,D_\sigma,0,\dotsc,0),
\]
with $D_{\sigma}$ on $j$-th
place.

It is clear that the Cartan $k$-forms can be identified with
multilinear
skew-symmetric \cd operators in $k$ arguments:
\[
\CLa{p}(\pi)=\CDiff_{(p)}^{\alt}(\varkappa,\F).
\]

Now consider the general case. Suppose that the equation $\E$ is given
by an operator $\Delta\colon\Gamma(\pi)\to\Gamma(\xi)$. Denote
$P=\F(\CE,\xi)$, so that $\ell_\E\colon\vk\to P$. From
\veqref{sec3:eq:evsym2} we get
\begin{proposition}\label{sec3:pr:fi}
\begin{enumerate}
\item The module $\Dv(\E)$ is isomorphic to the kernel of the
homomorphism $\psi_{\infty}^{\ell_\E}\colon\hJi(\varkappa)
\xra{}\hJi(P)$\textup{;}
\item the module $\CLa{p}(\E)$ is isomorphic to
$\CDiff_{(p)}^{\alt}(\varkappa,\F)$ modulo the submodule consisting of
the operators of the form $\nabla\circ\ell_\E$, where
$\nabla\in\CDiff(P,\CDiff_{(p-1)}^{\alt}(\varkappa,\F))$.
\end{enumerate}
\end{proposition}
\newpage

\section{Coverings and nonlocal symmetries}
\label{sec:nonloc}
The facts exposed in this section constitute a formal base to introduce
nonlocal variables to the differential setting, i.e., variables of the
type $\int\vf\,dx$, $\vf$ being a function on an infinitely prolonged
equation. A detailed exposition of this material can be found in
\cite{KrasVin1} and \cite{Symm}.
\subsection{Coverings}\label{sub:cov}
We start with fixing up the setting. To do it, extend the universum of
infinitely prolonged equations in the following way. Let $\CN$ be a chain of
smooth maps $\dotsb\xra{} N^{i+1}\xra{\tau_{i+1,i}}N^i\xra{}\dotsb$,
where $N^i$ are smooth finite-dimensional manifolds. Define the algebra $\CF(\CN)$ of
smooth functions on $\CN$ as the direct limit of the homomorphisms
$\dotsb\xra{}\Ci(N^i) \xra{\tau_{i+1,i}^*}\Ci(N^{i+1})\xra{}\dotsb$.
Then there exist natural homomorphisms
$\tau_{\infty,i}^*\colon \Ci(N^i)\to\CF(\CN)$ and the algebra $\CF(\CN)$ may
be considered to be filtered by the images of these maps.  Let us
consider calculus (see Section \ref{sec:calc}) over $\CF(\CN)$ agreed
with this filtration. Define the category\label{sec4:pg:catInf}
\textit{Inf} as follows:

\begin{enumerate}
\item The \emph{objects} of {\it
Inf} are the above introduced chains $\CN$ endowed with
\emph{integrable distributions} $\cD_\CN\sbs \Dr(\CF(\CN))$, where the
word ``integrable'' means that $[\cD_\CN,\cD_\CN] \sbs\cD_\CN$.
\item If $\CN_1=\{N_1^i,\tau^1_{i+1,i}\}$, $\CN_2=\{N_2^i,\tau^2_{i+1,i}\}$
are two objects of \textit{Inf}, then a \emph{morphism} $\vf\colon \CN_1\to\CN_2$ is a
system of smooth maps $\vf_i\colon N_1^{i+\al}\to N_2^i$, where
$\al\in\mathbb{Z}$ is independent of $i$, satisfying
$\tau_{i+1,i}^2\circ\vf_{i+1}=\vf_i\circ\tau_{i+\al+1,i+\al}^1$ and such
that $\vf_{*,\ta}(\cD_{\CN_1,\ta})\sbs\cD_{\CN_2,\vf(\ta)}$ for any point
$\ta\in\CN_1$.
\end{enumerate}
\bde\label{sec4:df:cov}
A morphism $\vf\colon \CN_1\to\CN_2$ is called a \emph{covering} in the category
\textit{Inf}, if $\vf_{*,\ta}\rest{\cD_{\CN_1,\ta}}\colon\cD_{\CN_1,\ta}\to
\cD_{\CN_2,\vf(\ta)}$ is an isomorphism for any point $\ta\in\CN_1$.
\ede
In particular, manifolds $\Ji(\pi)$ and $\Ei$ endowed with the corresponding
Cartan distributions are objects of \textit{Inf} and we can consider coverings
over these objects.
\begin{example}\label{sec4:ex:oper}
Let $\De\colon \Ga(\pi)\to\Ga(\pi')$ be a differential operator of order $\le
k$. Then the system of maps $\Phi_\De^{(l)}\colon J^{l+l}(\pi)\to J^l(\pi')$ (see
the proof of Lemma \vref{sec3:le:prol}) is a morphism of $\Ji(\pi)$ to
$\Ji(\pi')$. Under unrestrictive conditions of regularity, its image is of
the form $\Ei$ for some equation $\CE$ in the bundle $\pi'$ while the map
$\Ji(\pi)\to\Ei$ is a covering.
\end{example}
\bde\label{sec4:df:equviv}
Let $\vf'\colon \CN'\to\CN$ and $\vf''\colon \CN''\to\CN$ be two coverings.
\begin{enumerate}
\item A morphism $\psi\colon \CN'\to\CN''$ is said to be a \emph{morphism of
coverings}, if $\vf'=\vf''\circ\psi$.
\item
The coverings $\vf',\vf''$ are called \emph{equivalent}, if there exists
a morphism $\psi\colon \CN'\to\CN''$ which is a diffeomorphism.
\end{enumerate}
\ede
\bde\label{sec4:df:lincov}
A covering $\vf:\CN'\xra{}\CN$ is called \emph{linear}, if
\begin{enumerate}
\item $\vf$ is a linear bundle;
\item any element $X\in\cD(\CN')$ preserves the submodule of fiber-wise
linear (with respect to the projection $\vf$) functions in $\CF(\CN')$.
\end{enumerate}
\ede
Let $\CN$ be an object of \textit{Inf} and $W$ be a smooth manifold. Consider
the projection $\tau_W\colon \CN\times W\to\CN$ to the first factor. Then we
can make a covering of $\tau_W$ by lifting the distribution $\cD_\CN$ to
$\CN\times W$ in a trivial way.
\bde\label{sec4:df:triv}
A covering $\tau\colon \CN'\to\CN$ is called \emph{trivial}, if it is equivalent
to the covering $\tau_W$ for some $W$.
\ede
Let again $\vf'\colon\CN'\to\CN$, $\vf''\colon\CN''\to\CN$ be two coverings.
Consider the commutative diagram
$$\begin{CD}
\CN'\times_\CN\CN''@>{\vf''}^*(\vf')>>\CN''\\
@V{\vf'}^*(\vf'')VV@VV\vf''V\\
\CN'@>\vf'>>\CN
\end{CD}$$
where
$$\CN'\times_\CN\CN''=\{\,(\ta',\ta'')\in\CN'\times\CN''\mid\vf'(\ta')=
\vf''(\ta'')\,\}$$
while ${\vf'}^*(\vf'')$, ${\vf''}^*(\vf')$ are natural
projections. The manifold $\CN'\times_\CN\CN''$ is supplied with a natural
structure of an object of \textit{Inf} and the maps $(\vf')^*(\vf'')$,
$(\vf'')^*(\vf')$ become coverings.
\bde\label{sec4:df:Whitney}
The composition
$$
\vf'\times_\CN\vf''=\vf'\circ{\vf'}^*(\vf'')=\vf''\circ{\vf''}^*(\vf')\colon
\CN'\times_\CN\CN''\to\CN
$$
is called the \emph{Whitney product} of the coverings $\vf'$ and $\vf''$.
\ede
\bde\label{sec4:df:red}
A covering is said to be \emph{reducible}, if it is equivalent to a covering
of the form $\vf\times_\CN\tau$, where $\tau$ is a trivial covering.
Otherwise it is called \emph{irreducible}.
\ede
From now on, all coverings under consideration will be assumed to be smooth
fiber bundles. The fiber dimension is called the \emph{dimension of the
covering} $\vf$ under consideration and is denoted by $\dim\vf$.
\bpr\label{sec4:pr:covconn}
Let $\CE\sbs J^k(\pi)$ be an equation in the bundle $\pi\colon E\to M$ and $\vf\colon \CN
\to\Ei$ be a smooth fiber bundle. Then the following statements are
equivalent:
\begin{enumerate}
\item The bundle $\vf$ is equipped with a structure of a covering.
\item There exists a connection $\CC^\vf$ in the bundle $\pi_\infty\circ\vf
\colon \CN\to M$, $\CC^\vf\colon X\mapsto X^\vf$, $X\in\Dr(M)$, $X^\vf\in
\Dr(\CN)$, such that
\begin{enumerate}
\item[(a)] $[X^\vf,Y^\vf]=[X,Y]^\vf$, i.e., $\CC^\vf$ is
flat, and
\item[(b)] any vector field $X^\vf$ is projectible to $\Ei$ under
$\vf_*$ and $\vf_*(X^\vf)=\CC X$, where $\CC$ is the Cartan connection on
$\Ei$.
\end{enumerate}
\end{enumerate}
\epr
The proof reduces to the check of definitions.

Using this result, we shall now obtain coordinate description of coverings.
Namely, let $x_1,\dots,x_n,u^1,\dots,u^m$ be local coordinates in $J^0(\pi)$
and assume that internal coordinates in $\Ei$ are chosen. Suppose also that
over the neighborhood under consideration the bundle $\vf\colon \CN\to\Ei$ is
trivial with the fiber $W$ and $w^1,w^2,\dots,w^s,\dots$ are local
coordinates in $W$.
The functions $w^j$ are called \emph{nonlocal coordinates} in the covering
$\vf$. The connection $\CC^\vf$ puts into correspondence
to any partial derivative $\partial/\partial x_i$ the vector field
$\CC^\vf(\partial/\partial x_i)=\tilde{D_i}$. By Proposition
\vref{sec4:pr:covconn}, these vector fields are to be of the form
\begin{equation}\label{sec4:eq:totder}
\tilde{D_i}=D_i+X_i^v=D_i+\sum_\al X^\al_i\pdr{}{w^\al},\quad i=1,\dots,n,
\end{equation}
where $D_i$ are restrictions of total derivatives to $\Ei$, and satisfy
the conditions
\begin{align}\label{sec4:eq:totcomm}
\begin{split}
[\tilde{D_i},\tilde{D_i}]=[D_i,D_j]+[D_i,X_j^v]&+[X_i^v,D_j]+[X_i^v,X_j^v]\\
&=[D_i,X_j^v]+[X_i^v,D_j]+[X_i^v,X_j^v]=0
\end{split}
\end{align}
for all $i,j=1,\dots,n$.

We shall now prove a number of facts that simplify checking of triviality
and equivalence of coverings.
\bpr\label{sec4:pr:inv}
Let $\vf_1\colon \CN_1\to\Ei$ and $\vf_2\colon \CN_2\to\Ei$ be two coverings of the same
dimension $r<\infty$. They are equivalent if and only if there exists a
submanifold $X\sbs\CN_1\times_{\Ei}\CN_2$ such that
\begin{enumerate}
\item $\codim X=r$\textup{;}
\item The restrictions $\vf_1^*(\vf_2)\left|_X\right.$ and
$\vf_2^*(\vf_1)\left|_{X}\right.$ are surjections.
\item[(3)] $(\cD_{\CN_1\times_{\Ei}\CN_2})_\ta\sbs T_\ta X$ for any point
$\ta\in X$.
\end{enumerate}
\epr
\begin{proof}
In fact, if $\psi\colon \CN_1\to\CN_2$ is an equivalence, then its graph $G_\psi=
\{\,(y,\psi(y))\mid y\in\CN_1\,\}$ is the needed manifold $X$.
Conversely, if $X$ is a manifold satisfying (1)--(3), then the
correspondence
$$y\mapsto\vf_1^*(\vf_2)\big((\vf_1^*(\vf_2))^{-1}(y)\cap X\big)$$
is an equivalence.
\end{proof}
Submanifolds $X$ satisfying assumption (3) of the previous proposition are
called \emph{invariant}.
\bpr\label{sec4:pr:equiv}
Let $\vf_1\colon \CN_1\to\Ei$ and $\vf_2\colon \CN_2\to\Ei$ be two irreducible coverings
of the same dimension $r<\infty$. Assume that the Whitney product of $\vf_1$
and $\vf_2$ is reducible and there exists an invariant submanifold $X$ in
$\CN_1\times_{\Ei}\CN_2$ of codimension $r$. Then $\vf_1$ and $\vf_2$ are
equivalent almost everywhere.
\epr
\begin{proof}
Since $\vf_1,\vf_2$ are irreducible, $X$ is to be mapped surjectively almost
everywhere by $\vf_1^*(\vf_2)$ and $\vf_2^*(\vf_1)$ to $\CN_1$ and $\CN_2$
respectively (otherwise, their images would be invariant submanifolds).
Hence, the coverings are equivalent by Proposition \vref{sec4:pr:inv}.
\end{proof}
\bco\label{sec4:co:onedim}
If $\vf_1$ and $\vf_2$ are one-dimensional coverings over $\Ei$ and their
Whitney product is reducible, then they are equivalent.
\eco
\bpr\label{sec4:pr:kertot}
Let $\vf\colon \CN\to\Ei$ be a covering and $\CU\sbs\Ei$ be a domain such that the
the manifold $\tilde{\CU}=\vf^{-1}(\CU)$ is represented in the form $\CU
\times\BBR^r$, $r\le\infty$, while $\vf|_{\tilde{\CU}}$ is the
projection to the first factor. Then the covering $\vf$ is locally
irreducible if the system
\begin{equation}\label{sec4:eq:kertot}
D_1^\vf(f)=0,\dots,D_n^\vf(f)=0
\end{equation}
has only constant solutions.
\epr
\begin{proof}
Suppose that there exists a solution $f\neq\cnst$ of
\veqref{sec4:eq:kertot}.  Then, since the only solutions of the system
$$D_1(f)=0,\dots,D_n(f)=0,$$
where $D_i$ is the restriction of the $i$-th total derivative to $\Ei$, are
constants, $f$ depends on one nonlocal variable $w^\al$ at least. Without
loss of generality we may assume that $\partial f/\partial w^1\neq0$ in a
neighborhood $\CU'\times V$, $\CU'\sbs\CU$, $V\sbs\BBR^r$. Define the
diffeomorphism $\psi\colon \CU'\sbs\CU\to\psi(\CU'\sbs\CU)$ by setting
$$\psi(\dots,x_i,\dots,p_\si^j,\dots,w^\al,\dots)=
(\dots,x_i,\dots,p_\si^j,\dots,f,w^2,\dots,w^\al,\dots).
$$
Then $\psi_*(D_i^\vf)=D_i+\sum_{\al>1}X_i^\al\partial/\partial w^\al$ and
consequently $\vf$ is reducible.

Let now $\vf$ be a reducible covering, i.e., $\vf=\vf'\times_{\Ei}\tau$,
where $\tau$ is trivial. Then, if $f$ is a smooth function on the total
space of the covering $\tau$, the function $f^*=\big(\tau^*(\vf')\big)^*(f)$
is a solution of \veqref{sec4:eq:kertot}. Obviously, there exists an
$f$ such that $f^*\neq\cnst$.
\end{proof}

\subsection{Nonlocal symmetries and shadows}\label{sub:nsym}
Let $\CN$ be an object of
\textit{Inf} with the integrable distribution $\cD=\cD_\CN$. Define
$$\Dr_\cD(\CN)=\{\,X\in\Dr(\CN)\mid [X,\cD]\sbs\cD\,\}$$
and set $\sym\CN=\Dr_\cD(\CN)/\cD_\CN$. Obviously, $\Dr_\cD(\CN)$ is a Lie
$\BBR$-algebra and $\cD$ is its ideal. Elements of the Lie algebra $\sym\CN$
are called \emph{symmetries} of the object $\CN$.
\bde\label{sec4:df:nonloc}
Let $\vf\colon \CN\to\Ei$ be a covering. A \emph{nonlocal $\vf$-symmetry} of $\CE$
is an element of $\sym\CN$. The Lie algebra of such symmetries is denoted
by $\sym_\vf\CE$.
\ede
A base for computation of nonlocal symmetries is the following two results.
\bth\label{sec4:th:non-vert}
Let $\vf\colon \CN\to\Ei$ be a covering.
The algebra $\sym_\vf\CE$ is isomorphic to the Lie algebra of vector fields
$X$ on $\CN$ such that
\begin{enumerate}
\item $X$ is vertical, i.e., $X(\vf^*(f))=0$ for any function $f\in
\Ci(M)\sbs\CF(\CE)$\textup{;}
\item $[X,D_i^\vf]=0$, $i=1,\dots,n$.
\end{enumerate}
\ethm
\begin{proof}
Note that the first condition means that in coordinate representation the
coefficients of the field $X$ at all $\partial/\partial x_i$ vanish. Hence
the intersection of the set of vertical fields with $\cD$ vanish. On the
other hand, in any coset $[X]\in\sym_\vf\CE$ there exists one and only
one vertical element $X^v$. In fact, let $X$ be an arbitrary element of
$[X]$. Then $X^v=X-\sum_ia_iD_i^\vf$, where $a_i$ is the coefficient of $X$
at $\partial/\partial x_i$.
\end{proof}
\bth\label{sec4:th:non-ulin}
Let $\vf\colon \CN=\Ei\times\BBR^r\to\Ei$ be the covering locally determined by the
fields
$$
D_i^\vf=D_i+\sum_{\al=1}^rX_i^\al\pdr{}{w^\al},\quad i=1,\dots,n,
\quad X_i^\al\in
\CF(\CN),
$$
where $w^1,w^2,\dots$ are coordinates in $\BBR^r$ \emph{(}nonlocal
variables\emph{)}. Then any nonlocal $\vf$-symmetry of the equation $\CE=
\{F=0\}$ is of the form
\begin{equation}\label{sec4:eq:evnon}
\tilde{\re}_{\psi,a}=\tilde{\re}_\psi+\sum_{\al=1}^ra_\al\pdr{}{w^\al},
\end{equation}
where $\psi=\psi^1,\dots,\psi^m$, $a=(a^1,\dots,a^r)$, $\psi^i$, $a^\al\in
\CF(\CN)$ are functions satisfying the conditions
\begin{gather}
\tilde{\ell}_F(\psi)=0,\label{sec4:eq:linnon}\\
\tilde{D}_i(a^\al)=\tilde{\re}_{\psi,a}(X_i^\al)\label{sec4:eq:condnon}
\end{gather}
while
\begin{equation}\label{sec4:eq:evtil}
\tilde{\re}_\psi=\sum_{j,\si}\tilde{D}_\si(\psi)\pdr{}{u_\si^j}
\end{equation}
and $\tilde{\ell}_F$ is obtained from $\ell_F$ by changing total derivatives
$D_i$ for $D_i^\vf$.
\ethm
\begin{proof}
Let $X\in\sym_\vf\CE$. Using Theorem \vref{sec4:th:non-vert}, let us write
down the field $X$ in the form
\begin{equation}\label{sec4:eq:decomp}
X={\sum_{\si,j}^{}}'b_\si^j\pdr{}{u_\si^j}+\sum_{\al=1}^ra^\al\pdr{}{w^\al},
\end{equation}
where ``prime'' over the first sum means that the summation extends on
internal coordinates in $\Ei$ only. Then, equaling to zero the coefficient at
$\partial/\partial u_\si^j$ in the commutator $[X,D_i^\vf]$, we obtain
the following equations
$$
D_i^\vf(b_\si^j)=\left\{\begin{array}{ll}
b_{\si i}^j, &\text{ if $b_{\si i}^j$ is an internal coordinate,}\\[1mm]
X(u_{\si i}^j) &\text{ otherwise}.
\end{array}\right.
$$
Solving these equations, we obtain that the first summand in
\veqref{sec4:eq:decomp} is of the form $\tilde{\re}_\psi$, where $\psi$
satisfies \veqref{sec4:eq:linnon}.
\end{proof}

Comparing the result obtained with the description on local symmetries (see
Theorem \vref{sec3:th:highsym}), we see that in the nonlocal setting an
additional obstruction arises represented by equation
\veqref{sec4:eq:condnon}. Thus, in general, not every solution of
\eqref{sec4:eq:linnon} corresponds to a nonlocal $\vf$-symmetry. We
call vector fields $\tilde{\re}_\psi$ of the form
\eqref{sec4:eq:evtil}, where $\psi$ satisfies equation
\eqref{sec4:eq:linnon}, \emph{$\vf$-shadows}.  In the next subsection
it will be shown that for any $\vf$-shadow $\tilde{\re}_\psi$ there
exists a covering $\vf'\colon \CN'\to\CN$ and a nonlocal
$\vf\circ\vf'$-symmetry $S$ such that $\vf'_*(S)=\tilde{\re}_\psi$.

\subsection{Reconstruction theorems}\label{sub:recth}
Let $\CE\sbs J^k(\pi)$ be a differential equation. Let us first establish
relations between horizontal cohomology of $\CE$ (see Definition
\vref{sec3:df:hordR}) and coverings over $\Ei$. All constructions below are
realized in a local chart $\CU\sbs\Ei$.

Consider a horizontal 1-form $\om=\sum_{i=1}^nX_i\,dx_i\in\hL^1(\CE)$ and
define on the space $\Ei\times\BBR$ the vector fields
\begin{equation}\label{sec4:eq:formcov}
D_i^\om=D_i+X_i\partial/\partial w,\,X_i\in\CF(\CE),
\end{equation}
where $w$ is a coordinate along
$\BBR$. By direct computations, one can easily see that the conditions
$[D_i^\om,D_j^\om]=0$ fulfill if and only if $\hd\om=0$. Thus,
\eqref{sec4:eq:formcov} determines a covering structure in the bundle $\vf:
\Ei\times\BBR\to\Ei$ and this covering is denoted by $\vf^\om$. It is also
obvious that the covering $\vf^\om$ and $\vf^{\om'}$ are equivalent if and
only if the forms $\om$ and $\om'$ are cohomologous, i.e., if $\om-\om'=
\hd f$ for some $f\in\CF(\CE)$.

Let $[\om_1],\dots,[\om^\al],\dots$ be an $\BBR$-basis of the vector space
$\hH^1(\CE)$. Let us define the covering $\mathfrak{a}_{1,0}:\CA^1(\CE)\to\Ei$ as
the Whitney product of all $\vf^{\om_\al}$. It can be shown that the
equivalence class of $\mathfrak{a}_{1,0}$ does not depend on the basis choice.
Now, literary in the same manner as it was done in Definition
\vref{sec3:df:hordR} for $\Ei$, we can define horizontal cohomology for
$\CA^1(\CE)$ and construct the covering $\mathfrak{a}_{2,1}:\CA^2(\CE)\to
\CA^1(\CE)$, etc.
\bde\label{sec4:df:unab}
The inverse limit of the chain
\begin{equation}\label{sec4:eq:abtower}
\dots\to\CA^k(\CE)\xra{\mathfrak{a}_{k,k-1}}\CA^{k-1}(\CE)\to\dots\to
\CA^1(\CE)\xra{\mathfrak{a}_{1,0}}\Ei
\end{equation}
is called the \emph{universal Abelian covering} of the equation $\CE$ and
is denoted by $\mathfrak{a}:\CA(\CE)\to\Ei$.
\ede
Obviously, $\hH^1(\CA(\CE))=0$.
\bth[see \cite{Khorkova1}]\label{sec4:th:unab}
Let $\mathfrak{a}:\CA(\CE)\to\Ei$ be the universal Abelian covering of the
equation $\CE=\{F=0\}$. Then any $\mathfrak{a}$-shadow reconstructs up to a
nonlocal $\mathfrak{a}$-symmetry, i.e., for any solution $\psi=(\psi^1,\dots,
\psi^m)$, $\psi^j\in\CF(\CA(\CE))$, of the equation $\tilde{\ell}_F(\psi)=0$
there exists a set of functions $a=(a_{\al,i})$, where $a_{\al,i}\in
\CF(\CA(\CE))$ such that
$\tilde{\re}_{\psi,a}$ is a nonlocal $\mathfrak{a}$-symmetry.
\ethm
\begin{proof}
Let $w^{j,\al}$, $j\le k$, be nonlocal variables in $\CA^k(\CE)$ and
assume that the covering structure in $\mathfrak{a}$ is determined by the
vector fields $D_i^\mathfrak{a}=D_i+\sum_{j,\al}X_i^{j,\al}\partial/\partial
w^{j,\al}$, where, by construction, $X_i^{j,\al}\in\CF(\CA^{j-1}(\CE))$,
i.e., the functions $X_i^{j,\al}$ do not depend on $w^{k,\al}$ for all $k\ge
j$.

Our aim is to prove that the system
\begin{equation}\label{sec4:eq:recsh}
D_i^\mathfrak{a}(a_{j,\al})=\tilde{\re}_{\psi,a}(X_i^{j,\al})
\end{equation}
is solvable with respect to $a=(a_{j,\al})$ for any $\psi\in\ker
\tilde{\ell}_F$. We do this by induction on $j$. Note that
$$
[D_i^\mathfrak{a},\tilde{\re}_{\psi,a}]=\sum_{j,\al}
\big(D_i^\mathfrak{a}(a_{j,\al})-\tilde{\re}_{\psi,a}(X_i^{j,\al})\big)
\pdr{}{w^{j,\al}}
$$
for any set of functions $(a_{j,\al})$. Then for $j=1$ one has
$[D_i^\mathfrak{a},\tilde{\re}_{\psi,a}](X_k^{1,\al})=0$,
or
$$
D_i^\mathfrak{a}\big(\tilde{\re}_{\psi,a}(X_k^{1,\al})\big)=
\tilde{\re}_{\psi,a}\big(D_i^\mathfrak{a}(X_k^{1,\al})\big),
$$
since $X_k^{1,\al}$ are functions on $\Ei$.

But from the construction of the covering $\mathfrak a$ one has
$D_i^\mathfrak{a}(X_k^{1,\al})=D_k^\mathfrak{a}(X_i^{1,\al})$, and we finally
obtain
$$
D_i^\mathfrak{a}\big(\re_\psi(X_k^{1,\al})\big)=
D_k^\mathfrak{a}\big(\re_\psi(X_i^{1,\al})\big).
$$
Note now that the equality $\hH^1(\CA(\CE))=0$ implies existence of functions
$a_{1,\al}$ satisfying
$$D_i^\mathfrak{a}(a_{1,\al})=\re_\psi(X_i^{1,\al}),$$
i.e., equation \eqref{sec4:eq:recsh} is solvable for $j=1$.

Assume now that solvability of \eqref{sec4:eq:recsh} was proved for $j<s$ and
the functions $(a_{1,\al},\dots,a_{j-1,\al})$ are some solutions. Then, since
$[D_i^\mathfrak{a},\tilde{\re}_{\psi,a}]\left|_{\CA^{j-1}(\CE)}\right.=0$,
we obtain the needed $a_{j,\al}$ literally repeating the proof for the case
$j=1$.
\end{proof}

Let now $\vf\colon\CN\xra{}\Ei$ be an arbitrary covering. The next result
shows that any $\vf$-shadow is reconstructable.
\bth[see also \cite{Kiso1}]\label{sec4:th:recshad}
For any $\vf$-shadow, i.e., for a solution $\psi=(\psi^1,\dots,\psi^m)$,
$\psi^j\in\CF(\CN)$, of
the equation $\tilde{\ell}_F(\psi)=0$, there exists a covering $\vf_\psi
\colon\CN_\psi\xra{}\CN\xra{\vf}\Ei$ and a $\vf_\psi$-symmetry $S_\psi$,
such that $S_\psi\left|_{{\Ei}}\right.=\tilde{\re}_\psi\left|_{{\Ei}}
\right.$.
\ethm
\begin{proof}
Let locally the covering $\vf$ be represented by the vector fields
$$
D_i^\vf=D_i+\sum_{\al=1}^rX_i^\al\pdr{}{w^\al},
$$
$r\le\infty$ being the dimension of $\vf$. Consider the space
$\mathbb{R}^\infty$ with the coordinates $w_l^\al$, $\al=1,\dots,r$, $l=0,1,
2,\dots$, $w_0^\al=w^\al$, and set $\CN_\psi=\CN\times\mathbb{R}^\infty$ with
\begin{equation}\label{sec4:eq:totphi}
D_i^{\vf_\psi}=D_i+\sum_{l,\al}\left(\tilde{\re}_\psi+S_w\right)^l(X_i^\al)
\pdr{}{w_l^\al},
\end{equation}
where
\begin{equation}\label{sec4:eq:evphi}
\tilde{\re}_\psi={\sum_{\si,k}^{}}'D_\si^\vf(\psi^k)\pdr{}{u_\si^k},\quad
S_w=\sum_{\al,l}w_{l+1}^\al\pdr{}{w_l^\al}
\end{equation}
and ``prime'', as before, denotes summation over internal coordinates.

Set $S_\psi=\tilde{\re}_\psi+S_w$. Then
\begin{multline*}
[S_\psi,D_i^{\vf_\psi}]={\sum_{\si,k}^{}}'\tilde{\re}_\psi(\bar{u}_{\si i}^k)
\pdr{}{u_\si^k}+\sum_{l,\al}\left(\tilde{\re}_\psi+S_w\right)^{l+1}(X_i^\al)
\pdr{}{w_l^\al}\\
-{\sum_{\si,k}^{}}'D_i^{\vf_\psi}(D_\si^\vf(\psi^k))\pdr{}{u_\si^k}-
\sum_{l,\al}\left(\tilde{\re}_\psi+S_w\right)^{l+1}(X_i^\al)\pdr{}{w_l^\al}\\
={\sum_{\si,k}^{}}'\left(\tilde{\re}_\psi(\bar{u}_{\si i}^k)-
D_{\si i}^\vf(\psi^k)\right)\pdr{}{u_\si^k}=0.
\end{multline*}
Here, by definition, $\bar{u}_{\si i}^k=D_i^\vf(u_\si^k)\left|_{\CN}\right.$.

Now, using the above proved equality, one has
\begin{multline*}
[D_i^{\vf_\psi},D_j^{\vf_\psi}]=\sum_{l,\al}\left(
D_j^{\vf_\psi}\big(\tilde{\re}_\psi+S_w\big)^l(X_j^\al)-
D_j^{\vf_\psi}\big(\tilde{\re}_\psi+S_w\big)^l(X_i^\al)
\right)\pdr{}{w_l^\al}\\
=\sum_{l,\al}\big(\tilde{\re}_\psi+S_w\big)^l\big(
D_i^{\vf_\psi}(X_j^\al)-D_j^{\vf_\psi}(X_i^\al)\big)\pdr{}{w_l^\al}=0,
\end{multline*}
since $D_i^{\vf_\psi}(X_j^\al)-D_j^{\vf_\psi}(X_i^\al)=
D_i^\vf(X_j^\al)-D_j^\vf(X_i^\al)=0$.
\end{proof}

Let now $\vf\colon\CN\xra{}\Ei$ be a covering and $\vf'\colon\CN'\xra{}\CN
\xra{\vf}\Ei$ be another one. Then, by obvious reasons, any $\vf$-shadow
$\psi$ is a $\vf'$-shadow as well. Applying the construction of Theorem
\ref{sec4:th:recshad} to both $\vf$ and $\vf'$, we obtain two coverings,
$\vf_\psi$ and $\vf'_\psi$ respectively.
\ble\label{sec4:le:tworecs}
The following commutative diagram of coverings
$$\begin{CD}
\CN'_\psi@>>>\CN_\psi\\
@VVV@VVV\\
\CN'@>>>\CN@>>>\Ei
\end{CD}$$
takes place. Moreover, if $S_\psi$ and $S'_\psi$ are nonlocal symmetries
corresponding in $\CN_\psi$ and $\CN'_\psi$ constructed by Theorem
\vref{sec4:th:recshad}, then $S'_\psi\left|_{\CF(\CN_\psi)}\right.=S_\psi$.
\ele
\begin{proof}
It suffices to compare expressions \eqref{sec4:eq:totphi} and
\veqref{sec4:eq:evphi} for the coverings $\CN_\psi$ and $\CN'_\psi$.
\end{proof}

As a corollary of Theorem \ref{sec4:th:recshad} and of the previous lemma,
we obtain the following result.
\begin{theorem}\label{sec4:th:recshad1}
Let $\vf\colon\CN\xra{}\Ei$, $\CE=\{\,F=0\,\}$, be an arbitrary covering and
$\psi_1,\dots,\psi_s\in\CF(\CN)$, be solutions of the
equation $\tilde{\ell}_F(\psi)=0$. Then there exists a covering $\vf_\Psi
\colon\CN_\Psi\xra{}\CN\xra{\vf}\Ei$ and $\vf_\Psi$-symmetries $S_{\psi_1},
\dots,S_{\psi_s}$, such that $S_{\psi_s}\left|_{{\Ei}}\right.=
\tilde{\re}_{\psi_i}\left|_{{\Ei}}\right.$, $i=1,\dots,s$.
\end{theorem}
\begin{proof}
Consider the section $\psi_1$ and the covering $\vf_{\psi_1}\colon
\CN_{\psi_1}\xra{\bar{\vf}_{\psi_1}}\CN\xra{\vf}\Ei$ together with the symmetry $S_{\psi_1}$
constructed in Theorem \vref{sec4:th:recshad}. Then $\psi_2$ is a
$\vf_{\psi_1}$-shadow and we can construct the covering $\vf_{\psi_1,\psi_2}
\colon\CN_{\psi_1,\psi_2}\xra{\bar{\vf}_{\psi_1,\psi_2}}\CN_{\psi_1}\xra{\vf_{\psi_1}}\Ei$ with the
symmetry $S_{\psi_2}$. Applying this procedure step by step, we obtain the
series of coverings
$$\CN_{\psi_1,\dots,\psi_s}\xra{\bar{\vf}_{\psi_1,\dots,\psi_s}}
\CN_{\psi_1,\dots,\psi_{s-1}}\xra{\bar{\vf}_{\psi_1,\dots,\psi_{s-1}}}\dotsb
\xra{\bar{\vf}_{\psi_1,\psi_2}}\CN_{\psi_1}\xra{\bar{\vf}_{\psi_1}}\CN
\xra{\vf}\Ei.
$$
with the symmetries $S_{\psi_1},\dots,S_{\psi_s}$. But $\psi_1$ is a
$\vf_{\psi_1,\dots,\psi_s}$-shadow and we can construct the covering
$\vf_{\psi_1}\colon\CN_{\psi_1}^{(1)}\xra{}\CN_{\psi_1,\dots,\psi_s}\xra{}
\Ei$ with the symmetry $S_{\psi_1}^{(1)}$ satisfying
$S_{\psi_1}^{(1)}\left|_{\CF(\CN_{\psi_1})}\right.=S_{\psi_1}$ (see Lemma
\ref{sec4:le:tworecs}), etc. Passing to the inverse limit, we obtain the
covering $\CN_\Psi$ we need.
\end{proof}

\newpage

\section{\FN brackets and recursion operators}\label{sec:recop}
We return back to the general algebraic setting of Section \ref{sec:calc}
and extend standard constructions of calculus to form-valued derivations.
It allows us to define \emph{\FN brackets} and introduce a cohomology
theory ($\CC$-cohomologies) associated to commutative algebras with
\emph{flat connections}. Applying this theory to partial differential
equations, we obtain an algebraic description of recursion operators for
symmetries and describe efficient tools to compute these operators.
For technical details, examples and generalizations and we refer the reader
to the papers \cite{Kras1,Kras2,Kras3} and \cite{KrasKers1,KrasKers2,
KrasKers3}.

In Subsection \ref{hc.apkcc:subsec}, $\CC$-cohomologies will be discussed
again in the general framework of horizontal cohomologies with coefficients.

\subsection{Calculus in form-valued derivations}\label{sub:fvcalc}
Let $\Bbbk$ be a field of characteristic zero and $A$ be a commutative
unitary $\Bbbk$-algebra. Let us recall the basic notations from Section
\ref{sec:calc}:
\begin{itemize}
\item $\Dr(P)$ is the module of $P$-valued derivations $A\xra{}P$, where $P$
is an $A$-module;
\item $\Dr_i(P)$ is the module of $P$-valued skew-symmetric
$i$-de\-ri\-va\-tions. In particular, $\Dr_1(P)=\Dr(P)$;
\item $\La^i(A)$ is the module of differential $i$-forms of the algebra
$A$;
\item $d\colon \La^i(A)\xra{}\La^{i+1}(A)$ is the de Rham differential.
\end{itemize}
Recall also that the modules $\La^i(A)$ are representative objects for the
functors $\Dr_i\colon P\Rightarrow\Dr_i(P)$, i.e., $\Dr_i(P)=\Hom_A(\La^i(A),
P)$. The isomorphism $\Dr(P)=\Hom_A(\La^1(A),P)$ can be expressed in more
exact terms: for any derivation $X\colon A\xra{}P$, there exists a uniquely
defined homomorphism $\vf^X\colon\La^1(A)\xra{}P$ satisfying $X=\vf^X\circ
d$. Denote by $\langle Z,\om\rangle\in P$ the value of the derivation $Z\in
\Dr_i(P)$ at $\om\in\La^i(A)$.

Both $\La^*(A)=\bigoplus_{i\ge0}\La^i(A)$ and $\Dr_*(A)=\bigoplus_{i\ge0}\Dr_i(A)$
are endowed with the structures of superalgebras with respect to
the wedge product operation $\wg\colon\La^i(A)\ot\La^j(A)\xra{}\La^{i+j}(A)$
and $\wg:\Dr_i(A)\ot\Dr_j(A)\xra{}\Dr_{i+j}(A)$, the de Rham differential $d\colon
\La^*(A)\xra{}\La^*(A)$ becoming a derivation of $\La^*(A)$. Note also that
$\Dr_*(P)=\bigoplus_{i\ge0}\Dr_i(P)$ is a $\Dr_*(A)$-module.

Using the paring $\langle\cdot,\cdot\rangle$ and the wedge product, we
define the \emph{inner product} (or \emph{contraction}) $\ip_X\om\in
\La^{j-i}(A)$ of $X\in\Dr_i(A)$ and $\om\in\La^j(A)$, $i\le j$, by
setting
\begin{equation}\label{sec5:eq:intpr}
\langle Y,\ip_X\om\rangle=(-1)^{i(j-i)}\langle X\wg Y,\om\rangle,
\end{equation}
where $Y$ is an arbitrary element of $\Dr_{j-i}(P)$, $P$ being an $A$-module.
We formally set $\ip_X\om=0$ for $i>j$. When $i=1$, this definition coincides
with the one given in Section \ref{sec:calc}. Recall that the following duality
is valid:
\begin{equation}\label{sec5:eq:intdual}
\langle X,da\wg\om\rangle=\langle X(a),\om\rangle,
\end{equation}
where $\om\in\La^i(A)$, $X\in\Dr_{i+1}(P)$, and $a\in A$
(see Exercise \vref{sec1:es:dual}). Using the property
\eqref{sec5:eq:intdual}, one can show that
$$
\ip_X(\om\wg\ta)=\ip_X(\om)\wg\ta+(-1)^{X\om}\om\wg\ip_X(\om)
$$
for any $\om,\ta\in\La^*(A)$, where (as everywhere below) the symbol of a
graded object used as the exponent of $(-1)$ denotes the degree of that
object.

We now define the {\em Lie derivative} of $\om\in\La^*(A)$ along $X\in
\Dr_*(A)$ as
\begin{equation}\label{sec5:eq:Lie}
\Ld_X\om=\big(\ip_X\circ d-(-1)^X d\circ\ip_X\big)\om=[\ip_X, d]\om,
\end{equation}
where $[\cdot,\cdot]$ denotes the supercommutator: if $\De,\De':\La^*(A)
\xra{}\La^*(A)$ are
graded derivations, then $[\De,\De']=\De\circ\De'-(-1)^{\De\De'}\De'\circ
\De$. For $X\in\Dr(A)$ this definition coincides with the one given by
equality \veqref{sec1:eq:lie}.

Consider now the graded module $\Dr(\La^*(A))$ of $\La^*(A)$-valued
derivations $A\xra{}\La^*(A)$ (corresponding to form-valued vector
fields---or, which is the same---vector-valued differential forms on
a smooth manifold). Note that the graded structure in $\Dr(\La^*(A))$
is determined by the splitting
$\Dr(\La^*(A))=\bigoplus_{i\ge0}\Dr(\La^i(A))$ and thus elements of
grading $i$ are derivations $X$ such that $\im X\subset\La^i(A)$. We
shall need three algebraic structures associated to $\Dr(\La^*(A))$.
First note that $\Dr(\La^*(A))$ is a graded $\La^*(A)$-module: for any
$X\in\Dr (\La^*(A))$, $\om\in\La^*(A)$ and $a\in A$ we set $(\om\wg
X)a=\om\wg X(a)$.  Second, we can define the inner product
$\ip_X\om\in\La^{i+j-1}(A)$ of $X\in\Dr(\La^i(A))$ and $\om\in\La^j(A)$
in the following way. If $j=0$, we set $\ip_X\om=0$. Then, by induction
on $j$ and using the fact that $\La^*(A)$ as a graded $A$-algebra is
generated by the elements of the form $da$, $a\in A$, we set

\begin{equation}\label{sec5:eq:intpr1}
\ip_X(da\wg\om)=X(a)\wg\om-(-1)^Xda\wg\ip_X(\om),\quad a\in A.
\end{equation}
Finally, we can contract elements of $\Dr(\La^*(A))$ with each other in the
following way:
\begin{equation}\label{sec5:eq:intpr2}
(\ip_XY)a=\ip_X(Ya),\quad X,Y\in\Dr(\La^*(A)),\quad a\in A.
\end{equation}
Three properties of contractions are essential in the sequel.
\begin{proposition}\label{sec5:pr:cont}
Let $X,Y\in\Dr(\La^*(A))$ and $\om,\theta\in\La^*(A)$. Then
\begin{gather}
\ip_X(\om\wg\theta)=\ip_X(\om)\wg\theta+(-1)^{\om(X-1)}\om\wg\ip_X(\theta),
\label{sec5:eq:cont1}\\
\ip_X(\om\wg Y)=\ip_X(\om)\wg Y+(-1)^{\om(X-1)}\om\wg\ip_X(Y),
\label{sec5:eq:cont2}\\
[\ip_X,\ip_Y]=\ip_{\rnij{X}{Y}},\label{sec5:eq:cont3}
\end{gather}
where
\begin{eqnarray}
\rnij{X}{Y}=\ip_X(Y)-(-1)^{(X-1)(Y-1)}\ip_Y(X).\label{sec5:eq:RN}
\end{eqnarray}
\end{proposition}
\begin{proof}
Equality \eqref{sec5:eq:cont1} is a direct consequence of
\eqref{sec5:eq:intpr1}. To prove \eqref{sec5:eq:cont2}, it suffices to use
the definition and expressions \eqref{sec5:eq:intpr2} and \eqref{sec5:eq:cont1}.

Let us prove \eqref{sec5:eq:cont3} now. To do this, note first that due to
\eqref{sec5:eq:intpr2} the equality is sufficient to be checked on elements
$\om\in\La^j(A)$. Let us use induction on $j$. For $j=0$ it holds in a
trivial way. Let $a\in A$; then one has
\begin{multline*}
[\ip_X,\ip_Y](da\wg\om)=
\big(\ip_X\circ\ip_Y-(-1)^{(X-1)(Y-1)}\ip_Y\circ\ip_X\big)(da\wg\om)\\
=\ip_X(\ip_Y(da\wg\om))-(-1)^{(X-1)(Y-1)}\ip_Y(\ip_X(da\wg\om)).
\end{multline*}
But
\begin{multline*}
\ip_X(\ip_Y(da\wg\om))=\ip_X(Y(a)\wg\om-(-1)^Yda\wg\ip_Y\om)\\
=\ip_X(Y(a))\wg\om+(-1)^{(X-1)Y}Y(a)\wg\ip_X\om-
(-1)^Y(X(a)\wg\ip_Y\om\\
-(-1)^Xda\wg\ip_X(\ip_Y\om)),
\end{multline*}
while
\begin{multline*}
\ip_Y(\ip_X(da\wg\om)=\ip_Y(X(a)\wg\om-(-1)^Xda\wg\ip_X\om)\\
=\ip_Y(X(a))\wg\om+(-1)^{X(Y-1)}X(a)\wg\ip_Y\om-
(-1)^X(Y(a)\wg\om\\
-(-1)^Yda\wg\ip_Y(\ip_X\om)).
\end{multline*}
Hence,
\begin{multline*}
[\ip_X,\ip_Y](da\wg\om)=\big(\ip_X(Y(a))-(-1)^{(X-1)(Y-1)}\ip_Y(X(a))\big)
\wg\om\\
+(-1)^{X+Y}da\wg\big(\ip_X(\ip_Y\om)-(-1)^{(X-1)(Y-1)}\ip_Y(\ip_X\om)\big).
\end{multline*}
But, by definition,
\begin{multline*}
\ip_X(Y(a))-(-1)^{(X-1)(Y-1)}\ip_Y(X(a))\\
=(\ip_XY-(-1)^{(X-1)(Y-1)}\ip_YX)(a)=\rnij{X}{Y}(a),
\end{multline*}
whereas
$$
\ip_X(\ip_Y\om)-(-1)^{(X-1)(Y-1)}\ip_Y(\ip_X\om)=\ip_{\rnij{X}{Y}}(\om)
$$
by induction hypothesis.
\end{proof}
\begin{definition}\label{sec5:df:RN}
The element $\rnij{X}{Y}$ defined by equality \eqref{sec5:eq:RN} is called
the \emph{Richardson--Nijenhuis bracket} of elements $X$ and $Y$.
\end{definition}
Directly from Proposition \ref{sec5:pr:cont} we obtain the following
\begin{proposition}\label{sec5:pr:RN}
For any derivations $X,Y,Z\in\Dr(\La^*(A))$ and a form $\om\in\La^*(A)$ one
has
\begin{gather}
\rnij{X}{Y}+(-1)^{(X+1)(Y+1)}\rnij{Y}{X}=0,\label{sec5:eq:RN1}\\
\oint(-1)^{(Y+1)(X+Z)}\rnij{\rnij{X}{Y}}{Z}=0,\label{sec5:eq:RN2}\\
\rnij{X}{\om\wg Y}=\ip_X(\om)\wg Y+(-1)^{(X+1)\om}\om\wg\rnij{X}{Y}.
\label{sec5:eq:RN3}
\end{gather}
\end{proposition}
Here and below the symbol $\oint$ denotes the sum of cyclic permutations.
\bre
Note that Proposition \ref{sec5:pr:RN} means that $\Dr(\La^*(A))^\downarrow$
is a Gerstenhaber algebra with respect to the Richardson--Nijenhuis bracket
\cite{KosmannS}. Here the superscript $\downarrow$ denotes the shift of grading
by 1.
\ere
Similarly to \eqref{sec5:eq:Lie} define the Lie derivative of
$\om\in\La^*(A)$ along $X\in\Dr(\La^*(A))$ by
\begin{equation}\label{sec5:eq:Lied}
\Ld_X\om=(\ip_X\circ d+(-1)^X d\circ\ip_X)\om=[\ip_X,d]\om
\end{equation}
(the change of sign is due to the fact that $\deg(\ip_X)=\deg(X)-1$). From
the properties of $\ip_X$ and $d$ we obtain
\begin{proposition}\label{sec5:pr:Lie}
For any $X\in\Dr(\La^*(A))$ and $\om,\theta\in\La^*(A)$, one has the
following identities:
\begin{gather}
\Ld_X(\om\wg\theta)=\Ld_X(\om)\wg\theta+(-1)^{X\om}\om\wg\Ld_X(\theta),
\label{sec5:eq:Lied1}\\
\Ld_{\om\wg X}=\om\wg\Ld_X+(-1)^{\om+X} d(\om)\wg\ip_X,\label{sec5:eqLied2}\\
[\Ld_X,d]=0.\label{sec5:eq:Lied3}
\end{gather}
\end{proposition}

Our main concern now is to analyze the commutator $[\Ld_X,\Ld_Y]$ of
two Lie derivatives. It may be done efficiently for smooth
algebras (see Definition \vref{sec1:df:smootha}).
\begin{proposition}\label{sec5:pr:FNex}
Let $A$ be a smooth algebra. Then for any derivations $X,Y\in\Dr(\La^*(A))$
there exists a uniquely determined element $\fnij{X}{Y}\in\Dr(\La^*(A))$
such that
\begin{equation}\label{sec5:eq:FN}
[L_X,L_Y]=L_{\fnij{X}{Y}}.
\end{equation}
\end{proposition}
\begin{proof}
To prove existence, recall that for smooth algebras one has
$$
\Dr_i(P)=\Hom_A(\La^i(A),P)=P\ot_A\Hom_A(\La^i(A),A)=P\ot_A\Dr_i(A)
$$
for any $A$-module $P$ and integer $i\ge0$. Using this identification,
represent elements $X,Y\in\Dr(\La^*(A))$ in the form
$$X=\om\ot X'\mathrm{\ and\ }Y=\theta\ot Y'\mathrm{\ for\ }\om,\theta\in
\La^*(A),X',Y'\in\Dr(A).$$
Then it is easily checked that the element
\begin{align}\label{sec5:eq:FNexpl}
Z&=\om\wg\theta\ot[X',Y']+\om\wg\Ld_{X'}\theta\ot Y+(-1)^\om d\om\wg
\ip_{X'}\theta\ot Y'
\nonumber\\
&-(-1)^{\om\theta}\theta\wg\Ld_{Y'}\om\ot X'-(-1)^{(\om+1)\theta} d\theta
\wg\ip_{Y'}\om\ot X'\\
&=\om\wg\theta\ot[X',Y']+\Ld_X\theta\ot Y'-(-1)^{\om\theta}\Ld_Y\om\ot X'
\nonumber
\end{align}
satisfies \eqref{sec5:eq:FN}.

Uniqueness follows from the fact that $\Ld_X(a)=X(a)$ for any $a\in A$.
\end{proof}
\begin{definition}\label{sec5:df:FN}
The element $\fnij{X}{Y}\in\Dr^{i+j}(\La^*(A))$ defined by formula
\eqref{sec5:eq:FN} is called the \emph{\FN bracket} of elements
$X\in\Dr^i(\La^*(A))$ and $Y\in\Dr^j(\La^*(A))$.
\end{definition}
The basic properties of this bracket are summarized in the following
\begin{proposition}\label{sec5:pr:FN}
Let $A$ be a smooth algebra, $X,Y,Z\in\Dr(\La^*(A))$ be derivations and
$\om\in\La^*(A)$ be a differential form. Then the following identities are
valid:
\begin{gather}
\fnij{X}{Y}+(-1)^{XY}\fnij{Y}{X}=0,\label{sec5:eq:FN1}\\
\oint(-1)^{Y(X+Z)}\fnij{X}{\fnij{Y}{Z}}=0,\label{sec5:eq:FN2}\\
\ip_{\fnij{X}{Y}}=[\Ld_X,\ip_Y]+(-1)^{X(Y+1)}\Ld_{\ip_YX},
\label{sec5:eq:FN3}\\
\begin{split}
\ip_Z\fnij{X}{Y}=\fnij{\ip_ZX}{Y}&+(-1)^{X(Z+1)}\fnij{X}{\ip_ZY}\\
&+(-1)^X\ip_{\fnij{Z}{X}}Y-(-1)^{(X+1)Y}\ip_{\fnij{Z}{Y}}X,
\label{sec5:eq:FN4}\end{split}\\
\begin{split}
\fnij{X}{\om\wg Y}=\Ld_X\om\wg Y-(-1)^{(X+1)(Y+\om)}&d\om\wg\ip_YX\\
&+(-1)^{X\om}\om\wg\fnij{X}{Y}.\label{sec5:eq:FN5}
\end{split}
\end{gather}
\end{proposition}
Note that the first two equalities in the previous proposition mean that the
module $\Dr(\La^*(A))$ is a Lie superalgebra with respect to the \FN bracket.
\begin{remark}\label{sec5:rm:inf}
The above exposed algebraic scheme has a geometrical realization,
if one takes $A=C^\infty(M)$, $M$ being a smooth
finite-dimensional manifold.
The algebra $A=C^\infty(M)$ is smooth in this case. However, in the
geometrical theory of differential equations we have to work with
infinite-dimensional manifolds\footnote{Infinite jets, infinite prolongations
of differential equations, total spaces of coverings, etc.}
of the form $N=\projlim_{\{\pi_{k+1,k}\}}N_k$, where all the maps
$\pi_{k+1,k}\colon N_{k+1}\xra{}N_k$ are surjections of finite-dimensional smooth
manifolds. The corresponding algebraic object is a filtered algebra $A=
\bigcup_{k\in\BB{Z}}A_k,\ A_k\subset A_{k+1}$, where all $A_k$ are
subalgebras in $A$. As it was already noted, self-contained differential
calculus over $A$ is
constructed, if one considers the category of all filtered $A$-modules
with filtered homomorphisms for morphisms between them. Then all functors
of differential calculus in this category become filtered, as well as their
representative objects.

In particular, the $A$-modules $\La^i(A)$ are filtered by $A_k$-modules
$\La^i(A_k)$. We say that the algebra $A$ is {\em finitely smooth}, if
$\La^1(A_k)$ is a projective $A_k$-module of finite type for any $k\in
\BB{Z}$. For finitely smooth algebras, elements of $\Dr(P)$ may be represented
as formal infinite sums $\sum_kp_k\ot X_k$, such that any finite sum
$S_n=\sum_{k\le n}p_k\ot X_k$ is a derivation $A_n\xra{} P_{n+s}$ for some
fixed $s\in\BB{Z}$. Any derivation $X$ is completely determined by the
system $\{S_n\}$ and Proposition \ref{sec5:pr:FN} obviously remains valid.
\end{remark}

\subsection{Algebras with flat connections and cohomology}\label{sub:flatcon}
We now introduce the second object of our interest. Let $A$ be an
$\Bbbk$-algebra, $\Bbbk$ being a field of zero characteristic, and $B$ be an
algebra over $A$. We shall assume that the corresponding homomorphism
$\varphi\colon A\xra{}B$ is an embedding. Let $P$ be a $B$-module; then it is
an $A$-module as well and we can consider the $B$-module $\Dr(A,P)$ of
$P$-valued derivations $A\xra{}P$.
\begin{definition}\label{sec5:df:conn}
Let $\na^\bullet\colon\Dr(A,\cdot)\Ra\Dr(\cdot)$ be a natural
transformations of functors $\Dr(A,\cdot)\colon A\Ra\Dr(A,P)$ and
$\Dr(\cdot)\colon P\Ra\Dr(\cdot)$ in the category of $B$-modules, i.e., a
system of homomorphisms $\na^P\colon\Dr(A,P)\xra{}\Dr(P)$ such that the
diagram
$$
\begin{CD}
\Dr(A,P)@>\na^P>>\Dr(P)\\
@V\Dr(A,f)VV    @VV\Dr(f)V\\
\Dr(A,Q)@>\na^Q>>\Dr(Q)\\
\end{CD}
$$
is commutative for any $B$-homomorphism $f\colon P\xra{}Q$. We say that
$\na^\bullet$ is a \emph{connection} in the triad $(A,B,\varphi)$, if
$\na^P(X)\rest{A}=X$ for any $X\in \Dr(A,P)$.
\end{definition}
Here and below we use the notation $Y\rest{A}=Y\circ\varphi$ for any
$Y\in\Dr(P)$.
\begin{remark}\label{sec5:rm:equiv}
When $A=\Ci(M)$, $B=\Ci(E)$, $\varphi=\pi^*$, where $M$ and $E$ are smooth
manifolds and $\pi\colon E\xra{}M$ is a smooth fiber bundle, Definition
\ref{sec5:df:conn} reduces to the ordinary definition of a connection in the
bundle $\pi$. In fact, if we have a connection $\na^\bullet$ in the sense of
Definition \ref{sec5:df:conn}, then the correspondence
$$\Dr(A)\hookrightarrow\Dr(A,B)\xrightarrow{\na^B}\Dr(B)$$
allows one to lift any vector field on $M$ up to a $\pi$-projectible field on
$E$. Conversely, if $\na$ is such a correspondence, then we can construct a
natural transformation $\na^\bullet$ of the functors $\Dr(A,\cdot)$ and
$\Dr(\cdot)$ due to the fact that for smooth finite-dimensional manifolds
one has $\Dr(A,P)=
P\ot_A\Dr(A)$ and $\Dr(P)=P\ot_B\Dr(P)$ for an arbitrary $B$-module $P$. We
use the notation $\na=\na^B$ in the sequel.
\end{remark}
\begin{definition}\label{sec5:df:curv}
Let $\na^\bullet$ be a connection in $(A,B,\varphi)$ and $X,Y\in\Dr(A,B)$ be
two derivations. The \emph{curvature form} of the connection $\na^\bullet$
on the pair $X,Y$ is defined by
\begin{eqnarray}\label{sec5:eq:curv}
R_\na(X,Y)=[\na(X),\na(Y)]-\na(\na(X)\circ Y-\na(Y)\circ X).
\end{eqnarray}
Note that \eqref{sec5:eq:curv} makes sense, since $\na(X)\circ Y-\na(Y)\circ
X$ is a $B$-valued derivation of $A$.
\end{definition}

Consider now the de Rham differential $d=d_B\colon B\xra{}\La^1(B)$. Then the
composition $d_B\circ\varphi\colon A\xra{}B$ is a derivation. Consequently,
we may consider the derivation $\na(d_B\circ\varphi)\in\Dr(\La^1(B))$.
\begin{definition}\label{sec5:df:connf}
The element $U_\na\in\Dr(\La^1(B))$ defined by
\begin{equation}\label{sec5:eq:connf}
U_\na=\na(d_B\circ\varphi)-d_B
\end{equation}
is called the \emph{connection form} of $\na$.
\end{definition}
Directly from the definition we obtain the following
\begin{lemma}\label{sec5:le:connf1}
The equality
\begin{equation}\label{sec5:eq:connf1}
\ip_X(U_\na)=X-\na(X\rest{A})
\end{equation}
holds for any $X\in\Dr(B)$.
\end{lemma}
Using this result, we now prove
\begin{proposition}\label{sec5:le:curv1}
If $B$ is a smooth algebra, then
\begin{equation}\label{sec5:eq:curv1}
\ip_Y\ip_X\fnij{U_\na}{U_\na}=
2R_\na(X\rest{A},Y\rest{A})
\end{equation}
for any $X,Y\in\Dr(B)$.
\end{proposition}
\begin{proof}
First note that $\deg U_\na=1$. Then using \eqref{sec5:eq:FN4} and
\eqref{sec5:eq:FN1} we obtain
\begin{multline*}
\ip_X\fnij{U_\na}{U_\na}=\fnij{\ip_XU_\na}{U_\na}+
\fnij{U_\na}{\ip_XU_\na}-\ip_{\fnij{X}{U_\na}}U_\na
-\ip_{\fnij{X}{U_\na}}U_\na\\
=2\big(\fnij{\ip_XU_\na}{U_\na}-\ip_{\fnij{X}{U_\na}}U_\na\big).
\end{multline*}
Applying $\ip_Y$ to the last expression and using \eqref{sec5:eq:FN2} and
\eqref{sec5:eq:FN4}, we get now
$$
\ip_Y\ip_X\fnij{U_\na}{U_\na}=2\big(\fnij{\ip_XU_\na}{\ip_YU_\na}-
\ip_{\fnij{X}{Y}}U_\na\big).
$$
But $\fnij{V}{W}=[V,W]$ for any $V,W\in\Dr(\La^0(A))=\Dr(A)$. Hence, by
\eqref{sec5:eq:connf1}, we have
$$
\ip_Y\ip_X\fnij{U_\na}{U_\na}=
2\big([X-\na(X\rest{A}),Y-\na(Y\rest{A})]-
([X,Y]-\na([X,Y]\rest{A}))\big).
$$
It only remains to note now that $\na(X\rest{A})\rest{A}=X\rest{A}$ and
$[X,Y]\rest{A}=X\circ Y\rest{A}-Y\circ X\rest{A}$.
\end{proof}
\begin{definition}\label{sec5:df:flat}
A connection $\na$ in $(A,B,\varphi)$ is called \emph{flat}, if $R_\na=0$.
\end{definition}
Thus for flat connections we have
\begin{equation}\label{sec5:eq:integr}
\fnij{U_\na}{U_\na}=0.
\end{equation}

Let $U\in\Dr(\La^1(B))$ be an element satisfying \eqref{sec5:eq:integr}. Then
from the graded Jacobi identity \eqref{sec5:eq:FN2} we obtain $2\fnij{U}
{\fnij{U}{X}}=\fnij{\fnij{U}{U}}{X}=0$ for any $X\in\Dr(\La^*(A))$.
Consequently, the operator $\partial_U=\fnij{U}{\cdot}\colon\Dr(\La^i(B))
\xra{}\Dr(\La^{i+1}(B))$ defined by the equality $\partial_U(X)=\fnij{U}{X}$
satisfies the identity
$\partial_U\circ\partial_U=0$.

Consider now the case $U=U_\na$, where $\na$ is a flat connection.
\begin{definition}\label{sec5:df:vert}
An element $X\in\Dr(\La^*(B))$ is called \emph{vertical}, if $X(a)=0$ for any
$a\in A$. Denote the $B$-submodule of such elements by $\Dr^v(\La^*(B))$.
\end {definition}
\begin{lemma}\label{sec5:le:vert}
Let $\na$ be a connection in $(A,B,\varphi)$. Then
\begin{itemize}
\item[(1)] an element $X\in\Dr(\La^*(B))$ is vertical if and only if $\ip_X
U_\na=X$;
\item[(2)] the connection form $U_\na$ is vertical, $U_\na\in\Dr^v(\La^1(B))$;
\item[(3)] the map $\partial_{U_\na}$ preserves verticality,
$\partial_{U_\na}(\Dr^v(\La^i(B)))\subset\Dr^v(\La^{i+1}(B))$.
\end{itemize}
\end{lemma}
\begin{proof}
To prove (1), use Lemma \ref{sec5:le:connf1}: from \eqref{sec5:eq:connf1} it
follows that $\ip_XU_\na=X$ if and only if $\na(X\rest{A})=0$. But
$\na(X\rest{A})\rest{A}=X\rest{A}$. The second statements follows from the
same lemma and from the first one:
$$\ip_{U_\na}U_\na=U_\na-\na(U_\na\rest{A})=U_\na-\na\big((U_\na-
\na(U_\na\rest{A})\rest{A}\big)=U_\na.$$
Finally, (3) is a consequence of \eqref{sec5:eq:FN4}.
\end{proof}
\begin{definition}\label{sec5:df:Ucomp}
Denote the restriction $\partial_{U_\na}\rest{\Dr^v(\La^*(A))}$ by
$\partial_\na$ and call the complex
\begin{equation}\label{sec5:eq:Ucomp}
0\xra{}\Dr^v(B)\xra{\partial_\na}\Dr^v(\La^1(B))\xra{}\dotsb\xra{}
\Dr^v(\La^i(B))\xra{\partial_\na}\Dr^v(\La^{i+1}(B))\xra{}\dotsb
\end{equation}
the \emph{$\na$-complex} of the triple $(A,B,\varphi)$. The corresponding
cohomology is denoted by $H_\na^*(B;A,\varphi)=
\bigoplus_{i\ge0}H_\na^i(B;A,\varphi)$ and is called the
\emph{$\na$-cohomology} of the triple $(A,B,\varphi)$.
\end{definition}
Introduce the notation
\begin{eqnarray}\label{sec5:eq:dv}
 d_\na^v=\Ld_{U_\na}\colon \La^i(B)\xra{}\La^{i+1}(B).
\end{eqnarray}
\begin{proposition}\label{sec5:pr:Ucomp}
Let $\na$ be a flat connection in the triple $(A,B,\varphi)$ and $B$ be a
smooth \emph{(}or finitely smooth\emph{)} algebra. Then for any $X,Y\in
\Dr^v(\La^*(A))$ and $\om\in\La^*(A)$ one has
\begin{gather}
\partial_\na\fnij{X}{Y}=\fnij{\partial_\na X}{Y}+
(-1)^X\fnij{X}{\partial_\na Y},\label{sec5:eq:Ucomp1}\\
[\ip_X,\partial_\na]=(-1)^X\ip_{\partial_\na X},\label{sec5:eq:Ucomp2}\\
\partial_\na(\om\wg X)=( d_\na^v- d)(\om)\wg X+
(-1)^\om\om\wg\partial_\na X,\label{sec5:eq:Ucomp3}\\
[ d_\na^v,\ip_X]=\ip_{\partial_\na X}+(-1)^X\Ld_X.\label{sec5:eq:Ucomp4}
\end{gather}
\end{proposition}
\begin{proof}
Equality \eqref{sec5:eq:Ucomp1} is a direct consequence of
\eqref{sec5:eq:FN2}. Equality \eqref{sec5:eq:Ucomp2} follows from
\eqref{sec5:eq:FN4}. Equality \eqref{sec5:eq:Ucomp3} follows from
\eqref{sec5:eq:FN5} and \eqref{sec5:eq:connf1}. Finally,
\eqref{sec5:eq:Ucomp4} is obtained from \eqref{sec5:eq:FN3}.
\end{proof}
\begin{corollary}\label{sec5:co:homstr}
The cohomology module $H_\na^*(B;A,\varphi)$ inherits the graded Lie
algebra structure with respect to the \FN bracket $\fnij{\cdot}{\cdot}$, as
well as to the contraction operation.
\end{corollary}
\begin{proof}
Note that $\Dr^v(\La^*(A))$ is closed with respect to the \FN bracket: to
prove this fact, it suffices to apply \eqref{sec5:eq:FN4}. Then the first
statement follows from \eqref{sec5:eq:Ucomp1}. The second one is a
consequence of \eqref{sec5:eq:Ucomp2}.
\end{proof}
\begin{remark}\label{sec5:rm:struc}
We preserve the same notations for the inherited structures. Note, in
particular, that $H_\na^0(B;A,\varphi)$ is a Lie algebra with respect to the
\FN bracket (which reduces to the ordinary Lie bracket in this case).
Moreover, $H_\na^1(B;A,\varphi)$ is an associative algebra with respect to
the inherited contraction, while the action
$$\CR_\Om\colon X\mapsto\ip_X\Om,\quad X\in H_\na^0(B;A,\varphi),\quad
\Om\in H_\na^1(B;A,\varphi)$$
is a representation of this algebra as endomorphisms of
$H_\na^0(B;A,\varphi)$.
\end{remark}

Consider now the map $ d_\na^v\colon \La^*(B)\xra{}\La^*(B)$ defined by
\eqref{sec5:eq:dv} and define $ d_\na^h= d_B- d_\na^v$.
\begin{proposition}\label{sec5:pr:bicomp}
Let $B$ be a {\em(}finitely{\em)} smooth algebra and $\na$ be a smooth
connection in the triple $(B;A,\varphi)$. Then
\begin{itemize}
\item[(1)] The pair $( d_\na^h, d_\na^v)$ forms a bicomplex, i.e.
\begin{equation}\label{sec5:eq:bicomp1}
d_\na^v\circ d_\na^v=0,\quad d_\na^h\circ d_\na^h=0,\quad
d_\na^h\circ d_\na^v+ d_\na^v\circ d_\na^h=0.
\end{equation}
\item[(2)] The differential $ d_\na^h$ possesses the following
properties
\begin{gather}
[ d_\na^h,\ip_X]=-\ip_{\partial_\na X},\label{sec5:eq:bicomp2}\\
\partial_\na(\om\wg X)=- d_\na^h(\om)\wg X+(-1)^\om\om\wg\partial_\na X,
\label{sec5:eq:bicomp3}
\end{gather}
where $\om\in\La^*(B)$, $X\in\Dr^v(\La^*(B))$.
\end{itemize}
\end{proposition}
\begin{proof}
(1) Since $\deg d_\na^v=1$, we have
$$2 d_\na^v\circ d_\na^v=[ d_\na^v, d_\na^v]=[\Ld_{U_\na},\Ld_{U_\na}]=
\Ld_{\fnij{U_\na}{U_\na}}=0.$$
Since $ d_\na^v=\Ld_{U_\na}$, the identity $[ d_B, d_\na^v]=0$ holds
(see \eqref{sec5:eq:Lied3}), and it concludes the proof of the first part.

(2) To prove \eqref{sec5:eq:bicomp2}, note that
$$[ d_\na^h,\ip_X]=[d_B-d_\na^h,\ip_X]=(-1)^X\Ld_X-[d_\na^v,\ip_X],$$
and \eqref{sec5:eq:bicomp2}
holds due to \eqref{sec5:eq:Ucomp4}. Finally, \eqref{sec5:eq:bicomp3} is just
the other form of \eqref{sec5:eq:Ucomp3}.
\end{proof}
\begin{definition}\label{sec5:df:bicomp}
Let $\na$ be a connection in $(A,B,\varphi)$.
\begin{itemize}
\item[(1)] The bicomplex $(B, d_\na^h, d_\na^v)$ is called the {\em
variational bicomplex} associated to the connection $\na$.
\item[(2)] The corresponding spectral sequence is called the
 {\em $\na$-spectral sequence} of the triple $(A,B,\varphi)$.
\end{itemize}
\end{definition}
Obviously, the $\na$-spectral sequence converges to the de Rham cohomology
of $B$.

To finish this section, note the following. Since the module $\La^1(B)$
is generated by the image of the operator $d_B\colon B\xra{}\La^1(B)$ while
the graded algebra $\La^*(B)$ is generated by $\La^1(B)$, we have the
direct sum decomposition
$$
\La^*(B)=\bigoplus_{i\ge0}\bigoplus_{p+q=i}\La_v^p(B)\ot\La_h^q(B),
$$
where
$$
\La_v^p(B)=\underbrace{\La_v^1(B)\wg\dots\wg\La_v^1(B)}_{p\ {\rm times}},
\quad
\La_h^q(B)=\underbrace{\La_h^1(B)\wg\dots\wg\La_h^1(B)}_{q\ {\rm times}},
$$
while the submodules $\La_v^1(B)\subset\La^1(B)$, $\La_h^1(B)\subset\La^1(B)$
are spanned in $\La^1(B)$ by the images of the differentials $d_\na^v$ and
$d_\na^h$ respectively. Obviously, we have the following embeddings:
\begin{gather}
d_\na^h\left(\La_v^p(B)\ot\La_h^q(B)\right)
\subset\La_v^p(B)\ot\La_h^{q+1}(B),\nonumber\\
 d_\na^v\left(\La_v^p(B)\ot\La_h^q(B)\right)
\subset\La_v^{p+1}(B)\ot\La_h^q(B).\nonumber
\end{gather}

Denote by $\Dr^{p,q}(B)$ the module $\Dr^v(\La_v^p(B)\ot\La_h^q(B))$. Then,
obviously, $\Dr^v(B)=\bigoplus_{i\ge0}\bigoplus_{p+q=i}\Dr^{p,q}(B)$,
while from equalities \eqref{sec5:eq:bicomp2} and \eqref{sec5:eq:bicomp3} we
obtain
$$
\partial_\na\big(\Dr^{p,q}(B)\big)\subset\Dr^{p,q+1}(B).
$$
Consequently, the module $H_\na^*(B;A,\varphi)$ is split as
\begin{eqnarray}\label{sec5:eq:hsplit}
H_\na^*(B;A,\varphi)=\bigoplus_{i\ge0}\bigoplus_{p+q=i}
H_\na^{p,q}(B;A,\varphi)
\end{eqnarray}
with the obvious meaning of the notation $H_\na^{p,q}(B;A,\varphi)$.

\subsection{Applications to differential equations: recursion
operators}\label{sub:recop}
Now we apply the above exposed algebraic results to the case of
infinitely prolonged differential equations. Let us start with establishing
a correspondence between geometric constructions of Section \ref{sec:geom}
and algebraic ones presented in the previous two subsections.

Let $\CE\sbs J^k(\pi)$ be a formally integrable equation (see Definition
\vref{sec3:df:formint}) and $\Ei\sbs\Ji(\pi)$ be its infinite prolongations.
Then the bundle $\pi_\infty:\Ei\xra{}M$ is endowed with the Cartan connection
$\CC$ (Definition \vref{sec3:df:cartcon}) and this connection is flat
(Corollary \vref{sec3:co:cflat}). Thus the triple
$$
\big(A=\Ci(M),\, B=\CF(\CE),\, \vf=\pi_\infty^*\big)
$$
with $\na=\CC$ is an algebra with a flat connection, $A$ being a smooth and
$B$ being a finitely smooth algebra. The corresponding connection form is
exactly the structural element $U_\CC$ of the equation $\CE$ (see
Definition \vref{sec3:df:ue}).

Thus, to any formally integrable equation ${\CE}\subset J^k(\pi)$ we can
associate the complex
\begin{equation}\label{sec5:eq:Ccomp}
0\xra{}\Dv(\E)\xra{\partial_{\CC}}\Dv(\La^1(\E))\xra{}\dotsb\xra{}
\Dv(\La^i(\E))\xra{\partial_{\CC}}\Dv(\La^{i+1}(\E))\xra{}\dotsb
\end{equation}
and the cohomology theory determined by the Cartan connection. We denote
the corresponding cohomology modules by $H_{\CC}^*({\CE})=
\bigoplus_{i\ge0}H_{\CC}^i({\CE})$. In the case of the ``empty''
equation, we use the notation $H_{\CC}^*(\pi)=\bigoplus_{i\ge0}
H_{\CC}^i(\pi)$.
\begin{definition}\label{sec5:eq:Chom}
Let ${\CE}\subset J^k(\pi)$ be a formally integrable equation and ${\CC}$
be the Cartan connection in the bundle $\pi_\infty\colon \Ei\xra{} M$. Then
the module $H_{\CC}^*({\CE})$ is called the \emph{${\CC}$-cohomology} of
${\CE}$.
\end{definition}
\bre\label{sec5:re:coinc}
Let us also note that the above introduced modules $\La_h^q(B)$ are identical
to the modules $\hL^q(\CE)$ of horizontal $q$-forms on $\Ei$, the modules
$\La_v^p(B)$ coincide with the modules of Cartan forms $\CC^p\La(\CE)$, the
differential $d_\na^h$ is the extended horizontal de Rham differential $\hd$,
while $d_\na^v$ is the Cartan differential $d_\CC$ (cf.\ with constructions on
pp.~\pageref{sec3:p:start}--\pageref{sec3:p:end}). Thus we again obtain a
complete coincidence between algebraic and geometric approaches. In
particular, the $\na$-spectral sequence (Definition \vref{sec5:df:bicomp}
(2)) is the Vinogradov $\CC$-spectral sequence (see the Section
\ref{css:sec}).
\ere
The following result contains an interpretation of the first two
of ${\CC}$-cohomology groups.
\begin{theorem}\label{sec5:th:H01}
For any formally integrable equation ${\CE}\subset J^k(\pi)$, one has
the following identities:
\begin{itemize}
\item[(1)] The module $H_{\CC}^0({\CE})$ as a Lie algebra is
isomorphic to the Lie algebra $\sym{\CE}$ of higher
symmetries\footnote{See Definition \vref{sec3:df:symEi}.}
of the equation ${\CE}$.
\item[(2)] The module $H_{\CC}^1({\CE})$ is the set of the
equivalence classes of nontrivial vertical deformations of the equation
structure \emph{(}i.e., of the structural element\emph{)} on $\CE$.
\end{itemize}
\end{theorem}
\begin{proof}
To prove (1), take a vertical vector field $Y\in\Dr^v({\CE})$ and an
arbitrary field $Z\in\Dr({\CE})$. Then, due to \veqref{sec5:eq:FN4}, one has
\begin{multline*}
\ip_Z\partial_{\CC}Y=\ip_Z\fnij{U_{\CC}}{Y}
=[\ip_ZU_{\CC},Y]-\ip_{[Z,Y]}U_{\CC}\\
=[Z^v,Y]-[Z,Y]^v=[Z^v-Z,Y]^v,
\end{multline*}
where $Z^v=\ip_ZU_{\CC}$. Hence, $\partial_{\CC}Y=0$ if and only if
$[Z-Z^v,Y]^v=0$ for any $Z\in\Dr({\CE})$. But the last equality holds if
and only if $[{\CC}X,Y]=0$ for any $X\in\Dr(M)$ which means that
$$\ker\big(\partial_{\CC}\colon\Dr^v({\CE})\xra{}\Dr^v(\La^1({\CE}))\big)=
\sym{\CE}.$$

Consider the second statement now. Let $U(\varepsilon)\in\Dr^v(\La^1({\CE}))$ be
a deformation of the structural element satisfying the conditions
$\fnij{U(\varepsilon)}{U(\varepsilon)}=0$ and $U(0)=U_{\CC}$. Then
$U(\varepsilon)=U_{\CC}+U_1\varepsilon+O(\varepsilon^2)$. Consequently,
$$
\fnij{U(\varepsilon)}{U(\varepsilon)}=\fnij{U_{\CC}}{U_{\CC}}+
2\fnij{U_{\CC}}{U_1}\varepsilon+O(\varepsilon^2)=0,
$$
from which it follows that $\fnij{U_{\CC}}{U_1}=\partial_{\CC}U_1=0$.
Hence the linear part of the deformation $U(\varepsilon)$ determines an
element of $H_{\CC}^1({\CE})$ and vice versa. On the other hand, let $A\colon
\Ei\xra{}\Ei$ be
a diffeomorphism\footnote{%
Since $\Ei$ is, in general, infinite-dimensional,
vector fields on $\Ei$ do not usually possess
one-parameter groups of diffeomorphisms. Thus
the arguments below are of a heuristic nature.}
of $\Ei$. Define the action $A^*$ of $A$ on
the elements $\Om\in\Dr(\La^*({\CE}))$ in such a way that the diagram
$$
\begin{CD}
\La^*({\CE})@>\Ld_\Om>>\La^*({\CE})\\
@VA^*VV    @VVA^*V\\
\La^*({\CE})@>\Ld_\Om>>\La^*({\CE})\\
\end{CD}
$$
is commutative. Then, if $A_t$ is a one-parameter group of diffeomorphisms,
we have, obviously,
$$
\left.\frac{d}{dt}\right|_{t=0}A_{t,*}(\Ld_\Om)=
\left.\frac{d}{dt}\right|_{t=0}A_t^*\circ\Ld_\Om\circ(A_t^*)^{-1}=[\Ld_X,\Ld_\Om]=
\Ld_{\fnij{X}{\Om}}.
$$
Hence, the infinitesimal action is given by the \FN bracket. Taking $\Om=
U_{\CC}$ and $X\in\Dr^v({\CE})$, we see that $\im\partial_{\CC}$ consists of
infinitesimal deformations arising due to infinitesimal action of
diffeomorphisms on the structural element. Such deformations are naturally
called trivial.
\end{proof}
\begin{remark}\label{sec5:rm:obstr}
From the general theory \cite{GerstSch}, we obtain that the module
$H_{\CC}^2({\CE})$ consists of obstructions to prolongation of
infinitesimal deformations to formal ones. In the case under
consideration, elements $H_{\CC}^2({\CE})$ have another nice interpretation
discussed  later (see Remark \vref{sec5:rm:Lierec}).
\end{remark}

We shall now compute the modules $H_{\CC}^p(\pi),\ p\ge0$. To do this,
recall the splitting $\La^i({\CE})=\bigoplus_{p+q=i}\CC^p\La({\CE})\ot
\hL^q({\CE})$ (see Subsection \ref{sub:fvcalc}).
\begin{theorem}\label{th3.8}
For any $p\ge0$, one has
$$H_{\CC}^p(\pi)={\CF}(\pi,\pi)\ot_{{\CF}(\pi)}\CC^p\La({\pi}).$$
\end{theorem}
\begin{proof}
Define a filtration in $\Dr^v(\La^*(\pi))$ by setting
$$
F^l\Dr^v(\La^p(\pi))=\{X\in\Dr^v(\La^p(\pi))\mid X\rest{{\CF}_{l-p-1}}=0\}.
$$
Evidently,
$$F^l\Dr^v(\La^p(\pi))\subset F^{l+1}\Dr^v(\La^p(\pi)),\quad
\partial_\CC\big(F^l\Dr^v(\La^p(\pi))\big)\subset F^l\Dr^v(\La^{p+1}(\pi)).$$
Thus we obtain the spectral sequence associated to this filtration. To
compute the term $E_0$, choose local coordinates $x_1,\dots,x_n,u^1,\dots,
u^m$ in the bundle $\pi$ and consider the corresponding special coordinates
$u_\si^j$ in $\Ji(\pi)$. In these
coordinates, the structural element is represented as
\begin{eqnarray}\label{sec5:eq:Upi}
U_\CC=\sum_{|\si|\ge0}\sum_{j=1}^m\left( d u_\si^j-\sum_{i=1}^n
u_{\si i}^j\,dx_i\right)\ot\frac{\partial}{\partial u_\si^j},
\end{eqnarray}
while for $X=\sum_{\si,j}
\theta_\si^j\ot\partial/\partial u_\si^j,\ \theta\in\La^*(\pi)$, one has
\begin{eqnarray}\label{sec5:eq:dpi}
\partial_\CC(X)=\sum_{|\si|\ge0}\sum_{j=1}^m\sum_{i=1}^n d x_i\wg
\left(\theta_{\si i}^j-D_i(\theta_\si^j)\right)
\ot\frac{\partial}{\partial u_\si^j}.
\end{eqnarray}

Obviously, the term
$$E_0^{p,-q}=F^p\Dr^v(\La^{p-q}(\pi))/F^{p-1}\Dr^v(\La^{p-q}(\pi)),\quad
p\ge0,\quad 0\le q\le p,$$
is identified with the tensor product $\La^{p-q}(\pi)\ot_{{\CF}(\pi)}\Ga(\pi_
{\infty,q-1}^*(\pi_{q,q-1}))$. These modules can be locally represented as
${\CF}(\pi,\pi)\ot
\La^{p-q}(\pi)$-valued homogeneous polynomials of order $q$, while the
differential $\partial_0^{p,-q}\colon E_0^{p,-q}\xra{}E_0^{p,-q+1}$ acts as
the $\delta$-Spencer differential (or, which is the same, as the Koszul
differential; see Exercise \vref{sec1:es:Kosz}). Hence, all homology groups
are trivial except for those at the terms
$E_0^{p,0}$ and one has
$$\coker\partial_0^{p,0}={\CF}(\pi,\pi)\ot_{{\CF}(\pi)}\CC^p\La(\pi).$$
Consequently, only the $0$-th line survives in the term
$E_1$ and this line is of the form
\begin{multline*}
{\CF}(\pi,\pi)\xra{\partial_1^{0,0}}{\CF}(\pi,\pi)\ot_{{\CF}(\pi)}
\CC^1\La(\pi)\xra{}\dotsb\\
\dotsb\xra{}{\CF}(\pi,\pi)\ot_{{\CF}(\pi)}\CC^p\La(\pi)\xra{\partial_1^{p,0}}
{\CF}(\pi,\pi)\ot_{{\CF}(\pi)}\CC^{p+1}\La(\pi)\xra{}\dotsb
\end{multline*}
But the image of $\partial_{\CC}$ contains at least one horizontal
component (see equality \veqref{sec5:eq:Ucomp3}, where, by definition,
$d_\na^v-d=d_\CC-d=-\hd$). Therefore, all differentials
$\partial_1^{p,0}$ vanish.
\end{proof}

Let us now establish the correspondence between the last result (describing
${\CC}$-cohomology in terms of $\CC^*\La(\pi)$) and representation of
$H_{\CC}^*(\pi)$ as classes of derivations ${\CF}(\pi)\xra{}\La^*(\pi)$.
To do this, for any $\om=(\om^1,\dots,\om^m)\in{\CF}(\pi,\pi)
\ot_{{\CF}(\pi)}\CC^*\La(\pi)$ set
\begin{equation}\label{sec5:eq:supev}
\re_\om=\sum_{\si,j}D_\si(\om^j)\ot\frac{\partial}{\partial u_\si^j},
\end{equation}
where $D_\si=D_1^{\si_1}\circ\dots\circ D_n^{\si_n}$ for $\si=(\si_1,\dots
\si_n)$.
\begin{definition}\label{sec5:df:supev}
The element $\re_\om\in\Dr^v(\La^*(\pi))$ defined by
\eqref{sec5:eq:supev} is called the \emph{evolutionary superderivation}
with the \emph{generating section} $\om\in\CC^*\La(\pi)$.
\end{definition}
\begin{proposition}\label{sec5:pr:indep}
The definition of $\re_\om$ is independent of coordinate choice.
\end{proposition}
\begin{proof}
It is easily checked that
$$
\re_\om\big({\CF}(\pi)\big)\subset\La_v^*(\pi),\quad
\re_\om\in\ker\partial_{\CC}.
$$
But derivations possessing these properties are uniquely determined by
their restriction to ${\CF}_0(\pi)$ which coincides with the action of
the derivation $\om:{\CF}_0(\pi)\xra{}\CC^*\La(\pi)$. Let us prove this
fact.

Set $X=\re_\om$ and recall that the derivation $X$ is uniquely determined by the
corresponding Lie derivative $\Ld_X:\La^*(\pi)\xra{}\La^*(\pi)$. Further,
since $\Ld_Xd\ta=(-1)^Xd(\Ld_X\ta)$ (see \veqref{sec5:eq:Lied3}) for any $\ta
\in\La^*(\pi)$, the derivation $\Ld_X$ is determined by its restriction
to $\La^0(\pi)=\CF(\pi)$.

Now, from the identity $\pat_\CC X=0$ it follows that
\begin{equation}\label{sec5:eq:comm}
0=\fnij{U_\CC}{X}(f)=\Ld_{U_\CC}(\Ld_X(f))-(-1)^X\Ld_X(\Ld_{U_\CC}(f)),\quad
f\in\CF(\CE).
\end{equation}
Let now $X$ be such that $\Ld_X\rest{\CF_0(\pi)}=0$ and assume
that $\Ld_X\rest{\CF_r(\pi)}=0$ for some $r>0$. Then taking $f=u_\si^j$,
$|\si|=r$, and using \veqref{sec5:eq:comm} we obtain
$$
\Ld_X\left(du_\si^j-\sum_{i=1}^nu_{\si i}^jdx_i\right)=\Ld_Xd_\CC u_\si^j=
(-1)^Xd_\CC(L_X(u_\si^j))=0.
$$
In other words,
\begin{multline*}
\Ld_X\left(\sum_{i=1}^nu_{\si i}^jdx_i\right)=
\sum_{i=1}^n\Ld_X(u_{\si i}^jdx_i))=\sum_{i=1}^n\Ld_X(du_\si^j)\\
=(-1)^X\sum_{i=1}^nd(\Ld_Xu_\si^j)=0.
\end{multline*}
Hence, $\Ld_X(u_\si^j)=0$ and thus $\Ld_X\rest{\CF_{r+1}(\pi)}=0$.
\end{proof}

From this result and from Corollary \vref{sec5:co:homstr}, it follows
that if two evolutionary superderivations $\re_\om,\re_\theta$ are
given, the elements
\begin{itemize}
\item[(i)] $\quad\fnij{\re_\om}{\re_\theta}$,
\item[(ii)] $\quad\ip_{\re_\om}\left(\re_\theta\right)$
\end{itemize}
are evolutionary superderivations as well.

In the first case, the corresponding generating section is called the {\em
Jacobi superbracket} of elements $\om=(\om^1,\dots,\om^m)$ and $\theta=
(\theta^1,\dots,\theta^m)$ and is denoted by $\{\om,\theta\}$. The
components of this bracket are expressed by the formula
\begin{eqnarray}\label{sec5:eq:supJac}
\{\om,\theta\}^j=\Ld_{\re_\om}(\theta^j)-
(-1)^{\om\theta}\Ld_{\re_\theta}(\om^j),\ j=1,\dots,m.
\end{eqnarray}
Obviously, the module ${\CF}(\pi,\pi)\ot_{{\CF}(\pi)}\CC^*\La(\pi)$ is a
graded Lie algebra with respect to the Jacobi superbracket. The restriction
of $\{\cdot,\cdot\}$ to $\CF(\pi,\pi)\ot\CC^0\La(\pi)=\CF(\pi,\pi)$
coincides with the higher Jacobi bracket (see Definition
\vref{sec3:df:highj}).

In the case (ii), the generating section is $\ip_{\re_\om}(\theta)$. Note
now that any element $\rho\in\CC^1\La(\pi)$ is of the form $\rho=
\sum_{\si,\al}a_{\si,\al}\om_\si^\al$, where, as before,
$
\om_\si^\al= d_\CC u_\si^\al= d u_\si^\al-\sum_{i=1}^nu_{\si i}^\al d x_i
$
are the Cartan forms on $\Ji(\pi)$.
Hence, if $\theta\in{\CF}(\pi,\pi)\ot_{{\CF}(\pi)}\CC^1\La(\pi)$ and
$\theta^j=\sum_{\si,\al}a_{\si,\al}^j\om_\si^\al$, then
\begin{eqnarray}\label{sec5:eq:Cact}
\big(\ip_{\re_\om}(\theta)\big)^j=\sum_{\si,\al}a_{\si,\al}^jD_\si(\om^\al).
\end{eqnarray}
In particular, we see that \veqref{sec5:eq:Cact} establishes an isomorphism
between the modules ${\CF}(\pi,\pi)\ot_{{\CF}(\pi)}\CC^*\La(\pi)$ and
$\CDiff(\pi,\pi)$ and defines the action of ${\CC}$-differential
operators on elements of $\CC^*\La(\pi)$. This is a really well-defined
action because of the fact that $\ip_{{\CC}X}\om=0$ for any $X\in\Dr(M)$
and $\om\in\CC^*\La(\pi)$.

Consider now a formally integrable differential equation $\CE\subset
J^k(\pi)$ and assume that
it is determined by a differential operator $\De\in{\CF}(\pi,\xi)$. Denote,
as in Section \ref{sec:geom},
by $\ell_{\CE}$ the restriction of the operator of universal linearization
$\ell_\De$ to $\Ei$. Let $\ell_{\CE}^{[p]}$ be the extension of
$\ell_{\CE}$ to ${\CF}(\pi,\pi)\ot_{{\CF}(\pi)}\CC^p\La({\CE})$ which is
well defined due to what has been said above. Then the module $H_\CC^{p,0}
({\CE})$
is identified with the set of evolution superderivations $\re_\om$ whose
generating sections $\om\in{\CF}(\pi,\pi)\ot_{{\CF}(\pi)}\CC^p\La({\CE})$
satisfy the equation
\begin{eqnarray}\label{sec5:eq:ellp}
\ell_{\CE}^{[p]}(\om)=0
\end{eqnarray}
If, in addition, ${\CE}$ satisfies the assumptions of the two-line theorem,
then $H_{\CC}^{p,1}({\CE})$ is identified with the cokernel of
$\ell_{\CE}^{[p-1]}$ and thus
$$
H_{\CC}^i({\CE})=\ker\ell_{\CE}^{[i]}\oplus\coker\ell_{\CE}^{[i-1]}
$$
in this case. These two statements will be proved in Subsection
\ref{hc.apkcc:subsec}.

As it was noted in Remark \vref{sec5:rm:struc}, $H_{{\CC}}^1({\CE})$ is an
associative algebra with respect to contraction and is represented in the algebra
of endomorphisms of $H_{{\CC}}^0({\CE})$. It is easily seen that the
action of the $H_{{\CC}}^{0,1}({\CE})$ is trivial while
$H_{{\CC}}^{1,0}({\CE})$ acts on $H_{{\CC}}^0({\CE})=\sym{\CE}$ as
${\CC}$-differential operators (see above).
\begin{definition}\label{sec5:df:recop}
Elements of the module $H_{{\CC}}^{1,0}({\CE})$ are called \emph{recursion
operators} for symmetries of the equation ${\CE}$.
\end{definition}
We use the notation $\CR({\CE})$ for the algebra of recursion operators.
\begin{remark}\label{sec5:rm:recvoid}
The algebra $\CR({\CE})$ is always nonempty, since it contains the
structural element $U_{\CE}$ which is the unit of this algebra. ``Usually''
this is the only solution of \eqref{sec5:eq:ellp} for $p=1$
(see Example \ref{sec5:ex:Burg} below). This fact apparently
contradicts practical experience (cf.\ with well-known recursion operators
for the KdV and other integrable systems \cite{Olver}). The reason is that these
operators contain nonlocal terms like $D^{-1}$ or of a more complicated
form. An adequate framework to deal with such constructions will be
described in the next subsection.
\end{remark}
\bex\label{sec5:ex:Burg}
Let
\begin{equation}\label{sec5:eq:Burg}
u_t=uu_x+u_{xx}
\end {equation}
be the Burgers equation. For internal coordinates on $\Ei$ we choose
the functions $x,t,u=u_0,\dots, u_k,\dots$, where $u_k$ corresponds to the
partial derivative $\pat^ku/\pat x^k$.

We shall prove here that the only solution of the
equation $\ell_\CE^{[1]}(\om)=0$ for \eqref{sec5:eq:Burg} is $\om=\al\om_0$,
$\al=\cnst$, where
$$\om_k=d_\CC u_k=du_k-u_{k+1}dx-D_x^k(uu_1+u_2)dt.$$

Let $\om=\phi^0\om_0+\dots+\phi^r\om_r$. Then the equation \veqref{sec5:eq:ellp}
for $p=1$ transforms to
\begin{align}\label{sec5:eq:recBurg}
&u_0D_x(\phi^0)+D_x^2(\phi^0)=D_t(\phi^0)+\sum_{j=1}^ru_{j+1}\phi^j,
\nonumber\\
&u_0D_x(\phi^1)+D_x^2(\phi^1)+2D_x(\phi^0)=D_t(\phi^1)+
\sum_{j=2}^r(j+1)u_j\phi^j,
\nonumber\\
&\dots\\
&u_0D_x(\phi^i)+D_x^2(\phi^i)+2D_x(\phi^0)=D_t(\phi^i)+
\sum_{j=i+1}^r\binom{j+1}{i}u_{j-i+1}\phi^j,
\nonumber\\
&\dots
\nonumber\\
&u_0D_x(\phi^r)+D_x^2(\phi^r)+2D_x(\phi^{r-1})=D_t(\phi^r)+ru_1\phi^r,
\nonumber\\[3mm]
&D_x(\phi^r)=0.
\nonumber
\end{align}
To prove the result, we apply the scheme used in \cite{VinKras1} to
describe the symmetry algebra of the Burgers equation.

Denote by $\CK_r$ the set of solutions of \eqref{sec5:eq:recBurg}. Then a
direct computation shows that
\begin{equation}\label{sec5:eq:K1}
\CK^1=\{\al\om_0\mid\al\in\BBR\}
\end{equation}
and that any element $\om\in\CK_r$, $r>1$, is of the form
\begin{equation}\label{sec5:eq:Kr}
\om=\al_r+\left(\frac{r}{2}u_0\al_r+\frac12x\al_r^{(1)}+\al_{r-1}\right)+
\Om[r-2],
\end{equation}
where $\al_r=\al_r(t)$, $\al_{r-1}=\al_{r-1}(t)$, $\al^{(i)}$ denotes the
derivative $d^i\al/dt^i$, and $\Om[s]$ is an arbitrary linear combination
of $\om_0,\dots,\om_s$ with coefficients in $\CF(\CE)$.

Note now that for any evolution equation the embedding
$$\fnij{\sym\CE}{\ker\ell_\CE^{[1]}}\sbs\ker\ell_\CE^{[1]}$$
is valid. Consequently, if $\psi\in\sym\CE$ and $\om\in\ker\ell_\CE^{[1]}$,
then $\{\psi,\om\}\in\ker\ell_\CE^{[1]}$.

Since the function $u_1$ is a symmetry of the Burgers equation (translation
along $x$), one has
$$
\{u_1,\om\}=\left(\sum_ku_{k+1}\pdr{}{u_k}\right)\om-D_x\om=
-\pdr{}{x}\om.
$$
Hence, if $\om\in\CK_r$, then from \eqref{sec5:eq:Kr} we obtain that
$$\ad_{u_1}^{(r-1)}(\om)=\al_r^{(r-1)}\om_1+\Om[0]\in\CK_1,$$
where $\ad_\psi=\{\psi,\cdot\}$. Taking into account equation \eqref{sec5:eq:K1},
we get that $\al_r^{r-1}=0$, or
\begin{equation}\label{sec5:eq:polyt}
\al_r=a_0+a_1t+\dots+a_{r-2}t^{r-2},\quad a_i\in\BBR.
\end{equation}

We shall use now the fact that the element $\Phi=t^2u_2+(t^2u_0+tx)u_1+tu_0+1$ is a
symmetry of the Burgers equation (see \cite{VinKras1}). Then, since the
action of symmetries is permutable with the Cartan differential $d_\CC$,
we have
$$
\{\Phi,\phi^s\om_s\}=\re_\Phi(\phi^s\om_s)-\re_{\phi^s\om_s}(\Phi)=
\re_\Phi(\phi^s)\om_s+\phi^s\re_\Phi(\om_s)-\re_{\phi^s\om_s}(\Phi).
$$
But
\begin{multline*}
\re_\Phi(\om_s)=\re_\Phi d_\CC(u_s)=d_\CC\re_\Phi(u_s)=d_\CC D_x^s(\Phi)\\
=d_\CC\big(t^2u_{s+2}+(t^2u_0+tx)u_{s+1}+(s+1)(t^2u_1+t)u_s\big)+\Om[s-1].
\end{multline*}
On the other hand,
\begin{multline*}
\re_{\phi^s\om_s}(\Phi)=t^2\phi^s\om_{s+2}+
\big(2t^2D_x^2(\phi^s)+(t^2u_0+tx)\phi^s\big)\om_{s+1}\\
+\big(t^2D_x^2(\phi^s)+(t^2u_0+tx)D_x(\phi^s)+
(t^2u_1+t)\phi^s\big)\om_s.
\end{multline*}
Thus, we finally obtain
\begin{equation}\label{sec5:eq:Phi}
\begin{split}
\{\Phi,\phi^s\om_s\}=\{\Phi,\phi^s\}\om_s+(s+1)(t^2u_1&+t)\om_s\\
&-2t^2D_x(\phi^s)\om_{s+1}+\Om[s-1].
\end{split}
\end{equation}
Applying \eqref{sec5:eq:Phi} to \eqref{sec5:eq:Kr}, we get
\begin{equation}\label{sec5:eq:adPhi}
\ad_\Phi(\om)=(rt\al_r-t^2\al_r^{(1)})\om_r+\Om[r-1].
\end{equation}

Let now $\om\in\CK_r$ and assume that $\om$ has a nontrivial coefficient
$\al_r$ of the form \eqref{sec5:eq:polyt}, and $a_i$ be the first nontrivial
coefficient in $\al_r$. Then, by representation \eqref{sec5:eq:adPhi}, we have
$$\ad_\Phi^{r-i}(\om)=\al'_r\om_r+\Om[r-1]\in\CK_r,$$
where $\al'_r$ is a polynomial of degree $r-1$. This contradicts to
\eqref{sec5:eq:polyt} and thus concludes the proof.
\eex

\begin{remark}\label{sec5:rm:Lierec}
Let $\varphi\in\sym{\CE}$ be a symmetry and $R\in\CR({\CE})$ be a
recursion operator. Then we obtain a sequence of symmetries
$\varphi_0=\varphi$, $\varphi_1=R(\varphi)$, $\dots$, $\varphi_n=
R^n(\varphi)$, $\dots$. Using identity \veqref{sec5:eq:FN4}, one can compute
the commutators $[\varphi_m,\varphi_n]$ in
terms of $\fnij{\varphi}{R}\in H_{\CC}^{1,0}({\CE})$ and $\fnij{R}{R}\in
H_{\CC}^{2,0}({\CE})$. In particular, it can be shown that when both
$\fnij{\varphi}{R}$ and $\fnij{R}{R}$ vanish, all symmetries $\varphi_n$
mutually commute (see \cite{Kras3}).

For example, if ${\CE}$ is an evolution equation, $H_{\CC}^{p,0}({\CE})=
0$ for all $p\ge2$ (see Theorem \vref{sec6:th:Cvanish}).
Hence, if $\varphi$ is a symmetry and $R$ is a {\em
$\varphi$-invariant} recursion operator (i.e., such that $\fnij{\varphi}{R}=
0$), then $R$ generates a commutative sequence of symmetries. This is exactly
the case for the KdV and other integrable evolution equations.
\end{remark}

\subsection{Passing to nonlocalities}\label{sub:passnlc}
Let us now introduce nonlocal variables into the above described picture.
Namely, let $\CE$ be an equation and $\vf:\CN\xra{}\Ei$ be a covering over
its infinite prolongation. Then, due to Proposition \vref{sec4:pr:covconn},
the triad $\big(\CF(\CN),\Ci(M),(\pi_\infty\circ\vf)^*\big)$ is an algebra
with the flat connection $\CC^\vf$. Hence, we can apply the whole
machinery of Subsections \ref{sub:fvcalc}--\ref{sub:recop} to this
situation. To stress the fact that we are working over the covering $\vf$,
we shall add the symbol $\vf$ to all notations introduced in these
subsections. Denote by $U_\CC^\vf$ the connection form of the connection
$\CC^\vf$ (the structural element of the covering $\vf$).

In particular, on $\CN$ we have the $\CC^\vf$-differential
$\partial_\CC^\vf=\fnij{U_\CC^\vf}{\cdot}:
\Dr^v(\La^i(\CN))\to\Dr^v(\La^{i+1}(\CN))$, whose 0-cohomology
$H_\CC^0(\CE,\vf)$ coincides with the Lie algebra $\sym_\vf\CE$ of nonlocal
$\vf$-symmetries, while the module $H_\CC^{1,0}(\CE,\vf)$ identifies with
recursion operators acting on these symmetries and is denoted by
$\CR(\CE,\vf)$. We also have the horizontal and the Cartan differential
$\hd^\vf$ and $d_\CC^\vf$ on $\CN$ and the splitting $\La^i(\CN)=
\bigoplus_{p+q=i}\CC^p\La^p(\CN)\ot\hL^q(\CN)$.

Choose a trivialization of the bundle $\vf:\CN\xra{}\Ei$ and nonlocal
coordinates $w^1,w^2,\dots$ in the fiber. Then any derivation
$X\in\Dr^v(\La^i(\CN))$ splits to the sum $X=X_\CE+X^v$, where
$X_\CE(w^j)=0$ and $X^v$ is a $\vf$-vertical derivation.
\ble\label{sec5:le:dirs}
Let $\vf:\Ei\times\mathbb{R}^N\xra{}\Ei$, $N\le\infty$, be a covering. Then
$H_\CC^{p,0}(\CE,\vf)=\ker\partial_\CC^\vf\rest{\CC^p\La(\CN)}$. Thus
$H_\CC^{p,0}(\CE,\vf)$ consists of derivations $\Om:\CF(\CN)\xra{}
\CC^p\La(\CN)$ such that
\begin{equation}\label{sec5:eq:dirs}
\fnij{U_\CC^\vf}{\Om}_\CE=0,\quad\big(\fnij{U_\CC^\vf}{\Om}\big)^v=0.
\end{equation}
\ele
\begin{proof}
In fact, due to equality \veqref{sec5:eq:Ucomp3}, any element lying in the
image of $\partial_\CC^\vf$ contains at least one horizontal component, i.e.,
$$\partial_\CC^\vf\big(\Dr^v(\CC^p\La(\CN))\big)\sbs
\Dr^v(\CC^p\La(\CN)\ot\hL^1(\CN)).$$
Thus, equations \eqref{sec5:eq:dirs} should hold.
\end{proof}

We call the first equation in \veqref{sec5:eq:dirs} the \emph{shadow
equation} while the second one is called the \emph{relation equation}. This
is explained by the following result (cf.\ with Theorem
\vref{sec4:th:non-ulin}).
\bpr\label{sec5:pr:shadev}
Let $\CE$ be an evolution equation of the form
$$
u_t=f(x,t,u,\dots,\frac{\partial^ku}{\partial u^k})
$$
and $\vf:\CN=\Ei\times\mathbb{R}^N\xra{}\Ei$ be a covering given by the
vector fields\footnote{To simplify the notations of Section
\ref{sec:nonloc}, we denote the lifting of a \cd operator $\Delta$ to
$\CN$ by $\tilde{\Delta}$.}
$$
\tilde{D}_x=D_x+X,\quad\tilde{D}_t=D_t+T,
$$
where $[\tilde{D}_x,\tilde{D}_t]=0$ and
$$
X=\sum_sX^s\pdr{}{w^s},\quad T=\sum_sT^s\pdr{}{w^s},
$$
$w^1,\dots,w^s,\dots$ being nonlocal variables in $\vf$. Then the group
$H_\CC^{p,0}(\CE,\vf)$ consists of elements
$$\Psi=\sum_i\Psi_i\ot\pdr{}{u_i}+\sum_s\psi^s\pdr{}{w^s}\in
\Dr^v(\CC^p\La(\CN))$$
such that $\Psi_i=\tilde{D}_x^i\Psi_0$ and
\begin{gather}
\tilde{\ell}_\CE^{[p]}(\Psi_0)=0,\label{sec5:eq:hp01}\\
\sum_\al\pdr{X^s}{u_\al}\tilde{D}_x^\al(\Psi_0)+
\sum_\be\pdr{X^s}{w^\be}\psi^\be=\tilde{D}_x(\psi^s)\label{sec5:eq:hp02},\\
\sum_\al\pdr{T^s}{u_\al}\tilde{D}_x^\al(\Psi_0)+
\sum_\be\pdr{T^s}{w^\be}\psi^\be=\tilde{D}_t(\psi^s)\label{sec5:eq:hp03},
\end{gather}
$s=1,2,\dots$, where $\tilde{\ell}_\CE^{[p]}$ is the natural extension
of the operator $\ell_\CE^{[p]}$ to $\CN$.
\epr

\begin{proof}
Consider the Cartan forms
$$
\om_i=du_i-u_{i+1}dx-D_x^i(f)\,dt,\quad
\ta^s=dw^s-X^sdx-T^sdt
$$
on $\CN$. Then the derivation
$$
U_\CC^\vf=\sum_i\om_i\ot\pdr{}{u_i}+\sum_s\ta^s\ot\pdr{}{w^s}
$$
is the structural element of the covering $\vf$. Then, using representation
\veqref{sec5:eq:FNexpl}, we obtain
\begin{multline*}
\partial_\CC^\vf\Psi=dx\wg\sum_i\big(\Psi_{i+1}-\tilde{D}_x(\Psi_i)\big)\ot
\pdr{}{u_i}\\
+dt\wg\sum_i\Big(\sum_\al\pdr{(D_x^i f)}{u_\al}\Psi_\al-
\tilde{D}_t\Psi_i\Big)\ot\pdr{}{u_i}\\
+dx\wg\sum_s\Big(\sum_\al\pdr{X^s}{u_\al}\Psi_\al+
\sum_\be\pdr{X^s}{w^\be}\psi^\be-\tilde{D}_x(\psi^s)\Big)\ot\pdr{}{w^s}\\
+dt\wg\sum_s\Big(\sum_\al\pdr{T^s}{u_\al}\Psi_\al+
\sum_\be\pdr{T^s}{w^\be}\psi^\be-\tilde{D}_t(\psi^s)\Big)\ot\pdr{}{w^s},
\end{multline*}
which gives the needed result.
\end{proof}

Note that relations $\Psi_i=\tilde{D}_x^i(\Psi_0)$ together with equation
\eqref{sec5:eq:hp01} are equivalent to the shadow equations. In the case
$p=1$, we call the solutions of equation \eqref{sec5:eq:hp01} the
\emph{shadows of recursion operators} in the covering $\vf$. Equations
\eqref{sec5:eq:hp02} and \veqref{sec5:eq:hp03} are exactly the relation
equations on the case under consideration.

\begin{xca}\label{sec5:pr:shadgen}
Generalize the above result to general equations using the proof similar to
that of Theorem \vref{sec4:th:non-ulin}.
\end{xca}
Thus, any element of the group $H_\CC^{1,0}(\CE,\vf)$ is of the form
\begin{equation}\label{sec5:eq:phirec}
\Psi=\sum_i\tilde{D}_x^i(\psi)\ot\pdr{}{u_i}+\sum_s\psi^s\ot\pdr{}{w^s},
\end{equation}
where the forms $\psi=\Psi_0$, $\psi^s\in\CC^1\La(\CN)$ satisfy the
system of equations \eqref{sec5:eq:hp01}--\eqref{sec5:eq:hp03}.

As a direct consequence of the above said, we obtain the following
\bco\label{sec5:co:vertshad}
Let $\Psi$ be a derivation of the form \veqref{sec5:eq:phirec} with
$\psi,\psi^s\in\CC^p\La(\CN)$. Then $\psi$ is a solution of equation
\veqref{sec5:eq:hp01} in the covering $\vf$ if and only
if $\partial_\CC^\vf(\Psi)$ is a $\vf$-vertical derivation.
\eco
We can now formulate the main result of this subsection.
\bth\label{sec5:th:recexst}
Let $\vf:\CN\xra{}\Ei$ be a covering, $S\in\sym_\vf\CE$ be a $\vf$-symmetry,
and $\psi\in\CC^1\La(\CN)$ be a shadow of a recursion operator in the
covering $\vf$. Then $\psi'=\ip_S\psi$ is a shadow of a symmetry in $\vf$,
i.e., $\tilde{\ell}_\CE(\psi')=0$.
\ethm
\begin{proof}
In fact, let $\Psi$ be a derivation of the form \eqref{sec5:eq:phirec}. Then,
due to identity \veqref{sec5:eq:Ucomp2}, one has
$$
\partial_\CC^\vf(\ip_S\Psi)=\ip_{\partial_\CC^\vf S}-
\ip_S(\partial_\CC^\vf\Psi)=-\ip_S(\partial_\CC^\vf\Psi),
$$
since $S$ is a symmetry. But, by Corollary \vref{sec5:co:vertshad},
$\partial_\CC^\vf\Psi$ is a $\vf$-vertical derivation and consequently
$\partial_\CC^\vf(\ip_S\Psi)=-\ip_S(\partial_\CC^\vf\Psi)$ is $\vf$-vertical
as well. Hence, $\ip_S\Psi$ is a $\vf$-shadow by the same corollary.
\end{proof}

Using the last result together with Theorem \vref{sec4:th:recshad1}, we can
describe the process of generating a series of symmetries by shadows of
recursion operators. Namely, let $\psi$ be a symmetry and $\om\in\CC^1
\La(\CN)$ be a shadow of a recursion operator in a covering $\vf:\CN\xra{}
\Ei$. In particular, $\psi$ is a $\vf$-shadow. Then, by Theorem
\vref{sec4:th:recshad}, there exists a covering $\vf_\psi:\CN_\psi\xra{}\CN
\xra{\vf}\Ei$ where $\re_\psi$ can be lifted to as a
$\vf_\psi$-symmetry.  Obviously, $\om$ still remains a shadow in this
new covering. Therefore, we can act by $\om$ on $\psi$ and obtain a
shadow $\psi_1$ of a new symmetry on $\CN_\psi$. By Theorem
\vref{sec4:th:recshad1}, there exists a covering, where both $\psi$ and
$\psi_1$ are realized as nonlocal symmetries. Thus we can continue the
procedure applying $\om$ to $\psi_1$ and eventually arrive to a
covering in which the whole series $\{\psi_k\}$ is realized.
\bex\label{sec5:ex:RecBur}
Let $u_t=uu_x+u_{xx}$ be the Burgers equation. Consider the one-dimensional
covering $\vf:\Ei\times\mathbb{R}\xra{}\Ei$ with the nonlocal variable $w$
and defined by the vector fields
$$
D_x^\vf=D_x+u_0\pdr{}{w},\quad
D_t^\vf=D_t+\left(\frac{u_0^2}{2}+u_1\right)\pdr{}{w}.
$$
Then it easily checked that the form
$$\om=\om_1+\frac12\om_0+\frac12\ta,$$
where $\om_0$ and $\om_1$ are the Cartan forms $d_\CC u_0$ and $d_\CC u_1$
respectively and $\ta=dw-u_0dx-(u_0^2/2+u_1)dt$, is a solution of the
equation $\tilde{\ell}_\CE^{[1]}\om=0$. If $\re_\psi$ is a symmetry of
the Burgers equation, the corresponding action of $\om$ on $\psi$ is
$$D_x\psi+\frac12\psi+\frac12D_x^{-1}\psi$$
and thus coincides with the well-known recursion operator for this equation,
see \cite{Olver}.
\eex
\begin{xca}\label{sec5:es:RecKdV}
Let $u_t=uu_x+u_{xxx}$ be the KdV equation. Consider the one-dimensional
covering $\vf:\Ei\times\mathbb{R}\xra{}\Ei$ with the nonlocal variable $w$
and defined by the vector fields
$$
D_x^\vf=D_x+u_0\pdr{}{w},\quad
D_t^\vf=D_t+\left(\frac{u_0^2}{2}+u_2\right)\pdr{}{w}.
$$
Solve the equation $\tilde{\ell}_\CE^{[1]}\om=0$ in this covering and find
the corresponding recursion operator.
\end{xca}
\bre\label{sec5:re:super}
Recursion operators can be understood as supersymmetries (cf.\
Subsection \vref{sub:supeq}) of a certain superequation naturally
related to the initial one.  To such symmetries and equations one can
apply nonlocal theory of Section \ref{sec:nonloc} and prove the
corresponding reconstruction theorems, see \cite{Kras2,KrasKers3}.
\ere

\newpage

\section{Horizontal cohomology}
\label{hc:sec}
In this section we discuss the horizontal cohomology of differential
equations, i.e., the cohomology of the horizontal de~Rham complex (see
Definition \vref{sec3:df:hordR}). This cohomology has many physically
relevant applications.  To demonstrate this, let us start with the
notion of a conserved current. Consider a differential equation
$\E$. A \emph{conserved current} is a vector-function
$J=(J_1,\dots,J_n)$, where $J_k\in\F(\E)$, which is conserved modulo the
equation, i.e., that satisfies the equation
\eq{hc.cc:eq}
\sum_{k=1}^n D_k(J_k)=0,
\end{equation}
where $D_k$ are restrictions of total derivatives to $\Ei$.
For example, take the nonlinear Schr\"odinger equation\footnote{Here
$\psi$ is a complex function and \eqref{hc.Schr:eq} is to be understood
as a system of two equations.}
\eq{hc.Schr:eq}
i\psi_t=\Delta\psi+|\psi|^2\psi,\qquad
\Delta=\sum_{j=1}^{n-1}\frac{\partial^2}{\partial x_j^2}.
\end{equation}
Then it is straightforwardly verified, that the vector-function
\[
J=(|\psi|^2,i(\bar{\psi}\psi_{x_1}-\psi\bar{\psi}_{x_1}),\dots,
i(\bar{\psi}\psi_{x_{n-1}}-\psi\bar{\psi}_{x_{n-1}}))
\]
is a conserved current, i.e., that
\[
D_t(|\psi|^2)
+\sum_{k=1}^{n-1} D_k(i(\bar{\psi}\psi_{x_k}-\psi\bar{\psi}_{x_k}))
\]
vanishes by virtue of equation \veqref{hc.Schr:eq}.

A conserved current is called \emph{trivial}, if it has the form
\eq{hc.tcc:eq}
J_k=\sum_{l=1}^n D_l(\CL_{kl})
\end{equation}
for some skew-symmetric matrix, $\norm{\CL_{kl}}$, $\CL_{kl}=-\CL_{lk}$,
$\CL_{kl}\in\F(\E)$. The name ``trivial currents'' means that they are
trivially conserved regardless to the equation under consideration.
Two conserved currents are said to be \emph{equivalent} if they differ
by a trivial one. \emph{Conservation laws} are defined to be the
equivalent classes of conserved currents.

Let us assign the horizontal $(n-1)$-form $\om_J=\sum_{k=1}^n(-1)^{k-1}
J_k\,dx_1\wedge\dots\wedge\widehat{dx_k}\wedge\dots\wedge dx_n$ to each
conserved current $J=(J_1,\dots,J_n)$. Then equations \veqref{hc.cc:eq}
and \veqref{hc.tcc:eq} can be rewritten as $\hd\om_J=0$ and
$\om_J=\hd\eta$ respectively, where
$\eta=\sum_{k\,>\,l}(-1)^{k+l}\CL_{kl}\,dx_1
\wedge\dots\wedge\widehat{dx_l}\wedge\dots\wedge\widehat{dx_k}\wedge\dots\wedge
dx_n$. Thus, we see that the horizontal cohomology group in degree
$n-1$ of the equation $\E$ consists of conservation laws of $\E$.

In physical applications one also encounters the horizontal cohomology
in degree less than $n-1$. For instance, the Maxwell equations read
\[
\hd(*F)=0,
\]
where $F$ is the electromagnetic field strength tensor and $*$ is the
Hodge star operator. Clearly $*F$ is not exact. Another reason to
consider the low-dimensional horizontal cohomology is that it appears
as an auxiliary cohomology in calculation of the BRST cohomology
\cite{BarnBrandtHenn1}. Recently, by means of horizontal cohomology
the problem of consistent deformations and of candidate anomalies has
been completely solved in cases of Yang-Mills gauge theories and of
gravity \cite{BarnBrandtHenn2,BarnBrandtHenn3}.

The horizontal cohomology plays a central role in the Lagrangian
formalism as well. Really, it is easy to see that the horizontal
cohomology group in degree $n$ is exactly the space of actions of
variational problems constrained by equation $\E$.

For computing the horizontal cohomology there is a general method
based on \css. It can be outlined as follows. The horizontal cohomology
is the term $E_1^{0,\bullet}$ of \css and thereby related to the terms
$E_1^{p,\bullet}$, $p>0$. For each $p$, such a term is also a
horizontal cohomology but with some nontrivial coefficients. The
crucial observation is that the corresponding modules of coefficients
are supplied with filtrations such that the differentials of the
associated graded complexes are linear over the functions. Hence, the
cohomology can be computed algebraically. A detailed description of
these techniques is our main concern in this and the next sections.

\subsection{$\protect\mathcal{C}$-modules on differential equations}
\label{subsec:c-mod}

Let us begin with the definition of \cms, which are left differential
modules (see Definition \vref{sec1:df:lmod}) in \cd calculus and serve
as the modules of coefficients for horizontal de~Rham complexes.

\begin{proposition}
\label{hc.cmd:prop}
The following three definitions of a \cm are equivalent\textup{:}
\begin{enumerate}
\item An \fm $Q$ is called a \emph{$\mathcal{C}$-mod\-ule}, if $Q$ is
endowed with a left module structure over the ring $\CDiff(\F,\F)$,
i.e., for any scalar \cd operator $\Delta\in\CDiff_k(\F,\F)$ there
exists an operator $\Delta_Q\in \CDiff_k(Q,Q)$, with
\begin{enumerate}
\item[(1)] $(\sum_i f_i\Delta_i)_Q=\sum_i f_i(\Delta_i)_Q,\quad
f_i\in\F$,
\item[(2)] $(\id_{\F})_Q=\id_Q$,
\item[(3)] $(\Delta_1\circ\Delta_2)_Q=(\Delta_1)_Q\circ(\Delta_2)_Q$.
\end{enumerate}
\item A \cm is a module equipped with a \emph{flat horizontal
connection}, i.e., with an action on $Q$ of the module
$\CDer=\CDer(\E)$, $X\mapsto\nabla_X$, which is
$\F$-linear\textup{:}
\[
\nabla_{fX+gY}=f\nabla_X+g\nabla_Y,\quad f,g\in\F,\quad X,Y\in\CDer,
\]
satisfies the Leibniz rule\textup{:}
\[
\nabla_X(fq)=X(f)q+f\nabla_X(q),\quad
q\in Q,\quad X\in\CDer,\quad f\in\F,
\]
and is a Lie algebra homomorphism\textup{:}
\[
[\nabla_X,\nabla_Y]=\nabla_{[X,Y]}.
\]
\item A \cm is the module of sections of a linear covering, i.e., $Q$
is the module of sections of a vector bundle $\tau\colon W\to\Ei$,
$Q=\Gamma(\tau)$, equipped with a completely integrable $n$-dimensional
linear distribution \textup{(}see Definition
\vref{sec4:df:lincov}\textup{)} on $W$ which is projected onto the
Cartan distribution on $\Ei$.
\end{enumerate}
\end{proposition}

The proof is elementary.

\begin{xca}
Show that
\begin{enumerate}
\item in coordinates, the operator $( D_i)_Q=\norm{\Delta^k_j}$ is
a matrix operator of the form
\[
\Delta^k_j= D_i\delta^k_j+\Gamma^k_{ij},\quad\Gamma^k_{ij}\in\F,
\]
where $\delta^k_j$ is the Kronecker symbol;
\item the coordinate description of the corresponding flat horizontal
connection looks as
\[
\nabla_{ D_i}(s_j)=\sum_k\Gamma^k_{ij}s_k
\]
where $s_j$ are basis elements of $Q$;
\item the corresponding linear covering has the form
\[
\tilde{D_i}=D_i+\sum_{j,k}\Gamma^k_{ij}w^j\dd{}{w^k},
\]
where $w^i$ are fiber coordinates on $W$.
\end{enumerate}
\end{xca}

Here are basic examples of \cms.
\begin{example}
The simplest example of a \cm is $Q=\F$ with the usual action of \cd
operators.
\end{example}

\begin{example}
The module of vertical vector fields $Q=\Dv=\Dv(\E)$ with the
connection
\[
\nabla_X(Y)=[X,Y]^v,\quad X\in\CDer,\quad Y\in\Dv,
\]
where $Z^v=U_{\CC}(Z)$, is a \cm.
\end{example}

\begin{example}
\label{hc.cf:exmp}
Next example is the modules of Cartan forms
$Q=\CLa{k}=\CLa{k}(\E)$. A vector field $X\in\CDer$ acts on
$\CLa{k}$ as the Lie derivative $L_X$. It is easily seen that in
coordinates we have
\[
(D_i)_{\CLa{k}}(\omega_{\sigma}^j)=\omega_{\sigma i}^j.
\]
\end{example}

\begin{example}
The infinite jet module $Q=\hJi(P)$ of an \fm $P$ is a \cm via
\[
\Delta_{\hJi(P)}(f\ji(p))=\Delta(f)\ji(p),
\]
where $\Delta\in\CDiff(\F,\F)$, $f\in\F$, $p\in P$.
\end{example}

\begin{example}
Let us dualize the previous example. It is clear that for any \fm $P$
the module $Q=\CDiff(P,\F)$ is a \cm. The action of horizontal
operators is the composition:
\[
\Delta_Q(\nabla)=\Delta\circ\nabla,
\]
where $\Delta\in\CDiff(\F,\F)$, $\nabla\in Q=\CDiff(P,\F)$.
\end{example}

\begin{example}\label{hc.ker:exmp}
More generally, let $\Delta\colon P\to P_1$ be a \cd operator and
$\psi_\infty^\Delta\colon \hJi(P)\to\hJi(P_1)$ be the corresponding
prolongation of $\Delta$. Obviously, $\psi_\infty^\Delta$ is a morphism
of \cms, i.e., a homomorphism over the ring $\CDiff(\F,\F)$, so that
$\ker\psi_\infty^\Delta$ and $\coker\psi_\infty^\Delta$ are \cms.

On the other hand, the operator $\Delta$ gives rise to the morphism of
\cms $\CDiff(P_1,\F)\to\CDiff(P,\F)$, $\nabla\mapsto\nabla\circ\Delta$.
Thus the kernel and cokernel of this map are \cms as well.
\end{example}

\begin{example}
Given two \cms $Q_1$ and $Q_2$, we can define \cm structures on $Q_1\otimes_
{\F} Q_2$ and $\Hom_{\F}(Q_1,Q_2)$ by
\begin{align*}
\nabla_X(q_1\otimes q_2)&=\nabla_X(q_1)\otimes q_2+q_1\otimes\nabla_X(q_2),\\
\nabla_X(f)(q_1)&=\nabla_X(f(q_1))-f(\nabla_X(q_1)),
\end{align*}
where $X\in\CDer$, $q_1\in Q_1$, $q_2\in Q_2$, $f\in\Hom_{\F}(Q_1,Q_2)$.

For instance, one has \cm structures on
$Q=\hJi(P)\otimes_{\F}\CLa{k}$ and $Q=\CDiff(P,\CLa{k})$ for any \fm
$P$.
\end{example}

\begin{example}
Let $\mathfrak{g}$ be a Lie algebra and $\rho\colon\mathfrak{g}\to\gl(W)$ a
linear representation of $\mathfrak{g}$. Each $\mathfrak{g}$-valued
horizontal form $\omega\in\hL^1(\E){\otimes}_{\mathbb{R}}\,\mathfrak{g}$ that
satisfies the horizontal Maurer--Cartan condition
$\hd\omega+\frac12[\omega,\omega]=0$ defines on the module $Q$ of
sections of the trivial vector bundle $\Ei\times W\to\Ei$ the following
\cm structure:
\[
\nabla_X(q)_a=X(q)_a+\rho(\omega(X))(q_a),
\]
where $X\in\CDer$, $q\in Q$, $a\in\Ei$, and $X(q)$ means the
component-wise action.

\begin{xca}
Check that $Q$ is indeed a \cm.
\end{xca}
Such \cms are called \emph{zero-curvature representations} over $\Ei$.
Take the example of the KdV equation (in the form $u_t=uu_x+u_{xxx}$)
and $\mathfrak{g}=\sltwo(\mathbb{R})$. Then there exists a
one-parameter family of Maurer--Cartan forms
$\omega(\lambda)=A_1(\lambda)\,\hd x+A_2(\lambda) \,\hd t$, $\lambda$
being a parameter:
\[
A_1(\lambda)=\begin{pmatrix} 0& -(\lambda+u)\\
\frac16& 0
\end{pmatrix}
\]
and
\[
A_2(\lambda)=
\begin{pmatrix}
-\frac16 u_x&
-u_{xx}-\frac13 u^2+\frac13 \lambda u+\frac23 \lambda^2\\
\frac1{18}u-\frac19 \lambda&
\frac16 u_x
\end{pmatrix}.
\]
This is the zero-curvature representation used in the inverse scattering
method.
\end{example}

\begin{remark}
In parallel with left \cms one can consider \emph{right} \cms, i.e.,
right modules over the ring $\CDiff(\F,\F)$. There is a natural way to
pass from left \cms to right ones and back. Namely, for any left module
$Q$ set
\[
\mathrm{B}(Q)=Q\otimes_{\F}\hL^n(\E),
\]
with the right action of $\CDiff(\F,\F)$ on $\mathrm{B}(Q)$ given by
\begin{align*}
(q\otimes\omega)f&=fq\otimes\omega=q\otimes f\omega,\quad f\in\F,\\
(q\otimes\omega)X&=-\nabla_X(q)\otimes\omega-q\otimes L_X\omega,\quad
X\in\CDer.
\end{align*}
One can easily verify that $\mathrm{B}$ determines an equivalence
between the categories of left \cms and right \cms.
\end{remark}

By definition of a \cm, for a scalar \cd operator $\Delta\colon\F\to
\F$ there exists the extension $\Delta_Q\colon Q\to Q$ of $\Delta$ to
the \cm $Q$. Similarly to Lemma \vref{sec1:le:lmod} one has more: for
any \cd operator $\Delta\colon P\to S$ there exists the extension
$\Delta_Q\colon P\otimes_{\F}Q\to S\otimes_{\F}Q$.

\begin{proposition}
\label{hc.mq:prop}
Let $P,S$ be $\F$-modules. Then there exists a unique mapping
\[
\CDiff_k(P,S)\to\CDiff_k(P\otimes_{\F}Q,S\otimes_{\F}Q),\qquad
\Delta\mapsto\Delta_Q,
\]
such that the following conditions hold:
\begin{enumerate}
\item if $P=S=\F$ then the mapping is given by the \cm structure on $Q$,
\item $(\sum_i f_i\Delta_i)_Q=\sum_i f_i(\Delta_i)_Q,\quad f_i\in\F$,
\item if $\Delta\in\CDiff_0(P,S)=\Hom_{\F}(P,S)$ then
      $\Delta_Q=\Delta\otimes_{\F}\id_Q$,
\item if $R$ is another $\F$-module and $\Delta_1:P\to S,\Delta_2:S\to R$
are \cd operators, then $(\Delta_2\circ\Delta_1)_Q=(\Delta_2)_Q\circ
(\Delta_1)_Q$.
\end{enumerate}
\end{proposition}
\begin{proof}
The uniqueness is obvious. To prove the existence consider the family of
operators $\Delta(p,s^*)\colon\F\to\F$, $p\in P$, $s^*\in S^*=\Hom_{\F}(S,
\F)$, $\Delta(p,s^*)(f)=s^*(\Delta(fp))$, $f\in\F$. Clearly, the operator
$\Delta$ is defined by the family $\Delta(p,s^*)$. The following statement is
also obvious.
\begin{xca}
For the family of operators $\Delta[p,s^*]\in\CDiff_k(\F,\F)$, $p\in P$, $s^*
\in S^*$, we can find an operator $\Delta\in\CDiff_k(P,S)$ such that $\Delta
[p,s^*]=\Delta(p,s^*)$, if and only if
\begin{align*}
\Delta[p,\sum_i f_is_i^*]&=\sum_i f_i\Delta[p,s_i^*],\\
\Delta[\sum_i f_ip_i,s^*]&=\sum_i\Delta[p_i,s^*]f_i.
\end{align*}
\end{xca}
In view of this exercise, the family of operators
\[
\Delta_Q[p\otimes q,s^*\otimes q^*](f)=q^*(\Delta(p,s^*)_Q(fq))
\]
uniquely determines the operator $\Delta_Q$.
\end{proof}

\subsection{The horizontal de~Rham complex}
\label{hc.hdrc:subsec}
Consider a complex of \cd operators
$\dotsb\xra{}P_{i-1}\xra{\Delta_i}P_i\xra {\Delta_{i+1}} P_{i+1} \xra{}
\dotsb$. Multiplying it by a \cm $Q$ and taking into account
Proposition \vref{hc.mq:prop}, we obtain the complex
\[
\dotsb \xra{} P_{i-1}\otimes Q \xra{(\Delta_i)_Q} P_i\otimes Q
\xra{(\Delta_{i+1})_Q} P_{i+1}\otimes Q \xra{} \dotsb.
\]
Applying this construction to the horizontal de~Rham complex, we get
\emph{horizontal de~Rham complex with coefficients in $Q$}:
\[
0 \xra{} Q \xra{\hd_Q} \hL^1\otimes_{\F}Q \xra{\hd_Q} \dotsb
\xra{\hd_Q} \hL^n \otimes_{\F}Q \xra{} 0,
\]
where $\hL^i=\hL^i(\E)$.

The cohomology of the horizontal de~Rham complex with coefficients in
$Q$ is said to be \emph{horizontal cohomology} and is denoted by
$\hH^i(Q)$.

\begin{xca}
Proof that the differential $\hd=\hd_Q$ can also be defined by
\begin{align*}
(\hd q)(X)&=\nabla_X(q),\quad q\in Q,\\
\hd(\omega\otimes q)&=\hd\omega\otimes q+(-1)^p\omega\wedge\hd q,\quad
\omega\in\hL^p.
\end{align*}
\end{xca}
One easily sees that a morphism $f\colon Q_1\to Q_2$
of \cms gives rise to a cochain mapping of the de~Rham complexes:
\[
\begin{CD}
0 @>>> Q_1 @>\hd>> \hL^1\otimes_{\F}Q_1 @>\hd>> \dotsb @>\hd>>
\hL^n\otimes_{\F}Q_1 @>>> 0 \\
@.   @VVV   @VVV   @.   @VVV   @. \\
0 @>>>Q_2 @>\hd>> \hL^1\otimes_{\F}Q_2 @>\hd>> \dotsb @>\hd>>
\hL^n\otimes_{\F}Q_2 @>>> 0.
\end{CD}
\]

Let us discuss some examples of horizontal de~Rham complexes.

\begin{example}
The horizontal de~Rham complex with coefficients in the module $\hJi(P)$
\[
0 \xra{} \hJi(P) \xra{\hd} \hL^1\otimes\hJi(P) \xra{\hd}
\hL^2\otimes\hJi(P) \xra{\hd} \dotsb \xra{\hd} \hL^n\otimes\hJi(P)\xra{}0
\]
is the project limit of the \emph{horizontal Spencer
complexes}
\eq{hc.hsp:eq}
0 \xra{} \J^k(P) \xra{\bar{S}} \hL^1\otimes\J^{k-1}(P) \xra{\bar{S}}
\hL^2\otimes\J^{k-2}(P) \xra{\bar{S}}\dotsb,
\end{equation}
where $\bar{S}(\omega\otimes\hj_l(p))=\hd\omega\otimes\hj_{l-1}(p)$. As
usual Spencer complexes, they are exact in positive degrees and
\[
H^0(\hL^{\bullet}\otimes\J^{k-\bullet}(P))=P.
\]
Recall that one proves this fact by considering the commutative diagram
\[
\minCDarrowwidth=21pt
\begin{CD}
@. 0 @. 0 @. 0 @.  \\
@. @VVV @VVV @VVV @.\\
0 @>>> \hS^k\otimes P @>>> \J^k(P) @>>> \J^{k-1}(P)
@>>> 0 \\
@.  @VV\bar{\delta} V @VV \bar{S} V @VV \bar{S} V @. \\
0 @>>> \hL^1\otimes\hS^{k-1}\otimes P @>>>
\hL^1\otimes\J^{k-1}(P) @>>> \hL^1\otimes\J^{k-2}(P) @>>> 0\\
@.  @VV\bar{\delta} V @VV \bar{S} V @VV \bar{S} V @. \\
0 @>>> \hL^2\otimes\hS^{k-2}\otimes P @>>>
\hL^2\otimes\J^{k-2}(P) @>>> \hL^2\otimes\J^{k-3}(P) @>>> 0\\
@.  @VV\bar{\delta} V @VV \bar{S} V @VV \bar{S} V @. \\
@. \vdots @. \vdots @. \vdots @.
\end{CD}
\]
(see page \pageref{sec1:pg:jetsp}).
\begin{xca}
Multiply this diagram by a \cm $Q$ (possibly of
infinite rank) and prove that the complex
\[
0 \xra{} \hJi(P)\hatotimes Q \xra{\hd} \hL^1\otimes\hJi(P) \hatotimes Q
\xra{\hd} \dotsb \xra{\hd} \hL^n\otimes\hJi(P) \hatotimes Q \xra{} 0
\]
is exact in positive degrees and
\[
H^0(\hL^{\bullet}\otimes\hJi(P)\hatotimes Q)= P\otimes Q.
\]
Here
\[
\hJi(P)\hatotimes Q=\projlim\J^k(P)\otimes Q.
\]
\end{xca}
\end{example}

\begin{example}
The dualization of the previous example is as follows. The
coefficient module is $\CDiff(P,\F)$. The corresponding horizontal
de~Rham complex multiplied by a \cm $Q$ has the form
\begin{multline*}
0 \xra{} \CDiff(P,\F)\otimes Q \xra{\hd} \CDiff(P,\hL^1)\otimes Q
\xra{\hd}\dotsb\\
\dotsb \xra{\hd} \CDiff(P,\hL^n)\otimes Q \xra{} 0.
\end{multline*}
As in the previous example, it is easily shown that
\begin{align*}
H^i(\CDiff(P,\hL^{\bullet})\otimes Q)&=0\quad\text{for $i<n$,} \\
H^n(\CDiff(P,\hL^{\bullet})\otimes Q)&=\hat{P}\otimes Q,
\end{align*}
where $\hat{P}=\Hom_{\F}(P,\hL^n)$.

One can use this fact to define the notion of adjoint \cd operator
similarly to Definition \vref{sec2:df:adj}. The analog of
Proposition \vref{sec2:pr:adjprop} remains valid for \cd operators.
\end{example}

\begin{example}
\label{hc.kr:exmp}
Take the \cm
\[
Q=\bigoplus_p\Dv(\CLa{p})=\bigoplus_p\Hom_{\F}(\CLa{1},\CLa{p}).
\]
The horizontal de~Rham complex with coefficients in $Q$ can be written
as
\[
0\xra{}\Dv\xra{}\Dv(\Lambda^1)\xra{}\Dv(\Lambda^2)\xra{}\dotsb
\]
\begin{proposition}
The differential $d_{\Dv(\CLa{p})}$ of this complex is equal to
$-\partial_{\CC}$ \textup{(}see page~\pageref{sec5:eq:Ccomp}\textup{)},
so that the complex coincides up to sign with the complex
\veqref{sec5:eq:Ccomp}.
\end{proposition}
\begin{proof}
Take a vertical vector field $Y\in\Dv$ and an arbitrary vector field
$Z$. By \veqref{sec5:eq:FN4} we obtain (cf.\ the proof of Theorem
\vref{sec5:th:H01}) $\ip_Z\partial_{\CC}Y=[Z^v-Z,Y]^v$. Hence,
$\partial_{\CC}(\Dv)\subset\Dv\otimes\hL^1$ and
$\left.\partial_{\CC}\right|_{\Dv}=-d_{\Dv}$. This together with
formula \veqref{sec5:eq:bicomp3} and Remark \vref{sec5:re:coinc} yields
$\partial_{\CC}(\Dv(\CLa{p})\otimes\hL^q)
\subset\Dv(\CLa{p})\otimes\hL^{q+1}$ and $\left.\partial_{\CC}
\right|_{\Dv(\CLa{p})\otimes\hL^q}=-d_{\Dv(\CLa{p})}$.
\end{proof}
\end{example}

\subsection{Horizontal compatibility complex}
\label{hc.cc:subsec}
Consider a \cd operator $\Delta\colon P_0\to P_1$. It is clear that by
repeating word by word the construction of Subsection \vref{sub:comp}
one obtains the \emph{horizontal compatibility complex}

\eq{hc.fe:eq}
P_0 \xra{\Delta} P_1 \xra{\Delta_1} P_2 \xra{\Delta_2} P_3
\xra{\Delta_3} \dotsb,
\end{equation}
which is formally exact (see the end of Subsection \ref{sub:exact}
\vpageref{sec1:pg:ecc}).

Consider the \cm $\R_{\Delta}=\ker\psi^{\Delta}_\infty$ (cf.\ Example
\vref{hc.ker:exmp}). Then by Theorem \vref{sec1:th:jetcomp} the
cohomology of complex \veqref{hc.fe:eq} is isomorphic to the horizontal
cohomology with coefficients in $\R_{\Delta}$:
\begin{theorem}
\[
\hH^i(\R_{\Delta})=H^i(P_{\bullet}).
\]
\end{theorem}

Recall that this theorem follows from the spectral sequence arguments
applied to the commutative diagram
\[\minCDarrowwidth=20pt
\begin{CD}
@. \vdots @. \vdots @. \vdots @. \\
@.  @AAA  @AAA  @AAA  @. \\
0 @>>> \hL^2\otimes\hJi(P_0) @>>> \hL^2\otimes\hJi(P_1)
@>>> \hL^2\otimes\hJi(P_2) @>>> \dotsb \\
@.  @AA\hd A  @AA\hd A  @AA\hd A  @. \\
0 @>>> \hL^1\otimes\hJi(P_0) @>>> \hL^1\otimes\hJi(P_1)
@>>> \hL^1\otimes\hJi(P_2) @>>> \dotsb \\
@.  @AA\hd A  @AA\hd A  @AA\hd A  @. \\
0 @>>> \hJi(P_0) @>>> \hJi(P_1) @>>> \hJi(P_2) @>>> \dotsb \\
@.  @AAA  @AAA  @AAA  @. \\
@.  0  @. 0 @. 0  @.
\end{CD}
\]

Let us multiply this diagram by a \cm $Q$. This yields
\eq{hc.rd:eq}
\hH^i(\R_{\Delta}\hatotimes Q)=H^i(P_{\bullet}\otimes Q),
\end{equation}
where $\R_{\Delta}\hatotimes Q=\projlim\R_{\Delta}^l\otimes Q$, with
$\R_{\Delta}^l=\ker\psi^{\Delta}_{k+l}$, $\ord\Delta\le k$.

We can dualize our discussion. Namely, consider the commutative diagram
\[
\minCDarrowwidth=19pt
\begin{CD}
@. \vdots @. \vdots @. \vdots @. \\
@.  @VVV  @VVV  @VVV  @. \\
0 @<<< \CDiff(P_0,\hL^{n-2}) @<<< \CDiff(P_1,\hL^{n-2})
@<<< \CDiff(P_2,\hL^{n-2}) @<<< \dotsb \\
@.  @VV\hd V  @VV\hd V  @VV\hd V  @. \\
0 @<<< \CDiff(P_0,\hL^{n-1}) @<<< \CDiff(P_1,\hL^{n-1})
@<<< \CDiff(P_2,\hL^{n-1}) @<<< \dotsb \\
@.  @VV\hd V  @VV\hd V  @VV\hd V  @. \\
0 @<<< \CDiff(P_0,\hL^n) @<<< \CDiff(P_1,\hL^n) @<<< \CDiff(P_2,\hL^n)
@<<< \dotsb \\
@.  @VVV  @VVV  @VVV  @. \\
@.  0  @. 0 @. 0  @.
\end{CD}
\]

As above, we readily obtain
\[
\hH^i(\R_{\Delta}^*)=H_{n-i}(\hat{P}_{\bullet})
\]
and, more generally,
\eq{hc.rdst:eq}
\hH^i(\R_{\Delta}^*\otimes Q)=H_{n-i}(\hat{P}_{\bullet}\otimes Q),
\end{equation}
where $\R_{\Delta}^*=\Hom(\R_{\Delta},\F)$. The homology in the right-hand
side of these formulae is the homology of the complex
\[
\hat{P}_0 \xla{\Delta^*} \hat{P}_1 \xla{\Delta_1^*} \hat{P}_2
\xla{\Delta_2^*} \hat{P}_3 \xla{\Delta_3^*} \dotsb,
\]
dual to the complex \eqref{hc.fe:eq}.

\subsection{Applications to computing \kcc groups}
\label{hc.apkcc:subsec}

Let $\E$ be an equation,
\[
P_0=\varkappa \xra{\ell_{\E}} P_1 \xra{\Delta_1} P_2 \xra{\Delta_2} P_3
\xra{\Delta_3} P_4 \xra{\Delta_4} \dotsb
\]
the compatibility complex for the operator of universal linearization,
$\vk=\F(\E,\pi)$. Take a \cm $Q$.
\begin{theorem}
$\hH^i(\Dv(Q))=H^i(P_{\bullet}\otimes Q)$.
\end{theorem}
\begin{proof}
The statement follows immediately from \veqref{hc.rd:eq} and
Proposition \vref{sec3:pr:fi}.
\end{proof}

Let $Q=\CLa{p}$. The previous theorem gives a method for computing
of the cohomology groups $\hH^i(\Dv(\CLa{p}))$, which are \kcc groups (see
Example \vref{hc.kr:exmp}):

\begin{corollary}
\label{hc.cor:cor}
$\hH^i(\Dv(\CLa{p}))=H^i(P_{\bullet}\otimes\CLa{p})$.
\end{corollary}

Let us describe the isomorphisms given by this corollary in an explicit
form.

Consider an element
$\sum_{i\,\in\,I}\om_i^q\otimes\ji(s_i)\in\hL^q\otimes\Dv(\CLa{p})$, where
$\om_i^q\in\hL^q\otimes\CLa{p}$, $s_i\in\vk$, which is a horizontal
cocycle. This means that
\[
\sum_{i\,\in\,I}\om_i^q\otimes\ji(\ell_{\E}(s_i))=0 \text{ and }
\sum_{i\,\in\,I}\hd\om_i^q\otimes\ji(s_i)=0.
\]
From the second equality it easily follows that there exists an element
$\sum_{i\,\in\,I_1}\om_i^{q-1}\otimes\ji(s'_i)\in\hL^{q-1}\otimes\CLa{p}
\otimes\hJi(\vk)$, such that
$\sum_{i\,\in\,I_1}\hd\om_i^{q-1}\otimes\ji(s'_i)
=\sum_{i\,\in\,I}\om_i^q\otimes\ji(s_i)$. Denote
$s_i^1=\ell_{\E}(s'_i)$. The element
$\sum_{i\,\in\,I_1}\om_i^{q-1}\otimes\ji(s_i^1)\in\hL^{q-1}\otimes\CLa{p}
\otimes\hJi(P_1)$ satisfies
\[
\sum_{i\,\in\,I_1}\om_i^{q-1}\otimes\ji(\Delta_1(s_i^1))=0 \text{ and }
\sum_{i\,\in\,I_1}\hd\om_i^{q-1}\otimes\ji(s_i^1)=0.
\]
Continuing this process, we obtain elements
$\sum_{i\,\in\,I_l}\om_i^{q-l}\otimes\ji(s_i^l)\in\hL^{q-l}\otimes\CLa{p}
\otimes\hJi(P_l)$ such that
\[
\sum_{i\,\in\,I_l}\om_i^{q-l}\otimes\ji(\Delta_l(s_i^l))=0 \text{ and }
\sum_{i\,\in\,I_l}\hd\om_i^{q-l}\otimes\ji(s_i^l)=0.
\]
For $l=q$ these formulae mean that the element
$\sum_{i\,\in\,I_q}\om_i^0\otimes\ji(s_i^q)$ represents an element of the
module $P_q\otimes\CLa{p}$ that lies in the kernel of the operator
$\Delta_ {q+1}$. This is the element that gives rise to the
cohomology class in the group $H^q(P_{\bullet}\otimes\CLa{p})$
corresponding to the chosen element of $\hL^q\otimes\Dv(\CLa{p})$.

It follows from our results that if there is an integer $k$ such that
$P_k=P_{k+1}=P_{k+2}=\dots=0$, i.e., the compatibility complex has the
form
\[
P_0=\varkappa \xra{\ell_{\E}} P_1 \xra{\Delta_1} P_2 \xra{\Delta_2} P_3
\xra{\Delta_3} \dotsb \xra{\Delta_{k-2}} P_{k-1} \xra{} 0,
\]
then
\[
H^i(\Dv(\CLa{p}))=0\quad \text{for $i\ge k$}.
\]
This result is known as the \emph{$k$-line theorem} for \kcc.

What are the values of the integer $k$ for differential equations encountered
in mathematical physics? The existence of a compatibility operator $\Delta_1$
is usually due to the existence of dependencies between the equations under
consideration: $\Delta_1(F)=0$, $\E=\{F=0\}$. The majority of systems
that occur in practice consist of independent equations and for them
$k=2$.  Such systems of differential equations are said to be
\emph{$\ell$-normal}.  In the case of $\ell$-normal equations the
\emph{two-line theorem} for \kcc holds:

\begin{theorem}[the two-line theorem]
\label{hc.2lkcc:thm}
Let a differential equation $\E$ be $\ell$-normal. Then:
\begin{enumerate}
\item $H^i(\Dv(\CLa{p}))=0$\quad for $i\ge 2$,
\item $H^0(\Dv(\CLa{p}))=\ker(\ell_{\E})_{\CLa{p}}$,
\item $H^1(\Dv(\CLa{p}))=\coker(\ell_{\E})_{\CLa{p}}$.
\end{enumerate}
\end{theorem}
Further, we meet with the case $k>2$ in gauge theories, when the
dependencies $\Delta_1(F)=0$ are given by the second Noether theorem
(see page \pageref{css.snth:page}).
For usual irreducible gauge theories, like electromagnetism,
Yang\,-\,Mills models, and Einstein's gravity, the Noether identities
are independent, so that the operator $\Delta_2$ is trivial and, thus,
$k=3$. Finally, for an $L$-th stage reducible gauge theory, one has
$k=L+3$.

\begin{remark}
For the ``empty'' equation $\Ji(\pi)$ Corollary \vref{hc.cor:cor} yields
Theorem \vref{th3.8} (the \emph{one-line theorem}).
\end{remark}

\subsection{Example: Evolution equations}
\label{hc.exeveq:subsec}
Consider an evolution equation $\E=\{F=u_t-f(x,t,u_i)=0\}$, with
independent variables $x,t$ and dependent variable $u$; $u_i$ denotes
the set of variables corresponding to derivatives of $u$ with respect
to $x$.

Natural coordinates for $\Ei$ are $(x,t,u_i)$. The total derivatives
operators $D_x$ and $D_t$ on $\Ei$ have the form
\[
D_x=\dd{}{x}+\sum_i u_{i+1}\dd{}{u_i},\quad
D_t=\dd{}{t}+\sum_i D_x^i(f)\dd{}{u_i}.
\]
The operator of universal linearization is given by
\[
\ell_{\E}=D_t-\ell_f=D_t-\sum_i\dd{f}{u_i}D_x^i.
\]

Clearly, for an evolution equation the two-line theorem holds, hence
\kcc $\hH^q(\Dv(\CLa{p}))$ is trivial for $q\ge 2$.
Now, assume that the order of the equation $\E$ is greater than or
equal to $2$, i.e., $\ord\ell_f \ge 2$. Then one has more:
\begin{theorem}\label{sec6:th:Cvanish}
For any evolution equation of order $\ge2$, one has
\[
\hH^0(\Dv(\CLa{p}))=0\quad \text{for $p\ge 2$},
\]
\end{theorem}
\begin{proof}
It follows from Theorem \vref{hc.2lkcc:thm} that
$\hH^0(\Dv(\CLa{p}))=\ker( \ell_{\E})_{\CLa{p}}$. Hence to prove the
theorem it suffices to check that the equation
\begin{equation}
\label{hc.do:eq}
(D_t-\ell_f)(\omega)=0,
\end{equation}
with $\omega\in\CLa{p}$, has no nontrivial solutions for $p\ge 2$.

To this end consider the symbol of \veqref{hc.do:eq}. Denote
$\smbl(D_x)= \theta$. The symbol of $\ell_f$ has the form
$\smbl(\ell_f)=g\theta^k$, $k=\ord\ell_f \ge 2$, where
$g=\dd{f}{u_k}$. An element $\omega\in\CLa{p}$ can be identified with
a multilinear \cd operator, so the symbol of $\omega$ is a homogeneous
polynomial in $p$ variables
$\smbl(\omega)=\delta(\theta_1,\dots,\theta_p)$. Equation
\veqref{hc.do:eq} yields
\[
[g(\theta_1^k+\dots+\theta_p^k)-g(\theta_1+\dots+\theta_p)^k]\cdot
\delta(\theta_1,\dots,\theta_p)=0.
\]
The conditions $k\ge 2$ and $p\ge 2$ obviously imply that
$\delta(\theta_1,\dots,\theta_p)=0$. This completes the proof.
\end{proof}
\begin{remark}
This proof can be generalized for determined systems of evolution
equations with arbitrary number of independent variables (see
\cite{Gessler1}).
\end{remark}

\newpage

\section{Vinogradov's $\mathcal{C}$-spectral sequence}
\label{css:sec}

\subsection{Definition of \css}
\label{defcss:subsec}
Suppose $\E\subset J^k(\pi)$ is a formally integrable differential
equation. Consider the ideal $\CC\La^*=\CC\La^*(\E)$ of the exterior
algebra $\La^*(\E)$ of differential forms on $\Ei$ generated by the
Cartan submodule $\CLa1(\E)$ (see page \pageref{sec3:pg:catrsubm}):
$\CC\La^*=\CLa1(\E)\wedge\La^*(\E)$. Clearly, this ideal and all its
powers $(\CC\La^*)^{\wedge s}=\CLa{s}\wedge\La^*$, where
$\CLa{s}=\underbrace{\CLa1\wedge\dots\wedge\CLa1}_{\text{$s$ times}}$,
is stable with respect to the operator $d$, i.e.,
\[
d((\CC\La^*)^{\wedge s})\subset(\CC\La^*)^{\wedge s}.
\]
Thus, in the de~Rham complex on $\Ei$ we have the filtration
\[
\La^*\supset\CC\La^*\supset(\CC\La^*)^{\wedge
2}\supset\dotsb\supset(\CC\La^*)^{\wedge s}\supset\dotsb.
\]
The spectral sequence $(E_r^{p,q},d_r^{p,q})$ determined by this
filtration is said to be \emph{\css}of equation $\E$. As usual $p$ is
the filtration degree and $p+q$ is the total degree.

It follows from the direct sum decomposition \veqref{sec3:eq:split}
that $E_0^{p,q}$ can be identified with $\CLa{p}\otimes\hL^q$.

\begin{xca}
Prove that under this identification the operator $d_0^{p,q}$ coincides
with the horizontal de~Rham differential $\hd_{\CLa{p}}$ with
coefficients in $\CLa{p}$ (cf.\ Example \vref{hc.cf:exmp}).
\end{xca}

Thus, \css is one of two spectral sequences associated with the
variational bicomplex $(\CLa{p}\otimes\hL^q),\hd,d_\CC)$ constructed in
Subsection \ref{sub:basstr} \vpageref{sec3:pg:varbic}.

\begin{remark}
The second spectral sequences associated with the variational bicomplex
can be naturally identified with the Leray--Serre spectral sequence of
the de~Rham cohomology of the bundle $\Ei\to M$.  \end{remark}

\begin{remark}
The definition of \css given above remains valid for any object the
category \textit{Inf} (see page \pageref{sec4:pg:catInf}), whereas the
variational bicomplex exists only for an infinite prolonged equation.
\end{remark}

\begin{xca}
Prove that any morphism $F\colon\CN_1\to\CN_2$ in \textit{Inf} gives
rise to the homomorphism of \css for $\CN_2$ into \css for $\CN_1$.
\end{xca}

\subsection{The term $E_1$ for $\Ji(\pi)$}
\label{css.e1j:subsec}
Let us consider the term $E_1$ of \css for the ``empty'' equation
$\Ei=\Ji(\pi)$.

By definition the first term $E_1$ of a spectral sequence is the
cohomology of its zero term $E_0$. Thus, to describe the terms
$E_1^{p,q}(\pi)$ we must compute the cohomologies of complexes
\[
0\xra{}\CLa{p}(\pi) \xra{\hd}\CLa{p}(\pi) \otimes \hL^1(\pi) \xra{\hd}
\dotsb \xra{\hd}\CLa{p}(\pi) \otimes \hL^n(\pi) \xra{} 0.
\]

Using Proposition \vref{sec3:pr:fi}, this complex can be rewritten in
the form
\begin{multline*}
0\xra{}\CDiff_{(p)}^{\alt}(\vk(\pi),\F(\pi)) \xra{w}
\CDiff_{(p)}^{\alt}(\vk(\pi),\hL^1(\pi)) \xra{w} \dotsb \\
\xra{w}\CDiff_{(p)}^{\alt}(\vk(\pi),\hL^n(\pi)) \xra{} 0,
\end{multline*}
where $w(\Delta)=(-1)^p\hd\circ\Delta$.

Now from Theorem \vref{sec2:th:symskew} we obtain the following
description of the term $E_1$ for $\Ji(\pi)$:
\begin{theorem}
Let $\pi$ be a smooth vector bundle over a manifold
$M$, $\dim M=n$. Then\textup{:}
\begin{enumerate}
\item $E_1^{0,q}(\pi)=\hH^q(\pi)$ for all $q\ge 0$\textup{;}
\item $E_1^{p,q}(\pi)=0$ for $p>0$, $q\ne n$\textup{;}
\item $E_1^{p,n}(\pi)=L_p^{\alt}(\vk(\pi))$, $p>0$,
\end{enumerate}
where $L_p^{\alt}(\vk(\pi))$ was defined in Theorem
\vref{sec2:th:symskew}.
\end{theorem}

Since the \css converges to the de~Rham cohomology of the manifold
$J^{\infty}(\pi)$, this theorem has the following

\begin{corollary}
\label{css.h:cor}
For any smooth vector bundle $\pi$ over an $n$-dimensional smooth
manifold $M$ one has\textup{:}
\begin{enumerate}
\item  $E_r^{p,q}(\pi)=0$, $1\le r\le\infty$, if $p>0$, $q\ne n$ or
$p=0$, $q>n$\textup{;}
\item $E_1^{0,q}(\pi)=E_{\infty}^{0,q}(\pi)=H^q(J^{\infty}(\pi))
=H^q(J^0(\pi))$, $q<n$\textup{;}
\item $E_2^{p,n}(\pi)=E_{\infty}^{p,n}(\pi)=H^{p+n}(J^{\infty}(\pi))
=H^{p+n}(J^0(\pi))$, $p\ge0$.
\end{enumerate}
\end{corollary}

\begin{xca}
Prove that $H^q(\Ji(\pi))=H^q(J^0(\pi))$.
\end{xca}

We now turn our attention to the differentials $d_1^{p,n}$. They are
induced by the Cartan differential $d_\CC$. For $p=0$, we have
$d_\CC(\omega)=\ell_\omega$, $\omega\in\hL^n$. (Note that the
expression $\ell_\omega$ is correct, because $\omega$ is a horizontal
form, i.e., a nonlinear operator from $\Gamma(\pi)$ to $\La^n(M)$.)
Therefore the operator
\[
E_1^{0,n}(\pi)=\hH^n(\pi)\xra{d_1^{0,n}}E_1^{1,n}(\pi)=\hat{\vk}(\pi)
\]
is given by the formula $d_1^{0,n}([\omega])=\mu(\ell_\omega)
=\ell_{\omega}^*(1)$, where $\omega\in\hL^n(\pi)$, $[\omega]$ is the
horizontal cohomology class of $\omega$.

\begin{xca}
\label{css.Eu.ex}
Write down the coordinate expression for the operator $d_1^{0,n}$ and
show that it coincides with the standard Euler operator, i.e., with the
operator that takes a Lagrangian to the corresponding
Euler--Lagrange equation.
\end{xca}

Let us compute the operators  $d_1^{p,n}$, $p>0$.

Consider an element $\nabla\in L_p^{\alt}(\vk(\pi))$ and define the operator
$\square\in\CDiff_{(p+1)}(\vk(\pi),\hL^n(\pi))$ via
\begin{multline}
\label{css.Delta:eq}
\square(\chi_1,\dots,\chi_{p+1})=
\sum_{i=1}^{p+1}(-1)^{i+1}\re_{\chi_i}(\nabla(\chi_1,\dots,\hat{\chi}_i,
\dots,\chi_{p+1})) \\
+\sum_{1\le i<j\le p+1}
(-1)^{i+j}\nabla(\{\chi_i,\chi_j\},\chi_1,\dots,
\hat{\chi}_i,\dots,\hat{\chi}_j,\dots,\chi_{p+1}).
\end{multline}

\begin{xca}
Prove that $d_1^{p,n}(\nabla)=\mu_{(p+1)}(\square)$ (see page
\pageref{sec2:pg:mu} for the definition of $\mu_{(p+1)}$).
\end{xca}

\begin{remark}
Needless to say that this fact follows from the standard formula for
exterior differential. It needs however to be proved that one may use
this formula even though $\nabla$ as an element of
$\CDiff_{(p)}(\vk,\hL^n)$ is not skew-symmetric.
\end{remark}

From \veqref{css.Delta:eq} we get
\begin{multline*}
\square(\chi_1,\dots,\chi_{p+1})=\sum^p_{i=1}(-1)^{i+1}\re_{\chi_i}
(\nabla(\chi_1,\dots,\hat{\chi}_i,\dots,\chi_p))(\chi_{p+1}) \\
+\sum^p_{i=1}(-1)^{i+1}\nabla(\chi_1,\dots,\hat{\chi}_i,\dots,
\chi_p,\re_{\chi_i}(\chi_{p+1})) \\
+(-1)^p\re_{\chi_{p+1}}(\nabla(\chi_1,\dots,\chi_{p})) \\
+\sum_{1\le i<j\le p}(-1)^{i+j}\nabla(\{\chi_i,\chi_j\},\chi_1,\dots,
\hat{\chi}_i,\dots,\hat{\chi}_j,\dots,\chi_{p+1}) \\
+\sum_{i=1}^p(-1)^{i+p+1}\nabla(\{\chi_i,\chi_{p+1}\},\chi_1,\dots,
\hat{\chi}_i,\dots,\chi_p) \\
=\sum^p_{i=1}(-1)^{i+1}\re_{\chi_i}(\nabla(\chi_1,\dots,\hat{\chi}_i,
\dots,\chi_p))(\chi_{p+1}) \\
+\sum_{1\le i<j\le p}(-1)^{i+j}\nabla(\{\chi_i,\chi_j\},\chi_1,\dots,
\hat{\chi}_i,\dots,\hat{\chi}_j,\dots,\chi_{p+1}) \\
+\sum^p_{i=1}(-1)^{i+1}\nabla(\chi_1,\dots,\hat{\chi}_i,\dots,\chi_p,
\ell_{\chi_i}(\chi_{p+1}))+(-1)^p\ell_{\nabla(\chi_1,\dots,
\chi_{p})}(\chi_{p+1}).
\end{multline*}
Therefore
\begin{multline}
d_1^{p,n}(\nabla)(\chi_1,\dots,\chi_p)=
\mu_{(p+1)}(\square)(\chi_1,\dots,\chi_p)  \\
=\sum^p_{i=1}(-1)^{i+1}\re_{\chi_i}(\nabla(\chi_1,\dots,\hat{\chi}_i,
\dots,\chi_p)) \\
+\sum_{i<j}(-1)^{i+j}\nabla(\{\chi_i,\chi_j\},\chi_1,\dots,
\hat{\chi}_i,\dots,\hat{\chi}_j,\dots,\chi_{p}) \\
+\sum^p_{i=1}(-1)^{i+1}\ell^*_{\chi_i}(\nabla(\chi_1,\dots,
\hat{\chi}_i,\dots,\chi_p))+
(-1)^p\ell^*_{\nabla(\chi_1,\dots,\chi_p)}(1).
\label{css.chi:eq}
\end{multline}

\begin{xca}
\label{css.EL:ex}
Prove that
\[
\ell^*_{\psi(\vf)}(1)=\ell^*_{\psi}(\vf)+\ell^*_{\vf}(\psi),
\qquad \vf\in\vk(\pi), \quad \psi\in\hat{\vk}(\pi).
\]
\end{xca}
Using this formula, let us rewrite the last term of \veqref{css.chi:eq}
in the following way:
\begin{multline*}
(-1)^p\ell^*_{(\nabla(\chi_1,\dots,\chi_p))}(1)=
\frac1p\sum^p_{i=1}(-1)^i(\ell^*_{(\nabla(\chi_1,\dots,
\hat{\chi}_i,\dots,\chi_p,\chi_i))}(1) \\
=\frac1p\sum^p_{i=1}(-1)^i(\ell^*_{\nabla(\chi_1,\dots,
\hat{\chi}_i,\dots,\chi_p)}(\chi_i)+\ell^*_{\chi_i}
(\nabla(\chi_1,\dots,\hat{\chi}_i,\dots,\chi_p))).
\end{multline*}
Finally we obtain
\begin{multline*}
\left(d_1^{p,n}(\nabla)\right)(\chi_1,\dots,\chi_p)=
\sum^p_{i=1}(-1)^{i+1}\re_{\chi_i}(\nabla(\chi_1,\dots,\hat{\chi}_i,\dots,
\chi_p)) \\
+\sum_{i<j}(-1)^{i+j}\nabla(\{\chi_i,\chi_j\},\chi_1,\dots,
\hat{\chi}_i,\dots,\hat{\chi}_j,\dots,\chi_{p}) \\
+ \frac1p\sum^p_{i=1}(-1)^{i+1}((p-1)\ell^*_{\chi_i}(\nabla(\chi_1,
\dots,\hat{\chi}_i,\dots,\chi_p))-\ell^*_{\nabla(\chi_1,\dots,
\hat{\chi}_i,\dots,\chi_p)}(\chi_i)).
\end{multline*}

In particular, for $p=1$ we have
$d_1^{1,n}(\psi)(\vf)=\re_{\vf}(\psi)-\ell^*_{\psi}(\vf)=
\ell_{\psi}(\vf)-\ell^*_{\psi}(\vf)$,
$\psi\in\hat{\vk}(\pi)$, $\vf\in\vk(\pi)$,
that is
\[
d_1^{1,n}(\psi)=\ell_{\psi}-\ell^*_{\psi}.
\]

Consider the following complex, which is said to be the (global)
\emph{variational complex},
\[
0\xra{} \F(\pi) \xra{\hd}\hL^1(\pi)\xra{\hd}\dotsb
\xra{\hd}\hL^n(\pi)\xra{\bE}
E_1^{1,n}(\pi)\xra{d_1^{1,n}} E_1^{2,n}(\pi)\xra{d_1^{2,n}} \dotsb,
\]
where operator ${\bE}$ is equal to the composition of the natural
projection $\hL^n(\pi)\to \hH^n(\pi)$ and the operator $\hH^n(\pi)
\xra{d_1^{0,n}}E_1^{1,n}(\pi)$
\footnote{Below we use the notation $\bE$ for the operator
$d_1^{0,n}\colon\hH^n(\pi)\to E_1^{1,n}(\pi)$ as well.}.

In view of Corollary \vref{css.h:cor}, the cohomology of this complex
coincides with $H^*(J^0(\pi))$.

The operator $\bE$ is the Euler operator (see Exercise
\vref{css.Eu.ex}). It takes each Lagrangian density
$\omega\in\Lambda_0^n(\pi)$ to the left-hand part of the corresponding
Euler--Lagrange equation $\bE(\omega)=0$. Thus the action
functional
\[
s\mapsto\int_Mj_{\infty}(s)^*(\omega), \qquad s\in\Gamma(\pi),
\]
is stationary on the section $s$ if and only if $j_{\infty}(s)^*\left(
{{\bE}}(\omega)\right)=0$.

If the cohomology of the space $J^0(\pi)$ is trivial, then the
variational complex is exact. This immediately implies a number of
consequences. The three most important are:

\begin{enumerate}
\item $\ker\bE=\im\hd$ (``a Lagrangian with zero variational derivative
is a total divergence'');
\item $\hd\omega=0$ if and only if $\omega$ is of the form
$\omega=\hd\eta$, $\omega\in\hL^{n-1}(\pi)$ (``all zero total
divergence are total curls'');
\item $\ell_{\psi}=\ell_{\psi}^*$ if and only if $\psi$ is of
the form $\psi=\bE(\omega)$, $\psi\in\hat{\vk}(\pi)$ (this is the
solution of the inverse problem to the calculus of variations).
\end{enumerate}

Now suppose that we are given $\psi\in\hat{\vk}(\pi)$ such that
$\ell_{\psi}=\ell_{\psi}^*$. How one can find a Lagrangian $\omega$
such that $\psi=\bE(\omega)$? To this end take a one-parameter family
of fiberwise transformations $G_t\colon E\to E$, $0\le t\le1$, of the
space of the bundle $\pi\colon E\to M$, with $G_0=0$ and $G_1=\id_E$.
Consider the corresponding family of evolutionary vector field
$\re_{\vf_t}$, i.e.,
\[
\frac{d}{dt}G_t^{(\infty)*}=\re_{\vf_t}\circ G_t^{(\infty)*}
\]
for $t>0$. Let us compute the correspondent Lie derivative
$\re_{\vf_t}(\psi)$ (which is different from the usual
``component-wise'' derivative). Take $\Omega\in\hL^n\otimes\CLa{1}$,
$d\Omega=0$, that represents $\psi$. Then $\re_{\vf_t}(\Omega)
=d_{\CC}(\Omega(\re_{\vf_t}))$. Therefore
$\re_{\vf_t}(\psi)=\bE(\psi(\vf_t))$. Hence
\[
\frac{d}{dt}G_t^{(\infty)*}(\psi)=\bE(G_t^{(\infty)*}(\psi(\vf_t))),
\]
and integrating this with respect to $t$, we obtain the following
\emph{homotopy} (or \emph{inverse}) formula
\[
\psi=\bE\left(\int_0^1G_t^{(\infty)*}(\psi(\vf_t))\,dt\right).
\]
Take, for instance, $G_t(e_x)=te_x$, $e_x\in E_x=\pi^{-1}(x)$. Then
$\vf_t^i=\dfrac{u^i}{t}$ and we have
\[
\psi=\bE\left(\int_0^1\sum_i u^i\psi^i(x,tu_{\sigma}^j)\,dt\right).
\]

\begin{xca}
\label{css.uf:ex}
Let $\Delta\in\CDiff(P,\hL^n(\pi))$. Using the Green formula and
Exercise \vref{css.EL:ex}, prove that for any $p\in P$ one has
\[
\bE(\Delta(p))=\ell_p^*(\Delta^*(1))+\ell^*_{\Delta^*(1)}(p).
\]
Deduce from this formula that for any $\vf\in\vk(\pi)$ and
$\omega\in\hL^n(\pi)$ the following equality holds
\[
\bE(\re_{\vf}(\omega))=\re_{\vf}(\bE(\omega))+
\ell^*_{\vf}(\bE(\omega)).
\]
\end{xca}

\begin{xca}
Let $J=(J_0,J_1,\dots,J_n)$ be a conserved current for an evolution
equation $\E=\{u_t=f(t,x,u,u_x,u_{xx},\dotsc)\}$. Using the previous
exercise, prove that the vector-function $\psi=\bE(J_0)$, where $J_0$
is the $t$-component of the conserved current that is regarded as a
function of $(t,x,u,u_x,u_{xx},\dotsc)$, satisfies the equation
\[
D_t(\psi)+\ell_f^*(\psi)=0
\]
(cf.\ Theorem \vref{css.gf:thm}).
\end{xca}

\subsection{The term $E_1$ for an equation}
\label{css.e1eq:subsec}

Let $\E$ be an equation,
\[
P_0=\varkappa \xra{\ell_\E} P_1 \xra{\Delta_1} P_2 \xra{\Delta_2} P_3
\xra{\Delta_3} P_4 \xra{\Delta_4} \dotsb
\]
be the compatibility complex for the universal linearization operator,
and
\[
\hat{P_0}=\hat{\varkappa}\xla{\ell_E^*}\hat{P_1}\xla{\Delta_1^*}\hat{P}_2
\xla{\Delta_2^*}\hat{P}_3\xla{\Delta_3^*}\hat{P}_4\xla{\Delta_4^*}\dotsb
\]
be the dual complex. Take a \cm $Q$.
\begin{theorem}
\label{css.mth:thm}
For any equation $\E$ and a \cm $Q$ one has
\[
\hH^{n-i}(\CLa1\otimes Q)=H_i(\hat{P}_{\bullet}\otimes Q).
\]
\end{theorem}

\begin{proof}
The statement follows immediately from \veqref{hc.rdst:eq} and
Proposition \vref{sec3:pr:fi}.
\end{proof}

Let $Q=\CLa{p}$. The theorem gives a method for computing \css. Namely,
since the term $E_1^{p,q}=\hH^q(\CLa{p})$ of \css is a direct summand
in the cohomology group $\hH^q(\CLa1\otimes\CLa{p-1})$, we have a
description for the first term of \css. Thus:

\begin{corollary}
\label{css.cssc:cor}
The term $E_1^{p,q}$ of \css is the skew-symmetric part of the group
$H_{n-q}(\hat{P}_{\bullet}\otimes\CLa{p-1})$.
\end{corollary}

It is useful to describe the isomorphisms given by this corollary in
an explicit form.

Consider an operator $\nabla\in\CDiff(\varkappa,\hL^q\otimes\CLa{p-1})$
that represents an element of $E_1^{p,q}$. This means that
\[
\hd\circ\nabla=\nabla_1\circ\ell_\E
\]
for an operator $\nabla_1\in\CDiff(P_1,\hL^{q+1}\otimes\CLa{p-1})$.
Applying the operator $\hd$ to both sides of this formula and using
Exercise \vref{sec1:ex:comop}, we get
\[
\hd\circ\nabla_1=\nabla_2\circ\Delta_1
\]
for some operator $\nabla_2\in\CDiff(P_2,\hL^{q+2}\otimes\CLa{p-1})$.
Continuing this process, we obtain operators
$\nabla_i\in\CDiff(P_i,\hL^{q+i} \otimes\CLa{p-1})$, $i=1,2,\dots,n-q$,
such that
\[
\hd\circ\nabla_{i-1}=\nabla_i\circ\Delta_{i-1}.
\]
For $i=n-q$, this formula means that the operator $\nabla_{n-q}\in
\CDiff(P_{n-q},\hL^n\otimes\CLa{p-1})$ represents an element of the
module $\hat{P} _{n-q}\otimes\CLa{p-1}$ that lies in the kernel of the
operator $\Delta_ {n-q-1}^*$. This is the element that gives rise to
the homology class in $H_{n-q}(\hat{P}_{\bullet}\otimes\CLa{p-1})$
corresponding to the chosen element of $E_1^{p,q}$.

If the compatibility complex has the length $k$,
\[
P_0=\varkappa \xra{\ell_E} P_1 \xra{\Delta_1} P_2 \xra{\Delta_2} P_3
\xra{\Delta_3} \dotsb \xra{\Delta_{k-2}} P_{k-1} \xra{} 0,
\]
then $E_1^{p,q}=0$ for $p>0$ and $q\le n-k$.
This is the \emph{$k$-line theorem} for \css.

In the case $k=2$, i.e., for $\ell$-normal equations, the
\emph{two-line theorem} holds:

\begin{theorem}[the two-line theorem]
Let $\E$ be an $\ell$-normal differential equation. Then:
\begin{enumerate}
\item $E_1^{p,q}=0$\quad for $p>0$ and $q\le n-2$,
\item $E_1^{p,n-1}\subset\ker(\ell_E^*)_{\CLa{p-1}}$\quad for $p>0$,
\item $E_1^{p,n}\subset\coker(\ell_E^*)_{\CLa{p-1}}$\quad for $p>0$.
\end{enumerate}
\end{theorem}

This theorem has the following elementary

\begin{corollary}
The terms $E_r^{p,q}(\E)$ of \css satisfy the following:
\begin{enumerate}
\item $E_r^{p,q}(\E)=0$ if $p\ge1$, $q\ne n-1,n$, $1\le r\le\infty$\textup{;}
\item $E_3^{p,q}(\E)=E_{\infty}^{p,q}(\E)$\textup{;}
\item $E_1^{0,q}(\E)=E_{\infty}^{0,q}(\E)=H^q(\Ei)$, $q\le n-2$\textup{;}
\item $E_2^{0,n-1}(\E)=E_{\infty}^{0,n-1}(\E)=H^{n-1}(\Ei)$\textup{;}
\item $E_2^{1,n-1}(\E)=E_{\infty}^{1,n-1}(\E)$.
\end{enumerate}
\end{corollary}

\begin{example}
For an evolution equation $\E=\{F=u_t-f(x,t,u,u_x,u_{xx}\dotsc)=0\}$
the two-line theorem implies that \css is trivial for $q\ne1,2$, $p>0$,
and exactly as in Example \vref{hc.exeveq:subsec} one proves that
$E_1^{p,1}=0$ for $p\ge3$.
\end{example}

\subsection{Example: Abelian $p$-form theories}
\label{css.exapf:subsec}

Let $M$ be a (pseudo-)Rieman\-nian manifold and $\pi\colon E\to M$ the
$p$-th exterior power of the cotangent bundle over $M$, so that a
section of $\pi$ is a $p$-form on $M$. Evidently, on the jet space
$J^\infty(\pi)$ there exists a unique horizontal form
$A\in\hL^p(J^\infty(\pi))$ such that $j_ \infty^*(\omega)(A)=\omega$
for all $\omega\in\Lambda^p(M)$. Consider the equation $\E=\{F=0\}$,
with $F=\hd{*}\hd A$, where $*$ is the Hodge star operator. Our aim is
to calculate the terms of \css $E_1^{i,q}(\E)$ for $q\le n-2$. We shall
assume that $1\le p<n-1$ and that the manifold $M$ is topologically
trivial.

Obviously, we have $P_0=\varkappa=\hL^p$, $P_1=\hL^{n-p}$, and
$\ell_\E=\hd{*}\hd\colon\hL^p\to\hL^{n-p}$. Taking into account Example
\vref{sec1:ex:pf}, we see that the compatibility complex for $\ell_\E$
has the form
\begin{equation}
\label{eq:ccpf}
\begin{CD}
\hL^p @>\ell_\E>> \hL^{n-p} @>\hd>> \hL^{n-p+1} @>\hd>> \dotsb
@>\hd>> \hL^n @>>> 0 \\
@|   @|   @|   @.  @|  \\
P_0 @. P_1 @. P_2 @. @. P_{k-1}
\end{CD}
\end{equation}
Thus $k=p+2$ and the $k$-line theorem yields $E_1^{i,q}=0$ for $i>0$
and $q<n-p-1$. Since \css converges to the de~Rham cohomology of $\Ei$,
which is trivial, we also get $E_1^{0,q}=0$ for $0<q<n-p-1$, and $\dim
E_1^{0,0}=1$, i.e., $\hH^1=\hH^2=\dots=\hH^{n-p-2}=0$ and
$\dim\hH^0=1$. Next, consider the terms $E_1^{i,q}$ for $n-p-1\le
q<2(n-p-1)$ and $i>0$. In view of Corollary \vref{css.cssc:cor} one has
\[
E_1^{i,q}\subset\hH^{q-(n-p-1)}(\CLa{i-1})=E_1^{i-1,q-(n-p-1)},
\]
because the complex dual to the compatibility complex \eqref{eq:ccpf}
has the form
\[
\begin{CD}
\hL^{n-p} @<\ell_\E^*<< \hL^p @<\hd<< \hL^{p-1} @<\hd<< \dotsb
@<\hd<< \F @<<< 0. \\
@|   @|   @|   @.  @|  \\
\hat{P}_0 @. \hat{P}_1 @. \hat{P}_2 @. @. \hat{P}_{p+1}
\end{CD}
\]
(Throughout, it is assumed that $q\le n-2$.) Thus we obtain
$E_1^{i,q}=0$ for $n-p-1<q<2(n-p-1)$, $i>0$ and $\dim E_1^{1,n-p-1}=1$.
Again, taking into account that the spectral sequence converges to the
trivial cohomology, we get $E_1^{0,q}=0$ for $n-p-1<q<2(n-p-1)$ and
$\dim E_1^{0,n-p-1}=1$. In addition, the map $d_1^{0,n-p-1}\colon
E_1^{0,n-p-1}\to E_1^{1,n-p-1}$ is an isomorphism. Explicitly, one
readily obtains that the one-dimensional space $E_1^{0,n-p-1}$ is
generated by the element $*\hd A\in\hL^{n-p-1}$ and the map
$d_1^{0,n-p-1}$ takes this element to the operator
$*\hd\colon\varkappa=\hL^p\to\hL^{n-p-1}$, which generates the space
$E_1^{1,n-p-1}$.

Further, let us consider the terms $E_1^{i,q}$ for $2(n-p-1)\le q<3(n-p-1)$.
Arguing as before, we see that all these terms vanish unless $q=2(n-p-1)$ and
$i=0,1,2$, with $\dim E_1^{1,2(n-p-1)}=1$ and $\dim E_1^{i,2(n-p-1)}\le 1$,
$i=0,2$. To compute the terms $E_1^{i,2(n-p-1)}$ for $i=0$ and $i=2$, we have
to consider two cases: $n-p-1$ is even and $n-p-1$ is odd (see Diagram
\ref{css.pic}).

\begin{figure}
\begin{picture}(340,200)
\put(70,2){$n-p-1$ is even}
\put(250,2){$n-p-1$ is odd}
\put(60,15){\vector(1,0){100}}
\put(60,15){\vector(0,1){185}}
\put(240,15){\vector(1,0){100}}
\put(240,15){\vector(0,1){185}}
\put(53,195){$q$}
\put(156,20){$i$}
\put(233,195){$q$}
\put(336,20){$i$}
\multiput(78,15)(18,0){5}{\line(0,1){185}}
\multiput(258,15)(18,0){5}{\line(0,1){185}}
\put(55,180){\llap{$3(n-p-1)$}}
\put(55,125){\llap{$2(n-p-1)$}}
\put(55,70){\llap{$n-p-1$}}
\put(235,180){\llap{$3(n-p-1)$}}
\put(235,125){\llap{$2(n-p-1)$}}
\put(235,70){\llap{$n-p-1$}}
\put(69,19){\circle*{3}}
\multiput(69,72)(0,55){3}{\circle*{3}}
\multiput(87,72)(0,55){3}{\circle*{3}}
\multiput(72,72)(0,55){3}{\vector(1,0){12}}
\put(249,19){\circle*{3}}
\multiput(249,72)(18,55){3}{\circle*{3}}
\multiput(267,72)(18,55){3}{\circle*{3}}
\multiput(252,72)(18,55){3}{\vector(1,0){12}}
\end{picture}
\caption{}
\label{css.pic}
\end{figure}

In the first case, the map $d_1^{1,2(n-p-1)}\colon E_1^{1,2(n-p-1)}\to
E_1^ {2,2(n-p-1)}$ is trivial. Indeed, the operator $(*\hd
A)\wedge*\hd\colon \varkappa=\hL^p\to\hL^{2(n-p-1)}$, which generates
the space $E_1^{1,2(n-p- 1)}$, under the mapping $d_1^{1,2(n-p-1)}$ is
the antisymmetrization of the operator
$(\omega_1,\omega_2)\mapsto(*\hd\omega_1)\wedge(*\hd\omega_2)$,
$\omega_i\in\varkappa=\hL^p$. But this operator is symmetric, so that
$d_1^ {1,2(n-p-1)}=0$. Consequently, $E_1^{2,2(n-p-1)}=0$ and $\dim
E_1^{0,2(n-p- 1)}=1$. This settles the case when $n-p-1$ is even.

In the case when $n-p-1$ is odd, the operator $(\omega_1,\omega_2)
\mapsto(*\hd\omega_1)\wedge(*\hd\omega_2)$ is skew-symmetric, hence the
map $d_1^{1, 2(n-p-1)}$ is an isomorphism. Thus, $\dim
E_1^{2,2(n-p-1)}=1$ and $E_1^{0, 2(n-p-1)}=0$.

Continuing this line of reasoning, we obtain the following result.

\begin{theorem} For $i=q=0$ one has $\dim E_1^{0,0}=1$. If either or
both $i$ and $q$ are positive, there are two cases:

\begin{enumerate}
\item if $n-p-1$ is even then
\[
\dim E_1^{i,q}=
\begin{cases}
1&\text{for $i=l(n-p-1)$ and $q=0$, $1$}, \\
0&\text{otherwise};
\end{cases}
\]
\item if $n-p-1$ is odd then
\[
\dim E_1^{i,q}=
\begin{cases}
1&\text{for $i=l(n-p-1)$ and $q=l-1$, $l$}, \\
0&\text{otherwise}.
\end{cases}
\]
\end{enumerate}
Here $1\le l<\dfrac{n-1}{n-p-1}$.
\end{theorem}

In other words, let $\bar{\mathcal{A}}$ be the exterior algebra
generated by two forms: $\omega_1=*\hd A\in\hL^{n-p-1}$ and
$\omega_2=\hd_1(\omega_1)=*\hd \in\hL^{n-p-1}\otimes\CLa{1}$; then we
see that the space $\bigoplus_{i,q\le n-2}E_1^{i,q}$ is isomorphic to
the subspace of $\bar{\mathcal{A}}$ containing no forms of degree
$q>n-2$.

\subsection{Conservation laws and generating functions}
\label{css.genfun:subsec}
We start by describing the differentials $d_1^{0,n-1}$ and
$d_1^{1,n-1}$ for an $\ell$-normal equation since they directly relate
to the theory of conservation laws.

Suppose that an $\ell$-normal equation $\E\subset J^k(\pi)$ is given
by a section $F\in\F(\pi,\xi)=P$.

\begin{proposition}
The operator
\[
d_1^{0,n-1}\colon E_1^{0,n-1}(\E)=\hH^{n-1}(\E)\to E_1^{1,n-1}(\E)
=\ker\left(\ell_\E\right)^*\subset\hat P
\]
has the form
\[
d_1^{0,n-1}(h)=\square^*(1),
\]
where $h=[\omega]\in\hH^{n-1}(\E)$, $\omega\in\hL^{n-1}(\E)$
and $\square\in\CDiff(P,\hL^n(\E))$ is an operator satisfying
$\hd\omega=\square(F)$.
\end{proposition}

\begin{proof}
We have $\hd\circ\ell_\omega=\square\circ\ell_\E$. Thus $\square$
is an operator that represents the element $d_1^{0,n-1}(h)\in
E_1^{1,n-1}(\E)$. Hence $d_1^{0,n-1}(h)=\square^*(1)$.
\end{proof}

\begin{proposition}
\label{css.2nm1:prop}
The term $E_1^{2,n-1}(\E)$ can be described as the quotient
\[
\{\,\nabla\in\CDiff(\vk,\hat P)\mid
\ell_{\E}^*\circ\nabla=\nabla^*\circ\ell_{\E}\,\}
/\theta,
\]
where $\theta=\{\,\square\circ\ell_{\E}\mid\square\in\CDiff(P,\hat P),
\square^*=\square\,\}$.
\end{proposition}

\begin{proof}
Take a horizontal $(n-1)$-cocycle with coefficients in
$\CLa1\otimes\CLa1$. Let an operator $\Delta\in\CDiff(\vk,\hat P)$
corresponds to this cocycle by Theorem \vref{css.mth:thm}. Then there
exists an operator $A\in\CDiff(P,\hat{P})$ such that
$\ell_\E^*\circ\Delta=A\circ\ell_\E$. By the Green formula we have
\[
\langle\ell_\E^*(\Delta(\chi_1)),\chi_2\rangle
-\langle\Delta(\chi_1),\ell_\E(\chi_2)\rangle
=\hd(\Delta_1(\chi_1,\chi_2)),
\]
where $\chi_1,\chi_2\in\vk$, and
$\Delta_1\in\CDiff_{(2)}(\vk,\La^{n-1})$. The cocycle under
consideration belongs to $E_1^{2,n-1}$, if the operator $\Delta_1$ is
skew-symmetric:
\[
\Delta_1(\chi_1,\chi_2)=-\Delta_1(\chi_2,\chi_1)\bmod K,
\]
where $K\subset\CDiff_{(2)}(\vk,\hL^{n-1})$ is the submodule consisting
of the operators of the form
$\gamma(\chi_1,\chi_2)=\gamma_1(\ell_\E(\chi_1),\chi_2)
+\gamma_2(\ell_\E(\chi_2),\chi_1)$ for some operators
$\gamma_1,\gamma_2\in\CDiff(P,\CDiff(\vk,\hL^{n-1}))$. In this case
\begin{multline*}
\langle\ell_\E^*(\Delta(\chi_1)),\chi_2\rangle
-\langle\Delta(\chi_1),\ell_\E(\chi_2)\rangle
=-\langle\ell_\E^*(\Delta(\chi_2)),\chi_1\rangle
+\langle\Delta(\chi_2),\ell_\E(\chi_1)\rangle \\
=-\langle\ell_\E^*(\Delta(\chi_2)),\chi_1\rangle
+\langle\chi_2,\Delta^*(\ell_\E(\chi_1))\rangle\\
=-\langle A(\ell_\E(\chi_2)),\chi_1\rangle
+\langle\chi_2,\ell_\E^*(A^*(\chi_1))\rangle
\end{multline*}
modulo $\hd K$. This implies $\Delta=A^*+B\circ\ell_\E$ for an operator
$B\in\CDiff(P,\hat{P})$. One has $\ell_\E^*\circ
B\circ\ell_\E=\ell_\E^*\circ\Delta-\ell_\E^*\circ
A^*=\ell_\E^*\circ\Delta-\Delta^*\circ\ell_\E$, hence $B^*=-B$. Now we
see that the operator $\nabla=\Delta-\frac12B\circ\ell_\E$ satisfies
$\ell_{\E}^*\circ\nabla=\nabla^*\circ\ell_{\E}$. The
operator $\nabla$ is defined modulo the operators of the form
$\square\circ\ell_\E$. We have $\ell_{\E}^*\circ\square\circ\ell_\E
=\ell_{\E}^*\circ\square^*\circ\ell_\E$, so that $\square^*=\square$.
\end{proof}

\begin{proposition}
The operator $d_1^{1,n-1}\colon E_1^{1,n-1}(\E)=
\ker\ell_{\E}^*\to E_1^{2,n-1}(\E)$ is given by
\[
d_1^{1,n-1}(\psi)=(\ell_{\psi}+\Delta^*)\bmod\theta,
\]
where $\Delta\in\CDiff(P,\hat{\vk})$ is an operator satisfying
$\ell_F^*(\psi)=\Delta(F)$.
\end{proposition}

\begin{proof}
By Green's formula on $\Ji(\pi)$ we have
\[
\langle\psi,\ell_F(\chi)\rangle
-\langle\ell_F^*(\psi),\chi\rangle=\hd(\square(\chi)),
\]
where $\chi\in\vk(\pi)$, $\square\in\CDiff(\vk(\pi),\hL^{n-1}(\pi))
=\CLa1(\pi)\otimes\hL^{n-1}(\pi)$.
Let us compute $\hd\circ d_\CC(\square)\in\CLa2(\pi)\otimes\hL^n(\pi)$:
\begin{multline*}
\hd(d_\CC(\square)(\chi_1,\chi_2))=\re_{\chi_1}(\hd(\square(\chi_2)))
-\re_{\chi_2}(\hd(\square(\chi_1)))-\hd(\square(\{\chi_1,\chi_2\})) \\
=\re_{\chi_1}(\langle\psi,\ell_F(\chi_2)\rangle)
-\re_{\chi_2}(\langle\psi,\ell_F(\chi_1)\rangle)
-\langle\psi,\ell_F(\{\chi_1,\chi_2\})\rangle \\
-\re_{\chi_1}(\langle\ell_F^*(\psi),\chi_2\rangle)
+\re_{\chi_2}(\langle\ell_F^*(\psi),\chi_1\rangle)
+\langle\ell_F^*(\psi),\{\chi_1,\chi_2\}\rangle \\
=\langle\ell_\psi(\chi_1),\ell_F(\chi_2)\rangle
-\langle\ell_\psi(\chi_2),\ell_F(\chi_1)\rangle
-\langle\ell_{\Delta(F)}(\chi_1),\chi_2\rangle
+\langle\ell_{\Delta(F)}(\chi_2),\chi_1\rangle.
\end{multline*}
Therefore, the restriction of $\hd\circ d_\CC(\square)$ to $\Ei$ equals
to
\begin{multline*}
\left.\hd\circ d_\CC(\square)\right|_{\Ei}(\chi_1,\chi_2) \\
=\langle\ell_\psi(\chi_1),\ell_{\E}(\chi_2)\rangle
-\langle\ell_\psi(\chi_2),\ell_{\E}(\chi_1)\rangle
-\langle\Delta(\ell_{\E}(\chi_1)),\chi_2\rangle
+\langle\Delta(\ell_{\E}(\chi_2)),\chi_1\rangle \\
=\langle(\ell_\psi+\Delta^*)(\chi_1),\ell_{\E}(\chi_2)\rangle
-\langle(\ell_\psi+\Delta^*)(\chi_2),\ell_{\E}(\chi_1)\rangle
+\hd\gamma(\chi_1,\chi_2),
\end{multline*}
where $\gamma\in K$. This completes the proof.
\end{proof}

Now we apply these results to the problem of computing
conservation laws of an $\ell$-normal differential equation $\E$.

First, note that for a formally integrable equation $\E$ the
projections $\E^{(k+1)}\to\E^{(k)}$ are affine bundles, therefore
$\E^{(k+1)}$ and $\E^{(k)}$ are of the same homotopy type. Hence,
$H^*(\Ei)=H^*(\E)$.

Further, it follows from the two-line theorem that there exists the
following exact sequence:
\[
0\xra{} H^{n-1}(\E)\xra{} \hH^{n-1}(\E) \xra{d_1^{0,n-1}}
\ker\left(\ell_{\E}\right)^*.
\]
Recall that the group $\hH^{n-1}(\E)$ was interpreted as the group of
conservation laws of the equation $\E$ (see the beginning of Section
\vref{hc:sec}). Conservation laws $\omega\in H^{n-1}(\E)\subset{\bar
H}^{n-1}(\E)$ are called \emph{topological} (or \emph{rigid}), since
they are determined only by the topology of the equation $\E$. In
particular, the corresponding conserved quantities do not change under
deformations of solutions of the equation $\E$. Therefore topological
conservation laws are not very interesting for us and we consider the
quotient group $\cl(\E)=\hH^{n-1}(\E)/H^{n-1}(\E)$, called the group of
\emph{proper} conservation laws of the equation $\E$. The two-line
theorem implies immediately the following.

\begin{theorem}
\label{css.gf:thm}
If $\E$ is an $\ell$-normal equation, then
\[
\cl(\E)\subset\ker\ell_{\E}^*.
\]
If, moreover, $H^n(\E)\subset\hH^n(\E)$ \textup{(}in particular,
$H^n(\E)=0$\textup{)}, we have
\[
\cl(\E)=\ker d_1^{1,n-1}.
\]
\end{theorem}

Element $\psi\in\ker\ell_{\E}^*$ that corresponds to a
conservation law $[\omega]\in\cl(\E)$ is called its \emph{generating
function}.

Theorem \vref{css.gf:thm} gives an effective method for computing
conservation laws.

\begin{remark}
\label{css.dwd:rem}
In view of Proposition \vref{css.2nm1:prop}, elements of $E_1^{2,n-1}$
can be interpreted as mappings from $\ker\ell_\E^*$ to $\ker\ell_\E$,
i.e., from generating functions of conservation laws to symmetries of
$\E$.
\end{remark}

\begin{proposition}
\label{css.cleveq:prop}
Let $\E=\{u_t=f(t,x,u,u_x,u_{xx},\dotsc)\}$ be an evolution equation
and $J=(J_0,J_1,\dots,J_n)$ a conserved current for $\E$. Then the
generating function of $J$ is equal to $\psi=\bE(J_0)$, where $J_0$ is
the $t$-component of the conserved current that is regarded as a
function of $(t,x,u,u_x,u_{xx},\dotsc)$.
\end{proposition}

\begin{proof}
The restriction of the total derivative $D_t$ to the equation $\Ei$ has
the form: $D_t=\dd{}{t}+\re_f$. Hence
$\dd{J_0}{t}+\re_f(J_0)+\sum_{i=1}^nD_i(J_i)=0$. On the other hand,
$D_t=\dd{}{t}+\re_{u_t}$, therefore $D_t(J_0)+\sum_{i=1}^nD_i(J_i)
=\dd{J_0}{t}+\re_{u_t}(J_0)-\dd{J_0}{t}-\re_f(J_0)=\re_{u_t-f}(J_0)
=\ell_{J_0}(u_t-f)$. Thus $\psi=\ell_{J_0}^*(1)=\bE(J_0)$.
\end{proof}

Let $\vf\in\ker\ell_{\E}$ be a symmetry and
$[\omega]\in\hH^{n-1}(\E)$ a conservation law of the equation $\E$.
Then $\left[\re_{\vf}(\omega)\right]$ is a conservation law of $\E$
as well.

\begin{proposition}
If $\psi\in\ker\ell_{\E}^*$ is the generating function of a
conservation law $[\omega]$ of an $\ell$-normal equation $\E=\{F=0\}$,
then the generating function of the conservation law
$\left[\re_{\vf}(\omega) \right]$ has the form
$\re_{\vf}(\psi)+\Delta^*(\psi)$, where the operator
$\Delta\in\CDiff(P,P)$ is defined by $\re_{\vf}(F)=\Delta(F)$.
\end{proposition}

\begin{proof}
First, we have
\[
\langle\psi,\ell_\E(\chi)\rangle=\hd\ell_\omega(\chi)
+\hd\gamma(\ell_\E(\chi)),\qquad\chi\in\vk,
\]
where $\gamma\in\CDiff(P,\hL^{n-1})$.
Using the obvious formula
\[
\ell_{\re_{\chi_1}(\eta)}(\chi_2)=\re_{\chi_1}(\ell_\eta(\chi_2))
-\ell_\eta(\{\chi_1,\chi_2\}),\qquad\chi_1,\chi_1\in\vk,
\quad\eta\in\hL^n,
\]
where $\{\cdot,\cdot\}$ is the Jacobi bracket (see Definition
\vref{sec3:df:highj}, we obtain
\begin{multline*}
\hd\ell_{\re_\vf(\omega)}(\chi)
=\hd(\re_{\vf}(\ell_\omega(\chi)))-\hd(\ell_\omega(\{\vf,\chi\}))
=\re_{\vf}(\hd(\ell_\omega(\chi)))-\hd(\ell_\omega(\{\vf,\chi\})) \\
=\re_{\vf}(\langle\psi,\ell_\E(\chi)\rangle)
-\langle\psi,\ell_\E(\{\vf,\chi\})\rangle
-\re_{\vf}(\hd\gamma(\ell_\E(\chi)))
+\hd\gamma(\ell_\E(\{\vf,\chi\})) \\
=\langle\re_{\vf}(\psi),\ell_\E(\chi)\rangle
+\langle\psi,(\re_{\vf}(\ell_\E(\chi))
-\ell_\E(\{\vf,\chi\}))\rangle+\hd\gamma'(\ell_\E(\chi)) \\
=\langle\re_{\vf}(\psi),\ell_\E(\chi)\rangle
+\langle\psi,\left.\ell_{\re_\vf(F)}\right|_{\Ei}\rangle
+\hd\gamma'(\ell_\E(\chi)) \\
=\langle\re_{\vf}(\psi),\ell_\E(\chi)\rangle
+\langle\psi,\Delta(\ell_\E(\chi))\rangle+\hd\gamma'(\ell_\E(\chi)) \\
=\langle(\re_{\vf}+\Delta^*)(\psi),\ell_\E(\chi)\rangle
+\hd\gamma''(\ell_\E(\chi)),
\end{multline*}
where $\gamma',\gamma''\in\CDiff(P,\hL^{n-1})$.
This completes the proof.
\end{proof}

\subsection{Generating functions from the antifield-BRST standpoint}

A differential equation $\E=\{F=0\}$, $F\in P$, is called
\emph{normal}, if any \cd operator $\Delta$, such that $\Delta(F)=0$,
vanishes on $\Ei$. A normal equation is obviously $\ell$-normal.

Consider a normal equation $\E$ and the complex on $\Ji(\pi)$
\[
0\xla{}\F\xla{\delta}\CDiff(P,\F)\xla{\delta}\CDiff_{(2)}^{\alt}(P,\F)
\xla{\delta}\CDiff_{(3)}^{\alt}(P,\F)\xla{\delta}\dotsb,
\]
$\delta(\Delta)(p_1,\dots,p_k)=\Delta(F,p_1,\dots,p_k)$, $p_i\in P$.
This complex is exact in all terms except for the term $\F$. At points
$\theta\in\Ei$, the exactness follows immediately from the normality
condition. At points $\theta\notin\Ei$, this is a well known fact from
linear algebra (see Example \vref{aha.lv:exmp}). The homology in the
term $\F$ is clearly equal to $\F(\E)$.

In physics, this complex is said to be the \emph{Koszul--Tate
resolution}, and elements of $\CDiff_{(k)}^{\alt}(P,\F)$ are called
\emph{antifields}.

Consider the commutative diagram
\[
\begin{CD}
@. 0 @. 0 @. 0  \\
@. @AAA @AAA @AAA \\
0@<<<\hL^n @<\delta<< \CDiff(P,\hL^n) @<\delta<<
\CDiff_{(2)}^{\alt}(P,\hL^n) @<\delta<< \dotsb \\
@. @AA \hd A @AA \hd A @AA \hd A \\
0@<<<\hL^{n-1} @<\delta<< \CDiff(P,\hL^{n-1}) @<\delta<<
\CDiff_{(2)}^{\alt}(P,\hL^{n-1}) @<\delta<< \dotsb \\
@. @AA \hd A @AA \hd A @AA \hd A \\
0@<<<\hL^{n-2} @<\delta<< \CDiff(P,\hL^{n-2}) @<\delta<<
\CDiff_{(2)}^{\alt}(P,\hL^{n-2}) @<\delta<< \dotsb \\
@. @AA \hd A @AA \hd A @AA \hd A \\
@. \vdots @. \vdots @. \vdots
\end{CD}
\]
From the standard spectral sequence arguments (see the Appendix) and
Theorem \vref{sec2:th:symskew} it follows that
$\hH^q(\E)=H_{n-q}(L_{\bullet}^{\alt}(P),\delta)$. Since the complex
$(L_{\bullet}^{\alt}(P),\delta)$ is a direct summand in the complex
$(\CDiff_{(\bullet)}^{\alt}(P,\hat{P}),\delta)$, it is exact in all
degrees except for $0$ and $1$. This yields the \emph{two-line theorem}
for normal equations. We also get
\[
\hH^{n-1}(\E)=H_1(L_{\bullet}^{\alt}(P),\delta)
=\{\,\psi\in\hat{P}\bmod T\mid\langle\psi,F\rangle\in\hd\hL^{n-1}\,\},
\]
where $T=\{\,\psi\in\hat{P}\mid\psi=\square(F),
\ \square\in\CDiff(P,\hat{P}),\ \square^*=-\square\,\}$. The condition
$\langle\psi,F\rangle\in\hd\hL^{n-1}$ is equivalent to
$0=\bE\langle\psi,F\rangle=\ell_F^*(\psi)+\ell_{\psi}^*(F)$. So we
again obtain the correspondence between conservation laws and
generating functions together with the equation $\ell_\E^*(\psi)=0$.

\subsection{Euler--Lagrange equations}
\label{css.ele:subsec}

Consider the Euler--Lagrange equation $\E=\{{\bE}(\CL)=0\}$
corresponding to a Lagrangian $\CL=[\omega]\in\hH^n(\pi)$.
Let $\vf\in\vk(\pi)$ be a \emph{Noether symmetry} of $\CL$, i.e.,
$\re_{\vf}(\CL)=0$ on $\Ji(\pi)$.

\begin{xca}
Using Exercise \vref{css.uf:ex}, check that a Noether symmetry of $\CL$
is a symmetry of the corresponding equation $\E$ as well, i.e.,
$\sym(\CL)\subset\sym(\E)$.
\end{xca}

\begin{xca}
Show that if $E_2^{0,n}(\E)=0$, then finding of Noether symmetries of
the Lagrangian $\CL=[\omega]$ amounts to solution of the equation
${\bE}(\ell_{\omega}(\vf))
=\ell_{\bE(\CL)}(\vf)+\ell_{\vf}^*(\bE(\CL))=0$.
(Thus, to calculate the Noether symmetries of an Euler--Lagrange
equation one has no need to know the Lagrangian.)
\end{xca}

Let $\re_{\vf}(\omega)=\hd\nu$, where $\nu\in\hL^{n-1}(\pi)$.  By the
Green formula we have
\begin{multline*}
\re_{\vf}(\omega)-\hd\nu=\ell_{\omega}(\varphi)-\hd\nu
=\ell^*_{\omega}(1)(\varphi)+\hd\gamma(\varphi)-\hd\nu \\
=\bE(\CL)(\vf)+\hd(\gamma(\varphi)-\nu)=0.
\end{multline*}

Set
\[
\eta=\left.(\nu-\gamma(\varphi))\right|_{\Ei}\in\hL^{n-1}(\E).
\]
Thus, $\hd\left.\eta\right|_{\Ei}=0$, i.e., $[\eta]\in\hH^{n-1}(\E)$ is
a conservation law of the equation $\E$. The map
\[
\sym(\CL)\to\hH^{n-1}(\E),\qquad\varphi\mapsto[\eta],
\]
is said to be the \emph{Noether map}.

An arbitrariness in the choice of $\omega$ and $\nu$ leads to the
multivaluedness of the Noether map.

\begin{xca}
Check that the Noether map is well defined up to the image of the
natural homomorphism ${\hH}^{n-1}(\pi)\to{\hH}^{n-1}(\E)$.
\end{xca}

\begin{proposition}
If the Euler--Lagrange equation $\E$ corresponding to a Lagrangian
$\CL$ is $\ell$-normal, then the Noether map considered on the set of
Noether symmetries of $\CL$ is inverse to the differential
$d_1^{0,n-1}$.
\end{proposition}

\begin{proof}
On $\Ji(\pi)$ we have
\[
\hd\ell_{\eta}(\chi)=\ell_{\langle\bE(\CL),\vf\rangle}(\chi)
=\langle\ell_{\bE(\CL)}(\chi),\vf\rangle
+\langle\bE(\CL),\ell_{\vf}(\chi)\rangle.
\]
Therefore on $\Ei$ we obtain $\hd\ell_{\eta}(\chi)
=\langle\ell_\E(\chi),\vf\rangle$, i.e., $d_1^{0,n-1}([\eta])=\vf$.
\end{proof}

\begin{remark}
The Noether map can be understood as a procedure for finding a
conserved current corresponding to a given generating function.
\end{remark}

Thus, we see that if $\re_\vf$ is a Noether symmetry of a Lagrangian,
then $\vf$ is the generating function of a conservation law for the
corresponding Euler--Lagrange equation. This is the (first)
\emph{Noether theorem}. Note that since for Euler--Lagrange
equations one has $\ell_\E^*=\ell_\E$, the inverse Noether theorem is
obvious: if $\vf$ is the generating function of a conservation law for
an Euler--Lagrange equation, then $\vf$ is a symmetry for this
equation.

Let us discuss the Noether theorem from the antifield-BRST point of
view. Consider a $1$-cycle $\vf\in\vk$ of the complex
$L_{\bullet}^{\alt}(\hat{\vk})$. We have
$\langle\vf,\bE(\omega)\rangle\in\hd\hL^{n-1}$, where $\omega$ is a
density of the Lagrangian $\CL=[\omega]$. Hence
$\re_{\vf}(\omega)\in\hd\hL^{n-1}$ and, therefore, $\re_{\vf}(\CL)=0$,
i.e., $\vf$ is a Noether symmetry. Thus, the Koszul--Tate
resolution gives a homological interpretation of the Noether theorem.

Now, suppose that the Lagrangian has a \emph{gauge symmetry}, i.e.,
there exist an \fm $\mathfrak{a}$ and a \cd operator
$R\colon\mathfrak{a}\to\vk$ such that $R(\alpha)$ is a Noether symmetry
for any $\alpha\in\mathfrak{a}$.  This means that
$\re_{R(\alpha)}(\CL)=0$ or $\ell_\CL\circ R=0$. Hence
$R^*\circ\ell_\CL^*=0$ and, finally,
$R^*(\ell_\CL^*(1))=R^*(\bE(\CL))=0$. Thus, if the Lagrangian is
invariant under a gauge symmetry, then the\emph{Noether identities}
$R^*(\bE(\CL))=0$ between the Euler--Lagrange equations hold (the
second Noether theorem\label{css.snth:page}).

\subsection{The Hamiltonian formalism on $\Ji(\pi)$}
\label{css.hamform:subsec}

Let $A\in\CDiff(\hat{\vk}(\pi),\vk(\pi))$ be a \cd operator. Define the
\emph{Poisson bracket} on $\hH^n(\pi)$ corresponding to the operator
$A$ by the formula
\[
\{\omega_1,\omega_2\}_A
=\langle A(\bE(\omega_1)),\bE(\omega_2)\rangle,
\]
where $\langle\,,\,\rangle$ denotes the natural pairing
$\vk(\pi)\times\hat{\vk}(\pi)\to{\hH}^n(\pi)$.

The lemma below shows that the operator $A$ is uniquely determined by
the corresponding Poisson bracket.

\begin{lemma}
\label{css.sp:lem}
Let $\pi\colon E\to M$ be a vector bundle.
\begin{enumerate}
\item Consider an operator
$A\in\CDiff_{(l)}(\hat{\vk}(\pi),P)$, where $P$ is an $\F(\pi)$-module.
If for all $\omega_1,\dots,\omega_l\in{\hH}^n(\pi)$ one has
\[
A(\bE(\omega_1),\dots,\bE(\omega_l))=0,
\]
then $A=0$.
\item Consider an operator
$A\in\CDiff_{(l)}(\hat{\vk}(\pi),\hL^n(\pi))$. If for all cohomology classes
$\omega_1,\dots,\omega_l\in{\hH}^n(\pi)$ the element
$A(\bE(\omega_1),\dots,\bE(\omega_l))$ belongs to the image of
$\hd$, then $\im A\subset\im\hd$, i.e., $\mu_{(l-1)}(A)=0$ \textup{(}see
Subsection \ref{sub:Euler}\textup{)}.
\item Consider an operator
$A\in\CDiff_{(l-1)}(\hat{\vk}(\pi),\vk(\pi))$. If for all elements
$\omega_1,\dots,\omega_l\in{\hH}^n(\pi)$ one has $\langle A(\bE
(\omega_1),\dots,\bE(\omega_{l-1})),\bE(\omega_l)\rangle=0$, then $A=0$.
\end{enumerate}
\end{lemma}

\begin{proof}
$(1)$ It suffices to consider the case $l=1$. Obviously, on $\Ji(\pi)$
every element of $\hat{\vk}(\pi)=\hat{\F}(\pi,\pi)$ of the form
$\pi^*(f)$, with $f\in\hat{\Gamma}(\pi)$, (in other words, every
element of $\hat{\vk}(\pi)$ depending on base coordinates $x$ only) can
locally be presented in the form $\pi^*(f)=\bE(\omega)$ for some
$\omega\in\hL^n(\pi)$. Thus $A(\pi^*(f))=0$ for all $f$. Since $A$ is a
\cd operator, this implies $A=0$.

$(2)$ It is also sufficient to consider the case $l=1$. We have
$\bE(A(\bE(\omega)))=0$. Using Exercise \vref{css.uf:ex}, we get
\[
0=\bE(A(\bE(\omega)))=\ell_{\bE(\omega)}^*(A^*(1))
+\ell_{A^*(1)}^*(\bE(\omega))
\]
for all $\omega\in\hL^n(\pi)$. As above, we see that for any
$f\in\hat{\Gamma}(\pi)$ there exists $\omega\in\hL^n(\pi)$ such that
$\pi^*(f)=\bE(\omega)$. Since $\ell_{\pi^*(f)}=0$, we obtain
$\ell_{A^*(1)}^*(\pi^*(f))=0$. Hence $\ell_{A^*(1)}^*=0$, so that
$0=\bE(A(\bE(\omega))) =\ell_{\bE(\omega)}^*(A^*(1))$.

\begin{xca}
Prove that locally there exists a form $\omega\in\hL^n(\pi)$ such that
$\ell_{\bE(\omega)}$ is the identity operator.
\end{xca}
\noindent
Using this exercise, we get $0=A^*(1)=\mu(A)$, which is our claim.

$(3)$ The assertion follows immediately from (1) and (2) above.
\end{proof}

\begin{definition}
An operator $A\in\CDiff(\hat{\vk}(\pi),\vk(\pi))$ is called
\emph{Hamiltonian}, if its Poisson bracket defines a Lie algebra
structure on $\hH^n(\pi)$, i.e., if
\begin{gather}
\{\omega_1,\omega_2\}_A=-\{\omega_2,\omega_1\}_A, \label{css.ha:eq} \\
\{\{\omega_1,\omega_2\}_A,\omega_3\}_A+
\{\{\omega_2,\omega_3\}_A,\omega_1\}_A+
\{\{\omega_3,\omega_1\}_A,\omega_2\}_A=0.
\label{css.hj:eq}
\end{gather}
\end{definition}

The bracket $\{\,,\}_A$ is said to be a \emph{Hamiltonian structure}.

\begin{proposition}
The Poisson bracket $\{\,,\}_A$ is skew-symmetric, i.e., condition
\veqref{css.ha:eq} holds, if and only if the operator $A$ is
skew-adjoint, i.e., $A=-A^*$.
\end{proposition}

\begin{proof}
Since
\[
\{\omega_1,\omega_2\}_A+\{\omega_2,\omega_1\}_A
=\langle (A+A^*)(\bE(\omega_1)),\bE(\omega_2)\rangle,
\]
the claim follows immediately from the previous lemma.
\end{proof}

Now we shall prove criteria for checking an arbitrary skew-adjoint
operator $A\in\CDiff(\hat{\vk}(\pi),\vk(\pi))$ to be Hamiltonian.
For this, we need the following

\begin{lemma}
Consider an operator $A\in\CDiff(\hat{\vk}(\pi),\vk(\pi))$ and an
element $\psi\in\hat{\vk}(\pi)$. Define the operator
$\ell_{A,\psi}\in\CDiff(\vk(\pi),\vk(\pi))$ by
\[
\ell_{A,\psi}(\varphi)=(\ell_A(\varphi))(\psi)\qquad\vf\in\vk(\pi).
\]
Then
\eq{css.la:eq}
\ell^*_{A,\psi_1}(\psi_2)=\ell^*_{A^*,\psi_2}(\psi_1).
\end{equation}
\end{lemma}

\begin{proof}
By the Green formula,
\[
\langle A(\psi_1),\psi_2\rangle=\langle\psi_1,A^*(\psi_2)\rangle.
\]
Applying $\re_\vf$ to both sides, we get
\[
\langle\re_\vf(A)(\psi_1),\psi_2\rangle
=\langle\psi_1,\re_\vf(A^*)(\psi_2)\rangle,
\]
and so
\[
\langle\ell_{A,\psi_1}(\vf),\psi_2\rangle
=\langle\psi_1,\ell_{A^*,\psi_2}(\vf)\rangle.
\]
Again the Green formula yields
\[
\langle\vf,\ell_{A,\psi_1}^*(\psi_2)\rangle
=\langle\ell_{A^*,\psi_2}^*(\psi_1),\vf\rangle,
\]
and the lemma is proved.
\end{proof}

\begin{theorem}
Let $A\in\CDiff(\hat{\vk}(\pi),\vk(\pi))$ be a skew-adjoint
operator\textup{;} then the following conditions are
equivalent\textup{:}
\begin{enumerate}
\item $A$ is a Hamiltonian operator\textup{;}
\item $\langle\ell_A(A(\psi_1))(\psi_2),\psi_3\rangle
+\langle\ell_A(A(\psi_2))(\psi_3),\psi_1\rangle
+\langle\ell_A(A(\psi_3))(\psi_1),\psi_2\rangle=0$
for all $\psi_1,\psi_2,\psi_3\in\hat{\vk}(\pi)$\textup{;}
\item $\ell_{A,\psi_1}(A(\psi_2))-\ell_{A,\psi_2}(A(\psi_1))
=A(\ell^*_{A,\psi_2}(\psi_1))$ for all
$\psi_1,\psi_2\in\hat{\vk}(\pi)$\textup{;}
\item the expression
$\ell_{A,\psi_1}(A(\psi_2))+\frac12A(\ell^*_{A,\psi_1}(\psi_2))$
is symmetric with respect to $\psi_1,\psi_2\in\hat{\vk}(\pi)$\textup{;}
\item $[\re_{A(\psi)},A]=\ell_{A(\psi)}\circ A+A\circ\ell^*_{A(\psi)}$
for all $\psi\in\im\bE\subset\hat{\vk}(\pi)$.
\end{enumerate}
Moreover, it is sufficient to verify conditions \textup{(2)--(4)} for
elements $\psi_i\in\im\bE$ only.
\end{theorem}

\begin{proof}
Let $\omega_1,\omega_2,\omega_3\in{\hH}^n(\pi)$ and
$\psi_i=\bE(\omega_i)$. The Jacobi identity \veqref{css.hj:eq} yields

\begin{multline*}
\oint\{\{\omega_1,\omega_2\}_A,\omega_3\}_A
=\oint-\re_{A(\psi_3)}\langle A(\psi_1),\psi_2\rangle \\
=\oint-\langle \re_{A(\psi_3)}(A)(\psi_1),\psi_2\rangle
-\langle A(\ell_{\psi_1}(A(\psi_3))),\psi_2\rangle
-\langle A(\psi_1),\ell_{\psi_2}(A(\psi_3))\rangle \\
=\oint-\langle \ell_A(A(\psi_3))(\psi_1),\psi_2\rangle
+\langle A(\psi_2),\ell_{\psi_1}(A(\psi_3))\rangle
-\langle A(\psi_1),\ell_{\psi_2}(A(\psi_3))\rangle \\
=\oint-\langle \ell_A(A(\psi_3))(\psi_1),\psi_2\rangle=0,
\end{multline*}
where as above the symbol $\oint$ denotes the sum of cyclic
permutations. It follows from Lemma \vref{css.sp:lem} that this formula
holds for all $\psi_i\in\hat{\vk}(\pi)$. Criterion (2) is proved.

Rewrite the Jacobi identity in the form
\[
\langle\ell_{A,\psi_1}(A(\psi_2)),\psi_3\rangle
+\langle A(\psi_1),\ell^*_{A,\psi_3}(\psi_2)\rangle
-\langle A(\ell^*_{A,\psi_2}(\psi_1)),\psi_3\rangle=0.
\]
Using \veqref{css.la:eq}, we obtain
\[
\langle\ell_{A,\psi_1}(A(\psi_2)),\psi_3\rangle-
\langle\ell_{A,\psi_2}(A(\psi_1)),\psi_3\rangle-
\langle A(\ell^*_{A,\psi_2}(\psi_1)),\psi_3\rangle=0,
\]
which implies criterion (3).

The equivalence of criteria (3) and (4) follows from \veqref{css.la:eq}.

Finally, criterion (5) is equivalent to criterion (3) by virtue of the
following obvious equalities:
\begin{gather*}
[\re_{A(\psi_2)},A](\psi_1)=\ell_{A,\psi_1}(A(\psi_2)), \\
\ell_{A,\psi}\circ A=\ell_{A(\psi)}\circ A-A\circ\ell_{\psi}\circ A.
\end{gather*}
This concludes the proof.
\end{proof}

\begin{example}
Consider a skew-symmetric differential operator
$\Delta\colon\hat{\Gamma}(\pi)\to\Gamma(\pi)$. Then its lifting (see
Definition \vref{sec3:df:dolift})
$\CC\Delta\colon\hat{\vk}(\pi)\to\vk(\pi)$ is obviously a Hamiltonian
operator.
\end{example}

\begin{xca}
Check that in the case $n=\dim M=1$ and $m=\dim\pi=1$ operators of the
form $A=D_x^3+(\alpha+\beta u)D_x+\dfrac{\beta}{2}u_x$ are Hamiltonian.
\end{xca}

Let $A\colon\hat{\vk}(\pi)\to\vk(\pi)$ be a Hamiltonian operator. For
any $\omega\in{\hH}^n(\pi)$ the evolutionary vector field
$X_{\omega}=\re_{A({{\bE}}(\omega))}$ is called \emph{Hamiltonian}
vector field corresponding to the Hamiltonian~$\omega$.
Obviously,
\[
X_{\omega_1}(\omega_2)=\langle A{{\bE}}(\omega_1),
{\bE}(\omega_2)\rangle=\{\omega_1,\omega_2\}_A.
\]
This yields
\begin{multline*}
X_{\{\omega_1,\omega_2\}_A}(\omega)=\{\{\omega_1,\omega_2\}_A,\omega\}_A
=\{\omega_1,\{\omega_2,\omega\}_A\}_A
-\{\omega_2,\{\omega_1,\omega\}_A\}_A \\
= (X_{\omega_1}\circ X_{\omega_2}-X_{\omega_2}\circ X_{\omega_1})(\omega)
=[X_{\omega_1},X_{\omega_2}](\omega)
\end{multline*}
for all $\omega\in{\hH}^n(\pi)$. Thus
\eq{css.h:eq}
X_{\{\omega_1,\omega_2\}_A}=[X_{\omega_1},X_{\omega_2}].
\end{equation}

As with the finite dimensional Hamiltonian formalism, \vref{css.h:eq}
implies a result similar to the Noether theorem.

For each ${\CH}\in{\hH}^n(\pi)$, the evolution equation
\eq{css.he:eq}
u_t=A({{\bE}}({\CH})),
\end{equation}
corresponding to the Hamiltonian $\CH$ is called \emph{Hamiltonian}
evolution equation.

\begin{example}
The KdV equation $u_t=uu_x+u_{xxx}$ admits two Hamiltonian structures:
\begin{align*}
u_t&=D_x\left(\bE\left(\frac{u^3}{6}-\frac{u_x^2}{2}\right)\right) \\
\intertext{and}
u_t&=\left(D_x^3+\frac23uD_x+\frac13u_x\right)
\left(\bE\left(\frac{u^2}{2}\right)\right).
\end{align*}
\end{example}

\begin{theorem}
Hamiltonian operators take the generating function of a
conservation law of equation \veqref{css.he:eq} to the symmetry of this
equation.
\end{theorem}

\begin{proof}
Let $A$ be a Hamiltonian operator and
\[
\tilde\omega_0(t)+\tilde\omega_1(t)
\wedge dt\in\hL^n(\pi)\oplus\hL^{n-1}(\pi)\wedge dt
\]
be a conserved current of equation \veqref{css.he:eq}. This means that
$D_t(\omega_0(t))=0$, where $\omega_0(t)\in\hH^n(\pi)$ is the
horizontal cohomology class corresponding to the form
$\tilde\omega_0(t)$, and $D_t$ is the restriction of the total
derivative in $t$ to the equation.
Further,
\[
D_t(\omega_0)=\dd{\omega_0}{t}+
\re_{A({{\bE}}({\CH}))}(\omega_0)=\dd{\omega_0}{t}+\{{\CH},\omega_0\}.
\]
This yields
\[
\dd{}{t}X_{\omega_0}+[X_{\CH},X_{\omega_0}]=0.
\]
Hence $X_{\omega_0}=\re_{A({{\bE}}(\omega_0))}$ is  a symmetry
of~\veqref{css.he:eq}. It remains to recall that $\bE(\omega_0)$ is the
generating function of the conservation law under consideration (see
Proposition \vref{css.cleveq:prop}).
\end{proof}

\begin{remark}
Thus Hamiltonian operators are in a sense dual to elements of
$E_1^{2,n-1}$ (cf.\ Remark \vref{css.dwd:rem}).
\end{remark}

\subsection{On superequations}\label{sub:supeq}

The theory of this and preceding sections is based on the pure
algebraic considerations in Sections \ref{sec:calc} and
\ref{sec:lagr}. Therefore all results remain valid for the case of
differential superequations, provided one inserts the minus sign where
appropriate (detailed geometric definitions of superjets, super
Cartan distribution, and so on the reader can find, for example, in
\cite{HernRuipMuMasque1,HernRuipMuMasque2}). So we discuss here only a
couple of somewhat less obvious points and the coordinates formula.

Let $M$ be a supermanifold, $\dim M=n|m$, and $\pi$ be a superbundle
over $M$, $\dim\pi=s|t$. The following theorem is the superanalog of
theorem \vref{sec2:th:ber}.

\begin{theorem}
\begin{enumerate}
\item $\hat{A}_s=0$ for $s\ne n$.
\item $\hat{A}_n$ is the module of sections for the bundle $\Ber(M)$,
the latter being defined as follows: locally, sections of $\Ber(M)$
are written in the form $f(x)\bD(x)$, where  $f\in
\Ci(\CU)$ and $\bD$ is a basis local section that is
multiplied by the Berezin determinant of the Jacobi matrix under the
change of coordinates. The Berezin determinant of an {\text{even}} matrix
$\bigl(
\begin{smallmatrix}
A&B\\C&D
\end{smallmatrix}
\bigr)$ is equal to $\det(A-BD^{-1}C)(\det D)^{-1}$.
\end{enumerate}
\end{theorem}

\begin{proof} The assertion is local, so we can consider the domain
$\CU$ with local coordinates $x=(y_i,\xi_j)$, $i=1,\dots,n$,
$j=1,\dots,m$, and split the complex \veqref{sec2:eq:adjcom}
$\Diff^+(\Lambda^*)$ in the tensor product of complexes
$\Diff^+(\Lambda^*)_{\text{{\text{even}}}}
\otimes\Diff^+(\Lambda^*)_{\text{{\text{odd}}}}$,
where $\Diff^+(\Lambda^*)_{\text{{\text{even}}}}$ is complex
\eqref{sec2:eq:adjcom} on the underlying {\text{even}} domain of $\CU$
and $\Diff^+(\Lambda^*)_{\text{{\text{odd}}}}$ is the same complex for
the Grassmann algebra in variables $\xi_1,\ldots,\xi_m$.

We have $H^i(\Diff^+(\Lambda^*)_{{\text{even}}})=0$ for $i\ne n$
and $H^i(\Diff^+(\Lambda^*)_{{\text{even}}})=\Lambda_{\CU}^n$, where
$\Lambda_{\CU}^n$ is the module of $n$-form on the underlying even domain
of $\CU$. To compute the cohomology of $\Diff^+(\Lambda^*)_{{\text{odd}}}$
consider the quotient complexes
\[
0\xra{}\Smbl_k(A)_{{\text{odd}}}\xra{}\Smbl_{k+1}(\Lambda^1)_{{\text{odd}}}\xra{}\dotsb,
\]
where $\Smbl_k(P)_{{\text{odd}}}=\Diff_k^+(P)_{{\text{odd}}}/\Diff_{k-1}^+(P)_{{\text{odd}}}$.
Then an
easy calculation shows that these complexes are the Koszul complexes,
hence $H^i(\Diff^+(\Lambda^*))_{{\text{odd}}}=0$ for $i>0$ and
$H^0(\Diff^+(\Lambda^*))$ is a module of rank $1$. Therefore
$\hat{A}_i=H^i(\Diff^+(\Lambda^*))=0$ for $i\ne n$ and the only
operators that represent non-trivial cocycles have the form
$dy_1\wedge\dots\wedge dy_n\dd{^m}{\xi_1\cdots\partial\xi_m}f(y,\xi)$.

To complete the proof it remains to check that $\hat{A}_n$ is precisely
$\Ber(M)$, i.e., that changing coordinates we obtain:
\begin{multline*}
dy_1\wedge\dots\wedge dy_n\dd{^m}{\xi_1\dots\partial\xi_m}f\\
=dv_1\wedge\dots\wedge dv_n\dd{^m}{\eta_1\dots\partial\eta_m}
f\Ber\left(J\left(\frac{x}{z}\right)\right)+T,
\end{multline*}
where
$z=(v_i,\eta_j)$ is a new coordinate system on $\CU$, $\Ber$
denotes the Berezin determinant, $J\left(\dfrac{x}{z}\right)$ is the
Jacobi matrix, $T$ is cohomologous to zero. This is an immediate
consequence of the following well known formula for the Berezin
determinant:
$\Ber\bigl(
\begin{smallmatrix}
A&B\\C&D
\end{smallmatrix}
\bigr)=\det A\cdot\det\widetilde{D}$, where $\widetilde{D}$ is defined
by $\bigl(
\begin{smallmatrix}
A&B\\C&D
\end{smallmatrix}
\bigr)^{-1}=\bigl(
\begin{smallmatrix}
\widetilde{A}&\widetilde{B}\\\widetilde{C}&\widetilde{D}
\end{smallmatrix}
\bigr)$.
\end{proof}

The coordinate expression for the adjoint operator is as follows. Let
$\Delta\in\Diff(A,\CB)$ be a scalar operator
$\Delta=\sum_{\sigma}\bD a_{\sigma}\dd{^{|\sigma|}}
{x_{\sigma}}$.
Then
\[
\Delta^*=\sum_{\sigma}(-1)^{|\sigma|+a_{\sigma}x_{\sigma}}\bD
\dd{^{|\sigma|}}{x_{\sigma}}\circ a_{\sigma}.
\]
Here the symbol of an object used in exponent denotes the parity of the
object.

Now, consider a matrix operator $\Delta\colon P\to Q$,
$\Delta=\|\Delta_j^i\|$, where the matrix elements are defined by the
equalities $\Delta(\sum_{\alpha}e_{\alpha}f^{\alpha})
=\sum_{\alpha,\beta}e'_{\alpha}\Delta_{\beta}^{\alpha}(f^{\beta})$,
$\{e_i\}$ is a basis in $P$, $\{e'_i\}$ is a basis in $Q$. If $\bD$ is
even, then $\Delta^*$ has the form
\[
\bD(\Delta^*)_j^i=
(-1)^{(e_i+e'_j)(\Delta+e_i)}(\bD\Delta_i^j)^*.
\]
If $\bD $ is odd, then
\[
\bD((\Delta^*)^{\Pi})_j^i=
(-1)^{(e_i+\Delta)(e'_j+1)+\Delta e_i}(\bD\Delta_i^j)^*,
\]
where $\bigl(
\begin{smallmatrix}
A&B\\C&D
\end{smallmatrix}
\bigr)^{\Pi}=\bigl(
\begin{smallmatrix}
D&C\\B&A
\end{smallmatrix}
\bigr)$is the $\Pi$-transposition.

\begin{remark}
One has $(\Delta^{**})_j^i=(-1)^{e_i+e'_j}\Delta_j^i$.
\end{remark}

\begin{remark}
There is one point where we need to improve the algebraic theory of
differential operators to extend it to the supercase. This is the
definition of geometrical modules that should read:
\begin{definition}
A module $P$ over $\Ci(M)$ is called \emph{geometrical}, if
\[
\bigcap_{\substack{x\,\in\,M_{\text{rd}} \\ k\,\ge\,1}}\mu_x^kP=0,
\]
where $M_{\text{rd}}$ is the underlying even manifold of $M$ and
$\mu_x$ is the ideal in $\Ci(M)$ consisting of functions vanishing
at point $x\in M_{\text{rd}}$.
\end{definition}
\end{remark}

\newpage

\section*{Appendix: Homological algebra}
\label{sec:aha}
\refstepcounter{section}

In this appendix we sketch the basics of homological algebra. For an
extended discussion see, e.g.,
\cite{MacLane,HilStam,BottTu,McCleary,BourHomAlg}.

\subsection{Complexes}
\label{aha.es:subsec}

A sequence of vector spaces over a field $\Bbbk$ and linear mappings
\[
\dotsb\xra{}K^{i-1}\xra{d^{i-1}}K^i\xra{d^i}K^{i+1}\xra{d^{i+1}}\dotsb
\]
is said to be a \emph{complex} if the composition of any two
neighboring arrows is the zero map: $d^i\circ d^{i-1}=0$.

The maps $d^i$ are called \emph{differentials}. The index $i$ is
often omitted, so that the definition of a complex reads: $d^2=0$.

By definition, $\im d^{i-1}\subset\ker d^i$. The complex
$(K^{\bullet},d^{\bullet})$ is called \emph{exact} (or \emph{acyclic})
in degree $i$, if $\im d^{i-1}=\ker d^i$. A complex exact in all
degrees is called \emph{acyclic} (or \emph{exact}, or an \emph{exact
sequence}).

\begin{example}
The sequence $0\xra{}L\xra{f}K$ is always a complex. It is acyclic if
and only if $f$ is injection. The sequence $K\xra{g}M\xra{}0$ is always
a complex, as well. It is acyclic if and only if $g$ is surjection.

The sequence
\eq{aha.ses:eq}
0\xra{}L\xra{f}K\xra{g}M\xra{}0
\end{equation}
is a complex, if $g\circ
f=0$. It is exact, if and only if $f$ is injection, $g$ is surjection,
and $\im f=\ker g$. In this case we can identify $L$ with a subspace of
$K$ and $M$ with the quotient space $K/L$. Exact sequence
\veqref{aha.ses:eq} is called a \emph{short exact sequence} (or an
\emph{exact triple}).
\end{example}

\begin{example}
The \emph{de~Rham complex} is the complex of differential forms on a
smooth manifold $M$ with respect to the exterior derivation:
\[
\dotsb\xra{}\Lambda^{i-1}\xra{d}\Lambda^i\xra{d}\Lambda^{i+1}
\xra{d}\dotsb.
\]
\end{example}

The \emph{cohomology} of a complex $(K^{\bullet},d^{\bullet})$ is the
family of the spaces
\[
H^i(K^{\bullet},d^{\bullet})=\ker d^i/\im d^{i-1}.
\]
Thus, the equality $H^i(K^{\bullet},d^{\bullet})=0$ means that the
complex $(K^{\bullet},d^{\bullet})$ is acyclic in degree $i$. Note that
for the sake of brevity the cohomology is often denoted by
$H^i(K^{\bullet})$ or $H^i(d^{\bullet})$. Elements of $\ker
d^i\subset K^i$ are called $i$-dimensional \emph{cocycles}, elements of
$\im d^{i-1}\subset K^i$ are called $i$-dimensional
\emph{coboundaries}. Thus, the cohomology is the quotient space of the
space of all cocycles by the subspace of all coboundaries. Two cocycles
$k_1$ and $k_2$ from common cohomology coset, i.e., such that
$k_1-k_2\in\im d^{i-1}$, are called \emph{cohomologous}.

\begin{remark}
In the case of the complex of differential forms on a manifold cocycles
are called \emph{closed forms}, and coboundaries are called \emph{exact
forms}.
\end{remark}

\begin{remark}
It is clear that the definition of a complex can be immediately
generalized to modules over a ring instead of vector spaces.
\end{remark}

\begin{xca}
Prove that if
\[
\dotsb\xra{}Q^{i-1}\xra{d^{i-1}}Q^i\xra{d^i}Q^{i+1}\xra{d^{i+1}}\dotsb
\]
is a complex of modules (and $d^i$ are homomorphisms) and $P$ is a
projective module, then $H^i(Q^{\bullet}\otimes
P)=H^i(Q^{\bullet})\otimes P$.
\end{xca}

Complexes defined above are called \emph{cochain} to stress that the
differentials raise the dimension by 1. Inversion of arrows gives
\emph{chain} complexes
\[
\dotsb\xla{d_{i-1}}K_{i-1}\xla{d_i}K_i\xla{d_{i+1}}K_{i+1}\xla{}\dotsb,
\]
\emph{homology}, \emph{cycles}, \emph{boundaries}, etc. The difference
between these types of complex is pure terminological, so we shall
mainly restrict our considerations to cochain complexes.

A \emph{morphism} (or a \emph{cochain map}) of complexes $f\colon
K^{\bullet}\to L^{\bullet}$ is the family of linear mappings $f^i\colon
K^i\to L^i$ that commute with differentials, i.e., that make the
following diagram commutative:
\[
\begin{CD}
\dotsb@>>>K^{i-1}@>d_K^{i-1}>>K^i@>d_K^i>>K^{i+1}@>d_K^{i+1}>>\dotsb \\
@.     @VVf^{i-1}V         @VVf^iV     @VVf^{i+1}V            @. \\
\dotsb@>>>L^{i-1}@>d_L^{i-1}>>L^i@>d_L^i>>L^{i+1}@>d_L^{i+1}>>\dotsb.
\end{CD}
\]
Such a morphism induces the map $H^i(f)\colon H^i(K^{\bullet})\to
H^i(L^{\bullet})$, $[k]\mapsto [f(k)]$, where $k$ is a cocycle and
$[\;\cdot\;]$ denotes the cohomology coset. Clearly, $H^i(f\circ
g)=H^i(f)\circ H^i(g)$ (so that $H^i$ is a functor from the category of
complexes to the category of vector spaces). A morphism of complexes is
called \emph{quasiisomorphism} (or \emph{homologism}) if it induces an
isomorphism of cohomologies.

\begin{example}
A smooth map of manifolds $F\colon M_1\to M_2$ gives rise to the map of
differential forms
$F^*\colon\Lambda^{\bullet}(M_2)\to\Lambda^{\bullet}(M_1)$, such that
$d(F^*(\omega))=F^*(d(\omega))$. Thus $F^*$ is a cochain map and
induces the map of the de~Rham cohomologies $F^*\colon
H^{\bullet}(M_2)\to H^{\bullet}(M_1)$. In particular, if $M_1$ and
$M_2$ are diffeomorphic, then their de~Rham cohomologies are isomorphic.
\end{example}

\begin{xca}
Check that the wedge product on differential forms on $M$ induces a
well-defined multiplication on the de~Rham cohomology
$H^*(M)=\bigoplus_i H^i(M)$, which makes the de~Rham cohomology a
\emph{\textup{(}super\textup{)}algebra}, and not just a vector space.
Show that for diffeomorphic manifolds these algebras are isomorphic.
\end{xca}

Two morphisms of complexes $f^{\bullet},g^{\bullet}\colon
K^{\bullet}\to L^{\bullet}$ are called \emph{homotopic} if there exist
mappings $s^i\colon K^i\to L^{i-1}$, such that
\[
f^i-g^i=s^{i+1}d^i+d^{i-1}s^i.
\]
The mappings $s^i$ are called \emph{\textup{(}cochain\textup{)}
homotopy}.

\begin{proposition}
If morphisms $f^{\bullet}$ and $g^{\bullet}$ are homotopic, then
$H^i(f^{\bullet})=H^i(g^{\bullet})$ for all $i$.
\end{proposition}

\begin{proof}
Consider a cocycle $z\in K^i$, $dz=0$. Then
\[
f(z)-g(z)=(sd+ds)(z)=d(s(z)).
\]
Thus, $f(z)$ and $g(z)$ are cohomologous, and so
$H^i(f^{\bullet})=H^i(g^{\bullet})$.
\end{proof}

Two complexes $K^{\bullet}$ and $L^{\bullet}$ are said to be
\emph{cochain equivalent} if there exist morphisms $f^{\bullet}\colon
K^{\bullet}\to L^{\bullet}$ and $g^{\bullet}\colon L^{\bullet}\to
K^{\bullet}$ such that $g\circ f$ is homotopic to $\id_{K^{\bullet}}$
and $f\circ g$ is homotopic to $\id_{L^{\bullet}}$. Obviously, cochain
equivalent complexes have isomorphic cohomologies.

\begin{example}
Consider two maps of smooth manifolds $F_0,F_1\colon M_1\to M_2$ and
assume that they are homotopic (in the topological sense). Let us show
that the corresponding morphisms of the de~Rham complexes
$F_0^*,F_1^*\colon\Lambda^{\bullet}(M_2)\to\Lambda^{\bullet}(M_1)$ are
homotopic (in the above algebraic sense).

Let $F\colon M_1\times[0,1]\to M_2$ be the homotopy between $F_0$ and
$F_1$, $F_0(x)=F(x,0)$, $F_1(x)=F(x,1)$. Take a form
$\omega\in\Lambda^i(M_2)$. Then
\[
F^*(\omega)=\omega_1(t)+dt\wedge\omega_2(t),
\]
where $\omega_1(t)\in\Lambda^i(M_1)$,
$\omega_2(t)\in\Lambda^{i-1}(M_1)$ for each $t\in[0,1]$. In particular,
$F_0^*(\omega)=\omega_1(0)$ and $F_1^*(\omega)=\omega_1(1)$. Set
$s(\omega)=\int_0^1\omega_2(t)\,dt$. We have
$F^*(d\omega)=d(F^*(\omega))=d\omega_1(t)+dt\wedge\omega'_1(t)-dt\wedge
d\omega_2(t)$, where ${}'$ denotes the derivative in $t$. Hence,
$s(d(\omega))=\int_0^1(\omega'_1(t)-d\omega_2(t))\,dt
=\omega_1(1)-\omega_1(0)-d\int_0^1\omega_2(t)\,dt
=F_1^*(\omega)-F_0^*(\omega)-d(s(\omega))$, so $s$ is a homotopy
between $F_0^*$ and $F_1^*$.

\begin{xca}
Prove that if two manifolds $M_1$ and $M_2$ are homotopic (i.e., there
exist maps $f\colon M_1\to M_2$ and $g\colon M_2\to M_1$ such that the
maps $f\circ g$ and $g\circ f$ are homotopic to the identity maps),
then their cohomology are isomorphic.
\end{xca}

\begin{corollary}[Poincar\'e lemma]
Locally, every closed form $\omega\in\Lambda^i(M)$, $d\omega=0$,
$i\ge1$, is exact: $\omega=d\eta$.
\end{corollary}
\end{example}

A complex $K^{\bullet}$ is said to be \emph{homotopic to zero} if the
identity morphism $\id_{K^{\bullet}}$ homotopic to the zero morphism,
i.e., if there exist maps $s^i\colon K^i\to K^{i-1}$ such that
$\id_{K^{\bullet}}=sd+ds$. Obviously, a complex homotopic to zero has
the trivial cohomology.

\begin{example}
\label{aha.lv:exmp}
Let $V$ be a vector space. Take a nontrivial
linear functional $u\colon V\to\Bbbk$ and consider the complex
\[
0\xla{}\Bbbk\xla{d}V\xla{d}\Lambda^{2}(V)\xla{d}\dotsb\xla{d}
\Lambda^{n-1}(V)\xla{d}\Lambda^n(V)\xla{d}\dotsb,
\]
where $d$ is the inner product with $u$:
\[
d(v_1\wedge\dots\wedge v_k)=\sum_{i=1}^k(-1)^{i+1}u(v_i)
v_1\wedge\dots\wedge v_{i-1}\wedge v_{i+1}\wedge\dots\wedge v_k.
\]
Take also a nontrivial element $v\in V$ and consider the complex
\[
0\xra{}\Bbbk\xra{s}V\xra{s}\Lambda^{2}(V)\xra{s}\dotsb\xra{s}
\Lambda^{n-1}(V)\xra{s}\Lambda^n(V)\xra{s}\dotsb,
\]
where $s$ is the exterior product with $v$:
\[
s(v_1\wedge\dots\wedge v_k)=v\wedge v_1\wedge\dots\wedge v_k.
\]
Since $d$ is a derivation of the exterior algebra $\Lambda^*(V)$, we
have $(ds+sd)(w)=d(v\wedge w)+v\wedge dw=dv\wedge w=u(v)w$. This means
that both complexes under consideration are homotopic to zero and,
therefore, acyclic.
\end{example}

\begin{example}
Consider two complexes
\begin{gather}
\label{aha.kos:eq}
0\xla{}S^n(V)\xla{d}S^{n-1}(V)\otimes
V\xla{d}S^{n-2}(V)\otimes\Lambda^2(V)\xla{d}\dotsb, \\
\label{aha.derham:eq}
0\xra{}S^n(V)\xra{s}S^{n-1}(V)\otimes
V\xra{s}S^{n-2}(V)\otimes\Lambda^2(V)\xra{s}\dotsb,
\end{gather}
where
\begin{align*}
d(w\otimes v_1\wedge\dots\wedge v_q)
&=\sum_{i=1}^q(-1)^{i+1}v_iw\otimes v_1\wedge\dots\wedge v_{i-1}\wedge
v_{i+1}\wedge\dots\wedge v_q, \\
s(w_1\dotsm w_p\otimes v)
&=\sum_{i=1}^p w_1\dotsm w_{i-1}w_{i+1}\dotsm w_p\otimes w_i\wedge v.
\end{align*}
Both maps $d$ and $s$ are derivations of the algebra
$S^*(V)\otimes\Lambda^*(V)$, equipped with the grading induced from
$\Lambda^*(V)$, therefore their commutator is also a derivation. Noting
that on elements of $S^1(V)\otimes\Lambda^1(V)$ the commutator is
identical, we get the formula
\[
(ds+sd)(x)=(p+q)x,\qquad x\in S^p(V)\otimes\Lambda^q(V).
\]
Thus again both complexes under consideration are homotopic to zero
(for $n>0$). Complex \veqref{aha.kos:eq} is called the \emph{Koszul
complex}. Complex \veqref{aha.derham:eq} is the polynomial de~Rham
complex.
\end{example}

A complex $L^{\bullet}$ is called a \emph{subcomplex} of a complex
$K^{\bullet}$, if the spaces $L^i$ are subspaces of $K^i$, and the
differentials of $L^{\bullet}$ are restrictions of differentials of
$K^{\bullet}$, i.e., $d_K(L^{i-1})\subset L^i$. In this situation,
differentials of $K^{\bullet}$ induce differentials on quotient spaces
$M^i=K^i/L^i$ and we obtain the complex $M^{\bullet}$ called the
\emph{quotient complex} and denoted by $M^{\bullet}
=K^{\bullet}/L^{\bullet}$.

The cohomologies of complexes $K^{\bullet}$, $L^{\bullet}$, and
$M^{\bullet}=K^{\bullet}/L^{\bullet}$ are related to one another by the
following important mappings. First, the inclusion $\vf\colon
L^{\bullet}\to K^{\bullet}$ and the natural projection $\psi\colon
K^{\bullet}\to M^{\bullet}$ induce the cohomology mappings
$H^i(\vf)\colon H^i(L^{\bullet})\to H^i(K^{\bullet})$ and
$H^i(\psi)\colon H^i(K^{\bullet})\to H^i(M^{\bullet})$. There exists one
more somewhat less obvious mapping
\[
\partial^i\colon H^i(M^{\bullet})\to H^{i+1}(L^{\bullet})
\]
called the \emph{boundary} (or \emph{connecting}) mapping.

The map $\partial^i$ is defined as follows. Consider a cohomology class
$x\in H^i(M^{\bullet})$ represented by an element $y\in M^i$. Take an
element $z\in K^i$ such that $\psi(z)=y$. We have $\psi(dz)=d\psi(z)
=dy=0$, hence there exists an element $w\in L^{i+1}$ such that
$\vf(w)=dz$. Since $\vf(dw)=d\vf(w)=ddz=0$, we get $dw=0$, i.e., $w$ is
a cocycle. It can easily be checked that its cohomology class is
independent of the choice of $y$ and $z$. This class is the class
$\partial^i(x)$.

Thus, given a short exact sequence of complexes
\eq{aha.sesc:eq}
0\xra{}L^{\bullet}\xra{\vf}K^{\bullet}\xra{\psi}M^{\bullet}\xra{}0
\end{equation}
(this means that $\vf$ and $\psi$ are morphisms of complexes and for
each $i$ the sequences $0\xra{}L^i\xra{\vf^i}K^i\xra{\psi^i}M^i\xra{}0$
are exact), one has the following infinite sequence:
\begin{multline}
\label{aha.les:eq}
\dotsb\xra{H^{i-1}(\psi)}H^{i-1}(M^{\bullet})\xra{\partial^{i-1}}
H^i(L^{\bullet})\xra{H^i(\vf)}H^i(K^{\bullet})
\xra{H^i(\psi)}H^i(M^{\bullet}) \\
\xra{\partial^i}H^{i+1}(L^{\bullet})\xra{H^{i+1}(\vf)}\dotsb
\end{multline}

The main property of this sequence is the following.

\begin{theorem}
\label{aha.les:thm}
Sequence \veqref{aha.les:eq} is exact.
\end{theorem}

\begin{proof}
The proof is straightforward and is left to the reader.
\end{proof}

Sequence \veqref{aha.les:eq} is called the \emph{long exact sequence}
corresponding to short exact sequence of complexes \veqref{aha.sesc:eq}.

\begin{xca}
Consider the commutative diagram
\[
\begin{CD}
0 @>>> A_1 @>>> A_2 @>>> A_3 @>>> 0 \\
@. @VV f V @VV g V @VV h V  @. \\
0 @>>> B_1 @>>> B_2 @>>> B_3 @>>> 0.
\end{CD}
\]
Prove using Theorem \vref{aha.les:thm} that if $f$ and $h$ are
isomorphisms, then $g$ is also an isomorphism.
\end{xca}

\subsection{Spectral sequences}
\label{aha.ss:subsec}

Given a complex $K^{\bullet}$ and a subcomplex $L^{\bullet}\subset
K^{\bullet}$, the exact sequence \veqref{aha.les:eq} can tell something
about the cohomology of $K^{\bullet}$, if the cohomology of
$L^{\bullet}$ and $K^{\bullet}/L^{\bullet}$ are known. Now, suppose
that we are given a filtration of $K^{\bullet}$, that is a decreasing
sequence of subcomplexes
\[
K^{\bullet}\supset K_1^{\bullet}\supset K_2^{\bullet}
\supset K_3^{\bullet}\supset\dotsb.
\]
Then we obtain for each $p=0,1,2,\dotsc$ complexes
\[
\dotsb\xra{}E_0^{p,q-1}\xra{}E_0^{p,q}\xra{}E_0^{p,q+1}\xra{}\dotsb,
\]
where $E_0^{p,q}=K_p^{p+q}/K_{p+1}^{p+q}$.
The cohomologies $E_1^{p,q}=H^{p+q}(E_0^{p,\bullet})$ of these
complexes can be considered as the first approximation to the
cohomology of $K^{\bullet}$. The apparatus of spectral sequences
enables one to construct all successive approximations $E_r$, $r\ge1$.

\begin{definition}
A \emph{spectral sequence} is a sequence of vector spaces $E_r^{p,q}$,
$r\ge0$, and linear mappings $d_r^{p,q}\colon E_r^{p,q}\to
E_r^{p+r,q-r+1}$, such that $d_r^2=0$ (more precisely,
$d_r^{p+r,q-r+1}\circ d_r^{p,q}=0$) and the cohomology
$H^{p,q}(E_r^{\bullet,\bullet},d_r^{\bullet,\bullet})$ with respect to
the differential $d_r$ is isomorphic to $E_{r+1}^{p,q}$.
\end{definition}

Thus $E_r$ and $d_r$ determine $E_{r+1}$, but do not determine
$d_{r+1}$.

Usually, $p+q$, $p$, and $q$ are called respectively the degree,
the filtration degree, and the complementary degree.

It is convenient for each $r$ to picture the spaces $E_r^{p,q}$ as
integer points on the $(p,q)$-plane. The action of the differential
$d_r$ is shown as follows:

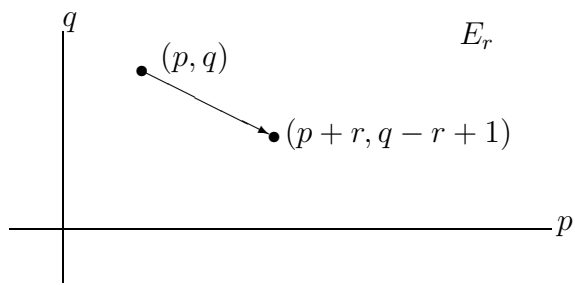
\begin{center}
\begin{picture}(220,110)
\put(0,20){\line(1,0){205}}
\put(20,0){\line(0,1){95}}
\put(207,20){$p$}
\put(20,97){$q$}
\put(50,80){\circle*{4}}
\put(100,55){\circle*{4}}
\put(50,80){\vector(2,-1){48}}
\put(57,82){$(p,q)$}
\put(104,53){$(p+r,q-r+1)$}
\put(170,90){$E_r$}
\end{picture}
\end{center}

Take an element $\alpha\in E_r^{p,q}$. If $d_r(\alpha)=0$ then $\alpha$
can be considered as an element of $E_{r+1}^{p,q}$. If again
$d_{r+1}(\alpha)=0$ then $\alpha$ can be considered as an element of
$E_{r+2}^{p,q}$ and so on. This allows us to define the following two
vector spaces:
\begin{align}
\label{aha.cbi:eq}
C_{\infty}^{p,q}&=\{\,\alpha\in E_0^{p,q}\mid\text{$d_0(\alpha)=0$,
$d_1(\alpha)=0$,\dots, $d_r(\alpha)=0$,\dots}\,\}, \\
B_{\infty}^{p,q}&=\{\,\alpha\in C_{\infty}^{p,q}\mid \text{there exists
an element $\beta\in E_r^{p,q}$ such that $\alpha=d_r(\beta)$}\,\}.
\notag
\end{align}
Set $E_{\infty}^{p,q}=C_{\infty}^{p,q}/B_{\infty}^{p,q}$. A spectral
sequence is called \emph{regular} if for any $p$ and $q$ there exists
$r_0$, such that $d_r^{p,q}=0$ for $r\ge r_0$. In this case there are
natural projections
\[
E_r^{p,q}\xra{}E_{r+1}^{p,q}\xra{}\dotsb\xra{}E_{\infty}^{p,q},\qquad
r\ge r_0,
\]
and $E_{\infty}^{p,q}=\injlim E_r^{p,q}$.

Let $E$ and ${}'\!E$ be two spectral sequences. A \emph{morphism}
$f\colon E\to {}'\!E$ is a family of mappings $f_r^{p,q}\colon
E_r^{p,q}\to {}'\!E_r^{p,q}$, such that $d_r\circ f_r=f_r\circ d_r$ and
$f_{r+1}=H(f_r)$. Obviously, a morphism $f\colon E\to {}'\!E$ induces
the maps $f_{\infty}^{p,q}\colon E_{\infty}^{p,q}\to
{}'\!E_{\infty}^{p,q}$. Further, it is clear that if $f_r$ is an
isomorphism, then $f_s$ are isomorphisms for all $s\ge r$. Moreover, if
the spectral sequences $E$ and ${}'\!E$ are regular, then $f_{\infty}$
is an isomorphism as well.

\begin{xca}
Assume that $E_r^{p,q}\ne0$ for $p\ge p_0$, $q\ge q_0$ only. Prove
that in this case there exists $r_0$ such that
$E_r^{p,q}=E_{r+1}^{p,q}=\dots=E_{\infty}^{p,q}$ for $r\ge r_0$.
\end{xca}

Consider a graded vector space $G=\bigoplus_{i\in\mathbb{Z}}G^i$ endowed
with a decreasing filtration $\dotsb\supset G_p\supset
G_{p+1}\supset\dotsb$, such that $\bigcap_pG_p=0$ and $\bigcup_pG_p=G$.
The filtration is called \emph{regular}, if for each $i$ there exists
$p$, such that $G_p^i=0$.

It is said that a spectral sequence $E$ \emph{converges} to $G$, if the
spectral sequence and the filtration of $G$ are regular and
$E_{\infty}^{p,q}$ is isomorphic to $G_p^{p+q}/G_{p+1}^{p+q}$.

\begin{xca}
Consider two spectral sequences $E$ and ${}'\!E$ that converge to $G$
and $G'$ respectively. Let $f\colon E\to {}'\!E$ be a morphism of
spectral sequences and $g\colon G\to G'$ be a map such that
$f_{\infty}^{p,q}\colon E_{\infty}^{p,q}\to {}'\!E_{\infty}^{p,q}$
coincides with the map induced by $g$. Prove that if the map
$f_r^{p,q}\colon E_r^{p,q}\to {}'\!E_r^{p,q}$ for some $r$ is an
isomorphism, then $g$ is an isomorphism too.
\end{xca}

Now we describe an important method for constructing spectral
sequences.

\begin{definition}
An \emph{exact couple} is a pair of vector spaces $(D,E)$ together with
mappings $i$, $j$, $k$, such that the diagram
\[
\begin{matrix}
D & \xra{\ i\ } & D \\
\rlap{${}_k$}\nwarrow\!\!\!\!\!\!\!\!\! &&\!\!\!\!\!\!\!\!\! \swarrow
\llap{${}_j$} \\
& E & \end{matrix}
\]
is exact in each vertex.
\end{definition}

Set $d=jk\colon E\to E$. Clearly, $d^2=0$, so that we can define
cohomology $H(E,d)$ with respect to $d$. Given an exact couple, one
defines the \emph{derived} couple
\[
\begin{matrix}
D' & \xra{\ i'\ } & D' \\
\rlap{${}_{k'}$}\nwarrow\!\!\!\!\!\!\!\!\! &&\!\!\!\!\!\!\!\!\! \swarrow
\llap{${}_{j'}$} \\
& E' & \end{matrix}
\]
as follows: $D'=\im i$, $E'=H(E,d)$, $i'$ is the restriction of $i$ to
$D'$, $j'(i(x))$ for $x\in D$ is the cohomology class of $j(x)$ in
$H(E)$, the map $k'$ takes a cohomology class $[y]$, $y\in E$, to the
element $k(y)\in D'$.

\begin{xca}
Check that mappings $i'$, $j'$, and $k'$ are well defined and that the
derived couple is an exact couple.
\end{xca}

Thus, starting from an exact couple $C_1=(D,E,i,j,k)$ we obtain the
sequence of exact couples $C_r=(D_r,E_r,i_r,j_r,k_r)$ such that
$C_{r+1}$ is the derived couple for $C_r$.

A direct description of $C_r$ in terms of $C_1$ is as follows.

\begin{proposition}
\label{aha.dd:prop}
The following isomorphisms hold for all $r$\textup{:}
\begin{align*}
D_r&=\im i^{r-1}, \\
E_r&=k^{-1}(\im i^{r-1})/j(\ker i^{r-1}).
\end{align*}
The map $i_r$ is the restriction of $i$ to $D_r$,
$j_r(i^{r-1}(x))=[j(x)]$, and $k_r([y])=k(y)$, where $[\;\cdot\;]$
denotes equivalence class modulo $j(\ker i^{r-1})$.
\end{proposition}

\begin{proof}
The proof is by induction on $r$ and is left to the reader.
\end{proof}

Now suppose that the exact couple $C_1$ is bigraded, i.e.,
$D=\bigoplus_{p,q}D^{p,q}$, $E=\bigoplus_{p,q}E^{p,q}$, and the
maps $i$, $j$, and $k$ have bidegrees $(-1,1)$, $(0,0)$,
$(1,0)$ respectively. In other words, one has:
\begin{align*}
i^{p,q}&\colon D^{p,q}\to D^{p-1,q+1}, \\
j^{p,q}&\colon D^{p,q}\to E^{p,q}, \\
k^{p,q}&\colon E^{p,q}\to D^{p+1,q}.
\end{align*}
It is clear that the derived couples $C_r$ are bigraded as well,
and the mappings $i_r$, $j_r$, and $k_r$ have bidegrees
$(-1,1)$, $(r-1,1-r)$, $(1,0)$ respectively. Therefore the differential
$d_r$ is a differential in $E_r$ and has bidegree $(r,1-r)$. Thus,
$(E_r^{p,q},d_r^{p,q})$ is a spectral sequence.

Now, suppose we are given a complex $K^{\bullet}$ with a decreasing
filtration $K_p^{\bullet}$. Each short exact sequence
\[
0\xra{}K_{p+1}^{\bullet}\xra{}K_p^{\bullet}\xra{}
K_p^{\bullet}/K_{p+1}^{\bullet}\xra{}0
\]
induces the corresponding long exact sequence:
\begin{multline*}
\dotsb\xra{k}H^{p+q}(K_{p+1}^{\bullet})\xra{i}H^{p+q}(K_p^{\bullet})
\xra{j}H^{p+q}(K_p^{\bullet}/K_{p+1}^{\bullet}) \\
\xra{k}H^{p+q+1}(K_{p+1}^{\bullet})\xra{i}\dotsb.
\end{multline*}
Hence, setting $D_1^{p,q}=H^{p+q}(K_p^{\bullet})$ and
$E_1^{p,q}=H^{p+q}(K_p^{\bullet}/K_{p+1}^{\bullet})$ we obtain a
bigraded exact couple, with mappings having bidegrees as above. Thus we
assign a spectral sequence to a complex with a filtration.

Let us compute the spaces $E_r^{p,q}$ in an explicit form. Consider
the upper term $k^{-1}(\im i^{r-1})$ from the expression for
$E_r^{p,q}$ (see Proposition \vref{aha.dd:prop}). An element of
$E_1^{p,q}$ is a class $[x]\in
H^{p+q}(K_p^{\bullet}/K_{p+1}^{\bullet})$, $x\in K_p^{p+q}$, $dx\in
K_{p+1}^{p+q}$. The class $[x]$lies in $k^{-1}(\im i^{r-1})$, if
$k([x])\in H^{p+q+1}(K_{p+r}^{\bullet})\subset
H^{p+q+1}(K_{p+1}^{\bullet})$. This is equivalent to $dx=y+dz$, with
$y\in K_{p+r}^{p+q}$, $z\in K_{p+1}^{p+q}$. Thus, we see that
$x=(x-z)+z$, with $d(x-z)\in K_{p+r}^{p+q}$. Denoting
\[
Z_r^{p,q}=\{\,w\in K_p^{p+q}\mid dw\in K_{p+r}^{p+q}\,\},
\]
we obtain $k^{-1}(\im i^{r-1})=Z_r^{p,q}+K_{p+1}^{p+q}$.

Further, consider the lower term $j(\ker i^{r-1})$ from the expression
for $E_r^{p,q}$. The kernel of the map $i^{r-1}\colon
H^{p+q}(K_p^{\bullet})\to H^{p+q}(K_{p-r+1}^{\bullet})$ consists
of cocycles $x\in K_p^{p+q}$ such that $x=dy$ for $y\in
K_{p-r+1}^{p+q-1}$. So $y\in Z_{r-1}^{p-r+1,q+r-2}$ and $\ker
i^{r-1}=dZ_{r-1}^{p-r+1,q+r-2}$. Then $j(\ker i^{r-1})
=dZ_{r-1}^{p-r+1,q+r-2}+K_{p+1}^{p+q}$.

Thus, we get
\[
E_r^{p,q}
=\frac{Z_r^{p,q}+K_{p+1}^{p+q}}{dZ_{r-1}^{p-r+1,q+r-2}+K_{p+1}^{p+q}}
=\frac{Z_r^{p,q}}{dZ_{r-1}^{p-r+1,q+r-2}+Z_{r-1}^{p+1,q-1}}.
\]

\begin{remark}
The last equality follows from the well known Noether modular
isomorphism
\[
\frac{M+N}{M_1+N}=\frac{M}{M_1+(M\cap N)},\qquad M_1\subset M.
\]
\end{remark}

\begin{theorem}
If the filtration of the complex $K^{\bullet}$ is regular, then the
spectral sequence of this complex converges to $H^{\bullet}(K^{\bullet})$
endowed with the filtration $H_p^k(K^{\bullet})=\im H^k(i_p)$, where
$i_p\colon K_p^{\bullet}\to K^{\bullet}$ is the natural inclusion.
\end{theorem}

\begin{proof}
Note first, that if the filtration of the complex $K^{\bullet}$ is
regular, then the spectral sequence of this complex is regular too.
Further, the spaces $C_{\infty}^{p,q}$ and $B_{\infty}^{p,q}$ (see
\veqref{aha.cbi:eq}) can easily be described by
\[
C_{\infty}^{p,q}=\frac{Z_{\infty}^{p,q}}{Z_{\infty}^{p+1,q-1}},\qquad
B_{\infty}^{p,q}=\frac{(K_p^{p+q}\cap d(K^{p+q-1}))
+Z_{\infty}^{p+1,q-1}}{Z_{\infty}^{p+1,q-1}},
\]
where $Z_{\infty}^{p,q}=\{\,w\in K_p^{p+q}\mid dw=0\,\}$, whence
\[
E_{\infty}^{p,q}=\frac{Z_{\infty}^{p,q}}{(K_p^{p+q}\cap
d(K^{p+q-1}))+Z_{\infty}^{p+1,q-1}}.
\]
Since $H_p^{p+q}(K^{\bullet})
=\dfrac{Z_{\infty}^{p,q}+d(K^{p+q-1})}{d(K^{p+q-1})}$, we have
\begin{multline*}
\frac{H_p^{p+q}(K^{\bullet})}{H_{p+1}^{p+q}(K^{\bullet})}
=\frac{Z_{\infty}^{p,q}+d(K^{p+q-1})}{Z_{\infty}^{p+1,q-1}
+d(K^{p+q-1})} \\
=\frac{Z_{\infty}^{p,q}}{Z_{\infty}^{p+1,q-1}+(K_p^{p+q}\cap
d(K^{p+q-1}))}=E_{\infty}^{p,q}.
\end{multline*}
This concludes the proof.
\end{proof}

\begin{definition}
A \emph{bicomplex} is a family of vector spaces $K^{\bullet,\bullet}$
and linear mappings $d'\colon K^{p,q}\to K^{p+1,q}$, $d''\colon
K^{p,q}\to K^{p,q+1}$, such that $(d')^2=0$, $(d'')^2=0$, and
$d'd''+d''d'=0$.
\end{definition}

Let $K^{\bullet}$ be the \emph{total} (or \emph{diagonal}) complex of
a bicomplex $K^{\bullet,\bullet}$, i.e., by definition,
$K^i=\bigoplus_{i=p+q}K^{p,q}$ and $d_K=d'+d''$. There are two obvious
filtration of $K^{\bullet}$:
\begin{align*}
\text{filtration I:}&\qquad{}'\!K_p^i=\bigoplus_{\substack{j+q=i\\ j\ge
p}}K^{j,q}, \\
\text{filtration II:}&\qquad{}''\!K_q^i=\bigoplus_{\substack{p+j=i\\
j\ge q}}K^{p,j}.
\end{align*}

These two filtrations yield two spectral sequences, denoted respectively
by ${}'\!E_r^{p,q}$ and ${}''\!E_r^{p,q}$.

It is easy to check that ${}'\!E_1^{p,q}={}''\!H^q(K^{p,\bullet})$ and
${}''\!E_1^{p,q}={}'\!H^q(K^{\bullet,p})$, where ${}'\!H$ (resp.,
${}''\!H$) denotes the cohomology with respect to $d'$ (resp., $d''$),
with the differential $d_1$ being induced respectively by $d'$ and
$d''$. Thus, we have:

\begin{proposition}
\label{aha.bic:prop}
${}'\!E_2^{p,q}={}'\!H^p({}''\!H^q(K^{\bullet,\bullet}))$ and
${}''\!E_2^{p,q}={}''\!H^p({}'\!H^q(K^{\bullet,\bullet}))$.
\end{proposition}

Now assume that both filtrations are regular.

\begin{xca}
Prove that
\begin{enumerate}
\item if $K^{p,q}=0$ for $q<q_0$ (resp., $p<p_0$), then the first
(resp., second) filtration is regular;
\item if $K^{p,q}=0$ for $q<q_0$ and $q>q_1$, then both filtration are
regular.
\end{enumerate}
\end{xca}

In this case both spectral sequences converge to the common limit
$H^{\bullet}(K^{\bullet})$.

\begin{remark}
This fact does not mean that both spectral sequences have a common
infinite term, because the two filtrations of $H^{\bullet}(K^{\bullet})$
are different.
\end{remark}

Let us illustrate Proposition \vref{aha.bic:prop}.

\begin{example}
\label{aha.ed:exmp}
Consider the commutative diagram
\[
\begin{CD}
@. \vdots @. \vdots @. \vdots @. \\
@.  @AAA  @AAA  @AAA  @. \\
0 @>>> K^{2,0} @>d_2>> K^{2,1} @>d_2>> K^{2,2} @>>> \dotsb \\
@.  @AA d_1 A  @AA d_1 A  @AA d_1 A  @. \\
0 @>>> K^{1,0} @>d_2>> K^{1,1} @>d_2>> K^{1,2} @>>> \dotsb \\
@.  @AA d_1 A  @AA d_1 A  @AA d_1 A  @. \\
0 @>>> K^{0,0} @>d_2>> K^{0,1} @>d_2>> K^{0,2} @>>> \dotsb \\
@.  @AAA  @AAA  @AAA  @. \\
@.  0  @. 0 @. 0  @.
\end{CD}
\]
and suppose that the differential $d_1$ is exact everywhere except for
the terms $K^{0,q}$ in the bottom row, and the differential $d_2$ is
exact everywhere except for the terms $K^{p,0}$ in the left column.
Thus, we have two complexes $L_1^{\bullet}$ and $L_2^{\bullet}$, where
$L_1^i=H^0(K^{i,\bullet},d_2)$, $L_2^i=H^0(K^{\bullet,i},d_1)$ and the
differential of $L_1$ (resp., $L_2$) is induced by $d_1$ (resp., $d_2$).
Consider the bicomplex $K^{\bullet,\bullet}$ with
$(d')^{p,q}=d_1^{p,q}$, $(d'')^{p,q}=(-1)^qd_2^{p,q}$. We easily get
\begin{align*}
{}'\!E_2^{p,q}&={}'\!E_3^{p,q}=\dots={}'\!E_{\infty}^{p,q}=
\begin{cases}
0& \text{if $q\ne0$}, \\
H^p(L_1^{\bullet})& \text{if $q=0$,}
\end{cases} \\
{}''\!E_2^{p,q}&={}''\!E_3^{p,q}=\dots={}''\!E_{\infty}^{p,q}=
\begin{cases}
0& \text{if $p\ne0$}, \\
H^q(L_2^{\bullet})& \text{if $p=0$.}
\end{cases}
\end{align*}
Since both spectral sequences converge to a common limit, we conclude
that $H^i(L_1^{\bullet})=H^i(L_2^{\bullet})$.

Let us describe this isomorphism in an explicit form. Consider a
cohomology class from $H^i(L_1^{\bullet})$. Choose an element
$k^{i,0}\in K^{i,0}$, $d_1(k^{i,0})=0$, $d_2(k^{i,0})=0$, that
represents this cohomology class. Since $d_1(k^{i,0})=0$, there exists
an element $x\in K^{i-1,0}$ such that $d_1(x)=k^{i,0}$. Set
$k^{i-1,1}=-d_2(x)\in K^{i-1,1}$. We have $d_2(k^{i-1,1})=0$ and
$d_1(k^{i-1,1})=-d_1(d_2(x))=-d_2(d_1(x))=-d_2(k^{i,0})=0$. Further,
the elements $k^{i,0}$ and $k^{i-1,1}$ are cohomologous in the total
complex $K^{\bullet}$: $k^{i,0}-k^{i-1,1}=d_1x+d_2x=(d'+d'')(x)$.
Continuing this process we obtain elements $k^{i-j,j}\in K^{i-j,j}$,
$d_1(k^{i-j,j})=0$, $d_2(k^{i-j,j})=0$, that are cohomologous in the
total complex $K^{\bullet}$. Thus, the above isomorphism takes the
cohomology class of $k^{i,0}$ to that of $k^{0,i}$.
\end{example}

\begin{xca}
Discuss an analog of Example \vref{aha.ed:exmp} for the commutative
diagram
\[
\begin{CD}
@. \vdots @. \vdots @. \vdots @. \\
@.  @VVV  @VVV  @VVV  @. \\
0 @<<< K^{2,0} @<d_2<< K^{2,1} @<d_2<< K^{2,2} @<<< \dotsb \\
@.  @VV d_1 V  @VV d_1 V  @VV d_1 V  @. \\
0 @<<< K^{1,0} @<d_2<< K^{1,1} @<d_2<< K^{1,2} @<<< \dotsb \\
@.  @VV d_1 V  @VV d_1 V  @VV d_1 V  @. \\
0 @<<< K^{0,0} @<d_2<< K^{0,1} @<d_2<< K^{0,2} @<<< \dotsb \\
@.  @VVV  @VVV  @VVV  @. \\
@.  0  @. 0 @. 0  @.
\end{CD}
\]
\end{xca}

\newpage


\end{document}